\input amstex

%

\def\next{AMS-SEKR}\ifx\styname\next \endinput\fi
\catcode`\@=11
\def\styname{AMS-SEKR}
\def\styversion{2.0}
{\W@{}\W@{\styname.STY - Version \styversion}\W@{}}
\hyphenation{acad-e-my acad-e-mies af-ter-thought anom-aly anom-alies
an-ti-deriv-a-tive an-tin-o-my an-tin-o-mies apoth-e-o-ses apoth-e-o-sis
ap-pen-dix ar-che-typ-al as-sign-a-ble as-sist-ant-ship as-ymp-tot-ic
asyn-chro-nous at-trib-uted at-trib-ut-able bank-rupt bank-rupt-cy
bi-dif-fer-en-tial blue-print busier busiest cat-a-stroph-ic
cat-a-stroph-i-cally con-gress cross-hatched data-base de-fin-i-tive
de-riv-a-tive dis-trib-ute dri-ver dri-vers eco-nom-ics econ-o-mist
elit-ist equi-vari-ant ex-quis-ite ex-tra-or-di-nary flow-chart
for-mi-da-ble forth-right friv-o-lous ge-o-des-ic ge-o-det-ic geo-met-ric
griev-ance griev-ous griev-ous-ly hexa-dec-i-mal ho-lo-no-my ho-mo-thetic
ideals idio-syn-crasy in-fin-ite-ly in-fin-i-tes-i-mal ir-rev-o-ca-ble
key-stroke lam-en-ta-ble light-weight mal-a-prop-ism man-u-script
mar-gin-al meta-bol-ic me-tab-o-lism meta-lan-guage me-trop-o-lis
met-ro-pol-i-tan mi-nut-est mol-e-cule mono-chrome mono-pole mo-nop-oly
mono-spline mo-not-o-nous mul-ti-fac-eted mul-ti-plic-able non-euclid-ean
non-iso-mor-phic non-smooth par-a-digm par-a-bol-ic pa-rab-o-loid
pa-ram-e-trize para-mount pen-ta-gon phe-nom-e-non post-script pre-am-ble
pro-ce-dur-al pro-hib-i-tive pro-hib-i-tive-ly pseu-do-dif-fer-en-tial
pseu-do-fi-nite pseu-do-nym qua-drat-ics quad-ra-ture qua-si-smooth
qua-si-sta-tion-ary qua-si-tri-an-gu-lar quin-tes-sence quin-tes-sen-tial
re-arrange-ment rec-tan-gle ret-ri-bu-tion retro-fit retro-fit-ted
right-eous right-eous-ness ro-bot ro-bot-ics sched-ul-ing se-mes-ter
semi-def-i-nite semi-ho-mo-thet-ic set-up se-vere-ly side-step sov-er-eign
spe-cious spher-oid spher-oid-al star-tling star-tling-ly
sta-tis-tics sto-chas-tic straight-est strange-ness strat-a-gem strong-hold
sum-ma-ble symp-to-matic syn-chro-nous topo-graph-i-cal tra-vers-a-ble
tra-ver-sal tra-ver-sals treach-ery turn-around un-at-tached un-err-ing-ly
white-space wide-spread wing-spread wretch-ed wretch-ed-ly Brown-ian
Eng-lish Euler-ian Feb-ru-ary Gauss-ian Grothen-dieck Hamil-ton-ian
Her-mit-ian Jan-u-ary Japan-ese Kor-te-weg Le-gendre Lip-schitz
Lip-schitz-ian Mar-kov-ian Noe-ther-ian No-vem-ber Rie-mann-ian
Schwarz-schild Sep-tem-ber
form per-iods Uni-ver-si-ty cri-ti-sism for-ma-lism}
\Invalid@\nofrills
\Invalid@\usualspace
\newif\ifnofrills@
\def\nofrills@#1#2{\relaxnext@
  \DN@{\ifx\next\nofrills
    \nofrills@true\let#2\relax\DN@\nofrills{\nextii@}%
  \else
    \nofrills@false\def#2{#1}\let\next@\nextii@\fi
\next@}}
\def\usualspace@#1{\ifnofrills@\def\usualspace{#1}\fi}
\def\addto#1#2{\csname \expandafter\eat@\string#1@\endcsname
  \expandafter{\the\csname \expandafter\eat@\string#1@\endcsname#2}}
\newdimen\bigsize@
\def\big@#1#2{{\hbox{$\left#2\vcenter to#1\bigsize@{}%
  \right.\nulldelimiterspace\z@\m@th$}}}
\def\big{\big@\@ne}
\def\Big{\big@{1.5}}
\def\bigg{\big@\tw@}
\def\Bigg{\big@{2.5}}
\def\raggedcenter@{\leftskip\z@ plus.4\hsize \rightskip\leftskip
 \parfillskip\z@ \parindent\z@ \spaceskip.3333em \xspaceskip.5em
 \pretolerance9999\tolerance9999 \exhyphenpenalty\@M
 \hyphenpenalty\@M \let\\\linebreak}
\def\upperspecialchars{\def\ss{SS}\let\i=I\let\j=J\let\ae\AE\let\oe\OE
  \let\o\O\let\aa\AA\let\l\L}
\def\uppercasetext@#1{%
  {\spaceskip1.2\fontdimen2\the\font plus1.2\fontdimen3\the\font
   \upperspecialchars\uctext@#1$\m@th\aftergroup\eat@$}}
\def\uctext@#1$#2${\endash@#1-\endash@$#2$\uctext@}
\def\endash@#1-#2\endash@{\uppercase{#1}\if\notempty{#2}--\endash@#2\endash@\fi}
\def\runaway@#1{\DN@{#1}\ifx\envir@\next@
  \Err@{You seem to have a missing or misspelled \string\end#1 ...}%
  \let\envir@\empty\fi}
\newif\iftemp@
\def\notempty#1{TT\fi\def\test@{#1}\ifx\test@\empty\temp@false
  \else\temp@true\fi \iftemp@}
\font@\tensmc=cmcsc10
\font@\sevenex=cmex7
\font@\sevenit=cmti7
\font@\eightrm=cmr8 
\font@\sixrm=cmr6 
\font@\eighti=cmmi8     \skewchar\eighti='177 
\font@\sixi=cmmi6       \skewchar\sixi='177   
\font@\eightsy=cmsy8    \skewchar\eightsy='60 
\font@\sixsy=cmsy6      \skewchar\sixsy='60   
\font@\eightex=cmex8
\font@\eightbf=cmbx8 
\font@\sixbf=cmbx6   
\font@\eightit=cmti8 
\font@\eightsl=cmsl8 
\font@\eightsmc=cmcsc8
\font@\eighttt=cmtt8 


\loadmsam
\loadmsbm
\loadeufm
\UseAMSsymbols
\newtoks\tenpoint@
\def\tenpoint{\normalbaselineskip12\p@
 \abovedisplayskip12\p@ plus3\p@ minus9\p@
 \belowdisplayskip\abovedisplayskip
 \abovedisplayshortskip\z@ plus3\p@
 \belowdisplayshortskip7\p@ plus3\p@ minus4\p@
 \textonlyfont@\rm\tenrm \textonlyfont@\it\tenit
 \textonlyfont@\sl\tensl \textonlyfont@\bf\tenbf
 \textonlyfont@\smc\tensmc \textonlyfont@\tt\tentt
 \textonlyfont@\bsmc\tenbsmc
 \ifsyntax@ \def\big##1{{\hbox{$\left##1\right.$}}}%
  \let\Big\big \let\bigg\big \let\Bigg\big
 \else
  \textfont\z@=\tenrm  \scriptfont\z@=\sevenrm  \scriptscriptfont\z@=\fiverm
  \textfont\@ne=\teni  \scriptfont\@ne=\seveni  \scriptscriptfont\@ne=\fivei
  \textfont\tw@=\tensy \scriptfont\tw@=\sevensy \scriptscriptfont\tw@=\fivesy
  \textfont\thr@@=\tenex \scriptfont\thr@@=\sevenex
        \scriptscriptfont\thr@@=\sevenex
  \textfont\itfam=\tenit \scriptfont\itfam=\sevenit
        \scriptscriptfont\itfam=\sevenit
  \textfont\bffam=\tenbf \scriptfont\bffam=\sevenbf
        \scriptscriptfont\bffam=\fivebf
  \setbox\strutbox\hbox{\vrule height8.5\p@ depth3.5\p@ width\z@}%
  \setbox\strutbox@\hbox{\lower.5\normallineskiplimit\vbox{%
        \kern-\normallineskiplimit\copy\strutbox}}%
 \setbox\z@\vbox{\hbox{$($}\kern\z@}\bigsize@=1.2\ht\z@
 \fi
 \normalbaselines\rm\ex@.2326ex\jot3\ex@\the\tenpoint@}
\newtoks\eightpoint@
\def\eightpoint{\normalbaselineskip10\p@
 \abovedisplayskip10\p@ plus2.4\p@ minus7.2\p@
 \belowdisplayskip\abovedisplayskip
 \abovedisplayshortskip\z@ plus2.4\p@
 \belowdisplayshortskip5.6\p@ plus2.4\p@ minus3.2\p@
 \textonlyfont@\rm\eightrm \textonlyfont@\it\eightit
 \textonlyfont@\sl\eightsl \textonlyfont@\bf\eightbf
 \textonlyfont@\smc\eightsmc \textonlyfont@\tt\eighttt
 \textonlyfont@\bsmc\eightbsmc
 \ifsyntax@\def\big##1{{\hbox{$\left##1\right.$}}}%
  \let\Big\big \let\bigg\big \let\Bigg\big
 \else
  \textfont\z@=\eightrm \scriptfont\z@=\sixrm \scriptscriptfont\z@=\fiverm
  \textfont\@ne=\eighti \scriptfont\@ne=\sixi \scriptscriptfont\@ne=\fivei
  \textfont\tw@=\eightsy \scriptfont\tw@=\sixsy \scriptscriptfont\tw@=\fivesy
  \textfont\thr@@=\eightex \scriptfont\thr@@=\sevenex
   \scriptscriptfont\thr@@=\sevenex
  \textfont\itfam=\eightit \scriptfont\itfam=\sevenit
   \scriptscriptfont\itfam=\sevenit
  \textfont\bffam=\eightbf \scriptfont\bffam=\sixbf
   \scriptscriptfont\bffam=\fivebf
 \setbox\strutbox\hbox{\vrule height7\p@ depth3\p@ width\z@}%
 \setbox\strutbox@\hbox{\raise.5\normallineskiplimit\vbox{%
   \kern-\normallineskiplimit\copy\strutbox}}%
 \setbox\z@\vbox{\hbox{$($}\kern\z@}\bigsize@=1.2\ht\z@
 \fi
 \normalbaselines\eightrm\ex@.2326ex\jot3\ex@\the\eightpoint@}

\font@\twelverm=cmr10 scaled\magstep1
\font@\twelveit=cmti10 scaled\magstep1
\font@\twelvesl=cmsl10 scaled\magstep1
\font@\twelvesmc=cmcsc10 scaled\magstep1
\font@\twelvett=cmtt10 scaled\magstep1
\font@\twelvebf=cmbx10 scaled\magstep1
\font@\twelvei=cmmi10 scaled\magstep1
\font@\twelvesy=cmsy10 scaled\magstep1
\font@\twelveex=cmex10 scaled\magstep1
\font@\twelvemsa=msam10 scaled\magstep1
\font@\twelveeufm=eufm10 scaled\magstep1
\font@\twelvemsb=msbm10 scaled\magstep1
\newtoks\twelvepoint@
\def\twelvepoint{\normalbaselineskip15\p@
 \abovedisplayskip15\p@ plus3.6\p@ minus10.8\p@
 \belowdisplayskip\abovedisplayskip
 \abovedisplayshortskip\z@ plus3.6\p@
 \belowdisplayshortskip8.4\p@ plus3.6\p@ minus4.8\p@
 \textonlyfont@\rm\twelverm \textonlyfont@\it\twelveit
 \textonlyfont@\sl\twelvesl \textonlyfont@\bf\twelvebf
 \textonlyfont@\smc\twelvesmc \textonlyfont@\tt\twelvett
 \textonlyfont@\bsmc\twelvebsmc
 \ifsyntax@ \def\big##1{{\hbox{$\left##1\right.$}}}%
  \let\Big\big \let\bigg\big \let\Bigg\big
 \else
  \textfont\z@=\twelverm  \scriptfont\z@=\tenrm  \scriptscriptfont\z@=\sevenrm
  \textfont\@ne=\twelvei  \scriptfont\@ne=\teni  \scriptscriptfont\@ne=\seveni
  \textfont\tw@=\twelvesy \scriptfont\tw@=\tensy \scriptscriptfont\tw@=\sevensy
  \textfont\thr@@=\twelveex \scriptfont\thr@@=\tenex
        \scriptscriptfont\thr@@=\tenex
  \textfont\itfam=\twelveit \scriptfont\itfam=\tenit
        \scriptscriptfont\itfam=\tenit
  \textfont\bffam=\twelvebf \scriptfont\bffam=\tenbf
        \scriptscriptfont\bffam=\sevenbf
  \setbox\strutbox\hbox{\vrule height10.2\p@ depth4.2\p@ width\z@}%
  \setbox\strutbox@\hbox{\lower.6\normallineskiplimit\vbox{%
        \kern-\normallineskiplimit\copy\strutbox}}%
 \setbox\z@\vbox{\hbox{$($}\kern\z@}\bigsize@=1.4\ht\z@
 \fi
 \normalbaselines\rm\ex@.2326ex\jot3.6\ex@\the\twelvepoint@}

\def\headfonts{\twelvepoint\bf}

\font@\fourteenrm=cmr10 scaled\magstep2
\font@\fourteenit=cmti10 scaled\magstep2
\font@\fourteensl=cmsl10 scaled\magstep2
\font@\fourteensmc=cmcsc10 scaled\magstep2
\font@\fourteentt=cmtt10 scaled\magstep2
\font@\fourteenbf=cmbx10 scaled\magstep2
\font@\fourteeni=cmmi10 scaled\magstep2
\font@\fourteensy=cmsy10 scaled\magstep2
\font@\fourteenex=cmex10 scaled\magstep2
\font@\fourteenmsa=msam10 scaled\magstep2
\font@\fourteeneufm=eufm10 scaled\magstep2
\font@\fourteenmsb=msbm10 scaled\magstep2
\newtoks\fourteenpoint@
\def\fourteenpoint{\normalbaselineskip15\p@
 \abovedisplayskip18\p@ plus4.3\p@ minus12.9\p@
 \belowdisplayskip\abovedisplayskip
 \abovedisplayshortskip\z@ plus4.3\p@
 \belowdisplayshortskip10.1\p@ plus4.3\p@ minus5.8\p@
 \textonlyfont@\rm\fourteenrm \textonlyfont@\it\fourteenit
 \textonlyfont@\sl\fourteensl \textonlyfont@\bf\fourteenbf
 \textonlyfont@\smc\fourteensmc \textonlyfont@\tt\fourteentt
 \textonlyfont@\bsmc\fourteenbsmc
 \ifsyntax@ \def\big##1{{\hbox{$\left##1\right.$}}}%
  \let\Big\big \let\bigg\big \let\Bigg\big
 \else
  \textfont\z@=\fourteenrm  \scriptfont\z@=\twelverm  \scriptscriptfont\z@=\tenrm
  \textfont\@ne=\fourteeni  \scriptfont\@ne=\twelvei  \scriptscriptfont\@ne=\teni
  \textfont\tw@=\fourteensy \scriptfont\tw@=\twelvesy \scriptscriptfont\tw@=\tensy
  \textfont\thr@@=\fourteenex \scriptfont\thr@@=\twelveex
        \scriptscriptfont\thr@@=\twelveex
  \textfont\itfam=\fourteenit \scriptfont\itfam=\twelveit
        \scriptscriptfont\itfam=\twelveit
  \textfont\bffam=\fourteenbf \scriptfont\bffam=\twelvebf
        \scriptscriptfont\bffam=\tenbf
  \setbox\strutbox\hbox{\vrule height12.2\p@ depth5\p@ width\z@}%
  \setbox\strutbox@\hbox{\lower.72\normallineskiplimit\vbox{%
        \kern-\normallineskiplimit\copy\strutbox}}%
 \setbox\z@\vbox{\hbox{$($}\kern\z@}\bigsize@=1.7\ht\z@
 \fi
 \normalbaselines\rm\ex@.2326ex\jot4.3\ex@\the\fourteenpoint@}

\def\chapheadfonts{\fourteenpoint\bf}

\font@\seventeenrm=cmr10 scaled\magstep3
\font@\seventeenit=cmti10 scaled\magstep3
\font@\seventeensl=cmsl10 scaled\magstep3
\font@\seventeensmc=cmcsc10 scaled\magstep3
\font@\seventeentt=cmtt10 scaled\magstep3
\font@\seventeenbf=cmbx10 scaled\magstep3
\font@\seventeeni=cmmi10 scaled\magstep3
\font@\seventeensy=cmsy10 scaled\magstep3
\font@\seventeenex=cmex10 scaled\magstep3
\font@\seventeenmsa=msam10 scaled\magstep3
\font@\seventeeneufm=eufm10 scaled\magstep3
\font@\seventeenmsb=msbm10 scaled\magstep3
\newtoks\seventeenpoint@
\def\seventeenpoint{\normalbaselineskip18\p@
 \abovedisplayskip21.6\p@ plus5.2\p@ minus15.4\p@
 \belowdisplayskip\abovedisplayskip
 \abovedisplayshortskip\z@ plus5.2\p@
 \belowdisplayshortskip12.1\p@ plus5.2\p@ minus7\p@
 \textonlyfont@\rm\seventeenrm \textonlyfont@\it\seventeenit
 \textonlyfont@\sl\seventeensl \textonlyfont@\bf\seventeenbf
 \textonlyfont@\smc\seventeensmc \textonlyfont@\tt\seventeentt
 \textonlyfont@\bsmc\seventeenbsmc
 \ifsyntax@ \def\big##1{{\hbox{$\left##1\right.$}}}%
  \let\Big\big \let\bigg\big \let\Bigg\big
 \else
  \textfont\z@=\seventeenrm  \scriptfont\z@=\fourteenrm  \scriptscriptfont\z@=\twelverm
  \textfont\@ne=\seventeeni  \scriptfont\@ne=\fourteeni  \scriptscriptfont\@ne=\twelvei
  \textfont\tw@=\seventeensy \scriptfont\tw@=\fourteensy \scriptscriptfont\tw@=\twelvesy
  \textfont\thr@@=\seventeenex \scriptfont\thr@@=\fourteenex
        \scriptscriptfont\thr@@=\fourteenex
  \textfont\itfam=\seventeenit \scriptfont\itfam=\fourteenit
        \scriptscriptfont\itfam=\fourteenit
  \textfont\bffam=\seventeenbf \scriptfont\bffam=\fourteenbf
        \scriptscriptfont\bffam=\twelvebf
  \setbox\strutbox\hbox{\vrule height14.6\p@ depth6\p@ width\z@}%
  \setbox\strutbox@\hbox{\lower.86\normallineskiplimit\vbox{%
        \kern-\normallineskiplimit\copy\strutbox}}%
 \setbox\z@\vbox{\hbox{$($}\kern\z@}\bigsize@=2\ht\z@
 \fi
 \normalbaselines\rm\ex@.2326ex\jot5.2\ex@\the\seventeenpoint@}

\font@\rrrrrm=cmr10 scaled\magstep4
\font@\bigtitlefont=cmbx10 scaled\magstep4

\parindent1pc
\normallineskiplimit\p@
\newdimen\indenti \indenti=2pc
\def\pageheight#1{\vsize#1}
\def\pagewidth#1{\hsize#1%
   \captionwidth@\hsize \advance\captionwidth@-2\indenti}
\pagewidth{30pc} \pageheight{47pc}
\def\topmatter{%
 \ifx\undefined\msafam
 \else\font@\eightmsa=msam8 \font@\sixmsa=msam6
   \ifsyntax@\else \addto\tenpoint{\textfont\msafam=\tenmsa
              \scriptfont\msafam=\sevenmsa \scriptscriptfont\msafam=\fivemsa}%
     \addto\eightpoint{\textfont\msafam=\eightmsa \scriptfont\msafam=\sixmsa
              \scriptscriptfont\msafam=\fivemsa}%
   \fi
 \fi
 \ifx\undefined\msbfam
 \else\font@\eightmsb=msbm8 \font@\sixmsb=msbm6
   \ifsyntax@\else \addto\tenpoint{\textfont\msbfam=\tenmsb
         \scriptfont\msbfam=\sevenmsb \scriptscriptfont\msbfam=\fivemsb}%
     \addto\eightpoint{\textfont\msbfam=\eightmsb \scriptfont\msbfam=\sixmsb
         \scriptscriptfont\msbfam=\fivemsb}%
   \fi
 \fi
 \ifx\undefined\eufmfam
 \else \font@\eighteufm=eufm8 \font@\sixeufm=eufm6
   \ifsyntax@\else \addto\tenpoint{\textfont\eufmfam=\teneufm
       \scriptfont\eufmfam=\seveneufm \scriptscriptfont\eufmfam=\fiveeufm}%
     \addto\eightpoint{\textfont\eufmfam=\eighteufm
       \scriptfont\eufmfam=\sixeufm \scriptscriptfont\eufmfam=\fiveeufm}%
   \fi
 \fi
 \ifx\undefined\eufbfam
 \else \font@\eighteufb=eufb8 \font@\sixeufb=eufb6
   \ifsyntax@\else \addto\tenpoint{\textfont\eufbfam=\teneufb
      \scriptfont\eufbfam=\seveneufb \scriptscriptfont\eufbfam=\fiveeufb}%
    \addto\eightpoint{\textfont\eufbfam=\eighteufb
      \scriptfont\eufbfam=\sixeufb \scriptscriptfont\eufbfam=\fiveeufb}%
   \fi
 \fi
 \ifx\undefined\eusmfam
 \else \font@\eighteusm=eusm8 \font@\sixeusm=eusm6
   \ifsyntax@\else \addto\tenpoint{\textfont\eusmfam=\teneusm
       \scriptfont\eusmfam=\seveneusm \scriptscriptfont\eusmfam=\fiveeusm}%
     \addto\eightpoint{\textfont\eusmfam=\eighteusm
       \scriptfont\eusmfam=\sixeusm \scriptscriptfont\eusmfam=\fiveeusm}%
   \fi
 \fi
 \ifx\undefined\eusbfam
 \else \font@\eighteusb=eusb8 \font@\sixeusb=eusb6
   \ifsyntax@\else \addto\tenpoint{\textfont\eusbfam=\teneusb
       \scriptfont\eusbfam=\seveneusb \scriptscriptfont\eusbfam=\fiveeusb}%
     \addto\eightpoint{\textfont\eusbfam=\eighteusb
       \scriptfont\eusbfam=\sixeusb \scriptscriptfont\eusbfam=\fiveeusb}%
   \fi
 \fi
 \ifx\undefined\eurmfam
 \else \font@\eighteurm=eurm8 \font@\sixeurm=eurm6
   \ifsyntax@\else \addto\tenpoint{\textfont\eurmfam=\teneurm
       \scriptfont\eurmfam=\seveneurm \scriptscriptfont\eurmfam=\fiveeurm}%
     \addto\eightpoint{\textfont\eurmfam=\eighteurm
       \scriptfont\eurmfam=\sixeurm \scriptscriptfont\eurmfam=\fiveeurm}%
   \fi
 \fi
 \ifx\undefined\eurbfam
 \else \font@\eighteurb=eurb8 \font@\sixeurb=eurb6
   \ifsyntax@\else \addto\tenpoint{\textfont\eurbfam=\teneurb
       \scriptfont\eurbfam=\seveneurb \scriptscriptfont\eurbfam=\fiveeurb}%
    \addto\eightpoint{\textfont\eurbfam=\eighteurb
       \scriptfont\eurbfam=\sixeurb \scriptscriptfont\eurbfam=\fiveeurb}%
   \fi
 \fi
 \ifx\undefined\cmmibfam
 \else \font@\eightcmmib=cmmib8 \font@\sixcmmib=cmmib6
   \ifsyntax@\else \addto\tenpoint{\textfont\cmmibfam=\tencmmib
       \scriptfont\cmmibfam=\sevencmmib \scriptscriptfont\cmmibfam=\fivecmmib}%
    \addto\eightpoint{\textfont\cmmibfam=\eightcmmib
       \scriptfont\cmmibfam=\sixcmmib \scriptscriptfont\cmmibfam=\fivecmmib}%
   \fi
 \fi
 \ifx\undefined\cmbsyfam
 \else \font@\eightcmbsy=cmbsy8 \font@\sixcmbsy=cmbsy6
   \ifsyntax@\else \addto\tenpoint{\textfont\cmbsyfam=\tencmbsy
      \scriptfont\cmbsyfam=\sevencmbsy \scriptscriptfont\cmbsyfam=\fivecmbsy}%
    \addto\eightpoint{\textfont\cmbsyfam=\eightcmbsy
      \scriptfont\cmbsyfam=\sixcmbsy \scriptscriptfont\cmbsyfam=\fivecmbsy}%
   \fi
 \fi
 \let\topmatter\relax}
\def\chapterno@{\uppercase\expandafter{\romannumeral\chaptercount@}}
\newcount\chaptercount@
\def\chapter{\nofrills@{\afterassignment\chapterno@
                        CHAPTER \global\chaptercount@=}\chapter@
 \DNii@##1{\leavevmode\hskip-\leftskip
   \rlap{\vbox to\z@{\vss\centerline{\eightpoint
   \chapter@##1\unskip}\baselineskip2pc\null}}\hskip\leftskip
   \nofrills@false}%
 \FN@\next@}
\newbox\titlebox@

\def\title{\nofrills@{\relax}\title@%
 \DNii@##1\endtitle{\global\setbox\titlebox@\vtop{\tenpoint\bf
 \raggedcenter@\ignorespaces
 \baselineskip1.3\baselineskip\title@{##1}\endgraf}%
 \ifmonograph@ \edef\next{\the\leftheadtoks}\ifx\next\empty
    \leftheadtext{##1}\fi
 \fi
 \edef\next{\the\rightheadtoks}\ifx\next\empty \rightheadtext{##1}\fi
 }\FN@\next@}
\newbox\authorbox@
\def\author#1\endauthor{\global\setbox\authorbox@
 \vbox{\tenpoint\smc\raggedcenter@\ignorespaces
 #1\endgraf}\relaxnext@ \edef\next{\the\leftheadtoks}%
 \ifx\next\empty\leftheadtext{#1}\fi}
\newbox\affilbox@
\def\affil#1\endaffil{\global\setbox\affilbox@
 \vbox{\tenpoint\raggedcenter@\ignorespaces#1\endgraf}}
\newcount\addresscount@
\addresscount@\z@
\def\address#1\endaddress{\global\advance\addresscount@\@ne
  \expandafter\gdef\csname address\number\addresscount@\endcsname
  {\vskip12\p@ minus6\p@\noindent\eightpoint\smc\ignorespaces#1\par}}
\def\email{\nofrills@{\eightpoint{\it E-mail\/}:\enspace}\email@
  \DNii@##1\endemail{%
  \expandafter\gdef\csname email\number\addresscount@\endcsname
  {\def\usualspace{{\it\enspace}}\smallskip\noindent\eightpoint\email@
  \ignorespaces##1\par}}%
 \FN@\next@}
\def\thedate@{}
\def\date#1\enddate{\gdef\thedate@{\tenpoint\ignorespaces#1\unskip}}
\def\thethanks@{}
\def\thanks#1\endthanks{\gdef\thethanks@{\eightpoint\ignorespaces#1.\unskip}}
\def\thekeywords@{}
\def\keywords{\nofrills@{{\it Key words and phrases.\enspace}}\keywords@
 \DNii@##1\endkeywords{\def\thekeywords@{\def\usualspace{{\it\enspace}}%
 \eightpoint\keywords@\ignorespaces##1\unskip.}}%
 \FN@\next@}
\def\thesubjclass@{}
\def\subjclass{\nofrills@{{\rm2020 {\it Mathematics Subject
   Classification\/}.\enspace}}\subjclass@
 \DNii@##1\endsubjclass{\def\thesubjclass@{\def\usualspace
  {{\rm\enspace}}\eightpoint\subjclass@\ignorespaces##1\unskip.}}%
 \FN@\next@}
\newbox\abstractbox@
\def\abstract{\nofrills@{{\smc Abstract.\enspace}}\abstract@
 \DNii@{\setbox\abstractbox@\vbox\bgroup\noindent$$\vbox\bgroup
  \def\envir@{abstract}\advance\hsize-2\indenti
  \usualspace@{{\enspace}}\eightpoint \noindent\abstract@\ignorespaces}%
 \FN@\next@}
\def\endabstract{\par\unskip\egroup$$\egroup}
\def\widestnumber#1#2{\begingroup\let\head\null\let\subhead\empty
   \let\subsubhead\subhead
   \ifx#1\head\global\setbox\tocheadbox@\hbox{#2.\enspace}%
   \else\ifx#1\subhead\global\setbox\tocsubheadbox@\hbox{#2.\enspace}%
   \else\ifx#1\key\bgroup\let\endrefitem@\egroup
        \key#2\endrefitem@\global\refindentwd\wd\keybox@
   \else\ifx#1\no\bgroup\let\endrefitem@\egroup
        \no#2\endrefitem@\global\refindentwd\wd\nobox@
   \else\ifx#1\page\global\setbox\pagesbox@\hbox{\quad\bf#2}%
   \else\ifx#1\item\setboxz@h{#2}\global\rosteritemwd\wdz@
        \global\advance\rosteritemwd by.5\parindent
   \else\message{\string\widestnumber is not defined for this option
   (\string#1)}%
\fi\fi\fi\fi\fi\fi\endgroup}
\newif\ifmonograph@
\def\Monograph{\monograph@true \let\headmark\rightheadtext
  \let\varindent@\indent \def\headfont@{\bf}\def\proclaimheadfont@{\smc}%
  \def\demofont@{\smc}}
\let\varindent@\indent

\newbox\tocheadbox@    \newbox\tocsubheadbox@
\newbox\tocbox@
\def\toc{\toc@{Contents}}
\def\newtocdefs{%
   \def \title##1\endtitle
       {\penaltyandskip@\z@\smallskipamount
        \hangindent\wd\tocheadbox@\noindent{\bf##1}}%
   \def \chapter##1{%
        Chapter \uppercase\expandafter{\romannumeral##1.\unskip}\enspace}%
   \def \specialhead##1\endspecialhead
       {\par\hangindent\wd\tocheadbox@ \noindent##1\par}%
   \def \head##1 ##2\endhead
       {\par\hangindent\wd\tocheadbox@ \noindent
        \if\notempty{##1}\hbox to\wd\tocheadbox@{\hfil##1\enspace}\fi
        ##2\par}%
   \def \subhead##1 ##2\endsubhead
       {\par\vskip-\parskip {\normalbaselines
        \advance\leftskip\wd\tocheadbox@
        \hangindent\wd\tocsubheadbox@ \noindent
        \if\notempty{##1}\hbox to\wd\tocsubheadbox@{##1\unskip\hfil}\fi
         ##2\par}}%
   \def \subsubhead##1 ##2\endsubsubhead
       {\par\vskip-\parskip {\normalbaselines
        \advance\leftskip\wd\tocheadbox@
        \hangindent\wd\tocsubheadbox@ \noindent
        \if\notempty{##1}\hbox to\wd\tocsubheadbox@{##1\unskip\hfil}\fi
        ##2\par}}}
\def\toc@#1{\relaxnext@
   \def\page##1%
       {\unskip\penalty0\null\hfil
        \rlap{\hbox to\wd\pagesbox@{\quad\hfil##1}}\hfilneg\penalty\@M}%
 \DN@{\ifx\next\nofrills\DN@\nofrills{\nextii@}%
      \else\DN@{\nextii@{{#1}}}\fi
      \next@}%
 \DNii@##1{%
\ifmonograph@\bgroup\else\setbox\tocbox@\vbox\bgroup
   \centerline{\headfont@\ignorespaces##1\unskip}\nobreak
   \vskip\belowheadskip \fi
   \setbox\tocheadbox@\hbox{0.\enspace}%
   \setbox\tocsubheadbox@\hbox{0.0.\enspace}%
   \leftskip\indenti \rightskip\leftskip
   \setbox\pagesbox@\hbox{\bf\quad000}\advance\rightskip\wd\pagesbox@
   \newtocdefs
 }%
 \FN@\next@}
\def\endtoc{\par\egroup}
\let\pretitle\relax
\let\preauthor\relax
\let\preaffil\relax
\let\predate\relax
\let\preabstract\relax
\let\prepaper\relax
\def\dedicatory #1\enddedicatory{\def\preabstract{{\medskip
  \eightpoint\it \raggedcenter@#1\endgraf}}}
\def\thetranslator@{}
\def\translator#1\endtranslator{\def\thetranslator@{\nobreak\medskip
 \line{\eightpoint\hfil Translated by \uppercase{#1}\qquad\qquad}\nobreak}}
\outer\def\endtopmatter{\runaway@{abstract}%
 \edef\next{\the\leftheadtoks}\ifx\next\empty
  \expandafter\leftheadtext\expandafter{\the\rightheadtoks}\fi
 \ifmonograph@\else
   \ifx\thesubjclass@\empty\else \makefootnote@{}{\thesubjclass@}\fi
   \ifx\thekeywords@\empty\else \makefootnote@{}{\thekeywords@}\fi
   \ifx\thethanks@\empty\else \makefootnote@{}{\thethanks@}\fi
 \fi
  \pretitle
  \ifmonograph@ \topskip7pc \else \topskip4pc \fi
  \box\titlebox@
  \topskip10pt
  \preauthor
  \ifvoid\authorbox@\else \vskip2.5pc plus1pc \unvbox\authorbox@\fi
  \preaffil
  \ifvoid\affilbox@\else \vskip1pc plus.5pc \unvbox\affilbox@\fi
  \predate
  \ifx\thedate@\empty\else \vskip1pc plus.5pc \line{\hfil\thedate@\hfil}\fi
  \preabstract
  \ifvoid\abstractbox@\else \vskip1.5pc plus.5pc \unvbox\abstractbox@ \fi
  \ifvoid\tocbox@\else\vskip1.5pc plus.5pc \unvbox\tocbox@\fi
  \prepaper
  \vskip2pc plus1pc
}
\def\document{\let\fontlist@\relax\let\alloclist@\relax
  \tenpoint}

\newskip\aboveheadskip       \aboveheadskip1.8\bigskipamount
\newdimen\belowheadskip      \belowheadskip1.8\medskipamount

\def\headfont@{\smc}
\def\penaltyandskip@#1#2{\relax\ifdim\lastskip<#2\relax\removelastskip
      \ifnum#1=\z@\else\penalty@#1\relax\fi\vskip#2%
  \else\ifnum#1=\z@\else\penalty@#1\relax\fi\fi}
\def\nobreak{\penalty\@M
  \ifvmode\def\penalty@{\let\penalty@\penalty\count@@@}%
  \everypar{\let\penalty@\penalty\everypar{}}\fi}
\let\penalty@\penalty
\def\heading#1\endheading{\head#1\endhead}

\def\specialheadfont@{\bf}
\outer\def\specialhead{\par\penaltyandskip@{-200}\aboveheadskip
  \begingroup\interlinepenalty\@M\rightskip\z@ plus\hsize \let\\\linebreak
  \specialheadfont@\noindent\ignorespaces}
\def\endspecialhead{\par\endgroup\nobreak\vskip\belowheadskip}
\let\headmark\eat@
\newskip\subheadskip       \subheadskip\medskipamount
\def\subheadfont@{\bf}
\outer\def\subhead{\nofrills@{.\enspace}\subhead@
 \DNii@##1\endsubhead{\par\penaltyandskip@{-100}\subheadskip
  \varindent@{\usualspace@{{\subheadfont@\enspace}}%
 \subheadfont@\ignorespaces##1\unskip\subhead@}\ignorespaces}%
 \FN@\next@}
\outer\def\subsubhead{\nofrills@{.\enspace}\subsubhead@
 \DNii@##1\endsubsubhead{\par\penaltyandskip@{-50}\medskipamount
      {\usualspace@{{\it\enspace}}%
  \it\ignorespaces##1\unskip\subsubhead@}\ignorespaces}%
 \FN@\next@}
\def\proclaimheadfont@{\bf}
\outer\def\proclaim{\runaway@{proclaim}\def\envir@{proclaim}%
  \nofrills@{.\enspace}\proclaim@
 \DNii@##1{\penaltyandskip@{-100}\medskipamount\varindent@
   \usualspace@{{\proclaimheadfont@\enspace}}\proclaimheadfont@
   \ignorespaces##1\unskip\proclaim@
  \sl\ignorespaces}%
 \FN@\next@}
\outer\def\endproclaim{\let\envir@\relax\par\rm
  \penaltyandskip@{55}\medskipamount}
\def\demoheadfont@{\it}
\def\demo{\runaway@{proclaim}\nofrills@{.\enspace}\demo@
     \DNii@##1{\par\penaltyandskip@\z@\medskipamount
  {\usualspace@{{\demoheadfont@\enspace}}%
  \varindent@\demoheadfont@\ignorespaces##1\unskip\demo@}\rm
  \ignorespaces}\FN@\next@}
\def\enddemo{\par\medskip}
\def\qed{\ifhmode\unskip\nobreak\fi\quad\ifmmode\square\else$\m@th\square$\fi}
\let\remark\demo
\let\endremark\enddemo
\def\definition{\runaway@{proclaim}%
  \nofrills@{.\demoheadfont@\enspace}\definition@
        \DNii@##1{\penaltyandskip@{-100}\medskipamount
        {\usualspace@{{\demoheadfont@\enspace}}%
        \varindent@\demoheadfont@\ignorespaces##1\unskip\definition@}%
        \rm \ignorespaces}\FN@\next@}


\newdimen\rosteritemwd
\newcount\rostercount@
\newif\iffirstitem@
\let\plainitem@\item
\newtoks\everypartoks@
\def\par@{\everypartoks@\expandafter{\the\everypar}\everypar{}}
\def\roster{\edef\leftskip@{\leftskip\the\leftskip}%
 \relaxnext@
 \rostercount@\z@  
 \def\item{\FN@\rosteritem@}%
 \DN@{\ifx\next\runinitem\let\next@\nextii@\else
  \let\next@\nextiii@\fi\next@}%
 \DNii@\runinitem  
  {\unskip  
   \DN@{\ifx\next[\let\next@\nextii@\else
    \ifx\next"\let\next@\nextiii@\else\let\next@\nextiv@\fi\fi\next@}%
   \DNii@[####1]{\rostercount@####1\relax
    \enspace{\rm(\number\rostercount@)}~\ignorespaces}%
   \def\nextiii@"####1"{\enspace{\rm####1}~\ignorespaces}%
   \def\nextiv@{\enspace{\rm(1)}\rostercount@\@ne~}%
   \par@\firstitem@false  
   \FN@\next@}%
 \def\nextiii@{\par\par@  
  \penalty\@m\smallskip\vskip-\parskip
  \firstitem@true}%
 \FN@\next@}
\def\rosteritem@{\iffirstitem@\firstitem@false\else\par\vskip-\parskip\fi
 \leftskip3\parindent\noindent  
 \DNii@[##1]{\rostercount@##1\relax
  \llap{\hbox to2.5\parindent{\hss\rm(\number\rostercount@)}%
   \hskip.5\parindent}\ignorespaces}%
 \def\nextiii@"##1"{%
  \llap{\hbox to2.5\parindent{\hss\rm##1}\hskip.5\parindent}\ignorespaces}%
 \def\nextiv@{\advance\rostercount@\@ne
  \llap{\hbox to2.5\parindent{\hss\rm(\number\rostercount@)}%
   \hskip.5\parindent}}%
 \ifx\next[\let\next@\nextii@\else\ifx\next"\let\next@\nextiii@\else
  \let\next@\nextiv@\fi\fi\next@}

\newif\ifnextRunin@
\def\endroster{\relaxnext@
 \par\leftskip@  
 \penalty-50 \vskip-\parskip\smallskip  
 \DN@{\ifx\next\Runinitem\let\next@\relax
  \else\nextRunin@false\let\item\plainitem@  
   \ifx\next\par 
    \DN@\par{\everypar\expandafter{\the\everypartoks@}}%
   \else  
    \DN@{\noindent\everypar\expandafter{\the\everypartoks@}}%
  \fi\fi\next@}%
 \FN@\next@}
\newcount\rosterhangafter@
\def\Runinitem#1\roster\runinitem{\relaxnext@
 \rostercount@\z@ 
 \def\item{\FN@\rosteritem@}%
 \def\runinitem@{#1}%
 \DN@{\ifx\next[\let\next\nextii@\else\ifx\next"\let\next\nextiii@
  \else\let\next\nextiv@\fi\fi\next}%
 \DNii@[##1]{\rostercount@##1\relax
  \def\item@{{\rm(\number\rostercount@)}}\nextv@}%
 \def\nextiii@"##1"{\def\item@{{\rm##1}}\nextv@}%
 \def\nextiv@{\advance\rostercount@\@ne
  \def\item@{{\rm(\number\rostercount@)}}\nextv@}%
 \def\nextv@{\setbox\z@\vbox  
  {\ifnextRunin@\noindent\fi  
  \runinitem@\unskip\enspace\item@~\par  
  \global\rosterhangafter@\prevgraf}%
  \firstitem@false  
  \ifnextRunin@\else\par\fi  
  \hangafter\rosterhangafter@\hangindent3\parindent
  \ifnextRunin@\noindent\fi  
  \runinitem@\unskip\enspace 
  \item@~\ifnextRunin@\else\par@\fi  
  \nextRunin@true\ignorespaces}%
 \FN@\next@}
\def\footmarkform@#1{$\m@th^{#1}$}
\let\thefootnotemark\footmarkform@
\def\makefootnote@#1#2{\insert\footins
 {\interlinepenalty\interfootnotelinepenalty
 \eightpoint\splittopskip\ht\strutbox\splitmaxdepth\dp\strutbox
 \floatingpenalty\@MM\leftskip\z@\rightskip\z@\spaceskip\z@\xspaceskip\z@
 \leavevmode{#1}\footstrut\ignorespaces#2\unskip\lower\dp\strutbox
 \vbox to\dp\strutbox{}}}
\newcount\footmarkcount@
\footmarkcount@\z@
\def\footnotemark{\let\@sf\empty\relaxnext@
 \ifhmode\edef\@sf{\spacefactor\the\spacefactor}\/\fi
 \DN@{\ifx[\next\let\next@\nextii@\else
  \ifx"\next\let\next@\nextiii@\else
  \let\next@\nextiv@\fi\fi\next@}%
 \DNii@[##1]{\footmarkform@{##1}\@sf}%
 \def\nextiii@"##1"{{##1}\@sf}%
 \def\nextiv@{\iffirstchoice@\global\advance\footmarkcount@\@ne\fi
  \footmarkform@{\number\footmarkcount@}\@sf}%
 \FN@\next@}
\def\footnotetext{\relaxnext@
 \DN@{\ifx[\next\let\next@\nextii@\else
  \ifx"\next\let\next@\nextiii@\else
  \let\next@\nextiv@\fi\fi\next@}%
 \DNii@[##1]##2{\makefootnote@{\footmarkform@{##1}}{##2}}%
 \def\nextiii@"##1"##2{\makefootnote@{##1}{##2}}%
 \def\nextiv@##1{\makefootnote@{\footmarkform@{\number\footmarkcount@}}{##1}}%
 \FN@\next@}
\def\footnote{\let\@sf\empty\relaxnext@
 \ifhmode\edef\@sf{\spacefactor\the\spacefactor}\/\fi
 \DN@{\ifx[\next\let\next@\nextii@\else
  \ifx"\next\let\next@\nextiii@\else
  \let\next@\nextiv@\fi\fi\next@}%
 \DNii@[##1]##2{\footnotemark[##1]\footnotetext[##1]{##2}}%
 \def\nextiii@"##1"##2{\footnotemark"##1"\footnotetext"##1"{##2}}%
 \def\nextiv@##1{\footnotemark\footnotetext{##1}}%
 \FN@\next@}
\def\adjustfootnotemark#1{\advance\footmarkcount@#1\relax}
\def\footnoterule{\kern-3\p@
  \hrule width 5pc\kern 2.6\p@} 
\def\captionfont@{\smc}
\def\topcaption#1#2\endcaption{%
  {\dimen@\hsize \advance\dimen@-\captionwidth@
   \rm\raggedcenter@ \advance\leftskip.5\dimen@ \rightskip\leftskip
  {\captionfont@#1}%
  \if\notempty{#2}.\enspace\ignorespaces#2\fi
  \endgraf}\nobreak\bigskip}
\def\botcaption#1#2\endcaption{%
  \nobreak\bigskip
  \setboxz@h{\captionfont@#1\if\notempty{#2}.\enspace\rm#2\fi}%
  {\dimen@\hsize \advance\dimen@-\captionwidth@
   \leftskip.5\dimen@ \rightskip\leftskip
   \noindent \ifdim\wdz@>\captionwidth@ 
   \else\hfil\fi 
  {\captionfont@#1}\if\notempty{#2}.\enspace\rm#2\fi\endgraf}}
\def\@ins{\par\begingroup\def\vspace##1{\vskip##1\relax}%
  \def\captionwidth##1{\captionwidth@##1\relax}%
  \setbox\z@\vbox\bgroup} 
\def\block{\RIfMIfI@\nondmatherr@\block\fi
       \else\ifvmode\vskip\abovedisplayskip\noindent\fi
        $$\def\endblock{\par\egroup$$}\fi
  \vbox\bgroup\advance\hsize-2\indenti\noindent}
\def\endblock{\par\egroup}
\def\cite#1{{\rm[{\citefont@\m@th#1}]}}
\def\citefont@{\rm}
\def\refsfont@{\eightpoint}
\outer\def\Refs{\runaway@{proclaim}%
 \relaxnext@ \DN@{\ifx\next\nofrills\DN@\nofrills{\nextii@}\else
  \DN@{\nextii@{References}}\fi\next@}%
 \DNii@##1{\penaltyandskip@{-200}\aboveheadskip
  \line{\hfil\headfont@\ignorespaces##1\unskip\hfil}\nobreak
  \vskip\belowheadskip
  \begingroup\refsfont@\sfcode`.=\@m}%
 \FN@\next@}
\def\endRefs{\par\endgroup}
\newbox\nobox@            \newbox\keybox@           \newbox\bybox@
\newbox\paperbox@         \newbox\paperinfobox@     \newbox\jourbox@
\newbox\volbox@           \newbox\issuebox@         \newbox\yrbox@
\newbox\pagesbox@         \newbox\bookbox@          \newbox\bookinfobox@
\newbox\publbox@          \newbox\publaddrbox@      \newbox\finalinfobox@
\newbox\edsbox@           \newbox\langbox@
\newif\iffirstref@        \newif\iflastref@
\newif\ifprevjour@        \newif\ifbook@            \newif\ifprevinbook@
\newif\ifquotes@          \newif\ifbookquotes@      \newif\ifpaperquotes@
\newdimen\bysamerulewd@
\setboxz@h{\refsfont@\kern3em}
\bysamerulewd@\wdz@
\newdimen\refindentwd
\setboxz@h{\refsfont@ 00. }
\refindentwd\wdz@
\outer\def\ref{\begingroup \noindent\hangindent\refindentwd
 \firstref@true \def\nofrills{\def\refkern@{\kern3sp}}%
 \ref@}
\def\ref@{\book@false \bgroup\let\endrefitem@\egroup \ignorespaces}
\def\moreref{\endrefitem@\endref@\firstref@false\ref@}%
\def\transl{\endrefitem@\endref@\firstref@false
  \book@false
  \prepunct@
  \setboxz@h\bgroup \aftergroup\unhbox\aftergroup\z@
    \def\endrefitem@{\unskip\refkern@\egroup}\ignorespaces}%
\def\emptyifempty@{\dimen@\wd\currbox@
  \advance\dimen@-\wd\z@ \advance\dimen@-.1\p@
  \ifdim\dimen@<\z@ \setbox\currbox@\copy\voidb@x \fi}
\let\refkern@\relax
\def\endrefitem@{\unskip\refkern@\egroup
  \setboxz@h{\refkern@}\emptyifempty@}\ignorespaces
\def\refdef@#1#2#3{\edef\next@{\noexpand\endrefitem@
  \let\noexpand\currbox@\csname\expandafter\eat@\string#1box@\endcsname
    \noexpand\setbox\noexpand\currbox@\hbox\bgroup}%
  \toks@\expandafter{\next@}%
  \if\notempty{#2#3}\toks@\expandafter{\the\toks@
  \def\endrefitem@{\unskip#3\refkern@\egroup
  \setboxz@h{#2#3\refkern@}\emptyifempty@}#2}\fi
  \toks@\expandafter{\the\toks@\ignorespaces}%
  \edef#1{\the\toks@}}
\refdef@\no{}{. }
\refdef@\key{[\m@th}{] }
\refdef@\by{}{}
\def\bysame{\by\hbox to\bysamerulewd@{\hrulefill}\thinspace
   \kern0sp}
\def\manyby{\message{\string\manyby is no longer necessary; \string\by
  can be used instead, starting with version 2.0 of \styname.STY}\by}
\refdef@\paper{\ifpaperquotes@``\fi\it}{}
\refdef@\paperinfo{}{}
\def\jour{\endrefitem@\let\currbox@\jourbox@
  \setbox\currbox@\hbox\bgroup
  \def\endrefitem@{\unskip\refkern@\egroup
    \setboxz@h{\refkern@}\emptyifempty@
    \ifvoid\jourbox@\else\prevjour@true\fi}%
\ignorespaces}
\refdef@\vol{\ifbook@\else\bf\fi}{}
\refdef@\issue{no. }{}
\refdef@\yr{}{}
\refdef@\pages{}{}
\def\page{\endrefitem@\def\pp@{\def\pp@{pp.~}p.~}\let\currbox@\pagesbox@
  \setbox\currbox@\hbox\bgroup\ignorespaces}
\def\pp@{pp.~}
\def\book{\endrefitem@ \let\currbox@\bookbox@
 \setbox\currbox@\hbox\bgroup\def\endrefitem@{\unskip\refkern@\egroup
  \setboxz@h{\ifbookquotes@``\fi}\emptyifempty@
  \ifvoid\bookbox@\else\book@true\fi}%
  \ifbookquotes@``\fi\it\ignorespaces}
\def\inbook{\endrefitem@
  \let\currbox@\bookbox@\setbox\currbox@\hbox\bgroup
  \def\endrefitem@{\unskip\refkern@\egroup
  \setboxz@h{\ifbookquotes@``\fi}\emptyifempty@
  \ifvoid\bookbox@\else\book@true\previnbook@true\fi}%
  \ifbookquotes@``\fi\ignorespaces}
\refdef@\eds{(}{, eds.)}
\def\ed{\endrefitem@\let\currbox@\edsbox@
 \setbox\currbox@\hbox\bgroup
 \def\endrefitem@{\unskip, ed.)\refkern@\egroup
  \setboxz@h{(, ed.)}\emptyifempty@}(\ignorespaces}
\refdef@\bookinfo{}{}
\refdef@\publ{}{}
\refdef@\publaddr{}{}
\refdef@\finalinfo{}{}
\refdef@\lang{(}{)}

\let\refdef@\relax 
\def\ppunbox@#1{\ifvoid#1\else\prepunct@\unhbox#1\fi}
\def\nocomma@#1{\ifvoid#1\else\changepunct@3\prepunct@\unhbox#1\fi}
\def\changepunct@#1{\ifnum\lastkern<3 \unkern\kern#1sp\fi}
\def\prepunct@{\count@\lastkern\unkern
  \ifnum\lastpenalty=0
    \let\penalty@\relax
  \else
    \edef\penalty@{\penalty\the\lastpenalty\relax}%
  \fi
  \unpenalty
  \let\refspace@\ \ifcase\count@,
\or;\or.\or 
  \or\let\refspace@\relax
  \else,\fi
  \ifquotes@''\quotes@false\fi \penalty@ \refspace@
}
\def\transferpenalty@#1{\dimen@\lastkern\unkern
  \ifnum\lastpenalty=0\unpenalty\let\penalty@\relax
  \else\edef\penalty@{\penalty\the\lastpenalty\relax}\unpenalty\fi
  #1\penalty@\kern\dimen@}
\def\endref{\endrefitem@\lastref@true\endref@
  \par\endgroup \prevjour@false \previnbook@false }
\def\endref@{%
\iffirstref@
  \ifvoid\nobox@\ifvoid\keybox@\indent\fi
  \else\hbox to\refindentwd{\hss\unhbox\nobox@}\fi
  \ifvoid\keybox@
  \else\ifdim\wd\keybox@>\refindentwd
         \box\keybox@
       \else\hbox to\refindentwd{\unhbox\keybox@\hfil}\fi\fi
  \kern4sp\ppunbox@\bybox@
\fi 
  \ifvoid\paperbox@
  \else\prepunct@\unhbox\paperbox@
    \ifpaperquotes@\quotes@true\fi\fi
  \ppunbox@\paperinfobox@
  \ifvoid\jourbox@
    \ifprevjour@ \nocomma@\volbox@
      \nocomma@\issuebox@
      \ifvoid\yrbox@\else\changepunct@3\prepunct@(\unhbox\yrbox@
        \transferpenalty@)\fi
      \ppunbox@\pagesbox@
    \fi 
  \else \prepunct@\unhbox\jourbox@
    \nocomma@\volbox@
    \nocomma@\issuebox@
    \ifvoid\yrbox@\else\changepunct@3\prepunct@(\unhbox\yrbox@
      \transferpenalty@)\fi
    \ppunbox@\pagesbox@
  \fi 
  \ifbook@\prepunct@\unhbox\bookbox@ \ifbookquotes@\quotes@true\fi \fi
  \nocomma@\edsbox@
  \ppunbox@\bookinfobox@
  \ifbook@\ifvoid\volbox@\else\prepunct@ vol.~\unhbox\volbox@
  \fi\fi
  \ppunbox@\publbox@ \ppunbox@\publaddrbox@
  \ifbook@ \ppunbox@\yrbox@
    \ifvoid\pagesbox@
    \else\prepunct@\pp@\unhbox\pagesbox@\fi
  \else
    \ifprevinbook@ \ppunbox@\yrbox@
      \ifvoid\pagesbox@\else\prepunct@\pp@\unhbox\pagesbox@\fi
    \fi \fi
  \ppunbox@\finalinfobox@
  \iflastref@
    \ifvoid\langbox@.\ifquotes@''\fi
    \else\changepunct@2\prepunct@\unhbox\langbox@\fi
  \else
    \ifvoid\langbox@\changepunct@1%
    \else\changepunct@3\prepunct@\unhbox\langbox@
      \changepunct@1\fi
  \fi
}
\outer\def\enddocument{%
 \runaway@{proclaim}%
\ifmonograph@ 
\else
 \nobreak
 \thetranslator@
 \count@\z@ \loop\ifnum\count@<\addresscount@\advance\count@\@ne
 \csname address\number\count@\endcsname
 \csname email\number\count@\endcsname
 \repeat
\fi
 \vfill\supereject\end}

\def\headfont@{\headfonts}
\def\proclaimheadfont@{\bf}
\def\specialheadfont@{\bf}
\def\subheadfont@{\bf}
\def\demoheadfont@{\smc}

\newif\ifThisToToc \ThisToTocfalse
\newif\iftocloaded \tocloadedfalse

\def\C@L{\noexpand\Cal}\def\B@B{\noexpand\Bbb}\def\fR@K{\noexpand\frak}
\def\S@{\noexpand\S}\def\P@P{\noexpand\"}
\def\xpar{\\}

\def\writetoc#1{\iftocloaded\ifThisToToc\begingroup\def\totoc{}
  \def\Cal{\noexpand\C@L}\def\Bbb{\noexpand\B@B}
  \def\frak{\noexpand\fR@K}\def\goth{\frak}\def\S{\noexpand\S@}
  \def\"{\noexpand\P@P}
  \def\xpar{\par\penalty100000 }\def\idx##1{##1}\def\\{\xpar}
  \edef\next@{\write\toc{\noindent#1\leaderfill\noexpand\folio\par}}%
  \next@\endgroup\global\ThisToTocfalse\fi\fi}
\def\leaderfill{\leaders\hbox to 1em{\hss.\hss}\hfill}

\newif\ifindexloaded \indexloadedfalse
\def\idx#1{\ifindexloaded\begingroup\def\ign{}\def\it{}\def\/{}%
 \def\smc{}\def\bf{}\def\tt{}%
 \def\Cal{\noexpand\C@L}\def\Bbb{\noexpand\B@B}%
 \def\frak{\noexpand\fR@K}\def\goth{\frak}\def\S{\noexpand\S@}%
  \def\"{\noexpand\P@P}%
 {\edef\next@{\write\index{#1, \noexpand\folio}}\next@}%
 \endgroup\fi{#1}}
\def\ign#1{}

\def\totoc{\global\ThisToToctrue}

\outer\def\head#1\endhead{\par\penaltyandskip@{-200}\aboveheadskip
 {\headfont@\raggedcenter@\interlinepenalty\@M
 \ignorespaces#1\endgraf}\nobreak
 \vskip\belowheadskip
 \headmark{#1}\writetoc{#1}}

\outer\def\chaphead#1\endchaphead{\par\penaltyandskip@{-200}\aboveheadskip
 {\chapheadfonts\raggedcenter@\interlinepenalty\@M
 \ignorespaces#1\endgraf}\nobreak
 \vskip3\belowheadskip
 \headmark{#1}\writetoc{#1}}

\def\folio{{\foliofont@\ifnum\pageno<\z@ \romannumeral-\pageno
 \else\number\pageno \fi}}
\newtoks\leftheadtoks
\newtoks\rightheadtoks

\def\leftheadtext{\nofrills@{\relax}\lht@
  \DNii@##1{\leftheadtoks\expandafter{\lht@{##1}}%
    \mark{\the\leftheadtoks\noexpand\else\the\rightheadtoks}
    \ifsyntax@\setboxz@h{\def\\{\unskip\space\ignorespaces}%
        \headlinefont@##1}\fi}%
  \FN@\next@}
\def\rightheadtext{\nofrills@{\relax}\rht@
  \DNii@##1{\rightheadtoks\expandafter{\rht@{##1}}%
    \mark{\the\leftheadtoks\noexpand\else\the\rightheadtoks}%
    \ifsyntax@\setboxz@h{\def\\{\unskip\space\ignorespaces}%
        \headlinefont@##1}\fi}%
  \FN@\next@}
\def\NoRunningHeads{\global\runheads@false\global\let\headmark\eat@}

\newif\iffirstpage@     \firstpage@true
\newif\ifrunheads@      \runheads@true

\newdimen\fullhsize \fullhsize=\hsize
\newdimen\fullvsize \fullvsize=\vsize
\def\fullline{\hbox to\fullhsize}

\def\pagenumbers{\gdef\folio{\folio@}}

\let\norunningheads\NoRunningHeads
\def\userunningheads{\global\runheads@true}
\norunningheads

\headline={\def\chapter#1{\chapterno@. }%
  \def\\{\unskip\space\ignorespaces}\ifrunheads@\headlinefont@
    \ifodd\pageno\rightheadline \else\leftheadline\fi
   \else\hfil\fi\ifNoRunHeadline\global\NoRunHeadlinefalse\fi}
\let\folio@\folio
\def\foliofont@{\foliofont}
\def\foliofont{\eightrm}
\def\headlinefont@{\headlinefont}
\def\headlinefont{\eightpoint\smc}
\def\leftheadline{\rlap{\folio}\hfill
   \ifNoRunHeadline\else\iftrue\topmark\fi\fi \hfill}
\def\rightheadline{\hfill\ifNoRunHeadline
   \else \expandafter\fi
  \hfill \llap{\folio}}
\footline={{\eightpoint\bottremark}%
   \ifrunheads@\else\hfil{\let\foliofont\tenrm\folio}\fi\hfil}
\def\bottremark{}
 
\newif\ifNoRunHeadline      
\def\norunninghead{\global\NoRunHeadlinetrue}
\norunninghead

\output={\output@}
%
\newif\ifoffset\offsetfalse
\output={\output@}
\def\output@{%
 \ifoffset 
  \ifodd\count0\advance\hoffset by0.5truecm
   \else\advance\hoffset by-0.5truecm\fi\fi
 \shipout\vbox{%
  \makeheadline \pagebody \makefootline }%
 \advancepageno \ifnum\outputpenalty>-\@MM\else\dosupereject\fi}

\def\indexoutput#1{%
  \ifoffset 
   \ifodd\count0\advance\hoffset by0.5truecm
    \else\advance\hoffset by-0.5truecm\fi\fi
  \shipout\vbox{\makeheadline
  \vbox to\fullvsize{\boxmaxdepth\maxdepth%
  \ifvoid\topins\else\unvbox\topins\fi%
  #1 %
  \ifvoid\footins\else 
    \vskip\skip\footins
    \footnoterule
    \unvbox\footins\fi
  \ifr@ggedbottom \kern-\dimen@ \vfil \fi}%
  \baselineskip2pc
  \makefootline}%
 \global\advance\pageno\@ne
 \ifnum\outputpenalty>-\@MM\else\dosupereject\fi}
 
 \newbox\partialpage \newdimen\halfsize \halfsize=0.5\fullhsize
 \advance\halfsize by-0.5em

 \def\begindoublecolumns{\output={\indexoutput{\unvbox255}}%
   \begingroup \def\line{\fullline}
   \output={\global\setbox\partialpage=\vbox{\unvbox255\bigskip}}\eject
   \output={\doublecolumnout}\hsize=\halfsize \vsize=2\fullvsize}
 \def\enddoublecolumns{\output={\balancecolumns}\eject
  \endgroup \pagegoal=\fullvsize%
  \output={\output@}}
\def\doublecolumnout{\splittopskip=\topskip \splitmaxdepth=\maxdepth
  \dimen@=\fullvsize \advance\dimen@ by-\ht\partialpage
  \setbox0=\vsplit255 to \dimen@ \setbox2=\vsplit255 to \dimen@
  \indexoutput{\pagesofar} \unvbox255 \penalty\outputpenalty}
\def\pagesofar{\unvbox\partialpage
  \wd0=\hsize \wd2=\hsize \hbox to\fullhsize{\box0\hfil\box2}}
\def\balancecolumns{\setbox0=\vbox{\unvbox255} \dimen@=\ht0
  \advance\dimen@ by\topskip \advance\dimen@ by-\baselineskip
  \divide\dimen@ by2 \splittopskip=\topskip
  {\vbadness=10000 \loop \global\setbox3=\copy0
    \global\setbox1=\vsplit3 to\dimen@
    \ifdim\ht3>\dimen@ \global\advance\dimen@ by1pt \repeat}
  \setbox0=\vbox to\dimen@{\unvbox1} \setbox2=\vbox to\dimen@{\unvbox3}
  \pagesofar}

\tenpoint
\catcode`\@=\active

\def\smallheadings{\let\chapheadfonts\tenpoint\let\headfonts\tenpoint}

\tenpoint
\catcode`\@=\active

\def\LL{\leavevmode\setbox0=\hbox{L}\hbox to\wd0{\hss\char'40L}}
\def\al{\alpha}
\def\be{\beta}

\def\th{\theta}

\def\la{\lambda}

\def\si{\sigma}

\def\ph{\varphi}

\def\om{\omega}


\def\FM{{\Cal F\Cal M}}

\def\V{{\Cal V}}

\def\today{\ifcase\month\or
 January\or February\or March\or April\or May\or June\or
 July\or August\or September\or October\or November\or December\fi
 \space\number\day, \number\year}

\def\({\left(}
\def\){\right)}
\def\[{\left[}
\def\]{\right]}

\def\sgn{\operatorname{sgn}}

\def\3{\ss}
\catcode`\@=11
\def\dddot#1{\vbox{\ialign{##\crcr
      .\hskip-.5pt.\hskip-.5pt.\crcr\noalign{\kern1.5\p@\nointerlineskip}
      $\hfil\displaystyle{#1}\hfil$\crcr}}}

\newif\iftab@\tab@false
\newif\ifvtab@\vtab@false
\def\tab{\bgroup\tab@true\vtab@false\vst@bfalse\Strich@false%
   \def\\{\global\hline@@false%
     \ifhline@\global\hline@false\global\hline@@true\fi\cr}
   \edef\l@{\the\leftskip}\ialign\bgroup\hskip\l@##\hfil&&##\hfil\cr}
\def\endtab{\cr\egroup\egroup}
\def\vtab{\vtop\bgroup\vst@bfalse\vtab@true\tab@true\Strich@false%
   \bgroup\def\\{\cr}\ialign\bgroup&##\hfil\cr}
\def\endvtab{\cr\egroup\egroup\egroup}
\def\stab{\D@cke0.5pt\null 
 \bgroup\tab@true\vtab@false\vst@bfalse\Strich@true\Let@@\vspace@
 \normalbaselines\offinterlineskip
  \openup\spreadmlines@
 \edef\l@{\the\leftskip}\ialign
 \bgroup\hskip\l@##\hfil&&##\hfil\crcr}
\def\endstab{\crcr\egroup
 \egroup}
\newif\ifvst@b\vst@bfalse
\def\vstab{\D@cke0.5pt\null
 \vtop\bgroup\tab@true\vtab@false\vst@btrue\Strich@true\bgroup\Let@@\vspace@
 \normalbaselines\offinterlineskip
  \openup\spreadmlines@\bgroup}
\def\endvstab{\crcr\egroup\egroup
 \egroup\tab@false\Strich@false}

\newdimen\htstrut@
\htstrut@8.5\p@
\newdimen\htStrut@
\htStrut@12\p@
\newdimen\dpstrut@
\dpstrut@3.5\p@
\newdimen\dpStrut@
\dpStrut@3.5\p@
\def\openup{\afterassignment\@penup\dimen@=}
\def\@penup{\advance\lineskip\dimen@
  \advance\baselineskip\dimen@
  \advance\lineskiplimit\dimen@
  \divide\dimen@ by2
  \advance\htstrut@\dimen@
  \advance\htStrut@\dimen@
  \advance\dpstrut@\dimen@
  \advance\dpStrut@\dimen@}
\def\Let@@{\relax%
    \def\\{\global\hline@@false%
     \ifhline@\global\hline@false\global\hline@@true\fi\cr}%
    \iffalse}\fi}
\def\matrix{\null\,\vcenter\bgroup
 \tab@false\vtab@false\vst@bfalse\Strich@false\Let@@\vspace@
 \normalbaselines\openup\spreadmlines@\ialign
 \bgroup\hfil$\m@th##$\hfil&&\quad\hfil$\m@th##$\hfil\crcr
 \Mathstrut@\crcr\noalign{\kern-\baselineskip}}
\def\endmatrix{\crcr\Mathstrut@\crcr\noalign{\kern-\baselineskip}\egroup
 \egroup\,}
\def\smatrix{\D@cke0.5pt\null\,
 \vcenter\bgroup\tab@false\vtab@false\vst@bfalse\Strich@true\Let@@\vspace@
 \normalbaselines\offinterlineskip
  \openup\spreadmlines@\ialign
 \bgroup\hfil$\m@th##$\hfil&&\quad\hfil$\m@th##$\hfil\crcr}
\def\endsmatrix{\crcr\egroup
 \egroup\,\Strich@false}
\newdimen\D@cke
\def\Dicke#1{\global\D@cke#1}
\newtoks\tabs@\tabs@{&}
\newif\ifStrich@\Strich@false
\newif\iff@rst

\def\Stricherr@{\iftab@\ifvtab@\errmessage{\noexpand\s not allowed
     here. Use \noexpand\vstab!}%
  \else\errmessage{\noexpand\s not allowed here. Use \noexpand\stab!}%
  \fi\else\errmessage{\noexpand\s not allowed
     here. Use \noexpand\smatrix!}\fi}
\def\format{\ifvst@b\else\crcr\fi\egroup\iffalse{\fi\ifnum`}=0 \fi\format@}
\def\format@#1\\{\def\preamble@{#1}%
 \def\Str@chfehlt##1{\ifx##1\s\Stricherr@\fi\ifx##1\\\let\Next\relax%
   \else\let\Next\Str@chfehlt\fi\Next}%
 \def\c{\hfil\noexpand\ifhline@@\hbox{\vrule height\htStrut@%
   depth\dpstrut@ width\z@}\noexpand\fi%
   \ifStrich@\hbox{\vrule height\htstrut@ depth\dpstrut@ width\z@}%
   \fi\iftab@\else$\m@th\fi\the\hashtoks@\iftab@\else$\fi\hfil}%
 \def\r{\hfil\noexpand\ifhline@@\hbox{\vrule height\htStrut@%
   depth\dpstrut@ width\z@}\noexpand\fi%
   \ifStrich@\hbox{\vrule height\htstrut@ depth\dpstrut@ width\z@}%
   \fi\iftab@\else$\m@th\fi\the\hashtoks@\iftab@\else$\fi}%
 \def\l{\noexpand\ifhline@@\hbox{\vrule height\htStrut@%
   depth\dpstrut@ width\z@}\noexpand\fi%
   \ifStrich@\hbox{\vrule height\htstrut@ depth\dpstrut@ width\z@}%
   \fi\iftab@\else$\m@th\fi\the\hashtoks@\iftab@\else$\fi\hfil}%
 \def\s{\ifStrich@\ \the\tabs@\vrule width\D@cke\the\hashtoks@%
          \fi\the\tabs@\ }%
 \def\sa{\ifStrich@\vrule width\D@cke\the\hashtoks@%
            \the\tabs@\ %
            \fi}%
 \def\se{\ifStrich@\ \the\tabs@\vrule width\D@cke\the\hashtoks@\fi}%
 \def\cd{\hfil\noexpand\ifhline@@\hbox{\vrule height\htStrut@%
   depth\dpstrut@ width\z@}\noexpand\fi%
   \ifStrich@\hbox{\vrule height\htstrut@ depth\dpstrut@ width\z@}%
   \fi$\dsize\m@th\the\hashtoks@$\hfil}%
 \def\rd{\hfil\noexpand\ifhline@@\hbox{\vrule height\htStrut@%
   depth\dpstrut@ width\z@}\noexpand\fi%
   \ifStrich@\hbox{\vrule height\htstrut@ depth\dpstrut@ width\z@}%
   \fi$\dsize\m@th\the\hashtoks@$}%
 \def\ld{\noexpand\ifhline@@\hbox{\vrule height\htStrut@%
   depth\dpstrut@ width\z@}\noexpand\fi%
   \ifStrich@\hbox{\vrule height\htstrut@ depth\dpstrut@ width\z@}%
   \fi$\dsize\m@th\the\hashtoks@$\hfil}%
 \ifStrich@\else\Str@chfehlt#1\\\fi%
 \setbox\z@\hbox{\xdef\Preamble@{\preamble@}}\ifnum`{=0 \fi\iffalse}\fi
 \ialign\bgroup\span\Preamble@\crcr}
\newif\ifhline@\hline@false
\newif\ifhline@@\hline@@false
\def\hlinefor#1{\multispan@{\strip@#1 }\leaders\hrule height\D@cke\hfill%
    \global\hline@true\ignorespaces}
\def\Item "#1"{\par\noindent\hangindent2\parindent%
  \hangafter1\setbox0\hbox{\rm#1\enspace}\ifdim\wd0>2\parindent%
  \box0\else\hbox to 2\parindent{\rm#1\hfil}\fi\ignorespaces}
\def\ITEM #1"#2"{\par\noindent\hangafter1\hangindent#1%
  \setbox0\hbox{\rm#2\enspace}\ifdim\wd0>#1%
  \box0\else\hbox to 0pt{\rm#2\hss}\hskip#1\fi\ignorespaces}
\def\item"#1"{\par\noindent\hang%
  \setbox0=\hbox{\rm#1\enspace}\ifdim\wd0>\the\parindent%
  \box0\else\hbox to \parindent{\rm#1\hfil}\enspace\fi\ignorespaces}
\let\plainitem@\item
\catcode`\@=13

\catcode`\@=11
\font\tenln    = line10
\font\tenlnw   = linew10

\newskip\Einheit \Einheit=0.5cm
\newcount\xcoord \newcount\ycoord
\newdimen\xdim \newdimen\ydim \newdimen\PfadD@cke \newdimen\Pfadd@cke

\newcount\@tempcnta
\newcount\@tempcntb

\newdimen\@tempdima
\newdimen\@tempdimb

\newdimen\@wholewidth
\newdimen\@halfwidth

\newcount\@xarg
\newcount\@yarg
\newcount\@yyarg
\newbox\@linechar
\newbox\@tempboxa
\newdimen\@linelen
\newdimen\@clnwd
\newdimen\@clnht

\newif\if@negarg

\def\@whilenoop#1{}
\def\@whiledim#1\do #2{\ifdim #1\relax#2\@iwhiledim{#1\relax#2}\fi}
\def\@iwhiledim#1{\ifdim #1\let\@nextwhile=\@iwhiledim
        \else\let\@nextwhile=\@whilenoop\fi\@nextwhile{#1}}

\def\@whileswnoop#1\fi{}
\def\@whilesw#1\fi#2{#1#2\@iwhilesw{#1#2}\fi\fi}
\def\@iwhilesw#1\fi{#1\let\@nextwhile=\@iwhilesw
         \else\let\@nextwhile=\@whileswnoop\fi\@nextwhile{#1}\fi}

\def\thinlines{\let\@linefnt\tenln \let\@circlefnt\tencirc
  \@wholewidth\fontdimen8\tenln \@halfwidth .5\@wholewidth}
\def\thicklines{\let\@linefnt\tenlnw \let\@circlefnt\tencircw
  \@wholewidth\fontdimen8\tenlnw \@halfwidth .5\@wholewidth}
\thinlines

\PfadD@cke1pt \Pfadd@cke0.5pt
\def\PfadDicke#1{\PfadD@cke#1 \divide\PfadD@cke by2 \Pfadd@cke\PfadD@cke \multiply\PfadD@cke by2}
\long\def\LOOP#1\REPEAT{\def\BODY{#1}\ITERATE}
\def\ITERATE{\BODY \let\next\ITERATE \else\let\next\relax\fi \next}
\let\REPEAT=\fi
\def\Punkt{\hbox{\raise-2pt\hbox to0pt{\hss$\ssize\bullet$\hss}}}
\def\DuennPunkt(#1,#2){\unskip
  \raise#2 \Einheit\hbox to0pt{\hskip#1 \Einheit
          \raise-2.5pt\hbox to0pt{\hss$\bullet$\hss}\hss}}
\def\NormalPunkt(#1,#2){\unskip
  \raise#2 \Einheit\hbox to0pt{\hskip#1 \Einheit
          \raise-3pt\hbox to0pt{\hss\twelvepoint$\bullet$\hss}\hss}}
\def\DickPunkt(#1,#2){\unskip
  \raise#2 \Einheit\hbox to0pt{\hskip#1 \Einheit
          \raise-4pt\hbox to0pt{\hss\fourteenpoint$\bullet$\hss}\hss}}
\def\Kreis(#1,#2){\unskip
  \raise#2 \Einheit\hbox to0pt{\hskip#1 \Einheit
          \raise-4pt\hbox to0pt{\hss\fourteenpoint$\circ$\hss}\hss}}

\def\Line@(#1,#2)#3{\@xarg #1\relax \@yarg #2\relax
\@linelen=#3\Einheit
\ifnum\@xarg =0 \@vline
  \else \ifnum\@yarg =0 \@hline \else \@sline\fi
\fi}

\def\@sline{\ifnum\@xarg< 0 \@negargtrue \@xarg -\@xarg \@yyarg -\@yarg
  \else \@negargfalse \@yyarg \@yarg \fi
\ifnum \@yyarg >0 \@tempcnta\@yyarg \else \@tempcnta -\@yyarg \fi
\ifnum\@tempcnta>6 \@badlinearg\@tempcnta0 \fi
\ifnum\@xarg>6 \@badlinearg\@xarg 1 \fi
\setbox\@linechar\hbox{\@linefnt\@getlinechar(\@xarg,\@yyarg)}%
\ifnum \@yarg >0 \let\@upordown\raise \@clnht\z@
   \else\let\@upordown\lower \@clnht \ht\@linechar\fi
\@clnwd=\wd\@linechar
\if@negarg \hskip -\wd\@linechar \def\@tempa{\hskip -2\wd\@linechar}\else
     \let\@tempa\relax \fi
\@whiledim \@clnwd <\@linelen \do
  {\@upordown\@clnht\copy\@linechar
   \@tempa
   \advance\@clnht \ht\@linechar
   \advance\@clnwd \wd\@linechar}%
\advance\@clnht -\ht\@linechar
\advance\@clnwd -\wd\@linechar
\@tempdima\@linelen\advance\@tempdima -\@clnwd
\@tempdimb\@tempdima\advance\@tempdimb -\wd\@linechar
\if@negarg \hskip -\@tempdimb \else \hskip \@tempdimb \fi
\multiply\@tempdima \@m
\@tempcnta \@tempdima \@tempdima \wd\@linechar \divide\@tempcnta \@tempdima
\@tempdima \ht\@linechar \multiply\@tempdima \@tempcnta
\divide\@tempdima \@m
\advance\@clnht \@tempdima
\ifdim \@linelen <\wd\@linechar
   \hskip \wd\@linechar
  \else\@upordown\@clnht\copy\@linechar\fi}

\def\@hline{\ifnum \@xarg <0 \hskip -\@linelen \fi
\vrule height\Pfadd@cke width \@linelen depth\Pfadd@cke
\ifnum \@xarg <0 \hskip -\@linelen \fi}

\def\@getlinechar(#1,#2){\@tempcnta#1\relax\multiply\@tempcnta 8
\advance\@tempcnta -9 \ifnum #2>0 \advance\@tempcnta #2\relax\else
\advance\@tempcnta -#2\relax\advance\@tempcnta 64 \fi
\char\@tempcnta}

\def\Vektor(#1,#2)#3(#4,#5){\unskip\leavevmode
  \xcoord#4\relax \ycoord#5\relax
      \raise\ycoord \Einheit\hbox to0pt{\hskip\xcoord \Einheit
         \Vector@(#1,#2){#3}\hss}}

\def\Vector@(#1,#2)#3{\@xarg #1\relax \@yarg #2\relax
\@tempcnta \ifnum\@xarg<0 -\@xarg\else\@xarg\fi
\ifnum\@tempcnta<5\relax
\@linelen=#3\Einheit
\ifnum\@xarg =0 \@vvector
  \else \ifnum\@yarg =0 \@hvector \else \@svector\fi
\fi
\else\@badlinearg\fi}

\def\@hvector{\@hline\hbox to 0pt{\@linefnt
\ifnum \@xarg <0 \@getlarrow(1,0)\hss\else
    \hss\@getrarrow(1,0)\fi}}

\def\@vvector{\ifnum \@yarg <0 \@downvector \else \@upvector \fi}

\def\@svector{\@sline
\@tempcnta\@yarg \ifnum\@tempcnta <0 \@tempcnta=-\@tempcnta\fi
\ifnum\@tempcnta <5
  \hskip -\wd\@linechar
  \@upordown\@clnht \hbox{\@linefnt  \if@negarg
  \@getlarrow(\@xarg,\@yyarg) \else \@getrarrow(\@xarg,\@yyarg) \fi}%
\else\@badlinearg\fi}

\def\@upline{\hbox to \z@{\hskip -.5\Pfadd@cke \vrule width \Pfadd@cke
   height \@linelen depth \z@\hss}}

\def\@downline{\hbox to \z@{\hskip -.5\Pfadd@cke \vrule width \Pfadd@cke
   height \z@ depth \@linelen \hss}}

\def\@upvector{\@upline\setbox\@tempboxa\hbox{\@linefnt\char'66}\raise
     \@linelen \hbox to\z@{\lower \ht\@tempboxa\box\@tempboxa\hss}}

\def\@downvector{\@downline\lower \@linelen
      \hbox to \z@{\@linefnt\char'77\hss}}

\def\@getlarrow(#1,#2){\ifnum #2 =\z@ \@tempcnta='33\else
\@tempcnta=#1\relax\multiply\@tempcnta \sixt@@n \advance\@tempcnta
-9 \@tempcntb=#2\relax\multiply\@tempcntb \tw@
\ifnum \@tempcntb >0 \advance\@tempcnta \@tempcntb\relax
\else\advance\@tempcnta -\@tempcntb\advance\@tempcnta 64
\fi\fi\char\@tempcnta}

\def\@getrarrow(#1,#2){\@tempcntb=#2\relax
\ifnum\@tempcntb < 0 \@tempcntb=-\@tempcntb\relax\fi
\ifcase \@tempcntb\relax \@tempcnta='55 \or
\ifnum #1<3 \@tempcnta=#1\relax\multiply\@tempcnta
24 \advance\@tempcnta -6 \else \ifnum #1=3 \@tempcnta=49
\else\@tempcnta=58 \fi\fi\or
\ifnum #1<3 \@tempcnta=#1\relax\multiply\@tempcnta
24 \advance\@tempcnta -3 \else \@tempcnta=51\fi\or
\@tempcnta=#1\relax\multiply\@tempcnta
\sixt@@n \advance\@tempcnta -\tw@ \else
\@tempcnta=#1\relax\multiply\@tempcnta
\sixt@@n \advance\@tempcnta 7 \fi\ifnum #2<0 \advance\@tempcnta 64 \fi
\char\@tempcnta}

\def\Diagonale(#1,#2)#3{\unskip\leavevmode
  \xcoord#1\relax \ycoord#2\relax
      \raise\ycoord \Einheit\hbox to0pt{\hskip\xcoord \Einheit
         \Line@(1,1){#3}\hss}}
\def\AntiDiagonale(#1,#2)#3{\unskip\leavevmode
  \xcoord#1\relax \ycoord#2\relax 
      \raise\ycoord \Einheit\hbox to0pt{\hskip\xcoord \Einheit
         \Line@(1,-1){#3}\hss}}
\def\Pfad(#1,#2),#3\endPfad{\unskip\leavevmode
  \xcoord#1 \ycoord#2 \thicklines\ZeichnePfad#3\endPfad\thinlines}
\def\ZeichnePfad#1{\ifx#1\endPfad\let\next\relax
  \else\let\next\ZeichnePfad
    \ifnum#1=1
      \raise\ycoord \Einheit\hbox to0pt{\hskip\xcoord \Einheit
         \vrule height\Pfadd@cke width1 \Einheit depth\Pfadd@cke\hss}%
      \advance\xcoord by 1
    \else\ifnum#1=2
      \raise\ycoord \Einheit\hbox to0pt{\hskip\xcoord \Einheit
        \hbox{\hskip-\PfadD@cke\vrule height1 \Einheit width\PfadD@cke depth0pt}\hss}%
      \advance\ycoord by 1
    \else\ifnum#1=3
      \raise\ycoord \Einheit\hbox to0pt{\hskip\xcoord \Einheit
         \Line@(1,1){1}\hss}
      \advance\xcoord by 1
      \advance\ycoord by 1
    \else\ifnum#1=4
      \raise\ycoord \Einheit\hbox to0pt{\hskip\xcoord \Einheit
         \Line@(1,-1){1}\hss}
      \advance\xcoord by 1
      \advance\ycoord by -1
    \else\ifnum#1=5
      \advance\xcoord by -1
      \raise\ycoord \Einheit\hbox to0pt{\hskip\xcoord \Einheit
         \vrule height\Pfadd@cke width1 \Einheit depth\Pfadd@cke\hss}%
    \else\ifnum#1=6
      \advance\ycoord by -1
      \raise\ycoord \Einheit\hbox to0pt{\hskip\xcoord \Einheit
        \hbox{\hskip-\PfadD@cke\vrule height1 \Einheit width\PfadD@cke depth0pt}\hss}%
    \else\ifnum#1=7
      \advance\xcoord by -1
      \advance\ycoord by -1
      \raise\ycoord \Einheit\hbox to0pt{\hskip\xcoord \Einheit
         \Line@(1,1){1}\hss}
    \else\ifnum#1=8
      \advance\xcoord by -1
      \advance\ycoord by +1
      \raise\ycoord \Einheit\hbox to0pt{\hskip\xcoord \Einheit
         \Line@(1,-1){1}\hss}
    \fi\fi\fi\fi
    \fi\fi\fi\fi
  \fi\next}
\def\hSSchritt{\leavevmode\raise-.4pt\hbox to0pt{\hss.\hss}\hskip.2\Einheit
  \raise-.4pt\hbox to0pt{\hss.\hss}\hskip.2\Einheit
  \raise-.4pt\hbox to0pt{\hss.\hss}\hskip.2\Einheit
  \raise-.4pt\hbox to0pt{\hss.\hss}\hskip.2\Einheit
  \raise-.4pt\hbox to0pt{\hss.\hss}\hskip.2\Einheit}
\def\vSSchritt{\vbox{\baselineskip.2\Einheit\lineskiplimit0pt
\hbox{.}\hbox{.}\hbox{.}\hbox{.}\hbox{.}}}
\def\DSSchritt{\leavevmode\raise-.4pt\hbox to0pt{%
  \hbox to0pt{\hss.\hss}\hskip.2\Einheit
  \raise.2\Einheit\hbox to0pt{\hss.\hss}\hskip.2\Einheit
  \raise.4\Einheit\hbox to0pt{\hss.\hss}\hskip.2\Einheit
  \raise.6\Einheit\hbox to0pt{\hss.\hss}\hskip.2\Einheit
  \raise.8\Einheit\hbox to0pt{\hss.\hss}\hss}}
\def\dSSchritt{\leavevmode\raise-.4pt\hbox to0pt{%
  \hbox to0pt{\hss.\hss}\hskip.2\Einheit
  \raise-.2\Einheit\hbox to0pt{\hss.\hss}\hskip.2\Einheit
  \raise-.4\Einheit\hbox to0pt{\hss.\hss}\hskip.2\Einheit
  \raise-.6\Einheit\hbox to0pt{\hss.\hss}\hskip.2\Einheit
  \raise-.8\Einheit\hbox to0pt{\hss.\hss}\hss}}
\def\SPfad(#1,#2),#3\endSPfad{\unskip\leavevmode
  \xcoord#1 \ycoord#2 \ZeichneSPfad#3\endSPfad}
\def\ZeichneSPfad#1{\ifx#1\endSPfad\let\next\relax
  \else\let\next\ZeichneSPfad
    \ifnum#1=1
      \raise\ycoord \Einheit\hbox to0pt{\hskip\xcoord \Einheit
         \hSSchritt\hss}%
      \advance\xcoord by 1
    \else\ifnum#1=2
      \raise\ycoord \Einheit\hbox to0pt{\hskip\xcoord \Einheit
        \hbox{\hskip-2pt \vSSchritt}\hss}%
      \advance\ycoord by 1
    \else\ifnum#1=3
      \raise\ycoord \Einheit\hbox to0pt{\hskip\xcoord \Einheit
         \DSSchritt\hss}
      \advance\xcoord by 1
      \advance\ycoord by 1
    \else\ifnum#1=4
      \raise\ycoord \Einheit\hbox to0pt{\hskip\xcoord \Einheit
         \dSSchritt\hss}
      \advance\xcoord by 1
      \advance\ycoord by -1
    \else\ifnum#1=5
      \advance\xcoord by -1
      \raise\ycoord \Einheit\hbox to0pt{\hskip\xcoord \Einheit
         \hSSchritt\hss}%
    \else\ifnum#1=6
      \advance\ycoord by -1
      \raise\ycoord \Einheit\hbox to0pt{\hskip\xcoord \Einheit
        \hbox{\hskip-2pt \vSSchritt}\hss}%
    \else\ifnum#1=7
      \advance\xcoord by -1
      \advance\ycoord by -1
      \raise\ycoord \Einheit\hbox to0pt{\hskip\xcoord \Einheit
         \DSSchritt\hss}
    \else\ifnum#1=8
      \advance\xcoord by -1
      \advance\ycoord by 1
      \raise\ycoord \Einheit\hbox to0pt{\hskip\xcoord \Einheit
         \dSSchritt\hss}
    \fi\fi\fi\fi
    \fi\fi\fi\fi
  \fi\next}
\def\Koordinatenachsen(#1,#2){\unskip
 \hbox to0pt{\hskip-.5pt\vrule height#2 \Einheit width.5pt depth1 \Einheit}%
 \hbox to0pt{\hskip-1 \Einheit \xcoord#1 \advance\xcoord by1
    \vrule height0.25pt width\xcoord \Einheit depth0.25pt\hss}}
\def\Koordinatenachsen(#1,#2)(#3,#4){\unskip
 \hbox to0pt{\hskip-.5pt \ycoord-#4 \advance\ycoord by1
    \vrule height#2 \Einheit width.5pt depth\ycoord \Einheit}%
 \hbox to0pt{\hskip-1 \Einheit \hskip#3\Einheit 
    \xcoord#1 \advance\xcoord by1 \advance\xcoord by-#3 
    \vrule height0.25pt width\xcoord \Einheit depth0.25pt\hss}}
\def\Gitter(#1,#2){\unskip \xcoord0 \ycoord0 \leavevmode
  \LOOP\ifnum\ycoord<#2
    \loop\ifnum\xcoord<#1
      \raise\ycoord \Einheit\hbox to0pt{\hskip\xcoord \Einheit\Punkt\hss}%
      \advance\xcoord by1
    \repeat
    \xcoord0
    \advance\ycoord by1
  \REPEAT}
\def\Gitter(#1,#2)(#3,#4){\unskip \xcoord#3 \ycoord#4 \leavevmode
  \LOOP\ifnum\ycoord<#2
    \loop\ifnum\xcoord<#1
      \raise\ycoord \Einheit\hbox to0pt{\hskip\xcoord \Einheit\Punkt\hss}%
      \advance\xcoord by1
    \repeat
    \xcoord#3
    \advance\ycoord by1
  \REPEAT}
\def\Label#1#2(#3,#4){\unskip \xdim#3 \Einheit \ydim#4 \Einheit
  \def\lo{\advance\xdim by-.5 \Einheit \advance\ydim by.5 \Einheit}%
  \def\llo{\advance\xdim by-.25cm \advance\ydim by.5 \Einheit}%
  \def\loo{\advance\xdim by-.5 \Einheit \advance\ydim by.25cm}%
  \def\o{\advance\ydim by.25cm}%
  \def\ro{\advance\xdim by.5 \Einheit \advance\ydim by.5 \Einheit}%
  \def\rro{\advance\xdim by.25cm \advance\ydim by.5 \Einheit}%
  \def\roo{\advance\xdim by.5 \Einheit \advance\ydim by.25cm}%
  \def\l{\advance\xdim by-.30cm}%
  \def\r{\advance\xdim by.30cm}%
  \def\lu{\advance\xdim by-.5 \Einheit \advance\ydim by-.6 \Einheit}%
  \def\llu{\advance\xdim by-.25cm \advance\ydim by-.6 \Einheit}%
  \def\luu{\advance\xdim by-.5 \Einheit \advance\ydim by-.30cm}%
  \def\u{\advance\ydim by-.30cm}%
  \def\ru{\advance\xdim by.5 \Einheit \advance\ydim by-.6 \Einheit}%
  \def\rru{\advance\xdim by.25cm \advance\ydim by-.6 \Einheit}%
  \def\ruu{\advance\xdim by.5 \Einheit \advance\ydim by-.30cm}%
  #1\raise\ydim\hbox to0pt{\hskip\xdim
     \vbox to0pt{\vss\hbox to0pt{\hss$#2$\hss}\vss}\hss}%
}
\catcode`\@=13

\hsize13cm
\vsize19cm
\newdimen\fullhsize
\newdimen\fullvsize
\newdimen\halfsize
\fullhsize13cm
\fullvsize19cm
\halfsize=0.5\fullhsize
\advance\halfsize by-0.5em

\magnification1200

\TagsOnRight
\loadbold

\catcode`\@=11
\def\iddots{\mathinner{\mkern1mu\raise\p@\hbox{.}\mkern2mu
    \raise4\p@\hbox{.}\mkern2mu\raise7\p@\vbox{\kern7\p@\hbox{.}}\mkern1mu}}
\catcode`\@=13

\catcode`\@=11
\font@\sixrm=cmr10 scaled 700
\font@\sixit=cmti10 scaled700
\font@\sixsl=cmsl10 scaled700
\font@\sixsmc=cmcsc10 scaled700
\font@\sixtt=cmtt10 scaled700
\font@\sixbf=cmbx10 scaled700
\font@\sixi=cmmi10 scaled700
\font@\sixsy=cmsy10 scaled700
\font@\sixex=cmex10 scaled700
\font@\sixmsa=msam10 scaled700
\font@\sixeufm=eufm10 scaled700
\font@\sixmsb=msbm10 scaled700
\font@\fiverm=cmr10 scaled 500
\font@\fiveit=cmti10 scaled500
\font@\fivesl=cmsl10 scaled500
\font@\fivesmc=cmcsc10 scaled500
\font@\fivett=cmtt10 scaled500
\font@\fivebf=cmbx10 scaled500
\font@\fivei=cmmi10 scaled500
\font@\fivesy=cmsy10 scaled500
\font@\fiveex=cmex10 scaled500
\font@\fivemsa=msam10 scaled500
\font@\fiveeufm=eufm10 scaled500
\font@\fivemsb=msbm10 scaled500
\newtoks\sixpoint@
\def\sixpoint{\normalbaselineskip10\p@
 \abovedisplayskip10\p@ plus2.4\p@ minus7.2\p@
 \belowdisplayskip\abovedisplayskip
 \abovedisplayshortskip\z@ plus2.4\p@
 \belowdisplayshortskip5.6\p@ plus2.4\p@ minus3.2\p@
 \textonlyfont@\rm\sixrm \textonlyfont@\it\sixit
 \textonlyfont@\sl\sixsl \textonlyfont@\bf\sixbf
 \textonlyfont@\smc\sixsmc \textonlyfont@\tt\sixtt
 \textonlyfont@\bsmc\sixbsmc
 \ifsyntax@\def\big##1{{\hbox{$\left##1\right.$}}}%
  \let\Big\big \let\bigg\big \let\Bigg\big
 \else
  \textfont\z@=\sixrm \scriptfont\z@=\fiverm \scriptscriptfont\z@=\fiverm
  \textfont\@ne=\sixi \scriptfont\@ne=\fivei \scriptscriptfont\@ne=\fivei
  \textfont\tw@=\sixsy \scriptfont\tw@=\fivesy \scriptscriptfont\tw@=\fivesy
  \textfont\thr@@=\sixex \scriptfont\thr@@=\fiveex
   \scriptscriptfont\thr@@=\fiveex
  \textfont\itfam=\sixit \scriptfont\itfam=\fiveit
   \scriptscriptfont\itfam=\fiveit
  \textfont\bffam=\sixbf \scriptfont\bffam=\fivebf
   \scriptscriptfont\bffam=\fivebf
 \setbox\strutbox\hbox{\vrule height7\p@ depth3\p@ width\z@}%
 \setbox\strutbox@\hbox{\raise.5\normallineskiplimit\vbox{%
   \kern-\normallineskiplimit\copy\strutbox}}%
 \setbox\z@\vbox{\hbox{$($}\kern\z@}\bigsize@=1.2\ht\z@
 \fi
 \normalbaselines\sixrm\ex@.2326ex\jot3\ex@\the\sixpoint@}
\catcode`\@=13

\def\AnBBAA{1}
\def\BaRoAA{2}
\def\BeSaAA{3}
\def\BoViAA{4}
\def\CaFoAA{5}
\def\CiglBD{6}
\def\FavaAA{7}
\def\FereAA{8}
\def\FoStAA{9}
\def\GaEgAA{10}
\def\GaRaAF{11}
\def\GeViAA{12}
\def\GeViAB{13}
\def\GoHaAA{14}
\def\HoZaAA{15}
\def\JaKKAA{16}
\def\JinYAA{17}
\def\KratAH{18}
\def\KratBP{19}
\def\KratAZ{20}
\def\KratBT{21}
\def\KratCL{22}
\def\KratCM{23}
\def\KrYaAA{24}
\def\LaPrAC{25}
\def\LindAA{26}
\def\LoebAA{27}
\def\OwPrAA{28}
\def\PropXA{29}
\def\SpeyXA{30}
\def\StanAY{31}
\def\StanAZ{32}
\def\StanXA{33}
\def\StanBI{34}
\def\StanAP{35}
\def\VienAE{36}
\def\VienAF{37}
\def\VienZZ{38}
\def\ZaimAA{39}

\def\BAc{1.1}
\def\BAb{2.1}
\def\BD{2.2}
\def\BGa{2.3}
\def\BGb{2.4}
\def\BGc{2.5}
\def\BG{3.1}
\def\BH{3.2}
\def\BI{3.3}
\def\BIb{3.4}
\def\BId{3.5}
\def\BIa{3.6}
\def\BIc{3.7}
\def\BIe{3.8}
\def\BEa{3.9}
\def\BFa{3.10}
\def\BA{3.11}
\def\BAa{3.12}
\def\BDb{4.1}
\def\BDa{4.2}
\def\AA{5.1}
\def\AB{5.2}
\def\AC{5.3}
\def\AD{5.4}
\def\AE{5.5}
\def\AEa{5.6}
\def\AEb{5.7}
\def\AF{5.8}
\def\AFb{5.9}
\def\CA{6.1}
\def\CB{6.2}
\def\CC{6.3}
\def\CD{6.4}
\def\CE{6.5}
\def\CEa{6.6}
\def\CF{6.7}

\def\DA{7.1}
\def\DB{7.2}
\def\DBa{7.3}
\def\DC{7.4}
\def\DD{7.5}
\def\DE{7.6}
\def\DF{7.7}
\def\EA{8.1}
\def\EAa{8.2}
\def\EAb{8.3}
\def\EAc{8.4}
\def\EB{8.5}
\def\EBa{8.6}
\def\EC{8.7}
\def\ECa{8.8}
\def\ED{8.9}
\def\EE{8.10}
\def\EF{8.11}
\def\EG{8.12}
\def\EH{8.13}
\def\EHa{8.14}
\def\EI{8.15}
\def\EJ{8.16}
\def\EJa{8.17}
\def\EJb{8.18}
\def\EJc{8.19}

\def\EK{9.1}
\def\EL{9.2}
\def\EM{9.3}
\def\EMa{9.4}
\def\EBb{9.5}
\def\EN{9.6}
\def\EO{9.7}
\def\EP{9.8}
\def\EQ{10.1}
\def\EQa{10.2}
\def\ER{10.3}
\def\ES{10.4}
\def\ESa{10.5}
\def\ET{10.6}
\def\EU{10.7}
\def\EV{10.8}
\def\EVa{10.9}
\def\EVb{10.10}
\def\EVc{10.11}
\def\EW{10.12}
\def\EX{10.13}
\def\EXa{10.14}
\def\EY{10.15}
\def\EZ{10.16}
\def\EZa{10.17}
\def\EZb{10.18}
\def\EZc{10.19}
\def\GAa{11.1}
\def\GA{11.2}
\def\GB{11.3}
\def\GC{11.4}
\def\GD{11.5}

\def\GF{11.7}

\def\HA{12.1}
\def\HB{12.2}
\def\HBa{12.3}
\def\HBb{12.4}
\def\HBd{12.5}
\def\HBf{12.6}
\def\HBg{12.7}
\def\HBh{12.8}
\def\HBi{12.9}
\def\HBj{12.10}
\def\HBk{12.11}
\def\HJ{12.12}
\def\HK{12.13}
\def\HL{12.14}
\def\HM{12.15}
\def\HN{12.16}
\def\HC{12.17}
\def\HD{12.18}
\def\HE{12.19}
\def\HEd{12.20}
\def\HEa{12.21}
\def\HEb{12.22}
\def\HEc{12.23}
\def\HF{12.24}
\def\HG{12.25}
\def\HH{12.26}
\def\HI{12.27}
\def\HBc{12.28}
\def\HBe{12.29}

\def\VA{1}
\def\VB{2}
\def\VC{3}
\def\UA{4}
\def\UB{5}
\def\UC{6}
\def\UD{7}
\def\UE{8}
\def\UG{9}
\def\UH{10}
\def\TA{11}
\def\UF{12}
\def\TB{13}
\def\TC{14}
\def\TD{15}
\def\TE{16}
\def\TF{17}
\def\TG{18}
\def\TH{19}
\def\TI{20}
\def\TJ{21}
\def\TK{22}
\def\TL{23}
\def\TM{24}
\def\TN{25}
\def\TO{26}
\def\TP{27}
\def\TQ{28}
\def\TR{29}
\def\TS{30}
\def\TT{31}
\def\TU{32}
\def\TV{33}
\def\TW{34}
\def\TX{35}
\def\TY{36}
\def\TZ{37}
\def\TTA{38}
\def\TTB{39}
\def\TTTTTT{40}
\def\TTC{41}
\def\TTD{42}
\def\TTE{43}
\def\TTF{44}
\def\TTG{45}
\def\TTH{46}
\def\TTI{47}
\def\TTJ{48}
\def\TTK{49}
\def\TTL{50}
\def\TTLa{51}
\def\TTLb{52}
\def\TTTT{53}
\def\TTTTT{54}
\def\TTM{55}
\def\TTN{56}
\def\TTNa{57}
\def\TTO{58}
\def\TTP{59}
\def\TTQ{60}
\def\TTR{61}
\def\TTS{62}
\def\TTT{63}
\def\TTU{64}

\def\FG{1}
\def\FA{2}
\def\FB{3}
\def\FC{4}
\def\FD{5}
\def\FDa{6}
\def\FDb{7}
\def\FE{8}
\def\FDc{9}
\def\FF{10}
\def\FDd{11}
\def\FDe{12}
\def\FI{13}
\def\FJ{14}
\def\FGa{15}
\def\FH{16}
\def\FK{17}
\def\FN{18}
\def\FL{19}
\def\FM{20}

\def\GF{\operatorname{GF}}

\def\fl#1{\lfloor#1\rfloor}
\def\cl#1{\lceil#1\rceil}

\topmatter 
\title Bounded Dyck paths, bounded alternating sequences, orthogonal polynomials,
and reciprocity
\endtitle 
\author J. Cigler and C.~Krattenthaler$^{\dagger}$
\endauthor 
\affil 
Fakult\"at f\"ur Mathematik, Universit\"at Wien,\\
Oskar-Morgenstern-Platz~1, A-1090 Vienna, Austria.\\
WWW: {\tt http://homepage.univie.ac.at/johann.cigler}\\
WWW: \tt http://www.mat.univie.ac.at/\~{}kratt
\endaffil
\address Fakult\"at f\"ur Mathematik, Universit\"at Wien,
Oskar-Morgenstern-Platz~1, A-1090 Vienna, Austria.\newline
WWW: \tt http://homepage.univie.ac.at/johann.cigler,
http://www.mat.univie.ac.at/\~{}kratt
\endaddress
\thanks $^\dagger$Research partially supported by the Austrian
Science Foundation FWF (grant S50-N15)
in the framework of the Special Research Program
``Algorithmic and Enumerative Combinatorics"%
\endthanks

\subjclass Primary 05A15; Secondary 05A19 05A30 11C20 15A15 33C45 42C05
\endsubjclass
\keywords Dyck paths, ballot paths, alternating sequences,
heaps of pieces, non-int\-er\-sect\-ing lattice paths, tableaux,
plane partitions, parallelogram polyominoes,
Chebyshev polynomials, orthogonal polynomials,
combinatorial reciprocity law
\endkeywords
\abstract 
The theme of this article is a ``reciprocity" between bounded up-down
paths and bounded alternating sequences. Roughly speaking, this
``reciprocity" manifests itself by the fact 
that the extension of the sequence of numbers
of paths of length~$n$,
consisting of diagonal up- and down-steps and being confined
to a strip of bounded width, to negative~$n$ produces numbers
of alternating sequences of integers that are bounded from below and
from above. We show that
this reciprocity extends to families of non-intersecting
bounded up-down paths and certain arrays of alternating sequences
which we call alternating tableaux. We provide as well weighted
versions of these results. Our proofs are based on Viennot's
theory of heaps of pieces and on the combinatorics of non-intersecting
lattice paths.
An unexpected application leads to a refinement of a result of
Bousquet-M\'elou and Viennot on the width-height-area generating
function of parallelogram polyominoes.
Finally, we exhibit the relation of
the arising alternating tableaux to plane partitions of strip shapes.
\endabstract
\endtopmatter
\document

\subhead 1. Introduction\endsubhead
{\it Reciprocity} is a much used term in mathematics. 
The immediate association is with the quadratic reciprocity law
of number theory and its numerous generalisations and variations.
In combinatorics, however, ``reciprocity" has a different meaning.
The term ``reciprocity law" was introduced by Richard Stanley
in \cite{\StanAY}. It refers to a situation where we encounter
numbers $a_n$, for $n\ge0$, with $a_n$ being the number of 
certain combinatorial
objects of ``size''~$n$. If these numbers $a_n$ satisfy a linear
recurrence with constant coefficients 
then we may extend the sequence $(a_n)_{n\ge0}$
to {\it negative} integers~$n$ by applying the recurrence
``backwards". If it should happen that $\vert a_n\vert$, for
$n<0$, has a combinatorial meaning as well, then we speak of a
({\it combinatorial}) ``{\it reciprocity law}". 
The simplest case of such a (combinatorial)
``reciprocity law" occurs when one considers the binomial coefficients
$\binom nk$ (for fixed~$k$), giving the number of subsets of
cardinality~$k$ of $\{1,2,\dots,n\}$ on the one hand, and
$\left\vert\binom {-n}k\right\vert=\binom {n+k-1}k$, giving the number
of sub{\it multi}sets (i.e., ``sets'' where repeated elements are allowed)
of cardinality~$k$ of $\{1,2,\dots,n\}$, on the other hand. 
The most well-known non-trivial instance of this phenomenon is {\it
Ehrhart--Macdonald reciprocity} for counting integer points in
polytopes, generalised by Stanley
in \cite{\StanAY} (see also \cite{\StanBI, Theorem~4.5.14}). 
As a matter of fact, a whole book has been devoted to this kind
of ``reciprocity", see \cite{\BeSaAA}. For further instances of combinatorial
reciprocity laws
see \cite{\AnBBAA, \FoStAA, \LoebAA, \PropXA, \SpeyXA, \StanXA
}.

The subject of this article is a reciprocity relation that seemingly
has not been noticed earlier. It is a reciprocity between
{\it up-down paths} whose height is bounded from below and from above
and {\it alternating sequences} whose elements are bounded from
below and from above. More precisely, let 
$C_{2n}^{(k)}$ denote the number of paths with steps $(1,1)$ and $(1,-1)$ 
starting at $(0,0)$ and ending at $(2n,0)$ never passing below the
$x$-axis and never passing above the horizontal line $y=k$. 
(These paths are also
known as bounded {\it Dyck paths}.) It is not difficult to see
--- and well-known ---
that, for fixed~$k$, the numbers $C_{2n}^{(k)}$ satisfy a linear
recurrence with constant coefficients in~$n$.\footnote{Equivalently,
the generating function $\sum_{n\ge0}C_{2n}^{(k)}x^{2n}$ is rational.
Theorem~\VA\ with $r=s=0$ shows that this is indeed the case.}
This recurrence can also be applied in the backwards direction, thus
defining  numbers $C_{2n}^{(k)}$ for {\it negative} integers~$n$.
With this notation, the first author observed that
$$
\det\left(C_{2n+2i+2j+2}^{(2k+1)}\right)_{0\le i,j\le k-1}
=
C_{-2n}^{(2k+1)}.
\tag\BAc
$$
He posted this conjecture on {\tt MathOverflow}.
In response, Stanley identified the number $C_{-2n}^{(2k+1)}$ as
the number of sequences $a_1\le a_2\ge a_3\le a_4\ge a_5\ge\dots\ge
a_{2n-1}$ of integers with $1\le a_i\le k+1$ for all~$i$.
(See Corollary~\TB\ with $k$ replaced by $k+1$ for a formal proof;
in his posting, Stanley referred to \cite{\StanAP, Ex.~3.66} and
\cite{\StanAZ, Ex.~3.2}.)
With this observation, and the observation that the determinant on the
left-hand side of (\BAc) can be interpreted as the number of families
of non-intersecting bounded Dyck paths,
it is then rather simple to confirm the
above reciprocity law, see the proof of Theorem~\TD\ with $m=1$.

As it turned out, the relation (\BAc) is just the peak of an iceberg.
In this article --- so-to-speak --- 
we ``unearth" the iceberg (or rather: ``lift
the iceberg out of the sea''). 

The first manifestation of this
``iceberg" arises if we lift the upper bound
$2k+1$ on the height of the paths to $2k+2m-1$, where $m$ is a
positive integer. Then there is still a reciprocity relation,
where on the right-hand side a {\it determinant\/} of numbers of paths
``with negative length" has to be placed; see Theorem~\TD\ in
Section~5.\footnote{We apologise to the reader that our presentation 
in this introduction follows the order in which things and ideas
developed, which is not entirely the same order as they appear in the
article. We will come back to Sections~2--4 in our presentation shortly;
in a sense, these sections contain preparatory material for the later
sections.} For the proof, a bijection is set up between families of
non-intersecting Dyck paths and certain arrays of integers of
trapezoidal shape in which each row is an alternating sequence of
odd length, see Proposition~\TE.
Due to their resemblance to semistandard tableaux
(cf\. e.g\. \cite{\StanBI, Sec.~7.10} for their definition),
we call these arrays {\it alternating tableaux}. The key for the proof
of Theorem~\TD\ is to show that these alternating tableaux of
trapezoidal shape are counted by the determinant on the right-hand
side of~(\AA), which is done in Theorem~\TF.

Section~6 presents a second manifestation of the mentioned ``iceberg",
with many parallels to the first manifestation:
it is a reciprocity between
bounded up-down paths starting at height~$0$ and {\it ending at the
maximally allowed height} and bounded alternating sequences of {\it
even} length. The main theorem of Section~6 is Theorem~\TG, which is
analogous to Theorem~\TD. It is based on a bijection between 
families of non-intersecting up-down paths of the described kind
and alternating tableaux of {\it rhomboidal\/} shape in which each
row has even length, see Proposition~\TH,
and the proof in Theorem~\TI\ 
that these rhomboidal alternating tableaux are counted
by the determinant on the right-hand side of~(\CA).

However, the reciprocity between bounded up-down paths and alternating
sequences occurs even at a finer level. It is known that generating
functions for the numbers of bounded up-down paths with specified
starting and ending height can be expressed as fractions involving
Chebyshev polynomials of the second kind, see Theorem~\VA\ in Section~2.
In Section~3, we show that the same is true for generating functions
for bounded alternating sequences with specified first and last
element, see Theorem~\UA. Although both theorems could be proved
using recurrence relations and matrix algebra, here we present proofs
that use Viennot's \cite{\VienAF} theory of heaps of pieces.
A comparison of the two results then entails further
reciprocity relations, see Section~4.

This finer reciprocity relation leads to a third manifestation of the
``iceberg", which is the subject of Section~7. The main results of this
section are Theorems~\TJ\ and~\TM\ presenting equalities between a determinant
of numbers of bounded up-down paths with positive lengths with specified
starting and ending heights and a similar determinant where however
the path lengths are negative. Here the proof is based on a
bijection between families of non-intersecting paths with given
starting and ending points and
{\it flagged\/} alternating tableaux of {\it rectangular} shape,
``flagged" meaning that the first and last entry in each row is
specified, see Propositions~\TK\ and \TN.
Theorems~\TL\ and~\TO\ then show that these alternating tableaux of
rectangular shape are counted by the determinants on the right-hand
sides of~(\DA) and (\DD), respectively.

Sections~8--11 are complements to the material in Sections~2--7.
First, the reader may wonder whether weights could be introduced
in the enumeration results in Sections~2--7. Sections~8--10 are devoted
to the derivation of the weighted generalisations of these results.
Here, the Chebyshev polynomials are replaced by more general
orthogonal polynomials $P_n(x)$ that are even for even~$n$ and odd for
odd~$n$. Section~8 addresses the weighted enumeration of bounded
up-down paths and bounded alternating sequences, 
the main results being Theorems
~\TP\ and \TQ. In Section~9, we
extend the enumeration results of Section~4 for paths ``with negative length" 
to the weighted setting. Section~10 then compiles the weighted
versions of the reciprocity laws from Sections 5--7, see Theorems~\TW, 
\TZ, \TTC, and~\TTF. 
Remarkably, one of these, namely Theorem~\TZ, implies a general reciprocity
relation between Hankel determinants of sequences that satisfy a linear recurrence
with constant coefficients, see Theorem~\TTTTTT.

Second, the enumeration of alternating tableaux plays an important role in
Sections~5--7. The purpose of Section~11 is to show that these
enumeration results are actually special cases of two fundamental
theorems for alternating tableaux, see Theorems~\TTI\ and~\TTJ.
In the second part of
that section, we explain that alternating tableaux are plane
partitions in disguise, and we discuss some insights and consequences
resulting from this observation. In particular, it turns out that
our determinantal formulae for generating functions for alternating
tableaux are closely related to (but not equivalent to)
the ribbon determinant formulae for Schur functions
due to Lascoux and Pragacz \cite{\LaPrAC}.

We close our article in Section~12 by a discussion and a presentation of some 
questions raised by the material presented in this article.
In particular, we discuss what happens for {\it even} upper bounds on the
up-down paths, we discuss the possibility of a further extension
to Motzkin paths, including a precise conjecture
(see Conjecture~\TTN) perhaps hinting at another iceberg, a ``coincidental''
application of Theorem~\UA\ to the enumeration of parallelogram polyominoes
which produces a seemingly new result 
(see Theorem~\TTO), the link being the heaps of segments that we use
in the proof of our enumeration results for bounded alternating sequences in
Section~3, and
several curious Hankel determinant evaluations involving numbers of bounded
up-down paths respectively of bounded alternating sequences 
(see Theorems~\TTLb\ and~\TTP\ and Corollaries~\TTT\ and~\TTU).

As a kind of ``afterthought", we added a ``Final Note" at the end of
the article pointing the reader to some interesting more recent developments.

\subhead 2. Enumeration of bounded up-down paths\endsubhead
In this section, we discuss the enumeration of bounded up-down paths.
The main result is
Theorem~\VA, in which the generating function for
bounded up-down paths with given starting and ending height is computed.
Although this is a known result, we present a proof 
that makes use of Viennot's \cite{\VienAF} theory
of heaps of pieces because it has not been recorded anywhere before,
and since it constitutes the inspiration for the proof of Theorem~\UA\ in the
next section.

\medskip
To start with, we fix notation. 
We write $C_n^{(k)}(r\to s)$ for the {\it number} 
of up-down paths from $(0,r)$ to
$(n,s)$ that do not pass below the $x$-axis and do not pass above the
horizontal line $y=k$.\footnote{It should be noted that, due to the
geometry of the paths, we need the
condition $r+s\equiv n$~(mod~$2$) for this number to be non-zero.} 
If $r=s=0$, i.e., if we are talking of Dyck paths, then, instead of
$C_n^{(k)}(0\to 0)$, we frequently write $C_{n}^{(k)}$ for short
(as we have done in the introduction).
There is another special case where we use an abridged notation:
we write $D_{2n}^{(k)}$ instead of $C_{2n+k}^{(k)}(0\to k)$ for short.

In the results of this section and later sections, 
the {\it Chebyshev polynomials of the
second kind\/} are ubiquitous.
As usual, we write $U_n(x)$ for the $n$-th Chebyshev polynomial
of the second kind, which is explicitly given by
$$
\align
U_n(\cos\th)&=\frac {\sin((n+1)\th)} {\sin\th},\\
  U_n(x)&=\sum_{j\ge0}(-1)^j\binom {n-j}{j}(2x)^{n-2j}.
\tag\BAb
\endalign
$$
It is well-known that these Chebyshev polynomials satisfy the
two-term recurrence
$$
2xU_n(x)=U_{n+1}(x)+U_{n-1}(x),
\tag\BD
$$
with initial conditions $U_0(x)=1$ and $U_1(x)=2x$.

\medskip
With the above notations, we are now ready to state the announced theorem,
giving closed form expressions for the generating function for
bounded up-down paths with given starting and ending height.

\proclaim{Theorem \VA}
For all non-negative integers $r,s,k$ with $0\le r,s\le k$, we have
$$
\sum_{n\ge0}C_n^{(k)}(r\to s)\, x^{n}
=\cases \dfrac {U_r(1/2x)\,U_{k-s}(1/2x)} {x\,U_{k+1}(1/2x)},
&\text{if }r\le s,\\
\dfrac {U_s(1/2x)\,U_{k-r}(1/2x)}
{x\,U_{k+1}(1/2x)},
&\text{if }r\ge s.
\endcases
\tag\BGa
$$
\endproclaim

As already mentioned, 
this is a known theorem, however not ``as known'' as it should be.
The oldest source that we are aware of is \cite{\VienAE, Ch.~V, Eq.~(27)},
where a combinatorial proof, based on the combinatorial theory developed
in \cite{\VienAE}, is sketched. The simplest proof is one that uses
the transfer matrix method, see \cite{\KratBP, proof of Theorem~A2} or
\cite{\KratCL, proof of Theorem~10.11.1}. (As a matter of fact, 
these three sources discuss a weighted generalisation for Motzkin
paths. See Section~12.(3).)
In our opinion, the most illuminating proof is one that is based on
Viennot's theory of
heaps of pieces \cite{\VienAF}, the latter being a geometric realisation of
the {\it``mono\"\i de partiellement commutatif"} of
Cartier and Foata \cite{\CaFoAA}, now known as the {\it
Cartier--Foata monoid} (see also~\cite{\KratAZ}).
This is the proof that we are going to present here since it cannot
be found anywhere else in the literature. It is inspired by ideas
that one finds in \cite{\GaEgAA, Ch.~4}. Moreover, it prepares for
the --- slightly more complicated --- 
proof of Theorem~\UA\ on generating functions for alternating
sequences.

The heaps that we need here are heaps of dimers that are contained
in the interval $[0,k]$. Here, a dimer 
is a vertical segment connecting the points $(a,b)$ and $(a,b+1)$, for
some positive integer~$a$ and non-negative integer~$b$. 
We write $d_b$ for such a dimer. It is intentional that we ignore
the abscissa~$a$ in this notation since we allow to move dimers
horizontally, thus modifying~$a$.
The ``rule of the game" is that
two dimers $d_i$ and $d_j$ can be (horizontally)
moved past one another if and only if they do not block each other 
``physically", that is, if and only if either $i>j+1$ or $j>i+1$.
In our picture, gravity pulls dimers (horizontally) to the left.
A {\it heap of dimers on~$[0,k]$} is then what the name suggests:
one piles dimers $d_i$ with $0\le i<k$ 
on each other, and gravity pulls them to their left-most possible
position according to the ``rule of the game".
An example of such a heap on $[0,k]$
is shown in Figure~\FG.b. We denote the set of all possible
dimers on $[0,k]$ by $\Cal D_k$, that is, $\Cal
D_k=\{d_0,d_1,\dots,d_{k-1}\}$.
It should be noted that a heap may contain several copies of a dimer
$d_i$. For example, the heap in Figure~\FG.b contains three copies
of the dimer $d_1$. 

We define a weight function $w$ on dimers
by $w(d_i)=x^2$, and, given a heap~$H$, we extend
the weight function to $H$ by declaring that $w(H)$ equals the product
of the weights of all dimers of~$H$.

A {\it maximal dimer} of a heap is a dimer that can be moved to the
far right without being blocked by any other dimer. 
Analogously, a {\it minimal dimer} of a heap is a dimer that can 
be moved to the
far left without being blocked by any other dimer. In our example in
Figure~\FG.b, the maximal dimers of the heap shown there are 
the right-most dimer $d_{r-1}$, the right-most dimer
$d_{r+1}$, and the dimer $d_s$, while the minimal dimers
are the dimer $d_0$, the left-most dimer $d_{r+1}$, and the dimer $d_{s+1}$.
Finally, a {\it trivial heap} is one in which all its dimers are maximal.

With these definitions, the main theorem in the theory of heaps
\cite{\VienAF, Prop.~5.3} (see also \cite{\KratAZ, Theorem~4.1}) implies
that the generating function $\sum_Hw(H)$ for all heaps $H$ whose
maximal dimers are contained in a given subset $\Cal M$ of $\Cal D_k$
is given by
$$
\underset \text{maximal dimers}\subseteq\Cal M
\to{\sum_{H\text{ heap of dimers on }[0,k]}}
\kern-.3cm w(H)
=\frac {\dsize\underset \text{dimers}\subseteq \Cal D_k\backslash\Cal M
\to{\sum_{T\text{ trivial heap}}}(-1)^{\vert T\vert}w(T)} 
{\dsize \underset \text{dimers}\subseteq \Cal D_k
\to{\sum_{T\text{ trivial heap}}}(-1)^{\vert T\vert}w(T)} ,
\tag\BGb
$$
where $\vert T\vert$ denotes the number of dimers of~$T$.

The purpose of the next lemma is to transfer the problem of
enumeration of bounded up-down paths to a problem of enumeration
of heaps of dimers.

\proclaim{Lemma \VB}
Let $n,k,r,s$ be non-negative integers with $0\le r\le s\le k$
and\linebreak $n\equiv r+s$~{\rm(mod~$2$)}.
There is a bijection between up-down paths from $(0,r)$ to $(n,s)$
that do not pass below the $x$-axis and do not pass above the
horizontal line $y=k$
and heaps $H$ of $(n-(s-r))/2$ 
dimers on $[0,k]$ whose maximal dimers are contained
in $[r-1,s+1]$.
\endproclaim

\demo{Proof}
Let $P$ be an up-down path $P$ as in the statement of the lemma.
See Figure~\FG.a for an example. Out of~$P$, we are going to,
step-by-step, build a heap~$H$ of dimers. In the beginning,
$H$ is empty.

\midinsert
$$
\Einheit.38cm
\Koordinatenachsen(34,10)(0,0)
\Pfad(0,3),33444443334334443433333443333344\endPfad
\PfadDicke{.1pt}
\Pfad(0,3),111111111111111111111111111111111\endPfad
\Pfad(0,7),111111111111111111111111111111111\endPfad
\PfadDicke{1pt}
\Pfad(0,9),111111111111111111111111111111111\endPfad
\DickPunkt(0,3)
\DuennPunkt(1,4)
\DuennPunkt(2,5)
\DuennPunkt(3,4)
\Kreis(4,3)
\DuennPunkt(5,2)
\DuennPunkt(6,1)
\DuennPunkt(7,0)
\DuennPunkt(8,1)
\DuennPunkt(9,2)
\Kreis(10,3)
\DuennPunkt(11,2)
\Kreis(12,3)
\DuennPunkt(13,4)
\Kreis(14,3)
\DuennPunkt(15,2)
\DuennPunkt(16,1)
\DuennPunkt(17,2)
\DuennPunkt(18,1)
\DuennPunkt(19,2)
\Kreis(20,3)
\DuennPunkt(21,4)
\DuennPunkt(22,5)
\DuennPunkt(23,6)
\DuennPunkt(24,5)
\DuennPunkt(25,4)
\DuennPunkt(26,5)
\DuennPunkt(27,6)
\DuennPunkt(28,7)
\DuennPunkt(29,8)
\DuennPunkt(30,9)
\DuennPunkt(31,8)
\DickPunkt(32,7)
\Label\l{0}(0,0)
\Label\l{r}(0,3)
\Label\l{s}(0,7)
\Label\l{k}(0,9)
\hskip11.8cm
$$
\centerline{\eightpoint a. An up-down path from height $r$ to height $s$}
$$
\Pfad(0,-1),22222222222\endPfad
\PfadDicke{.1pt}
\Pfad(0,0),11111111\endPfad
\Pfad(0,1),11111111\endPfad
\Pfad(0,2),11111111\endPfad
\Pfad(0,3),11111111\endPfad
\Pfad(0,4),11111111\endPfad
\Pfad(0,5),11111111\endPfad
\Pfad(0,6),11111111\endPfad
\Pfad(0,7),11111111\endPfad
\Pfad(0,8),11111111\endPfad
\Pfad(0,9),11111111\endPfad
\Kreis(1,0)
\Kreis(1,1)
\Kreis(1,8)
\Kreis(1,9)
\Kreis(1,4)
\Kreis(1,5)
\Kreis(2,1)
\Kreis(2,2)
\Kreis(2,3)
\Kreis(2,4)
\Kreis(2,7)
\Kreis(2,8)
\Kreis(2,5)
\Kreis(2,6)
\Kreis(3,2)
\Kreis(3,3)
\Kreis(4,2)
\Kreis(4,3)
\Kreis(5,1)
\Kreis(5,2)
\Kreis(6,1)
\Kreis(6,2)
\Kreis(7,2)
\Kreis(7,3)
\Kreis(5,3)
\Kreis(5,4)
\Kreis(6,4)
\Kreis(6,5)
\PfadDicke{2pt}
\Pfad(1,0),2\endPfad
\Pfad(2,1),2\endPfad
\Pfad(1,4),2\endPfad
\Pfad(2,3),2\endPfad
\Pfad(1,8),2\endPfad
\Pfad(3,2),2\endPfad
\Pfad(4,2),2\endPfad
\Pfad(5,1),2\endPfad
\Pfad(6,1),2\endPfad
\Pfad(7,2),2\endPfad
\Pfad(2,5),2\endPfad
\Pfad(2,7),2\endPfad
\Pfad(5,3),2\endPfad
\Pfad(6,4),2\endPfad
\Label\l{0}(0,0)
\Label\l{r}(0,3)
\Label\l{s}(0,7)
\Label\l{k}(0,9)
\hbox{\hskip.3cm}
\raise.3cm\hbox{%
\SPfad(0,-1),12121225656566\endSPfad
\SPfad(3,1),122566\endSPfad
\SPfad(4,1),61121225655\endSPfad
\Label\ru{H_3}(1,0)
\Label\ro{H_2}(3,3)
\Label\ru{H_1}(6,1)
}%
\hskip3.2cm
$$
\centerline{\eightpoint b. The corresponding heap of dimers}
\vskip7pt
\centerline{\eightpoint Figure \FG}
\endinsert

We cut $P$ into path portions at the touching points with the
horizontal line $y=r$. In the figure, these touching points
are marked by circles. 
Path portions between touching points with $y=r$ are treated
differently depending on whether they stay above or below
that line. We call the path portions which stay above {\it positive},
and those which stay below {\it negative}.

We now scan the path from left to right.
While scanning a positive path portion, 
we ignore up-steps, and for each down-step
from height $h+1$ to~$h$ we put a dimer $d_h$ on~$H$, in the
order we encounter the down-steps. On the other hand,
while scanning a negative path portion,
we ignore down-steps, and for each up-step
from height $h$ to~$h+1$ we put a dimer $d_h$ on~$H$, in the
order we encounter the up-steps.
So, for instance, our example path in Figure~\FG.a starts with
a positive portion $UUDD$ (here, $U$ stands for an up-step and
$D$ for a down-step). The first down-step in this portion is from
height $r+2$ to $r+1$, therefore we put the dimer $d_{r+1}$ on
the (at this point empty) heap~$H$. The second down-step is from
height $r+1$ to height~$r$, therefore we next put a dimer $d_r$
on~$H$. On the other hand, the second path portion, $DDDUUU$ 
is negative. We ignore the down-steps, and the up-steps, 
read from left to right, correspond
to the dimers $d_0,d_1,d_2$, in this order. 
They are put on~$H$, in this order. 

This process is continued until all steps of~$P$ have been taken
into account.

It is not difficult to see that a heap obtained in this way has indeed
the property that its maximal dimers are contained in $[r-1,s+1]$.

\medskip
In order to see that this is a bijection, we have to describe the
inverse mapping. Given a heap~$H$ as in the statement of the lemma,
we have to construct an up-down path~$P$. Again, we do this
step-by-step, but from the ``back to the front". We start with the path
portion that just consists of the point~$(n,s)$.
We look for the top-most maximal dimer of~$H$, say $d_i$.
By the condition on maximal dimers in the lemma, we must have $i\le
s$. If also $i\ge r$,
then we prepend a down-step from height $i+1$
to height~$i$ followed by the appropriate number of up-steps to
connect to the previously constructed path portion. Subsequently, we remove
this dimer $d_i$ from~$H$. 

We repeat this procedure until we encounter a top-most maximal dimer $d_i$
with $i=r-1$. (Because of the condition on maximal
dimers, $i$ cannot be smaller than that.)
Let us illustrate the process described so far by considering our
example heap in Figure~\FG.b. There, the top-most maximal dimer is
$d_s$. Hence the last step of the up-down path to be determined is
a down-step from height~$s+1$ to height~$s$ (with zero up-steps to be
inserted in order to connect with $(n,s)$). See Figure~\FG.a.
We remove $d_s$ from the heap. Now the top-most maximal dimer
is~$d_{k-1}$. Hence, we prepend a down-step from height~$k$ to
height~$k-1$ to the path portion constructed so far. See again
Figure~\FG.a. We remove $d_{k-1}$ from the heap. Now the top-most
maximal dimer is~$d_{r+1}$. Consequently, we prepend a down-step
from height $r+2$ to height $r+1$ followed by five up-steps to
the already constructed path portion. After removal of $d_{r+1}$,
the top-most maximal dimer is $d_{r+2}$. We prepend a down-step
from height $r+3$ to height $r+2$ to our path portion, and we
remove $d_{r+2}$ from the heap. Now the top-most maximal dimer
is a dimer $d_{r-1}$.

We now continue the description of the general case.
If, during the process, we encounter a dimer $d_{r-1}$, 
then we push this dimer to the left. Thereby, 
due to the ``rule of the game" that dimers cannot be moved
past one another 
if they ``physically" touch each other, other dimers may also be
moved to the left. We are only interested into those which lie
strictly below the horizontal line $y=r$
and do not have another dimer $d_{r-1}$
``in between". In that way, a subheap is determined. 
In our example in Figure~\FG.b, subheaps arising in this way are
surrounded by dotted paths. The subheap that we would encounter
next is the subheap labelled~$H_1$. 

Such a subheap is now treated in the following way:
in the beginning, the top-most maximal dimer is $d_{r-1}$.
We prepend an up-step (!) from height~$r-1$ to height~$r$ followed by 
the appropriate number of up-steps to
connect to the previously constructed path portion. Subsequently, we
remove $d_{r-1}$ from the subheap. Next we consider the --- now ---
top-most maximal dimer {\it in the subheap}, $d_j$ say. We prepend
an up-step (!) from height~$j$ to height~$j+1$ followed by 
the appropriate number of down-steps (!) to connect to the so far
constructed path portion, and we remove $d_j$ from the subheap. 
We repeat this latter procedure until the subheap is emptied. 

This process is iterated, that is, again the maximal dimer is
determined and it is treated as above depending on whether it is
a dimer $d_i$ with $i\ge r$ or with $i=r-1$.

It is not too difficult to see that this algorithm is the inverse
of the map from up-down paths to heaps that we described in the first
part of this proof.\quad \quad \qed
\enddemo

\remark{Remark}
The obvious extension of the above bijection between up-down paths
and heaps of dimers to a bijection between paths with allowed
steps being $(1,1)$, {\it($\text{\it 1},\text{\it 0}$)}, $(1,-1)$ 
and heaps of dimers {\it and monomers}
also proves the earlier mentioned more general result in
\cite{\VienAE, Ch.~V, Eq.~(27)}, \cite{\KratBP, Theorem~A2}, and
\cite{\KratCL, Theorem~10.11.1}.
\endremark

If we combine Lemma~\VB\ with (\BGb), then we see that, for the proof
of Theorem~\VA, we need to
compute the generating function $\sum_T(-1)^{\vert T\vert}w(T)$
for all trivial heaps with dimers in $\Cal D_k$, that is, all
trivial heaps on $[0,k]$, and the analogous generating function
for trivial heaps whose dimers are in $[0,r-1]\cup[s+1,k]$.
The former is accomplished in the lemma below.
The computation of the latter is then an easy consequence as, due to the
shift-invariance of our weight~$w$, it decomposes into the product
of the generating function for trivial heaps on $[0,s-1]$ and
the generating function for trivial heaps on $[0,k-s-1]$.

\proclaim{Lemma \VC}
Let $k$ be a non-negative integer.
The generating function $\sum_T (-1)^{\vert T\vert}w(T)$, 
where the sum is over all
trivial heaps $T$ of dimers on $[0,k]$, is given by $x^{k+1}U_{k+1}(1/2x)$.
\endproclaim

\demo{Proof}
We prove the claim by induction on $k$. The claim is certainly true for
$k=0$ --- in this case there is only the empty heap with weight~$1$;
indeed, $x^1U_1(1/2x)=1$ ---, and for $k=1$
--- in this case there exist two heaps: the empty heap with weight~$1$,
and the heap consisting of the dimer~$d_0$ that has weight $x^2$; indeed,
$x^2U_2(1/2x)=1-x^2$.

For the induction step, 
consider trivial heaps of dimers on $[0,k+1]$. There are two possibilities:
\roster
\item"(a)"the heap contains the dimer $d_k$;
\item"(b)"the heap has all its dimers in $[0,k]$.
\endroster
By induction hypothesis,
the contribution to the generating function of the heaps in
Case~(a) is $(-x^2)x^{k-1}U_{k-1}(1/2x)$. Similarly, the contribution of
the heaps in Case~(b) is $x^kU_{k}(1/2x)$. 
Hence, in total we obtain
$$
(-x^2)x^{k-1}U_{k-1}(1/2x)+x^kU_{k}(1/2x)=x^{k+1}U_{k+1}(1/2x),
\tag\BGc
$$
the equality following from the two-term recurrence (\BD) for 
Chebyshev polynomials.\quad \quad \qed
\enddemo

We are now ready for the proof of Theorem~\VA.

\demo{Proof of Theorem \VA}
We treat the case $r\le s$. The other case follows by a reflection
of paths in a vertical line which implies a switch of the roles of~$r$
and~$s$.

We use Lemma~\VB\ to see that the up-down paths that we are interested
in are in bijection with heaps that are described in the statement
of the lemma. Thus, by (\BGb) with $\Cal M$ the set of dimers
contained in~$[r-1,s+1]$, we see that
$$
\sum_{n\ge0}C_n^{(k)}(r\to s)\, x^{n}
=x^{s-r}
\frac {\dsize\underset \text{dimers}\subseteq [0,r-1]\cup[s+1,k]
\to{\sum_{T\text{ trivial heap on }[0,k]}}(-1)^{\vert T\vert}w(T)} 
{\dsize {\sum_{T\text{ trivial heap on }[0,k]}}(-1)^{\vert T\vert}w(T)} .
$$
Here, the factor $x^{s-r}$ must be inserted since an up-down path from
height~$r$ to height~$s$ has $s-r$ more up-steps than down steps.
The claimed expression in the first line on the right-hand side of
(\BGa) now follows from Lemma~\VC.\quad \quad \qed
\enddemo

\subhead 3. Enumeration of bounded alternating sequences\endsubhead
In this section, we discuss the enumeration of 
bounded alternating sequences. The main result is
Theorem~\UA, in which we compute the generating function for
bounded alternating sequences with given first and last element.
Also here,  we present a proof that makes 
use of Viennot's \cite{\VienAF} theory
of heaps of pieces. 
Specialisations of Theorem~\UA\ can be used to find the
generating functions for the numbers of alternating sequences in which
first and last element are not specified.
These are recorded here in Corollaries~\UG--\TA.
In particular, the final result in this section, Corollary~\TA, sets up
the connection to Stanley's observation mentioned in the introduction.

\medskip
We start again by fixing notation. 
We let $\Cal A_n^{(k)}(r\to s)$ denote the {\it set\/} of alternating sequences
$$r\le a_2\ge a_3\le a_4\ge \cdots {}\diamond\,a_{n-1}{}\square\,s,$$
where $\diamond=\,\ge$ and $\square=\,\le$ if $n$ is even and
$\diamond=\,\le$ and $\square=\,\ge$ if $n$ is odd,
in which all $a_i$'s are integers between~$1$ and~$k$.
If we do not want to specify the first and last element of an
alternating sequence then we omit the corresponding indication. More
precisely, we write $\Cal A_{n}^{(k)}$ for the union of the sets
$\Cal A_n^{(k)}(r\to s)$ over all~$r$ and ~$s$.\footnote{The reader
should be aware that the meaning of the shorthand in the abridged
notations $C_n^{(k)}$ from Section~2 
and $\Cal A_n^{(k)}$ is {\it not\/}
 the same: while
in $C_n^{(k)}$ the omission of ``$(r\to s)$" means $r=s=0$,
the omission of ``$(r\to s)$" in $\Cal A_n^{(k)}$ means that we
allow {\it all\/} possible~$r$ and~$s$.}

\medskip
The theorem below
gives closed form expressions for the generating function for
bounded alternating sequences with given first and last element.
The reader should note the similarities with the result in
Theorem~\VA. These will be subsequently exploited in Section~4.

\proclaim{Theorem \UA}
For all positive integers $r,s,k$ with $1\le r,s\le k$, we have
$$
\sum_{n\ge0}\big\vert\Cal A_{2n+1}^{(k)}(r\to s)\big\vert x^{2n}=
\cases 
\dsize
(-1)^{r+s+1}\frac {xU_{2r-2}(x/2)U_{2k+1-2s}(x/2)} 
{U_{2k}(x/2)},&\text{if }r< s,\\
\dsize
1-\frac {xU_{2r-2}(x/2)U_{2k+1-2r}(x/2)} 
{U_{2k}(x/2)},&\text{if }r=s,\\
\dsize
(-1)^{r+s+1}\frac {xU_{2s-2}(x/2)U_{2k+1-2r}(x/2)} 
{U_{2k}(x/2)},&\text{if }r> s,
\endcases
\tag\BG
$$
and
$$
\sum_{n\ge0}\big\vert\Cal A_{2n+2}^{(k)}(r\to s)\big\vert x^{2n+1}=
\cases 
\dsize
(-1)^{r+s+1}\frac {xU_{2r-2}(x/2)U_{2k-2s}(x/2)}
{U_{2k}(x/2)},&\text{if }r\le s,\\
\dsize
(-1)^{r+s+1}\frac {xU_{2s-1}(x/2)U_{2k-2r+1}(x/2)}
{U_{2k}(x/2)},&\text{if }r> s.
\endcases
\tag\BH
$$
\endproclaim


As already announced,
for the proof of Theorem~\UA\ we use again Viennot's theory of
heaps of pieces \cite{\VienAF}.
The heaps that we need here are heaps of segments that are contained
in the interval $[1,k]$. Here, a segment $j${}$\pmb-${}$i$, with 
$1\le i\le j\le k$,
is a vertical segment connecting the points $(a,i)$ and $(a,j)$, for
some positive integer~$a$. 
The ``rule of the game" is that
two segments $j_1${}$\pmb-${}$i_1$ and $j_2${}$\pmb-${}$i_2$ can be (horizontally)
moved past one another if and only if they do not block each other 
``physically", that is, if and only if either $i_1>j_2$ or $i_2>j_1$.
As in Section~2,
in our picture, gravity pulls segments (horizontally) to the left.
An example of such a {\it heap of segments on $[1,k]$}
is shown in the bottom of Figure~\FA. We denote the set of all possible
segments in $[1,k]$ by $\Cal S_k$. Also here
it should be noted that a heap can contain several copies of a segment
$j${}$\pmb-${}$i$. For example, the heap in Figure~\FA\ contains two copies
of the segment $3${}$\pmb-${}$2$. We point out that this kind of heaps of
segments (without the upper bound imposed by~$k$) appeared earlier
in \cite{\BoViAA} in the context of enumeration of directed
convex polyominoes, see Section~12.(5) for more details on this
coincidence, including an implied generating function result for
parallelogram polyominoes that seems to be new.

We define a weight function $w$ on the segments
by $w(\text{$j${}$\pmb-${}$i$})=x^2$, and, given a heap~$H$, we extend
the weight function to $H$ by declaring that $w(H)$ equals the product
of the weights of all segments of~$H$.

We define {\it maximal\/} and {\it minimal segments} of a heap 
in the same way as in Section~2. In our example in
Figure~\FA, the maximal segments of the heap shown there are the segments
$1${}$\pmb-${}$1$, $2${}$\pmb-${}$2$ and $8${}$\pmb-${}$6$, while the minimal segments
are $4${}$\pmb-${}$3$ and $6${}$\pmb-${}$6$.
Finally, a {\it trivial heap} is one in which all its segments are maximal.

With these definitions, the main theorem in the theory of heaps
\cite{\VienAF, Prop.~5.3} (see also \cite{\KratAZ, Theorem~4.1}) implies
that the generating function $\sum_Hw(H)$ for all heaps $H$ whose
maximal segments are contained in a given subset $\Cal M$ of $\Cal S_k$
is given by
$$
\underset \text{maximal segments}\subseteq\Cal M
\to{\sum_{H\text{ heap of segments on }[1,k]}}
\kern-.3cm w(H)
=\frac {\dsize\underset \text{segments}\subseteq \Cal S_k\backslash\Cal M
\to{\sum_{T\text{ trivial heap}}}(-1)^{\vert T\vert}w(T)} 
{\dsize \underset \text{segments}\subseteq \Cal S_k
\to{\sum_{T\text{ trivial heap}}}(-1)^{\vert T\vert}w(T)} ,
\tag\BI
$$
where, as before, $\vert T\vert$ denotes the number of segments of~$T$.

The next lemma says that an alternating sequence can be converted into
a heap of segments in a rather straightforward way, and that this conversion
is actually a bijection. However, the obtained heaps do not form a set of
heaps to which we could apply~(\BI). We need
to modify this construction. This is what we do in the proof of Lemma~\UC.
The construction in the proof of Lemma~\UB\ will however be a part
of it.

\proclaim{Lemma \UB}
Let $n,k,r,s$ be non-negative integers with $1\le r\le s\le k$.
There is a bijection between $\Cal A_{2n+1}^{(k)}(r\to s)$ and heaps $H$ of 
$n$~segments
on $[1,k]$ with the following three properties:

\roster 
\item $H$ has a maximal segment of the form $j${}$\pmb-${}$s$ with $j\ge s$.
\item $H$ does not have any maximal segments that are contained in $[s+1,k]$.
\item $H$ does not have any minimal segments that are contained in $[1,r-1]$.
\endroster
\endproclaim

\demo{Proof}
Given an alternating sequence 
$r\le a_2\ge a_3\le a_4\ge\cdots\le a_{2n}\ge s$,
we build a heap by piling the segments
$a_2${}$\pmb-${}$a_3$,
$a_4${}$\pmb-${}$a_5$, \dots,
$a_{2n}${}$\pmb-${}$s$ on each other, in this order.
See Figure~\FA\ for an example in which $n=14$, $k=8$, $r=4$, and $s=6$.
We also mark the point $(0,r)$
and the lower point of the segment $a_{2n}${}$\pmb-${}$s$
in order to transfer the information about the first and
last entry in the alternating sequence to the heap picture.
Thus, in Figure~\FA\ the points $(0,4)$ and $(8,6)$ are marked.

\midinsert
$$
\multline
4\le 4\ge 3\le 3\ge 1\le 1\ge 1\le 3\ge 2\le 6\ge 6\le 8\ge 4\le 7\ge 4\le 4\ge 2\\
\le 3\ge 2\le 2\ge 2\le 5\ge 5\le 6\ge 3\le 6\ge 5\le 8\ge 6
\endmultline
$$
$$\updownarrow$$
$$
\Einheit.6cm
\Koordinatenachsen(9,9)(0,0)
\PfadDicke{.1pt}
\Pfad(0,1),111111111\endPfad
\Pfad(0,2),111111111\endPfad
\Pfad(0,3),111111111\endPfad
\Pfad(0,4),111111111\endPfad
\Pfad(0,5),111111111\endPfad
\Pfad(0,6),111111111\endPfad
\Pfad(0,7),111111111\endPfad
\Pfad(0,8),111111111\endPfad
\PfadDicke{2pt}
\Pfad(1,3),2\endPfad
\Pfad(2,1),22\endPfad
\Pfad(3,2),2\endPfad
\Pfad(2,4),2222\endPfad
\Pfad(3,4),222\endPfad
\Pfad(4,2),22\endPfad
\Pfad(5,2),2\endPfad
\Pfad(6,3),222\endPfad
\Pfad(7,5),2\endPfad
\Pfad(8,6),22\endPfad
\Kreis(1,3)
\Kreis(1,4)
\Kreis(1,6)
\Kreis(2,1)
\Kreis(2,3)
\Kreis(2,4)
\Kreis(2,8)
\Kreis(3,1)
\Kreis(3,2)
\Kreis(3,3)
\Kreis(3,4)
\Kreis(3,7)
\Kreis(4,2)
\Kreis(4,4)
\Kreis(4,5)
\Kreis(5,2)
\Kreis(5,3)
\Kreis(6,2)
\Kreis(6,3)
\Kreis(6,6)
\Kreis(7,5)
\Kreis(7,6)
\Kreis(8,6)
\Kreis(8,8)
\DickPunkt(0,4)
\DickPunkt(8,6)
\Label\l{1}(0,1)
\Label\l{2}(0,2)
\Label\l{3}(0,3)
\Label\l{4}(0,4)
\Label\l{5}(0,5)
\Label\l{6}(0,6)
\Label\l{7}(0,7)
\Label\l{8}(0,8)
\Label\u3(1,3)
\Label\o4(1,4)
\Label\o6(1,6)
\Label\u1(2,1)
\Label\l{\kern5pt3}(2,3)
\Label\r{4\kern5pt}(2,4)
\Label\o8(2,8)
\Label\r{4\kern5pt}(3,4)
\Label\r{3\kern5pt}(3,3)
\Label\u1(3,1)
\Label\l{\kern5pt2}(3,2)
\Label\o7(3,7)
\Label\u2(4,2)
\Label\r4(4,4)
\Label\o5(4,5)
\Label\u2(5,2)
\Label\o3(5,3)
\Label\u2(6,2)
\Label\r{3\kern5pt}(6,3)
\Label\o6(6,6)
\Label\u5(7,5)
\Label\o6(7,6)
\Label\u6(8,6)
\Label\o8(8,8)
\hbox{\hskip5.4cm}
$$
\centerline{\eightpoint Illustration of the bijection of Lemma~\UB\ for
$n=14$, $k=8$, $r=4$, $s=6$}
\vskip7pt
\centerline{\eightpoint Figure \FA}
\endinsert

It is not difficult to see that a heap obtained in this way satisfies
the properties (1)--(3).

\medskip
In order to see that this is indeed a bijection we describe the
inverse mapping. Again, we work from the ``back to the front". 
Given a heap as described in the statement of the lemma,
we look for the top-most maximal segment, say $j${}$\pmb-${}$i$.
Then $j\ge i$ are the last two entries in the corresponding alternating
sequence. Subsequently, we remove this segment from the heap and proceed iteratively
to reconstruct the remaining entries of the alternating 
sequence.\quad \quad \qed
\enddemo

As announced earlier, the heaps in the statement of Lemma~\UB\ are not suited for application
of (\BI) since Condition~(3) above is one on {\it minimal\/} segments. In order to achieve
our goal --- the proof of Theorem~\UA --- we need to modify the above bijection.

\proclaim{Lemma \UC}
Let $n,k,r,s$ be non-negative integers with $1\le r\le s\le k$.
There is a bijection between $\Cal A_{2n+1}^{(k)}(r\to s)$
and heaps $H'$ of $n$~segments
on $[1,k]$ with the following two properties:

\roster 
\item"(1')" $H'$ has a maximal segment of the form $j${}$\pmb-${}$s$ with
$j\ge s$.
\item"(2')" $H'$ does not have any maximal segments that are contained in $[1,r-1]$ or $[s+1,k]$.
\endroster

\endproclaim

\demo{Proof}
Given an alternating sequence 
$r\le a_2\ge a_3\le a_4\ge\cdots\le a_{2n}\ge s$,
we again group the entries in pairs:
$a_2${}$\pmb-${}$a_3$,
$a_4${}$\pmb-${}$a_5$, \dots,
$a_{2n}${}$\pmb-${}$s$. However, before piling them on each other, we first
modify the order of these pairs (which will become segments in the heap picture).

\midinsert
$$
\multline
4\le 4\ge 3\le 3\ge 1\le 1\ge 1\le 3\ge 2\le 6\ge 6\le 8\ge 4\le 7\ge 4\le 4\ge 2\\
\le 3\ge 2\le 2\ge 2\le 5\ge 5\le 6\ge 3\le 6\ge 5\le 8\ge 6
\endmultline
$$
$$\downarrow$$
$$
\multline
\text{$3${}$\pmb-${}$2$},\
\text{$1${}$\pmb-${}$1$},\
\text{$3${}$\pmb-${}$1$},\
\text{$4${}$\pmb-${}$3$},\ \big\vert\
\text{$6${}$\pmb-${}$6$},\ \big\vert\
\text{$8${}$\pmb-${}$4$},\ \big\vert\
\text{$7${}$\pmb-${}$4$},\ \big\vert\ \\
\text{$2${}$\pmb-${}$2$},\
\text{$3${}$\pmb-${}$2$},\
\text{$4${}$\pmb-${}$2$},\ \big\vert\
\text{$5${}$\pmb-${}$5$},\ \big\vert\
\text{$6${}$\pmb-${}$3$},\ \big\vert\
\text{$6${}$\pmb-${}$5$},\ \big\vert\
\text{$8${}$\pmb-${}$6$}
\endmultline
$$
$$\downarrow$$
$$
\Einheit.58cm
\Koordinatenachsen(10,9)(0,0)
\PfadDicke{.1pt}
\Pfad(0,1),1111111111\endPfad
\Pfad(0,2),1111111111\endPfad
\Pfad(0,3),1111111111\endPfad
\Pfad(0,4),1111111111\endPfad
\Pfad(0,5),1111111111\endPfad
\Pfad(0,6),1111111111\endPfad
\Pfad(0,7),1111111111\endPfad
\Pfad(0,8),1111111111\endPfad
\PfadDicke{2pt}
\Pfad(3,3),2\endPfad
\Pfad(2,1),22\endPfad
\Pfad(1,2),2\endPfad
\Pfad(4,4),2222\endPfad
\Pfad(5,4),222\endPfad
\Pfad(6,2),22\endPfad
\Pfad(4,2),2\endPfad
\Pfad(7,3),222\endPfad
\Pfad(8,5),2\endPfad
\Pfad(9,6),22\endPfad
\Kreis(3,3)
\Kreis(3,4)
\Kreis(1,6)
\Kreis(2,1)
\Kreis(2,3)
\Kreis(3,2)
\Kreis(4,4)
\Kreis(4,8)
\Kreis(1,1)
\Kreis(1,2)
\Kreis(1,3)
\Kreis(5,4)
\Kreis(5,7)
\Kreis(6,2)
\Kreis(6,4)
\Kreis(6,5)
\Kreis(4,2)
\Kreis(4,3)
\Kreis(7,3)
\Kreis(7,6)
\Kreis(8,5)
\Kreis(8,6)
\Kreis(9,6)
\Kreis(9,8)
\DickPunkt(0,4)
\DickPunkt(9,6)
\Label\l{1}(0,1)
\Label\l{2}(0,2)
\Label\l{3}(0,3)
\Label\l{4}(0,4)
\Label\l{5}(0,5)
\Label\l{6}(0,6)
\Label\l{7}(0,7)
\Label\l{8}(0,8)
\Label\r{3\kern5pt}(3,3)
\Label\o4(3,4)
\Label\o6(1,6)
\Label\u1(2,1)
\Label\o{3}(2,3)
\Label\u{2}(3,2)
\Label\r{4\kern5pt}(4,4)
\Label\o8(4,8)
\Label\r{4\kern5pt}(5,4)
\Label\o{3}(1,3)
\Label\u1(1,1)
\Label\l{\kern5pt2}(1,2)
\Label\o7(5,7)
\Label\u2(6,2)
\Label\r{4\kern5pt}(6,4)
\Label\o5(6,5)
\Label\u2(4,2)
\Label\r{3\kern5pt}(4,3)
\Label\u3(7,3)
\Label\o6(7,6)
\Label\u5(8,5)
\Label\o6(8,6)
\Label\u6(9,6)
\Label\o8(9,8)
\hbox{\hskip.3cm}
\raise.3cm\hbox{%
}%
\hbox{\hskip4.9cm}
$$
\centerline{\eightpoint Illustration of the bijection of Lemma~\UC\ for
$n=14$, $k=8$, $r=4$, $s=6$}
\vskip7pt
\centerline{\eightpoint Figure \FB}
\endinsert

The reordering is accomplished by 
the following algorithm, in which we construct a sequence, $S$ say, of pairs.
Initialise $\ell:=1$ and $S:=\emptyset$, the empty sequence. 

\medskip
\roster
\item""{\smc Procedure.}
If $\ell=n+1$ then the algorithm terminates and its output is $S$.

Otherwise, consider the pair $a_{2\ell}${}$\pmb-${}$a_{2\ell+1}$.
From the construction, whenever we arrive at this point, 
we have $a_{2\ell}\ge r$. For $\ell=1$, this is true since the
first element of our alternating sequence is~$r$.
If also $a_{2\ell+1}\ge r$, then we let 
$S:=S,${}$a_{2\ell}${}$\pmb-${}$a_{2\ell+1}$ and $\ell:=\ell+1$,
and we repeat the Procedure.

If, on the other hand, $a_{2\ell+1}<r$, then we continue reading the 
pairs $a_{2i}${}$\pmb-${}$a_{2i+1}$,
for $i=\ell+1,\ell+2,\dots$,
until we meet again a pair $a_{2j}${}$\pmb-${}$a_{2j+1}$ in which $a_{2j}\ge r$.
In the alternating sequence in Figure~\FB, we have
$a_2${}$\pmb-${}$a_3${}${}={}$$4${}$\pmb-${}$3$,
$a_4${}$\pmb-${}$a_5${}${}={}$$3${}$\pmb-${}$1$,
$a_6${}$\pmb-${}$a_7${}${}={}$$1${}$\pmb-${}$1$, 
$a_8${}$\pmb-${}$a_9${}${}={}$$3${}$\pmb-${}$2$, 
$a_{10}${}$\pmb-${}$a_{11}${}${}={}$$6${}$\pmb-${}$6$, that is, $j=5$.
We reverse the order of the pairs $a_{2\ell}${}$\pmb-${}$a_{2\ell+1}$,
\dots $a_{2j-2}${}$\pmb-${}$a_{2j-1}$ and
append this reversed sequence to $S$, that is, we let
$$
S:=S,\text{$a_{2j-2}${}$\pmb-${}$a_{2j-1}$, \dots, $a_{2\ell}${}$\pmb-${}$a_{2\ell+1}$}.
$$
Next we let $\ell:=j$ and we repeat the Procedure.
\endroster

\medskip
In the centre of Figure~\FB, the sequence of pairs that is produced
by the above algorithm when applied to the alternating sequence on the top
of Figure~\FB\ is displayed. 
The bars indicate when a new loop in the Procedure was
started.

\medskip
Finally, these pairs are interpreted as segments as in the proof of Lemma~\UB, and
they are piled on each other as described in that proof. The bottom of Figure~\FB\
shows the resulting heap in our example. 

\medskip
It is clear from the construction that a heap constructed in the above
manner will 
satisfy Conditions (1') and (2'). The question is how one goes back
from the heap 
to the corresponding alternating sequence. In order to get an idea that this is
non-trivial, one should observe that the inverse mapping of the proof
of Lemma~\UB\ 
applied to the heap in Figure~\FB\ {\it does not\/} produce the
sequence of pairs 
in the centre of the figure.

\midinsert
$$
\Einheit.58cm
\Koordinatenachsen(10,9)(0,0)
\PfadDicke{.1pt}
\Pfad(0,1),1111111111\endPfad
\Pfad(0,2),1111111111\endPfad
\Pfad(0,3),1111111111\endPfad
\Pfad(0,4),1111111111\endPfad
\Pfad(0,5),1111111111\endPfad
\Pfad(0,6),1111111111\endPfad
\Pfad(0,7),1111111111\endPfad
\Pfad(0,8),1111111111\endPfad
\PfadDicke{2pt}
\Pfad(3,3),2\endPfad
\Pfad(2,1),22\endPfad
\Pfad(1,2),2\endPfad
\Pfad(4,4),2222\endPfad
\Pfad(5,4),222\endPfad
\Pfad(6,2),22\endPfad
\Pfad(4,2),2\endPfad
\Pfad(7,3),222\endPfad
\Pfad(8,5),2\endPfad
\Pfad(9,6),22\endPfad
\Kreis(3,3)
\Kreis(3,4)
\Kreis(1,6)
\Kreis(2,1)
\Kreis(2,3)
\Kreis(3,2)
\Kreis(4,4)
\Kreis(4,8)
\Kreis(1,1)
\Kreis(1,2)
\Kreis(1,3)
\Kreis(5,4)
\Kreis(5,7)
\Kreis(6,2)
\Kreis(6,4)
\Kreis(6,5)
\Kreis(4,2)
\Kreis(4,3)
\Kreis(7,3)
\Kreis(7,6)
\Kreis(8,5)
\Kreis(8,6)
\Kreis(9,6)
\Kreis(9,8)
\DickPunkt(0,4)
\DickPunkt(9,6)
\Label\l{1}(0,1)
\Label\l{2}(0,2)
\Label\l{3}(0,3)
\Label\l{4}(0,4)
\Label\l{5}(0,5)
\Label\l{6}(0,6)
\Label\l{7}(0,7)
\Label\l{8}(0,8)
\Label\r{3\kern5pt}(3,3)
\Label\o4(3,4)
\Label\o6(1,6)
\Label\u1(2,1)
\Label\o{3}(2,3)
\Label\u{2}(3,2)
\Label\r{4\kern5pt}(4,4)
\Label\o8(4,8)
\Label\r{4\kern5pt}(5,4)
\Label\o{3}(1,3)
\Label\u1(1,1)
\Label\l{\kern5pt2}(1,2)
\Label\o7(5,7)
\Label\u2(6,2)
\Label\r{4\kern5pt}(6,4)
\Label\o5(6,5)
\Label\u2(4,2)
\Label\r{3\kern5pt}(4,3)
\Label\u3(7,3)
\Label\o6(7,6)
\Label\u5(8,5)
\Label\o6(8,6)
\Label\u6(9,6)
\Label\o8(9,8)
\hbox{\hskip.3cm}
\raise.3cm\hbox{%
\SPfad(5,4),65565611112225\endSPfad
\SPfad(0,0),11221225655666\endSPfad
\SPfad(0,5),1256\endSPfad
\SPfad(4,3),2222256666\endSPfad
\SPfad(5,4),2225\endSPfad
\SPfad(6,4),25\endSPfad
\SPfad(6,2),1222256\endSPfad
\SPfad(7,4),1225\endSPfad
\SPfad(8,5),1222566\endSPfad
\Label\r{H_1}(1,6)
\Label\r{H_2}(2,0)
\Label\r{H_3}(2,7)
\Label\r{H_4}(5,7)
\Label\ru{H_5}(5,6)
\Label\u{H_6}(6,1)
\Label\u{H_7}(7,2)
\Label\u{H_8}(8,4)
\Label\u{H_9}(9,5)
}%
\hbox{\hskip6cm$\longrightarrow$\hskip.9cm}
\Koordinatenachsen(9,9)(0,0)
\PfadDicke{.1pt}
\Pfad(0,1),111111111\endPfad
\Pfad(0,2),111111111\endPfad
\Pfad(0,3),111111111\endPfad
\Pfad(0,4),111111111\endPfad
\Pfad(0,5),111111111\endPfad
\Pfad(0,6),111111111\endPfad
\Pfad(0,7),111111111\endPfad
\Pfad(0,8),111111111\endPfad
\PfadDicke{2pt}
\Pfad(1,3),2\endPfad
\Pfad(2,1),22\endPfad
\Pfad(3,2),2\endPfad
\Pfad(2,4),2222\endPfad
\Pfad(3,4),222\endPfad
\Pfad(4,2),22\endPfad
\Pfad(5,2),2\endPfad
\Pfad(6,3),222\endPfad
\Pfad(7,5),2\endPfad
\Pfad(8,6),22\endPfad
\Kreis(1,3)
\Kreis(1,4)
\Kreis(1,6)
\Kreis(2,1)
\Kreis(2,3)
\Kreis(2,4)
\Kreis(2,8)
\Kreis(3,1)
\Kreis(3,2)
\Kreis(3,3)
\Kreis(3,4)
\Kreis(3,7)
\Kreis(4,2)
\Kreis(4,4)
\Kreis(4,5)
\Kreis(5,2)
\Kreis(5,3)
\Kreis(6,2)
\Kreis(6,3)
\Kreis(6,6)
\Kreis(7,5)
\Kreis(7,6)
\Kreis(8,6)
\Kreis(8,8)
\DickPunkt(0,4)
\DickPunkt(8,6)
\Label\l{1}(0,1)
\Label\l{2}(0,2)
\Label\l{3}(0,3)
\Label\l{4}(0,4)
\Label\l{5}(0,5)
\Label\l{6}(0,6)
\Label\l{7}(0,7)
\Label\l{8}(0,8)
\Label\u3(1,3)
\Label\o4(1,4)
\Label\o6(1,6)
\Label\u1(2,1)
\Label\l{\kern5pt3}(2,3)
\Label\r{4\kern5pt}(2,4)
\Label\o8(2,8)
\Label\r{4\kern5pt}(3,4)
\Label\r{3\kern5pt}(3,3)
\Label\u1(3,1)
\Label\l{\kern5pt2}(3,2)
\Label\o7(3,7)
\Label\u2(4,2)
\Label\r4(4,4)
\Label\o5(4,5)
\Label\u2(5,2)
\Label\o3(5,3)
\Label\r{3\kern5pt}(6,3)
\Label\u2(6,2)
\Label\o6(6,6)
\Label\u5(7,5)
\Label\o6(7,6)
\Label\u6(8,6)
\Label\o8(8,8)
\Label\r{H_1}(0,7)
\Label\r{H_2}(3,0)
\Label\o{H_3}(1,8)
\Label\r{H_4}(3,8)
\Label\r{H_6}(6,1)
\hbox{\hskip.3cm}
\raise.3cm\hbox{%
\SPfad(3,4),666111252525\endSPfad
\SPfad(3,3),552566166112\endSPfad
\SPfad(0,5),1256\endSPfad
\SPfad(2,3),2222256666\endSPfad
\SPfad(3,4),2225\endSPfad
\SPfad(4,4),25\endSPfad
\SPfad(5,2),122225666\endSPfad
\SPfad(6,4),1225\endSPfad
\SPfad(7,5),1222566\endSPfad
\Label\o{H_5}(4,5)
\Label\ro{H_7}(6,2)
\Label\u{H_8}(7,4)
\Label\u{H_9}(8,5)
}%
\hbox{\hskip4.4cm}
$$
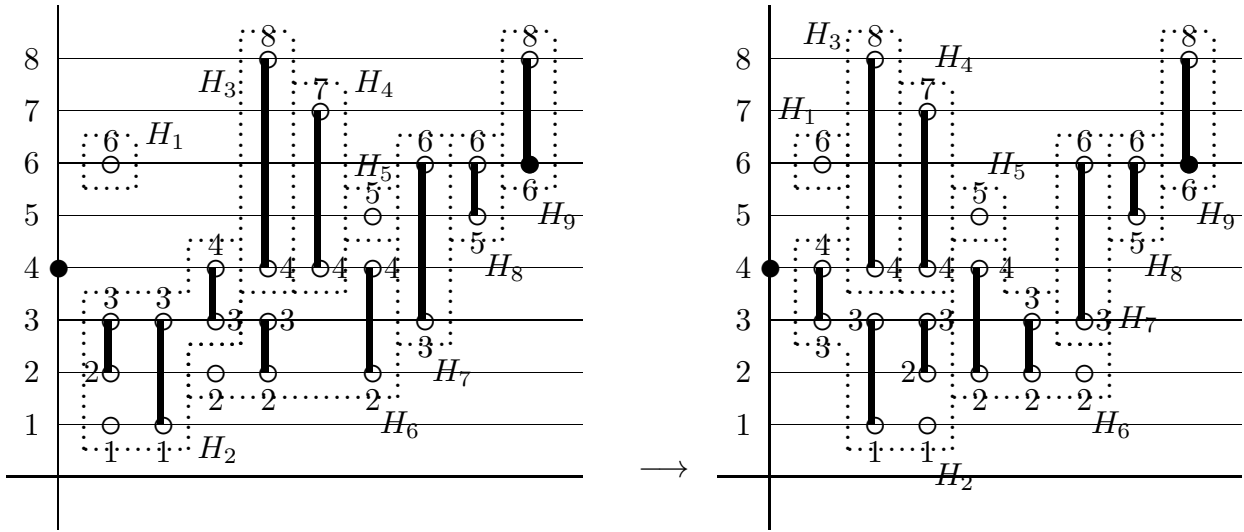
\centerline{\eightpoint Local modification of the heap}
\vskip7pt
\centerline{\eightpoint Figure \FC}
\endinsert

\medskip
Instead, given a heap $H'$ as described in the statement of the lemma,
we consider the 
top-most {\it minimal\/} segment, $\al$ say. If it should not be
contained in $[1,r-1]$, then 
we consider it as a subheap (consisting of a single segment), 
remove it from $H'$, and put it on a heap $H''$ that we are going to
build up. On the other hand,
if it is contained in $[1,r-1]$, then we push it to the right.
Thereby, due to the ``rule of the game" that segments cannot be moved
past one another 
if they ``physically" touch each other, other segments may also be
moved to the right. 
We collect all segments that are between $\al$ and the left-most moved
segment not contained 
in $[1,r-1]$, $\om$ say, into a subheap which we denote by $G$. 
We illustrate the construction in Figure~\FC. In the left part of the
figure, the heap from Figure~\FB\ is shown. (At this point, dotted
paths and labels should be ignored.) 
The top-most minimal
segment is $6${}$\pmb-${}$6$. It is not contained in $[1,r-1]=[1,3]$ and
therefore it is considered as a subheap consisting of a single
segment, indicated by dotted paths that surround it, and labelled by~$H_1$
in the left part of Figure~\FC.
As a subheap, it is subsequently removed, and it is put on a new heap~$H''$,
which we show in the right part of Figure~\FC. (At this point, the reader
should ignore the subheaps $H_2$, \dots, $H_9$; 
by now, $H''$ just consists of~$H_1$.)
After removal of $6${}$\pmb-${}$6$, the top-most minimal segment of~$H'$
is $3${}$\pmb-${}$2$. This is contained in $[1,r-1]=[1,3]$.
If we push it to the right then also the segments
$3${}$\pmb-${}$1$ and $4${}$\pmb-${}$3$ are moved to the right, where
$4${}$\pmb-${}$3$ is the left-most moved segment that is not contained in
$[1,r-1]=[1,3]$, that is, $\om${}${}={}${}$4${}$\pmb-${}$3$ in the example. 
Thus, the subheap~$G$ consists of the segments 
$3${}$\pmb-${}$2$, $3${}$\pmb-${}$1$, and $4${}$\pmb-${}$3$ at this point.
(It is true that also the segment $2${}$\pmb-${}$2$ is moved to the right
when $3${}$\pmb-${}$2$ is pushed to the right. However, $2${}$\pmb-${}$2$ is
not {\it between} $\al=3${}$\pmb-${}$2$ and $\om=4${}$\pmb-${}$3$ since
$4${}$\pmb-${}$3$ stays put when $2${}$\pmb-${}$2$ is pushed to the right.)
For the moment, we leave $G$ at its place.

We now consider the next-to-the-top-most minimal element, $\be$ say, 
and push it
to the right. If the left-most moved segment not contained in
$[1,r-1]$ is again~$\om$, then we put all moved elements between 
$\be$ and $\om$ also into the subheap $G$. In our running example, the
next-to-the-top-most minimal segment is $1${}$\pmb-${}$1$, that is, 
we have $\be=1${}$\pmb-${}$1$. If we move it to the right, then also
$1${}$\pmb-${}$1$, $3${}$\pmb-${}$1$, $4${}$\pmb-${}$3$ are moved to the
right. Since
$\om${}${}={}${}$4${}$\pmb-${}$3$ is among these elements, we augment
our subheap $G$ to $3${}$\pmb-${}$2$, $1${}$\pmb-${}$1$, 
$3${}$\pmb-${}$1$, $4${}$\pmb-${}$3$. See the left part of Figure~\FC, where
the subheap $G$ that we found so far is the subheap labelled by~$H_2$
on the left of Figure~\FC.
We continue in this manner until there are no minimal segments left, or
until a minimal segment is found that
does not move~$\om$. We then remove the subheap $G$ that we 
built up to this point from $H'$, reflect it into a vertical line, 
and put it in this reflected order on $H''$. In our example,
there are no more minimal segments to be considered, therefore
$G=H_2$ is removed from the heap $H'$ on the left of Figure~\FC, it
is reflected, and then put in this reflected order 
on the heap~$H''$ (at this point just
consisting of the segment $6${}$\pmb-${}$6$), see the right part of
Figure~\FC.

This process is repeated until $H'$ is emptied and a modified
heap $H''$ has been built up. In Figure~\FC, the subheaps
that were found during the construction and moved from $H'$
(on the left) to $H''$ (on the right) are surrounded by closed
dotted paths and labelled $H_1$, $H_2$, \dots, $H_9$.

\medskip
The final step then consists in applying the inverse mapping from the
proof of Lemma~\UB\ to~$H''$ in order to obtain an alternating
sequence with $2n+1$ elements, all of which between $1$ and~$k$,
and which starts with~$r$ and ends with~$s$. 
It can be seen that this yields the inverse mapping of the map
from alternating sequences to heaps that was described at the
beginning of this proof and illustrated in Figure~\FB.\quad \quad \qed
\enddemo

If we combine Lemma~\UC\ with (\BI), then we see that, for the proof
of Theorem~\UA, we need to
compute the generating function $\sum_T(-1)^{\vert T\vert}w(T)$
for all trivial heaps with segments in~$\Cal S_k$, that is, all
trivial heaps on $[1,k]$, and the analogous generating function
for trivial heaps whose segments are in $[1,r-1]\cup[s+1,k]$ and
with a maximal segment of the form $j-s$ for some $j\ge s$.
The former is accomplished in the lemma below, while the latter
is done in Lemma~\UE.

\proclaim{Lemma \UD}
Let $k$ be a non-negative integer.
The generating function $\sum_T (-1)^{\vert T\vert}w(T)$, 
where the sum is over all
trivial heaps $T$ of segments on $[1,k]$, is given by $(-1)^{k}U_{2k}(x/2)$.
\endproclaim

\demo{Proof}
We prove the claim by induction on $k$. The claim is certainly true for
$k=0$ --- in this case there is only the empty heap --- and for $k=1$
--- in this case there exist two heaps: the empty heap with weight~$1$,
and the heap consisting of the segment 1$-$1 contributing $-x^2$
to the signed sum of weights; indeed,
$(-1)^1U_2(x/2)=1-x^2$.

For the induction step, 
consider trivial heaps on $[1,k+1]$. There are three possibilities:
\roster
\item"(a)"the heap contains the segment $(k+1)${}$\pmb-${}$(k+1)$;
\item"(b)"the heap contains a segment $(k+1)${}$\pmb-${}$j$ with $j<k+1$;
\item"(c)"the heap does not contain a segment $(k+1)${}$\pmb-${}$j$ for
any~$j$.
\endroster
By induction hypothesis,
the contribution to the generating function of the heaps in
Case~(a) is $(-x^2)(-1)^kU_{2k}(x/2)$. Similarly, the contribution of
the heaps in Case~(c) is $(-1)^kU_{2k}(x/2)$. 
Finally, if we have a heap in Case~(b), then we may replace the
segment $(k+1)${}$\pmb-${}$j$ by the segment $k${}$\pmb-${}$j$. Thus, we
obtain a heap on $[1,k]$, and the correspondence is a bijection between
heaps in Case~(b) and heaps on $[1,k]$ with one segment of the form
$k${}$\pmb-${}$j$. Again by induction, the contribution of these heaps
to the generating function is $(-1)^kU_{2k}(x/2)-(-1)^{k-1}U_{2k-2}(x/2)$.
Hence, in total we obtain
$$
\align
(-x^2)(-1)^kU_{2k}(x/2)&+2(-1)^kU_{2k}(x/2)+(-1)^kU_{2k-2}(x/2)\\
&=
(-1)^{k+1}\big(
x^2U_{2k}(x/2)-2U_{2k}(x/2)-U_{2k-2}(x/2)
\big)\\
&=
(-1)^{k+1}U_{2k+2}(x/2),
\endalign
$$
where in the last line we used the recurrence formula (\BD) for 
Chebyshev polynomials iteratively.\quad \quad \qed
\enddemo

As we said above Lemma~\UD, we must also compute the generating function for
trivial heaps whose segments are
in $[1,r-1]\cup[s+1,k]$ and
with a maximal segment of the form $j-s$ for some $j\ge s$.
In such heaps, there is in particular no
maximal segment in $[1,r-1]\cup[j+1,k]$. Obviously, the same applies
when we remove the segment $j${}$\pmb-${}$s$, which has weight~$x^2$, 
from the top of this heap. The correspoonding generating function is given
in the lemma below.

\proclaim{Lemma \UE}
Let $r,s,k$ be positive integers with $1\le r\le s\le k$.
The sum of generating functions 
$$\sum_{j=s}^k x^2
\underset T\subseteq[1,r-1]\cup[j+1,k]
\to{\sum_{T\text{ trivial heap}}}
(-1)^{\vert T\vert}w(T)
\tag\BIb$$ 
is given by 
$$
(-1)^{k+r+s+1} x{U_{2r-2}(x/2)U_{2k+1-2s}(x/2)} .
$$
\endproclaim

\demo{Proof}
We may rewrite the sum in (\BIb) as
$$
-\underset T\not\subseteq[1,r-1]\cup[s+1,k]
\to{\underset T\subseteq[1,r-1]\cup[s,k]
\to{\sum_{T\text{ trivial heap}}}}
(-1)^{\vert T\vert}w(T),
\tag\BId$$
that is, we sum over all trivial heaps on $[1,r-1]\cup[s,k]$ that
contain a segment $j${}$\pmb-${}$s$ for some $j$ with $j\ge s$.

The parts of the trivial heap $T$ on $[1,r-1]$ and on $[s,k]$
are independent. Therefore the generating function is the product
of the corresponding generating functions. Thus, by Lemma~\UD,
we obtain
$$\multline
-(-1)^{r-1}U_{2r-2}(x/2)\left((-1)^{k-s+1}U_{2k-2s+2}(x/2)
-(-1)^{k-s}U_{2k-2s}(x/2)\right)\\
=
-(-1)^{r-1}U_{2r-2}(x/2)(-1)^{k-s+1}
xU_{2k-2s+1}(x/2),
\endmultline
\tag\BIa
$$
which implies the claimed expression.\quad \quad \qed
\enddemo

We have now collected all ingredients in order to establish Theorem~\UA.

\demo{Proof of Theorem \UA}
We start with the proof of (\BG) in the case where $r\le s$.

As we already announced before Lemma~\UD, we use Lemma~\UC\ to see
that the alternating sequences that we are interested in are in
bijection with the heaps that are described in the statement of
Lemma~\UC. If we give an alternating sequence with $2n+1$ elements
the weight $x^{2n}$, then this bijection is weight-preserving.

Let us now consider heaps as described in Lemma~\UC\ with a maximal
segment $j${}$\pmb-${}$s$. By Condition~(2') in Lemma~\UC, there is no
maximal segment in $[1,r-1]\cup[j+1,k]$. Obviously, the same applies
when we remove the segment $j${}$\pmb-${}$s$, which has weight~$x^2$, 
from the top of this heap. Consequently, by (\BI), the generating
function for these heaps is equal to
$$
x^2\times
\frac {\dsize\underset \text{segments}\subseteq [1,r-1]\cup[j+1,k]
\to{\sum_{T\text{ trivial heap}}}(-1)^{\vert T\vert}w(T)} 
{\dsize 
{\sum_{T\text{ trivial heap on }[1,k]}}(-1)^{\vert T\vert}w(T)} .
$$
In the end, we have to sum these expressions over $j=s,s+1,\dots,k$.
Lemma~\UD\ provides a closed form expression for the denominator in
the above fraction, while in Lemma~\UE\ the sum of the numerators is
computed. This yields the first expression on the right-hand side of
(\BG). In the case where $r=s$ one needs to add~1, the weight of the
empty heap. Symmetry in $r$ and~$s$ then also proves the third
expression on the right-hand side of~(\BG).

\medskip
The proof of (\BH) amounts to an exercise in summing Chebyshev polynomials.
We start with the case where $r>s$. 
By straightforward reasoning, we have
$$
\big\vert\Cal A_{2n+2}^{(k)}(r\to s)\big\vert
=\sum_{j=1}^s\big\vert\Cal A_{2n+1}^{(k)}(r\to j)\big\vert,
\tag\BIc$$
and therefore
$$
\sum_{n\ge0}\big\vert\Cal A_{2n+2}^{(k)}(r\to s)\big\vert x^{2n+1}=
\sum_{j=1}^s\sum_{n\ge0}\big\vert\Cal A_{2n+1}^{(k)}(r\to j)\big\vert
x^{2n+1}.$$
By the third case of (\BG), we obtain
$$
\align
\sum_{n\ge0}\big\vert\Cal A_{2n+2}^{(k)}(r\to s)\big\vert x^{2n+1}
&=
\sum_{j=1}^sx(-1)^{r+j+1}\frac {xU_{2j-2}(x/2)U_{2k+1-2r}(x/2)} 
{U_{2k}(x/2)}\\
&=(-1)^{r}\frac {xU_{2k+1-2r}(x/2)} 
{U_{2k}(x/2)}
\sum_{j=1}^s(-1)^{j-1}xU_{2j-2}(x/2)\\
&=(-1)^{r}\frac {xU_{2k+1-2r}(x/2)} 
{U_{2k}(x/2)}
\sum_{j=1}^s(-1)^{j-1}\left(U_{2j-1}(x/2)+U_{2j-3}(x/2)\right)\\
&=(-1)^{r}\frac {xU_{2k+1-2r}(x/2)} 
{U_{2k}(x/2)}
(-1)^{s-1}U_{2s-1}(x/2).
\endalign$$
Here we used the recurrence (\BD) to obtain the next-to-last line,
and a telescoping argument to arrive at the last line.
Thus, we obtain exactly the second expression on the right-hand side
of~(\BH). 

In the case where $r\le s$, we could do an analogous (but more
complicated) computation. However, there is a simpler argument,
using again heaps of segments. Namely, given a sequence
$r\le a_2\ge a_3\le a_4\ge \cdots \ge a_{2n+1}\le s$,
we apply the bijections from the proofs of Lemmas~\UB\ and~\UC\
to map the sequence to a heap of $n$~segments (ignoring the last
element~$s$ of the sequence). As is easy to see, this sets up a
bijection between $\Cal A_{2n+2}^{(k)}(r\to s)$ and heaps of $n$~segments on
$[1,k]$ whose maximal segments are not contained in $[1,r-1]\cup[s+1,k]$.
Hence, according to~(\BI), the generating function on the left-hand
side of~(\BH) equals
$$
\frac {\dsize\underset \text{segments}\subseteq [1,r-1]\cup[s+1,k]
\to{\sum_{T\text{ trivial heap on }[1,k]}}(-1)^{\vert T\vert}w(T)} 
{\dsize {\sum_{T\text{ trivial heap on }[1,k]}}(-1)^{\vert T\vert}w(T)} .
\tag\BIe
$$
Now application of Lemma~\UD\ immediately yields the first expression
on the right-hand side of~(\BH).

\medskip
This completes the proof of
the theorem.\quad \quad \qed
\enddemo

Next we show that
specialisations of Theorem~\UA\ may be used to find the
generating functions for the numbers of alternating sequences in which
first and last element are not specified.

\proclaim{Corollary \UG}
For all positive integers $k$, we have
$$
\sum_{n\ge1}\big\vert\Cal A_{2n-1}^{(k)}\big\vert x^{2n}
=-\frac {xU_{2k-1}(x/2)} {U_{2k}(x/2)}.
\tag\BEa
$$
\endproclaim

\demo{Proof}
The special case $r=s=1$ of (\BG) is equivalent with (\BEa).
Indeed, in any sequence
$$1\le a_2\ge a_3\le a_4\ge\dots\le a_{2n}\ge1$$
with entries between $1$ and $k$,
we may skip the $1$'s at the beginning and at the end, and then
replace every $a_i$ by $k+1-a_i$ in order to obtain
$$k+1-a_2\le k+1-a_3\ge k+1-a_4\le\dots\ge k+1-a_{2n},$$
which is an element of $\Cal A_{2n-1}^ {(k)}$.

\medskip
Alternatively, 
by using that $U_{2k-2}(x/2)=xU_{2k-1}(x/2)-U_{2k}(x/2)$ (cf\. (\BD)),
it can also be seen that the special case $r=s=k$ of (\BH)
is also equivalent with (\BEa).\quad \quad \qed
\enddemo

\proclaim{Corollary \UH}
For all positive integers $k$, we have
$$
\sum_{n\ge0}\big\vert\Cal A_{2n}^{(k)}\big\vert x^{2n+1}
=(-1)^k\frac {x} {U_{2k}(x/2)}.
\tag\BFa
$$
\endproclaim

\demo{Proof}
The special case $r=1$ and $s=k$ of (\BH) is equivalent
with (\BFa). Indeed, similarly to the first proof of Corollary~\UG, this is
seen by skipping the~$1$ at the beginning and the~$k$ at the end
of a sequence in $\Cal A_{2n+2}^{(k)}(1\to k)$ and then 
replacing each element $a_i$ in the sequence by $k+1-a_i$.

\medskip
Alternatively, 
with some additional work one can see that the special case $r=k$
and $s=1$ of (\BG) is also equivalent with (\BFa).\quad \quad \qed
\enddemo

In Exercise~3.66 of \cite{\StanAP}, Stanley considers 
the cumulative generating function
$$G_k(x)=1+\sum_{n\ge0}\big\vert\Cal A_{n}^{(k)}\big\vert x^{n+1}$$
for {\it all\/} bounded alternating sequences, regardless whether they
have odd or even length. In the solution section in \cite{\StanAP}, 
it is shown how to
derive a recursion formula for $G_k(x)$. Furthermore, it is pointed
out that \cite{\StanAZ, Ex.~3.2} gives an explicit formula for
$G_k(x)$. This formula is not very illuminating. Only specialists
may notice that, in the background, there lurk Chebyshev polynomials
of the second kind, which the corollary below reveals.
Our proof is independent of the above mentioned results as it is
based on our findings in Corollaries~\UG\ and~\UH.

\proclaim{Corollary \TA}
For $k\ge1$, we have
$$
G_{k}(x)=-\frac {U_{k-2}(x/2)+(-1)^kU_{k-3}(x/2)} 
{U_{k}(x/2)+(-1)^kU_{k-1}(x/2)}.
\tag\BA
$$
\endproclaim

\demo{Proof}
According to Corollaries~\UG\ and~\UH, we have
$$
G_k(x)=1+(-1)^k\frac {x} {U_{2k}(x/2)}
-\frac {xU_{2k-1}(x/2)} {U_{2k}(x/2)}.
\tag\BAa
$$

Let first $k$ be even.
We remember that, according to the definition (\BAb) of Chebyshev
polynomials of the second kind, we have 
$$U_n(X)=\frac {z^{n+1}-z^{-n-1}} {z-z^ {-1}},$$
where $z=e^{it}$ and $X=\cos t=\frac {1} {2}(z+z^{-1})$. 
We substitute these alternative expressions in~$z$
in the right-hand side of (\BAa) and simplify (using a computer).
The result is
$$
G_{k}(x)=-\frac {z^{2 k}-z^3} {z (z^{2 k+1}-1)}.
$$
If one does the same substitutions on the right-hand side of~(\BA)
(with $k$ even), then one obtains the same result.

The case of odd~$k$ can be treated similarly.\quad \quad \qed
\enddemo

\subhead 4. Numbers of bounded up-down paths ``with negative length'' and 
bounded alternating sequences\endsubhead
The strong similarity between the expressions in Theorems~\VA\ and~\UA\
has interesting consequences, which we reveal in this section.
It leads to first reciprocity relations of the kind that 
numbers of certain bounded up-down paths ``{\it of negative length}"
give numbers of certain bounded alternating sequences;
see Corollaries~\UF--\TC. 

The reader should recall that, given a rational power series 
$$f(x)=\frac {p(x)} {q(x)}=\sum_{n\ge0}f_nx^n,
\tag\BDb
$$
where the degree in $x$ of $p(x)$ is less than the degree of $q(x)$,
we have (cf\. \cite{\StanAP, Prop.~4.2.3})
$$
\sum_{n\ge1}f_{-n}x^n=-f(1/x),
\tag\BDa
$$
where the $f_{-n}$'s are defined by ``running" the linear recurrrence
for the sequence $\big(f_n\big)_{n\ge0}$ that results from~(\BDb) ``backwards".

\proclaim{Corollary \UF}
Let $n,k,r,s$ be positive integers with $1\le r,s\le k$.
The number $(-1)^{r+s}C_{-2n}^{(2k-1)}(2r-2\to 2s-2)$ equals the number of
sequences $r\le a_2\ge a_3\le a_4\ge \cdots \le a_{2n}\ge s$,
in which all $a_i$'s are integers between~$1$ and~$k$. Furthermore,
the number $(-1)^{r+s}C_{-2n+1}^{(2k-1)}(2r-2\to 2s-1)$ equals the number of
sequences $r\le a_2\ge a_3\le a_4\ge \cdots \ge a_{2n-1}\le s$,
in which all $a_i$'s are integers between~$1$ and~$k$.
\endproclaim

\demo{Proof}
By (\BGa) and (\BDa), we have
$$
\sum_{n\ge1}C_{-2n}^{(2k-1)}(2r-2\to 2s-2)\, x^{2n}
=\cases -\dfrac {xU_{2r-2}(x/2)\,U_{2k+1-2s}(x/2)} {U_{2k}(x/2)},
&\text{if }r\le s,\\
-\dfrac {xU_{2s-2}(x/2)\,U_{2k+1-2r}(x/2)}
{U_{2k}(x/2)},
&\text{if }r\ge s.
\endcases
$$
Comparison with Theorem~\UA\ completes the proof of
the first claim.

The second claim is established completely analogously, here using
(\BH) instead of (\BG).\quad \quad \qed
\enddemo

By comparison with appropriate special cases of Theorem~\VA, we obtain
more ``reciprocity relations''. 

\proclaim{Corollary \TB}
For positive integers $n$ and $k$,
the number $C_{-2n}^{(2k-1)}$ equals the number of
sequences $a_1\le a_2\ge a_3\le a_4\ge \dots \ge a_{2n-1}$
of positive integers with $a_i\le k$ for all~$i$.
\endproclaim

\demo{Proof}
This follows from applying (\BDa) to Theorem~\VA\ with $r=s=0$ and
$k$ replaced by $2k-1$, and
comparing the result with Corollary~\UG.\quad \quad \qed
\enddemo

For the statement of the next reciprocity relation, let us recall from the
beginning of Section~2 that
we write $D_{2n}^{(k)}$ instead of $C_{2n+k}^{(k)}(0\to k)$ for short.

\proclaim{Corollary \TC}
For integers $n$ and $k$ with $n\ge0$ and $k\ge1$,
the number $(-1)^{k+1}D_
{-2n-2k}^{(2k-1)}$ equals the number of
sequences $a_1\le a_2\ge a_3\le a_4\ge \dots \le a_{2n}$
of positive integers with $a_i\le k$ for all~$i$.
\endproclaim

\demo{Proof}
By Theorem~\VA\ with
$k$ replaced by $2k-1$, $r=0$, and $s=2k-1$, we have
$$
\sum_{n\ge0}D_{2n}^{(2k-1)} x^{2n+2k-1}
=\sum_{n\ge0}C_n^{(2k-1)}(0\to 2k-1)\, x^{n}
=\frac {1} {x\,U_{2k}(1/2x)},
$$
or, equivalently,
$$
\sum_{n\ge0}D_{2n}^ {(2k-1)}x^{2n}=\frac {1} {x^{2k}U_{2k}(1/2x)}.
$$
We now apply (\BDa) to get
$$
\sum_{n\ge1}D_{-2n}^ {(2k-1)}x^{2n}=-\frac {x^{2k}} {U_{2k}(x/2)}.
$$
Comparison with Corollary~\UH\ completes the proof.\quad \quad \qed
\enddemo

\subhead 5. A first reciprocity theorem\endsubhead
The main result in this section is a reciprocity law for
Hankel determinants of numbers of bounded Dyck paths, see Theorem~\TD.
As we show in Proposition~\TE\ and Theorem~\TF, in combinatorial terms
this reciprocity sets up a correspondence between families of
non-intersecting bounded Dyck paths and alternating tableaux
(see Section~11 for their formal definition)
of trapezoidal shape. 

\proclaim{Theorem \TD}
For all non-negative integers $n,k,m$, we have
$$
\det\left(C_{2n+2i+2j+4m-2}^{(2k+2m-1)}\right)_{0\le i,j\le k-1}
=
\det\left(C_{-2n-2i-2j}^{(2k+2m-1)}\right)_{0\le i,j\le m-1}.
\tag\AA
$$
\endproclaim

\remark{Remark}
Equation~(\BAc) is the special case $m=1$ of this theorem.
\endremark

Our strategy to prove Theorem~\TD\ is as follows.
By the Lindstr\"om--Gessel--Viennot theorem \cite{\LindAA, Lemma~1}
(see also \cite{\GeViAA, \GeViAB}), the determinant
on the left-hand side of (\AA) equals the number of 
families $(P_0,P_1,\dots,P_{k-1})$ of non-intersecting Dyck
paths of height at most $2k+2m-1$, where $P_i$ runs from
$(-2i,0)$ to $(2n+4m+2i-2,0)$, $i=0,1,\dots,k-1$. 
An example of such a family of
non-intersecting Dyck paths for $k=4$, $m=3$, and $n=3$ is 
shown in Figure~\FD. (At this point, circled points and attached
labels, as well as dotted lines should be ignored.)

\midinsert
$$
\Einheit.4cm
\Gitter(24,15)(-7,-1)
\Koordinatenachsen(24,15)(-7,-1)
\Pfad(-6,0),3333333343333344434444444444\endPfad
\Pfad(-4,0),333334333444334343444444\endPfad
\Pfad(-2,0),33343343443434334444\endPfad
\Pfad(0,0),3434334434343344\endPfad
\SPfad(1,-1),222222222222222\endSPfad
\SPfad(15,-1),222222222222222\endSPfad
\SPfad(2,8),33333\endSPfad
\SPfad(9,13),444\endSPfad
\PfadDicke{1pt}
\Pfad(-7,13),111111111111111111111111111111\endPfad
\DickPunkt(-6,0)
\DickPunkt(-4,0)
\DickPunkt(-2,0)
\DickPunkt(0,0)
\DickPunkt(16,0)
\DickPunkt(18,0)
\DickPunkt(20,0)
\DickPunkt(22,0)
\Kreis(2,6)
\Kreis(3,9)
\Kreis(4,2)
\Kreis(4,10)
\Kreis(5,5)
\Kreis(5,11)
\Kreis(6,0)
\Kreis(6,8)
\Kreis(6,12)
\Kreis(7,7)
\Kreis(7,9)
\Kreis(7,13)
\Kreis(8,6)
\Kreis(8,8)
\Kreis(8,10)
\Kreis(9,7)
\Kreis(9,9)
\Kreis(9,13)
\Kreis(10,4)
\Kreis(10,8)
\Kreis(10,12)
\Kreis(11,7)
\Kreis(11,11)
\Kreis(12,4)
\Kreis(12,8)
\Kreis(13,7)
\Kreis(14,0)
\Label\l{4\kern-5pt}(2,6)
\Label\l{5\kern-5pt}(3,9)
\Label\l{2\kern-5pt}(4,2)
\Label\l{6\kern-5pt}(4,10)
\Label\l{3\kern-5pt}(5,5)
\Label\l{6\kern-5pt}(5,11)
\Label\l{1\kern-5pt}(6,0)
\Label\l{5\kern-5pt}(6,8)
\Label\l{7\kern-5pt}(6,12)
\Label\l{4\kern-5pt}(7,7)
\Label\l{5\kern-5pt}(7,9)
\Label\l{7\kern-5pt}(7,13)
\Label\l{4\kern-5pt}(8,6)
\Label\l{5\kern-5pt}(8,8)
\Label\l{6\kern-5pt}(8,10)
\Label\l{4\kern-5pt}(9,7)
\Label\l{5\kern-5pt}(9,9)
\Label\l{7\kern-5pt}(9,13)
\Label\l{3\kern-5pt}(10,4)
\Label\l{5\kern-5pt}(10,8)
\Label\l{7\kern-5pt}(10,12)
\Label\l{4\kern-5pt}(11,7)
\Label\l{6\kern-5pt}(11,11)
\Label\l{3\kern-5pt}(12,4)
\Label\l{5\kern-5pt}(12,8)
\Label\l{4\kern-5pt}(13,7)
\Label\l{1\kern-5pt}(14,0)
%
\hskip6.3cm
$$
\centerline{\eightpoint A family of non-intersecting bounded Dyck paths}
\vskip7pt
\centerline{\eightpoint Figure \FD}
\endinsert

In Proposition~\TE\ below, we shall show that there is a bijection
between these families of non-intersecting Dyck paths and certain
trapezoidal arrays of alternating sequences. Subsequently, in Theorem~\TF,
it is shown that the number of these trapezoidal arrays is given by the
determinant on the right-hand side of~(\AA).
Hence, the combination of the above
observation with Proposition~\TE\ and Theorem~\TF\ will establish Theorem~\TD.

\proclaim{Proposition \TE}
Let $n$ be a non-negative integer and $k,m$ be positive integers.
There is a bijection between 
families $(P_0,P_1,\dots,P_{k-1})$ of non-intersecting Dyck
paths of height at most $2k+2m-1$, where $P_i$ runs from
$(-2i,0)$ to $(2n+4m+2i-2,0)$, $i=0,1,\dots,k-1$,
and trapezoidal arrays of integers of the form
{\sixpoint
$$
\matrix 
&&&&&&a_{1,2m-1}&\dots&a_{1,2n+2m-3}\\
&&&&&\iddots&&&&\hskip-.5cm\ddots\hskip1cm\\
&&&&a_{m-2,5}&\innerhdotsfor4\after\quad &a_{m-2,M-5}\\
&&a_{m-1,3}&a_{m-1,4}&a_{m-1,5}
&\innerhdotsfor4\after\quad &a_{m-1,M-5}&a_{m-1,M-4}&a_{m-1,M-3}\\
a_{m,1}&a_{m,2}&a_{m,3}&a_{m,4}&a_{m,5}
&\innerhdotsfor4\after\quad &
a_{m,M-5}&a_{m,M-4}&a_{m,M-3}&a_{m,M-2}&a_{m,M-1}\vphantom{\vrule depth6pt}\\
\endmatrix
\tag{\tenpoint\AB}
$$}%
where $M=2n+4m-4$, in which each row is alternating, that is,
$$
a_{i,2m-2i+1}\le a_{i,2m-2i+2}\ge a_{i,2m-2i+3}\le
a_{i,2m-2i+4}\ge\dots\ge a_{i,2n+2m+2i-5}
\tag\AC
$$
for all $i$,
and in which we have 
$$
1\le a_{i,j}\le k+m
\tag\AD
$$
and
$$
a_{i+1,2j}<a_{i,2j+1}>a_{i+1,2j+2}
\tag\AE
$$
for all $i$ and $j$.
In the particular case where $n=0$, the above must be
additionally complemented by the following interpretation:
row~$1$ is empty and all three entries in row~$2$
are bounded above by $k+m-1$.
\endproclaim

\demo{Proof}
The reader should consult Figure~\FD\ while reading the following
explanations. As mentioned earlier, the figure shows an example
of a family of non-intersecting Dyck paths as described in the
statement of the proposition
for $k=4$, $m=3$, and $n=3$.

Since the Dyck paths are non-intersecting, the path $P_i$
has to start with $2i+1$ up-steps, and it has to end with
$2i+1$ down-steps, for $i=0,1,\dots,k-1$. 
In other words, everything is uniquely
determined left of the abscissa $x=1$ and right of the
abscissa $x=2n+4m-3$. (These abscissa are marked by dotted lines
in the figure.)
Only ``in between" there is (some) ``freedom".

Now, for all $i$ and $j$, 
we mark all lattice points $(2i-1,2j-1)$ and all lattice
points $(2i,2j)$ in the region
$$
\{(x,y):0\le y\le 2k+2m-1,\ 2\le x\le 2n+4m-4,\ y\le x+2k-2,\
x+y\le 2n+4m+2k-4\}
\tag\AEa
$$
that are not occupied by any of the Dyck paths.
(This region is the trapezoidal region bounded below by the $x$-axis,
bounded above by the horizontal line $y=2k+2m-1$, bounded on the
left by the abscissa $x=2$ and
the diagonal line starting in the left-most starting point
of the Dyck paths, and bounded on the right by the abscissa $2n+2m-4$
and the anti-diagonal line
ending in the right-most end point of the Dyck paths.)
In the figure, these lattice points are marked by circles.

Each marked point $X$ acquires a label. The label is equal to 
$1$ plus the
number of paths and marked points that are below $X$ and have the
same abscissa as~$X$. Equivalently, the label
of a marked point $X=(x,y)$ equals $\cl{(y+1)/2}$.
The figure also
shows the labels of the marked points.

Finally, for $i=2,3,\dots,2n+4m-4$, we read the labels of the marked
points with abscissa~$i$ and put them into a column, and then
we concatenate the columns into a bottom-justified array of
trapezoidal shape as in (\AB). In this way, for the example in Figure~\FD,
we obtain the array
$$
\matrix 
 & & & &7&7&6&7&7\\
 & &6&6&5&5&5&5&5&6&5\\
4&5&2&3&1&4&4&4&3&4&3&4&1
\endmatrix
\tag\AEb$$

It is not difficult to see that (in general) the obtained array satisfies
the constraints in (\AC)--(\AE). Moreover, in the case where $n=0$,
the obtained array satisfies the additional properties that are
described at the end of the statement of the proposition.

Conversely, given an array of the form (\AB) which satisfies
(\AC)--(\AE), if one reverses the previously described construction
then one sees that this array determines a unique family 
 of non-intersecting Dyck paths as described in the
statement of the proposition.

To complete the proof, it should be noted that, if $n=0$,
then the construction described here indeed goes with the additional
restrictions mentioned at the end of the statement of the
proposition.\quad \quad \qed
\enddemo

\proclaim{Theorem \TF}
Let $n$ be a non-negative integer and $k,m$ be positive integers.
The number of trapezoidal arrays of integers of the form {\rm(\AB)}
that satisfy {\rm(\AC)--(\AE)}, and the additional restrictions
mentioned in the statement of Proposition~{\rm\TE} in the case where
$n=0$, is equal to
$$
\det\left(C_{-2n-2i-2j}^{(2k+2m-1)}\right)_{0\le i,j\le m-1}.
\tag\AF
$$
\endproclaim

\demo{Proof}
It should be noted that, by Corollary~\TB, the determinant in
(\AF) can be alternatively written as
$$
\det\left(\big\vert\Cal A_{2n+2i+2j-1}^{(k+m)}\big\vert\right)_{0\le i,j\le m-1}.
\tag\AFb
$$
To be in line with $C_0^{(2k+2m-1)}=1$, 
we have to {\it define} that $\Cal A_{-1}^{(k+m)}$
consists of just one element, namely the empty sequence.

We shall interpret the determinant in (\AFb) in terms of non-intersecting
lattice paths. The underlying directed graph for these paths is
the graph $\Cal G_{k+m}$, which by definition is the graph
with vertices being the integer lattice points $(x,y)$,
where $x\ge0$ and $1\le y\le k+m$, and with directed edges in the set
$$\multline
\{(2i,j+1)\to(2i,j):i\ge0\text{ and }1\le j\le k+m-1\}\\
\cup
\{(2i+1,j)\to(2i+1,j+1):i\ge0\text{ and }1\le j\le k+m-1\}\\
\cup
\{(i,j)\to(i+1,j):i\ge0\text{ and }1\le j\le k+m\}.
\endmultline$$
See Figure~\FDa\ for a portion of the directed graph $\Cal G_{k+m}$
for $k=4$ and $m=3$.

\midinsert
$$
\Gitter(10,8)(0,0)
\Koordinatenachsen(10,8)(0,0)
\Vektor(0,-1)1(0,7)
\Vektor(0,-1)1(0,6)
\Vektor(0,-1)1(0,5)
\Vektor(0,-1)1(0,4)
\Vektor(0,-1)1(0,3)
\Vektor(0,-1)1(0,2)
%
\Vektor(0,-1)1(2,7)
\Vektor(0,-1)1(2,6)
\Vektor(0,-1)1(2,5)
\Vektor(0,-1)1(2,4)
\Vektor(0,-1)1(2,3)
\Vektor(0,-1)1(2,2)
%
\Vektor(0,-1)1(4,7)
\Vektor(0,-1)1(4,6)
\Vektor(0,-1)1(4,5)
\Vektor(0,-1)1(4,4)
\Vektor(0,-1)1(4,3)
\Vektor(0,-1)1(4,2)
%
\Vektor(0,-1)1(6,7)
\Vektor(0,-1)1(6,6)
\Vektor(0,-1)1(6,5)
\Vektor(0,-1)1(6,4)
\Vektor(0,-1)1(6,3)
\Vektor(0,-1)1(6,2)
%
\Vektor(0,-1)1(8,7)
\Vektor(0,-1)1(8,6)
\Vektor(0,-1)1(8,5)
\Vektor(0,-1)1(8,4)
\Vektor(0,-1)1(8,3)
\Vektor(0,-1)1(8,2)
%
\Vektor(0,1)1(1,6)
\Vektor(0,1)1(1,5)
\Vektor(0,1)1(1,4)
\Vektor(0,1)1(1,3)
\Vektor(0,1)1(1,2)
\Vektor(0,1)1(1,1)
%
\Vektor(0,1)1(3,6)
\Vektor(0,1)1(3,5)
\Vektor(0,1)1(3,4)
\Vektor(0,1)1(3,3)
\Vektor(0,1)1(3,2)
\Vektor(0,1)1(3,1)
%
\Vektor(0,1)1(5,6)
\Vektor(0,1)1(5,5)
\Vektor(0,1)1(5,4)
\Vektor(0,1)1(5,3)
\Vektor(0,1)1(5,2)
\Vektor(0,1)1(5,1)
%
\Vektor(0,1)1(7,6)
\Vektor(0,1)1(7,5)
\Vektor(0,1)1(7,4)
\Vektor(0,1)1(7,3)
\Vektor(0,1)1(7,2)
\Vektor(0,1)1(7,1)
%
\Vektor(0,1)1(9,6)
\Vektor(0,1)1(9,5)
\Vektor(0,1)1(9,4)
\Vektor(0,1)1(9,3)
\Vektor(0,1)1(9,2)
\Vektor(0,1)1(9,1)
%
\Vektor(1,0)1(1,7)
\Vektor(1,0)1(1,6)
\Vektor(1,0)1(1,5)
\Vektor(1,0)1(1,4)
\Vektor(1,0)1(1,3)
\Vektor(1,0)1(1,2)
\Vektor(1,0)1(1,1)
\Vektor(1,0)1(3,7)
\Vektor(1,0)1(3,6)
\Vektor(1,0)1(3,5)
\Vektor(1,0)1(3,4)
\Vektor(1,0)1(3,3)
\Vektor(1,0)1(3,2)
\Vektor(1,0)1(3,1)
\Vektor(1,0)1(5,7)
\Vektor(1,0)1(5,6)
\Vektor(1,0)1(5,5)
\Vektor(1,0)1(5,4)
\Vektor(1,0)1(5,3)
\Vektor(1,0)1(5,2)
\Vektor(1,0)1(5,1)
\Vektor(1,0)1(7,7)
\Vektor(1,0)1(7,6)
\Vektor(1,0)1(7,5)
\Vektor(1,0)1(7,4)
\Vektor(1,0)1(7,3)
\Vektor(1,0)1(7,2)
\Vektor(1,0)1(7,1)
\Vektor(1,0)1(0,7)
\Vektor(1,0)1(0,6)
\Vektor(1,0)1(0,5)
\Vektor(1,0)1(0,4)
\Vektor(1,0)1(0,3)
\Vektor(1,0)1(0,2)
\Vektor(1,0)1(0,1)
\Vektor(1,0)1(2,7)
\Vektor(1,0)1(2,6)
\Vektor(1,0)1(2,5)
\Vektor(1,0)1(2,4)
\Vektor(1,0)1(2,3)
\Vektor(1,0)1(2,2)
\Vektor(1,0)1(2,1)
\Vektor(1,0)1(4,7)
\Vektor(1,0)1(4,6)
\Vektor(1,0)1(4,5)
\Vektor(1,0)1(4,4)
\Vektor(1,0)1(4,3)
\Vektor(1,0)1(4,2)
\Vektor(1,0)1(4,1)
\Vektor(1,0)1(6,7)
\Vektor(1,0)1(6,6)
\Vektor(1,0)1(6,5)
\Vektor(1,0)1(6,4)
\Vektor(1,0)1(6,3)
\Vektor(1,0)1(6,2)
\Vektor(1,0)1(6,1)
\Vektor(1,0)1(8,7)
\Vektor(1,0)1(8,6)
\Vektor(1,0)1(8,5)
\Vektor(1,0)1(8,4)
\Vektor(1,0)1(8,3)
\Vektor(1,0)1(8,2)
\Vektor(1,0)1(8,1)
\hskip4.5cm
$$
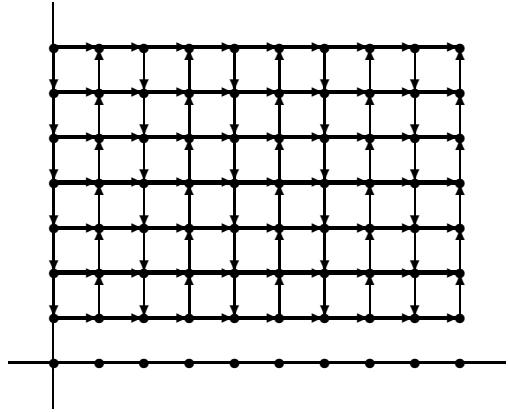
\centerline{\eightpoint The directed graph $\Cal G_7$}
\vskip7pt
\centerline{\eightpoint Figure \FDa}
\endinsert

Now let first $n\ge1$.
We claim that the determinant in (\AFb) equals the number of families
$(P_0,P_1,\dots,P_{m-1})$ of non-intersecting lattice paths in the
directed graph $\Cal G_{k+m}$, where $P_i$ runs from $(2m-2i-2,k+m)$ to
$(2n+2m+2i-3,k+m)$. Figure~\FDb\ shows an example of such a path family
for $k=4$, $m=3$, and $n=3$. (The numbers should be ignored
at this point.)

\midinsert
$$
\Gitter(15,8)(0,0)
\Koordinatenachsen(15,8)(0,0)
\Pfad(0,7),666121666121661222111612161216661222222\endPfad
\Pfad(2,7),611611111216122\endPfad
\Pfad(4,7),1161211\endPfad
\DickPunkt(0,7)
\DickPunkt(2,7)
\DickPunkt(4,7)
\DickPunkt(9,7)
\DickPunkt(11,7)
\DickPunkt(13,7)
\Label\r{P_0}(9,7)
\Label\r{P_1}(11,6)
\Label\r{P_2}(13,6)
\Label\lo{\text{\it4}}(1,4)
\Label\ro{\text{\it5}}(1,5)
\Label\lo{\text{\it2}}(3,2)
\Label\lo{\text{\it6}}(3,6)
\Label\ro{\text{\it3}}(3,3)
\Label\ro{\text{\it6}}(3,6)
\Label\lo{\text{\it1}}(5,1)
\Label\lo{\text{\it5}}(5,5)
\Label\lo{\text{\it7}}(5,7)
\Label\ro{\text{\it4}}(5,4)
\Label\ro{\text{\it5}}(5,5)
\Label\ro{\text{\it7}}(5,7)
\Label\lo{\text{\it4}}(7,4)
\Label\lo{\text{\it5}}(7,5)
\Label\lo{\text{\it6}}(7,6)
\Label\ro{\text{\it4}}(7,4)
\Label\ro{\text{\it5}}(7,5)
\Label\ro{\text{\it7}}(7,7)
\Label\lo{\text{\it3}}(9,3)
\Label\lo{\text{\it5}}(9,5)
\Label\lo{\text{\it7}}(9,7)
\Label\ro{\text{\it4}}(9,4)
\Label\ro{\text{\it6}}(9,6)
\Label\lo{\text{\it3}}(11,3)
\Label\lo{\text{\it5}}(11,5)
\Label\ro{\text{\it4}}(11,4)
\Label\lo{\text{\it1}}(13,1)
\hskip7cm
$$
\centerline{\eightpoint A family of non-intersecting paths in the
directed graph $\Cal G_7
$}
\vskip7pt
\centerline{\eightpoint Figure \FDb}
\endinsert

By the Lindstr\"om--Gessel--Viennot theorem \cite{\LindAA, Lemma~1},
to verify the claim it suffices to show that the number of paths
from $(2x,k+m)$ to $(2x+2s-1,k+m)$ in the graph $\Cal G_{k+m}$
is equal to $\vert\Cal A_{2s-1}\vert$, the number of alternating
sequences of integers $a_1\le a_2\ge a_3\le\dots \ge a_{2s-1}$
with $1\le a_i\le k+m$ for all~$i$. This is indeed easy to see
if one labels each horizontal step from $(i,j)\to(i+1,j)$
by~$j$, for all~$i$ and~$j$.
See Figure~\FDb\ for these labels.
If one then reads labels along a path from left to right, then one
reads an alternating sequence, and this correspondence sets up a
bijection between paths from $(2x,k+m)$ to $(2x+2s-1,k+m)$ in the
graph $\Cal G_{k+m}$ and alternating sequences in 
$\Cal A^{(k+m)}_{2s-1}$. For example, the path~$P_2$ in Figure~\FDb\
corresponds to the alternating sequence
$4\le 5\ge2\le3\ge1\le4\ge4\le4\ge3\le4\ge3\le4\ge1$.

\medskip
In view of the above considerations, the proof would be complete if
there is a bijection between families
$(P_0,P_1,\dots,P_{m-1})$ of non-intersecting lattice paths in 
$\Cal G_{k+m}$, where $P_i$ runs from $(2m-2i-2,k+m)$ to
$(2n+2m+2i-3,k+m)$, and arrays of integers as described in
Proposition~\TE. This bijection is easily set up: the $i$-th row in an array
of the form~(\AB) is translated into a path $P_{i-1}$ in $\Cal
G_{k+m}$, $i=1,2,\dots,m$; more precisely, for $i=1,2,\dots,m$,
the entry $a_{i,j}$ (if it exists)
is translated into the horizontal step
$(j-1,a_{i,j})\to(j,a_{i,j})$, for all $j\ge1$, and these
horizontal steps are then connected by vertical up- and down-steps
to form a path in $\Cal G_{k+m}$. This correspondence is 
illustrated in Figure~\FDb, which shows in fact the family of
non-intersecting paths corresponding to the array in~(\AEb).

\medskip
Finally we address the special case where $n=0$. In that case, 
when we attempt to carry through the same programme, we face the
problem that there is no family 
$(P_0,P_1,\dots,P_{m-1})$ of non-intersecting lattice paths in 
$\Cal G_{k+m}$, where $P_i$ runs from $(2m-2i-2,k+m)$ to
$(2m+2j-3,k+m)$, because $P_0$ would have to run from $(2m-2,k+m)$
to $(2m-3,k+m)$. We remedy the situation by extending $\Cal G_{k+m}$
by the additional ``artificial" (backwards) edge $(2m-2,k+m)\to(2m-3,k+m)$.
It is not difficult to see that, with this modified graph, 
the earlier arguments now also apply to this case.

This completes the proof of the theorem.\quad \quad \qed
\enddemo

\subhead 6. A second reciprocity theorem\endsubhead
In this section, the main result is a reciprocity law for
Toeplitz determinants of numbers of bounded up-down paths connecting the
$x$-axis with the upper boundary, see Theorem~\TG.
As we show in Proposition~\TH\ and Theorem~\TI, in combinatorial terms
this reciprocity sets up a correspondence between families of
non-intersecting bounded up-down paths connecting the $x$-axis with
the upper boundary and alternating tableaux
(see Section~11 for their formal definition)
of rhomboidal shape. 

\proclaim{Theorem \TG}
For all non-negative integers $n$ and positive integers $k,m$, we have
$$
\det\left(D_{2n+2j-2i}^{(2k+2m-1)}\right)_{0\le i,j\le k-1}
=
(-1)^{km}\det\left(D_{-2n-2j+2i-2k-2m}^{(2k+2m-1)}\right)_{0\le i,j\le m-1}.
\tag\CA
$$
\endproclaim

For the proof of the theorem, we proceed in a similar manner as for
the proof of Theorem~\TD. Namely,
by the Lindstr\"om--Gessel--Viennot theorem \cite{\LindAA, Lemma~1}, 
the determinant
on the left-hand side of (\CA) equals the number of 
families $(P_0,P_1,\dots,P_{k-1})$ of non-intersecting up-down
paths of height at most $2k+2m-1$ that do not pass below the $x$-axis, 
where $P_i$ runs from
$(2i,0)$ to $(2n+2m+2k+2i-1,2k+2m-1)$, $i=0,1,\dots,k-1$. 
An example of such a family of
non-intersecting up-down paths for $k=4$, $m=3$, and $n=5$ is 
shown in Figure~\FE. (At this point, circled points and attached
labels, as well as dotted lines should be ignored.)

\midinsert
$$
\Einheit.4cm
\Gitter(31,15)(-1,-1)
\Koordinatenachsen(31,15)(-1,-1)
\Pfad(0,0),33333333433333443433343\endPfad
\Pfad(2,0),33333433344433333343333\endPfad
\Pfad(4,0),33343343443333343333333\endPfad
\Pfad(6,0),34343344333433333333333\endPfad
\SPfad(0,0),3333333333333\endSPfad
\SPfad(16,0),3333333333333\endSPfad
\SPfad(7,-1),222222222222222\endSPfad
\SPfad(22,-1),222222222222222\endSPfad
\PfadDicke{1pt}
\Pfad(-1,13),1111111111111111111111111111111\endPfad
\DickPunkt(0,0)
\DickPunkt(2,0)
\DickPunkt(4,0)
\DickPunkt(6,0)
\DickPunkt(23,13)
\DickPunkt(25,13)
\DickPunkt(27,13)
\DickPunkt(29,13)
\Kreis(8,6)
\Kreis(9,9)
\Kreis(10,2)
\Kreis(10,10)
\Kreis(11,5)
\Kreis(11,11)
\Kreis(12,0)
\Kreis(12,8)
\Kreis(12,12)
\Kreis(13,7)
\Kreis(13,9)
\Kreis(13,13)
\Kreis(14,6)
\Kreis(14,8)
\Kreis(14,10)
\Kreis(15,7)
\Kreis(15,9)
\Kreis(15,13)
\Kreis(16,0)
\Kreis(16,8)
\Kreis(16,12)
\Kreis(17,1)
\Kreis(17,9)
\Kreis(17,13)
\Kreis(18,4)
\Kreis(18,12)
\Kreis(19,5)
\Kreis(19,13)
\Kreis(20,8)
\Kreis(21,11)
\Label\l{4\kern-5pt}(8,6)
\Label\l{5\kern-5pt}(9,9)
\Label\l{2\kern-5pt}(10,2)
\Label\l{6\kern-5pt}(10,10)
\Label\l{3\kern-5pt}(11,5)
\Label\l{6\kern-5pt}(11,11)
\Label\l{1\kern-5pt}(12,0)
\Label\l{5\kern-5pt}(12,8)
\Label\l{7\kern-5pt}(12,12)
\Label\l{4\kern-5pt}(13,7)
\Label\l{5\kern-5pt}(13,9)
\Label\l{7\kern-5pt}(13,13)
\Label\l{4\kern-5pt}(14,6)
\Label\l{5\kern-5pt}(14,8)
\Label\l{6\kern-5pt}(14,10)
\Label\l{4\kern-5pt}(15,7)
\Label\l{5\kern-5pt}(15,9)
\Label\l{7\kern-5pt}(15,13)
\Label\l{1\kern-5pt}(16,0)
\Label\l{5\kern-5pt}(16,8)
\Label\l{7\kern-5pt}(16,12)
\Label\l{1\kern-5pt}(17,1)
\Label\l{5\kern-5pt}(17,9)
\Label\l{7\kern-5pt}(17,13)
\Label\l{3\kern-5pt}(18,4)
\Label\l{7\kern-5pt}(18,12)
\Label\l{3\kern-5pt}(19,5)
\Label\l{7\kern-5pt}(19,13)
\Label\l{5\kern-5pt}(20,8)
\Label\l{6\kern-5pt}(21,11)
%
\hskip11.8cm
$$
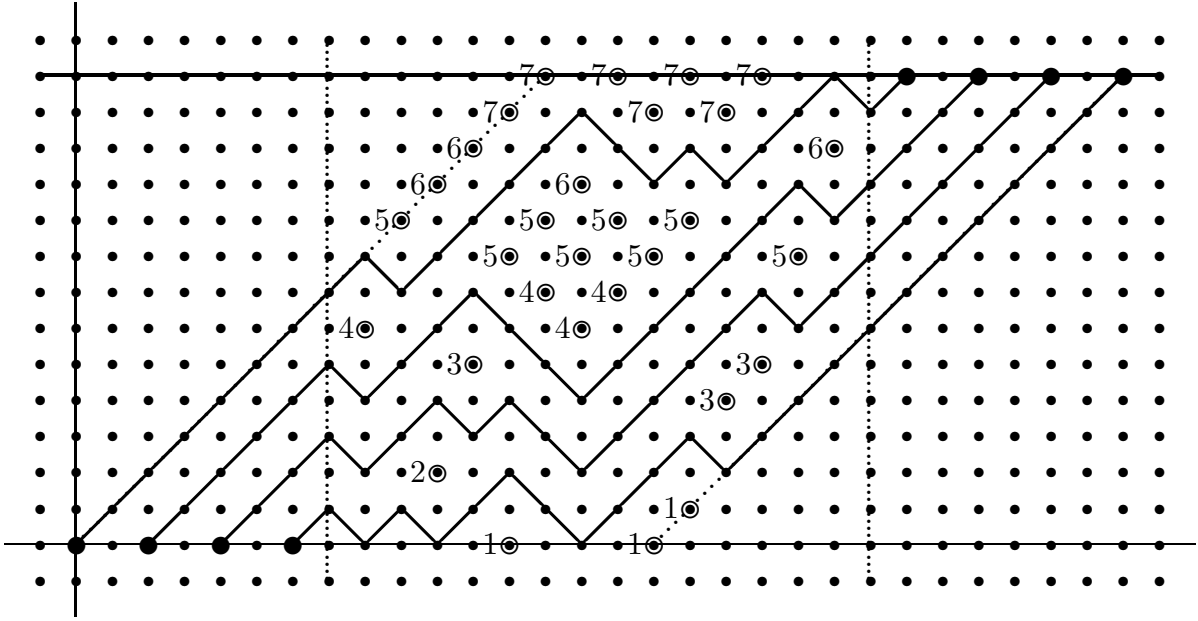
\centerline{\eightpoint A family of non-intersecting up-down paths
starting at the $x$-axis and ending at the upper bound}
\vskip7pt
\centerline{\eightpoint Figure \FE}
\endinsert

In Proposition~\TH\ below, we shall show that there is a bijection
between these families of non-intersecting up-down paths and certain
rhomboidal arrays of alternating sequences. Then, in Theorem~\TI,
it is shown that the number of these rhomboidal arrays is given by the
determinant on the right-hand side of~(\CA).
Hence, the combination of the above
observation with Proposition~\TH\ and Theorem~\TI\ will establish Theorem~\TG.

\proclaim{Proposition \TH}
Let $n$ be a non-negative integer and $k,m$ be positive integers.
There is a bijection between 
families $(P_0,P_1,\dots,P_{k-1})$ of non-intersecting up-down
paths of height at most $2k+2m-1$ that do not pass below the $x$-axis, 
where $P_i$ runs from
$(2i,0)$ to $(2n+2m+2k+2i-1,2k+2m-1)$, $i=0,1,\dots, k-1$, 
and rhomboidal arrays of integers of the form
{\sixpoint
$$
\matrix 
&&&&&&a_{1,2m-1}&\innerhdotsfor3\after\quad 
&a_{1,2m+2n-4}&a_{1,2m+2n-3}&a_{1,2m+2n-2}\\
&&&&a_{2,2m-3}&a_{2,2m-2}&a_{2,2m-1}&\dots
&a_{2,2m+2n-6}&a_{2,2m+2n-5}&a_{2,2m+2n-4}\\
&&a_{3,2m-5}&a_{3,2m-4}&a_{3,2m-3}&\innerhdotsfor3\after\quad &a_{3,2m+2n-6}\\
&\iddots&&&&&&\iddots\hskip-1cm\\
a_{m,1}&a_{m,2}&a_{m,3}
&\innerhdotsfor4\after\quad &
a_{m,2n}&
\endmatrix
\tag{\tenpoint\CB}
$$}%
\vskip6pt
\noindent
in which each row is alternating, that is,
$$
a_{i,2m-2i+1}\le a_{i,2m-2i+2}\ge a_{i,2m-2i+3}\le
a_{i,2m-2i+4}\ge\dots\le a_{i,2n+2m-2i}
\tag\CC
$$
for all $i$,
and in which we have 
$$
1\le a_{i,j}\le k+m
\tag\CD
$$
and
$$
a_{i+1,2j}<a_{i,2j+1}>a_{i+1,2j+2}
\tag\CE
$$
for all $i$ and $j$. In the case where $n=0$, the array in 
{\rm(\CB)} as to be interpreted as the empty array.
\endproclaim

\demo{Proof}
The reader should consult Figure~\FE\ while reading the following
explanations. As mentioned earlier, the figure shows an example
of a family of non-intersecting up-down paths as described in the
statement of the proposition
for $k=4$, $m=3$, and $n=5$.

Since the paths are non-intersecting, the path $P_i$
has to start with $2k-2i-1$ up-steps, and it has to end with
$2i+1$ up-steps, for $i=0,1,\dots,k-1$. 
In other words, everything is uniquely
determined left of the abscissa $x=2k-1$ and right of the
abscissa $x=2n+2m+2k-2$. (These abscissa are marked by dotted vertical lines
in the figure.)
Only ``in between" there is (some) ``freedom".

Now, for all $i$ and $j$, 
we mark all lattice points $(2i-1,2j-1)$ and all lattice
points $(2i,2j)$ in the region
$$
\{(x,y):0\le y\le 2k+2m-1,\ 2k\le x\le 2n+2m+2k-3,\ y\le x\le
y+2n+2k-2\}
$$
that are not occupied by any of the paths.
(This region is the rhomboidal region bounded below by the $x$-axis,
bounded above by the horizontal line $y=2k+2m-1$, bounded on the
left by the abcissa $x=2k$ and 
the diagonal line starting in the left-most starting point
of the paths, and bounded on the right by the abscissa $x=2n+2m+2k-3$ and
the diagonal line
ending in the right-most end point of the paths.)
In the figure, these lattice points are marked by circles.

As before, each marked point $X$  acquires a label, given by
$1$ plus the
number of paths and marked points that are below $X$ and have the
same abscissa as~$X$. Again, equivalently, the label
of a marked point $X=(x,y)$ equals $\cl{(y+1)/2}$. The figure 
shows the labels of these marked points. 

Finally, for $i=2k,2k+1,\dots,2n+2m+2k-3$, we read the labels of the marked
points with abscissa~$i$ and put them into a column, and then
we concatenate the columns into an array of
rhomboidal shape as in (\CB). In this way, for the example in Figure~\FE,
we obtain the array
$$
\matrix 
 & & & &7&7&6&7&7&7&7&7&5&6\\
 & &6&6&5&5&5&5&5&5&3&3\\
4&5&2&3&1&4&4&4&1&1
\endmatrix
\tag\CEa$$

It is not difficult to see that (in general) the obtained array satisfies
the constraints in (\CC)--(\CE). 

Conversely, given an array of the form (\CB) which satisfies
(\CC)--(\CE), if one reverses the previously described construction
then one sees that this array determines a unique family 
 of non-intersecting up-down paths as described in the
statement of the proposition.\quad \quad \qed
\enddemo

\proclaim{Theorem \TI}
Let $n$ be a non-negative integer and $k,m$ be positive integers.
The number of rhomboidal arrays of integers of the form {\rm(\CB)}
that satisfy {\rm(\CC)--(\CE)} is equal to
$$
\align
(-1)^{km}\det\left(D_{-2n-2i+2j-2k-2m}^{(2k+2m-1)}\right)_{0\le i,j\le m-1}\\
&\kern-1cm
=(-1)^{km+\binom {m}2}
\det\left(D_{-2n-2j-2i-2k-2}^{(2k+2m-1)}\right)_{0\le i,j\le m-1}\\
&\kern-1cm
=(-1)^{km}\det\left(D_{-2n-2j+2i-2k-2m}^{(2k+2m-1)}\right)_{0\le i,j\le m-1}.\\
\tag\CF
\endalign
$$
\endproclaim

\demo{Sketch of proof of Theorem \TI}
Obviously, the equality between the above three determinants 
results from reversing the order of rows and/or columns.

Next, it should be observed that, by Corollary~\TC, the first determinant in
(\CF) can be alternatively written as
$$
\det\left(\big\vert\Cal A_{2n+2i-2j}^{(k+m)}\big\vert\right)_{0\le i,j\le m-1}.
$$

We claim that this determinant equals the number of families
$(P_0,P_1,\dots,P_{m-1})$ of non-intersecting lattice paths in the
directed graph $\Cal G_{k+m}$ (see the proof of Theorem~\TF\ and
in particular Figure~\FDa), where $P_i$ runs from $(2i,k+m)$ to
$(2n+2i,0)$. Figure~\FDc\ shows an example of such a path family
for $k=4$, $m=3$, and $n=5$. (The numbers should be ignored
at this point.)

\midinsert
$$
\Gitter(15,8)(0,0)
\Koordinatenachsen(15,8)(0,0)
\Pfad(0,7),666121666121661222111666116\endPfad
\Pfad(2,7),61161111116611666\endPfad
\Pfad(4,7),116121111166121666666\endPfad
\DickPunkt(0,7)
\DickPunkt(2,7)
\DickPunkt(4,7)
\DickPunkt(10,0)
\DickPunkt(12,0)
\DickPunkt(14,0)
\Label\ro
{P_0}(10,0)
\Label\ro{P_1}(12,0)
\Label\ro{P_2}(14,0)
\Label\lo{\text{\it4}}(1,4)
\Label\ro{\text{\it5}}(1,5)
\Label\lo{\text{\it2}}(3,2)
\Label\lo{\text{\it6}}(3,6)
\Label\ro{\text{\it3}}(3,3)
\Label\ro{\text{\it6}}(3,6)
\Label\lo{\text{\it1}}(5,1)
\Label\lo{\text{\it5}}(5,5)
\Label\lo{\text{\it7}}(5,7)
\Label\ro{\text{\it4}}(5,4)
\Label\ro{\text{\it5}}(5,5)
\Label\ro{\text{\it7}}(5,7)
\Label\lo{\text{\it4}}(7,4)
\Label\lo{\text{\it5}}(7,5)
\Label\lo{\text{\it6}}(7,6)
\Label\ro{\text{\it4}}(7,4)
\Label\ro{\text{\it5}}(7,5)
\Label\ro{\text{\it7}}(7,7)
\Label\lo{\text{\it1}}(9,1)
\Label\lo{\text{\it5}}(9,5)
\Label\lo{\text{\it7}}(9,7)
\Label\ro{\text{\it1}}(9,1)
\Label\ro{\text{\it5}}(9,5)
\Label\ro{\text{\it7}}(9,7)
\Label\lo{\text{\it3}}(11,3)
\Label\lo{\text{\it7}}(11,7)
\Label\ro{\text{\it3}}(11,3)
\Label\ro{\text{\it7}}(11,7)
\Label\lo{\text{\it5}}(13,5)
\Label\ro{\text{\it6}}(13,6)
\hskip7cm
$$
\centerline{\eightpoint A family of non-intersecting paths in the
directed graph $\Cal G_7
$}
\vskip7pt
\centerline{\eightpoint Figure \FDc}
\endinsert

From here on, everything 
just works completely analogously to the corresponding arguments in the
proof of Theorem~\TF. In particular, under the correspondence between
families of non-intersecting paths and arrays of alternating sequences
described there, the path family in Figure~\FDc\ corresponds to
the array in (\CEa).
We leave the details to the reader.\quad \quad \qed
\enddemo

\subhead 7. A third and fourth reciprocity theorem\endsubhead
The main results in this section are reciprocity laws for
determinants of numbers of bounded up-down paths with specified starting
and ending heights, see Theorems~\TJ\ and~\TM.
As we show in Propositions~\TK, \TN\ and Theorems~\TL, \TO, 
in combinatorial terms
these reciprocity laws set up correspondences between families of
non-intersecting bounded up-down paths with specified starting and
ending heights and flagged alternating tableaux
(see Section~11 for their formal definition)
of rectangular shape. 

The first set of results in this section concerns determinants of numbers of
paths of even length and flagged alternating tableaux of rectangular
shape with an odd number of columns. 

\proclaim{Theorem \TJ}
Let $n$ be a non-negative integer and $k,m$ be positive integers,
and let $r_0<r_1<\dots<r_{k-1}$
and $s_0<s_1<\dots<s_{k-1}$ be sequences of positive integers with
$1\le r_i,s_i\le k+m$ for all~$i$. Then
$$
\multline
\det\left(C_{2n}^{(2k+2m-1)}(2r_i-2\to 2s_j-2)
\right)_{0\le i,j\le k-1}\\
=(-1)^{\sum_{i=0}^{m-1}(\bar r_i+\bar s_i)}
\det\left(C_{-2n}^{(2k+2m-1)}(2\bar r_i-2\to 2\bar s_j-2)\right)_{0\le i,j\le m-1},
\endmultline
\tag\DA
$$
where 
$$
\align
\{\bar r_0,\bar r_1,\dots,\bar r_{m-1}\}&=
\{1,2,\dots,k+m\}\setminus\{r_0,r_1,\dots,r_{k-1}\},\\
\{\bar s_0,\bar s_1,\dots,\bar s_{m-1}\}&=
\{1,2,\dots,k+m\}\setminus\{s_0,s_1,\dots,r_{k-1}\},
\endalign
$$
and where we assume that 
$\bar r_0<\bar r_1<\dots<\bar r_{m-1}$ and
$\bar s_0<\bar s_1<\dots<\bar s_{m-1}$.
\endproclaim

Our line of argument to prove the above theorem is completely analogous to
the ones proving Theorems~\TD\ and~\TG. Namely,
by the Lindstr\"om--Gessel--Viennot theorem \cite{\LindAA, Lemma~1}, 
the determinant
on the left-hand side of (\DA) equals the number of 
families $(P_0,P_1,\dots,P_{k-1})$ of non-intersecting up-down
paths of height at most $2k+2m-1$ that do not pass below the $x$-axis, 
where $P_i$ runs from
$(0,2r_i-2)$ to $(2n,2s_i-2)$, $i=0,1,\dots,k-1$. 
An example of such a family of
non-intersecting up-down paths for $k=4$, $m=3$, $n=10$,
$r_0=1$, $r_1=4$, $r_2=5$, $r_3=7$, 
$s_0=2$, $s_1=3$, $s_2=5$, and $s_3=7$
is 
shown in Figure~\FF. (At this point, circled points and attached
labels should be ignored.) 

\midinsert
$$
\Einheit.4cm
\Gitter(22,15)(-1,-1)
\Koordinatenachsen(22,15)(-1,-1)
\Pfad(0,12),43433444433433444333\endPfad
\Pfad(0,8),43343443443444333433\endPfad
\Pfad(0,6),44334444334344343334\endPfad
\Pfad(0,0),33344344343344343334\endPfad
\SPfad(20,-2),22222222222222222\endSPfad
\PfadDicke{1pt}
\Pfad(-1,13),1111111111111111111111\endPfad
\DickPunkt(0,0)
\DickPunkt(0,6)
\DickPunkt(0,8)
\DickPunkt(0,12)
\DickPunkt(20,2)
\DickPunkt(20,4)
\DickPunkt(20,8)
\DickPunkt(20,12)
\Kreis(0,2)
\Kreis(0,4)
\Kreis(0,10)
\Kreis(1,3)
\Kreis(1,9)
\Kreis(1,13)
\Kreis(2,0)
\Kreis(2,6)
\Kreis(2,10)
\Kreis(3,1)
\Kreis(3,7)
\Kreis(3,13)
\Kreis(4,0)
\Kreis(4,4)
\Kreis(4,10)
\Kreis(5,3)
\Kreis(5,7)
\Kreis(5,11)
\Kreis(6,0)
\Kreis(6,6)
\Kreis(6,10)
\Kreis(7,5)
\Kreis(7,9)
\Kreis(7,13)
\Kreis(8,4)
\Kreis(8,6)
\Kreis(8,12)
\Kreis(9,5)
\Kreis(9,11)
\Kreis(9,13)
\Kreis(10,2)
\Kreis(10,8)
\Kreis(10,12)
\Kreis(11,5)
\Kreis(11,9)
\Kreis(11,13)
\Kreis(12,0)
\Kreis(12,8)
\Kreis(12,12)
\Kreis(13,7)
\Kreis(13,9)
\Kreis(13,13)
\Kreis(14,6)
\Kreis(14,8)
\Kreis(14,10)
\Kreis(15,7)
\Kreis(15,9)
\Kreis(15,13)
\Kreis(16,4)
\Kreis(16,8)
\Kreis(16,12)
\Kreis(17,5)
\Kreis(17,11)
\Kreis(17,13)
\Kreis(18,0)
\Kreis(18,8)
\Kreis(18,12)
\Kreis(19,1)
\Kreis(19,9)
\Kreis(19,13)
\Kreis(20,0)
\Kreis(20,6)
\Kreis(20,10)
\Label\l{2\kern-4pt}(0,2)
\Label\l{3\kern-4pt}(0,4)
\Label\l{6\kern-4pt}(0,10)
\Label\l{2\kern-4pt}(1,3)
\Label\l{5\kern-4pt}(1,9)
\Label\l{7\kern-4pt}(1,13)
\Label\l{1\kern-4pt}(2,0)
\Label\l{4\kern-4pt}(2,6)
\Label\l{6\kern-4pt}(2,10)
\Label\l{1\kern-4pt}(3,1)
\Label\l{4\kern-4pt}(3,7)
\Label\l{7\kern-4pt}(3,13)
\Label\l{1\kern-4pt}(4,0)
\Label\l{3\kern-4pt}(4,4)
\Label\l{6\kern-4pt}(4,10)
\Label\l{2\kern-4pt}(5,3)
\Label\l{4\kern-4pt}(5,7)
\Label\l{6\kern-4pt}(5,11)
\Label\l{1\kern-4pt}(6,0)
\Label\l{4\kern-4pt}(6,6)
\Label\l{6\kern-4pt}(6,10)
\Label\l{3\kern-4pt}(7,5)
\Label\l{5\kern-4pt}(7,9)
\Label\l{7\kern-4pt}(7,13)
\Label\l{3\kern-4pt}(8,4)
\Label\l{4\kern-4pt}(8,6)
\Label\l{7\kern-4pt}(8,12)
\Label\l{3\kern-4pt}(9,5)
\Label\l{6\kern-4pt}(9,11)
\Label\l{7\kern-4pt}(9,13)
\Label\l{2\kern-4pt}(10,2)
\Label\l{5\kern-4pt}(10,8)
\Label\l{7\kern-4pt}(10,12)
\Label\l{3\kern-4pt}(11,5)
\Label\l{5\kern-4pt}(11,9)
\Label\l{7\kern-4pt}(11,13)
\Label\l{1\kern-4pt}(12,0)
\Label\l{5\kern-4pt}(12,8)
\Label\l{7\kern-4pt}(12,12)
\Label\l{4\kern-4pt}(13,7)
\Label\l{5\kern-4pt}(13,9)
\Label\l{7\kern-4pt}(13,13)
\Label\l{4\kern-4pt}(14,6)
\Label\l{5\kern-4pt}(14,8)
\Label\l{6\kern-4pt}(14,10)
\Label\l{4\kern-4pt}(15,7)
\Label\l{5\kern-4pt}(15,9)
\Label\l{7\kern-4pt}(15,13)
\Label\l{3\kern-4pt}(16,4)
\Label\l{5\kern-4pt}(16,8)
\Label\l{7\kern-4pt}(16,12)
\Label\l{3\kern-4pt}(17,5)
\Label\l{6\kern-4pt}(17,11)
\Label\l{7\kern-4pt}(17,13)
\Label\l{1\kern-4pt}(18,0)
\Label\l{5\kern-4pt}(18,8)
\Label\l{7\kern-4pt}(18,12)
\Label\l{1\kern-4pt}(19,1)
\Label\l{5\kern-4pt}(19,9)
\Label\l{7\kern-4pt}(19,13)
\Label\l{1\kern-4pt}(20,0)
\Label\l{4\kern-4pt}(20,6)
\Label\l{6\kern-4pt}(20,10)
%
\hskip8.8cm
$$
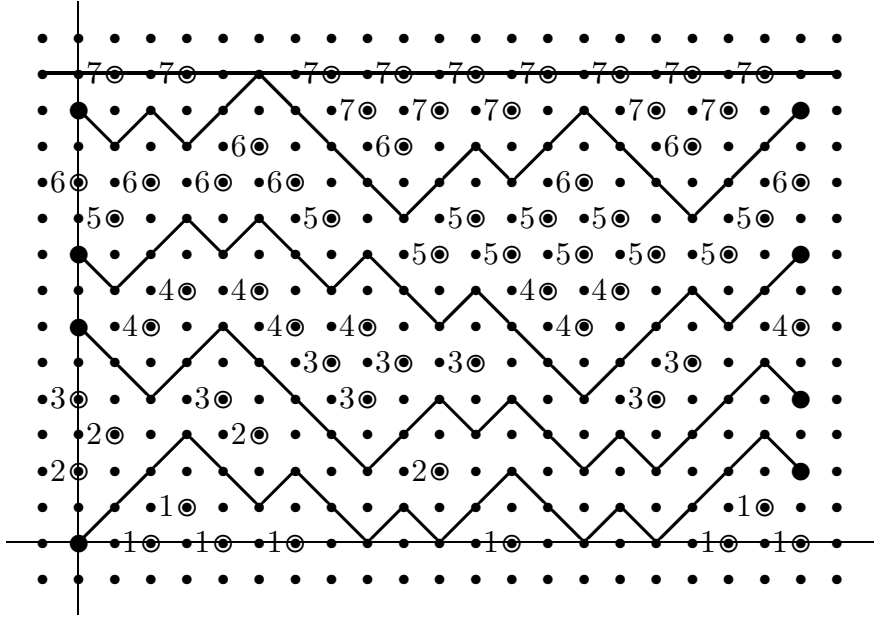
\centerline{\eightpoint A family of non-intersecting bounded up-down
paths with arbitrary starting and ending points}
\vskip7pt
\centerline{\eightpoint Figure \FF}
\endinsert

In Proposition~\TK\ below, we shall show that there is a bijection
between these families of non-intersecting up-down paths and certain
rectangular arrays of alternating sequences. Subsequently, in Theorem~\TL,
it is shown that the number of these rectangular arrays is given by the
determinant on the right-hand side of~(\DA).
Hence, the combination of the above
observation with Proposition~\TK\ and Theorem~\TL\ will establish Theorem~\TJ.

\proclaim{Proposition \TK}
With the assumptions and notation of Theorem~{\rm\TJ},
there is a bijection between 
families $(P_0,P_1,\dots,P_{k-1})$ of non-intersecting up-down
paths of height at most $2k+2m-1$ that do not pass below the $x$-axis, 
where $P_i$ runs from
$(0,2r_i-2)$ to $(2n,2s_i-2)$ and rectangular arrays of integers of the form
$$
\matrix 
a_{1,1}&a_{1,2}&\dots&a_{1,2n+1}\\
a_{2,1}&a_{2,2}&\dots&a_{2,2n+1}\\
\hdotsfor4\\
a_{m,1}&a_{m,2}&\dots&a_{m,2n+1}\\
\endmatrix
\tag\DB
$$
where $a_{i,1}=\bar r_{m-i}$ and $a_{i,2n+1}=\bar s_{m-i}$ for all~$i$,
in which each row is alternating, that is,
$$
a_{i,1}\le a_{i,2}\ge a_{i,3}\le
a_{i,4}\ge\dots\ge a_{i,2n+1}
$$
for all $i$,
and in which we have 
$$
1\le a_{i,j}\le k+m
$$
and
$$
a_{i+1,2j}<a_{i,2j+1}>a_{i+1,2j+2}
$$
for all $i$ and $j$.
\endproclaim

\demo{Proof}
The reader should consult Figure~\FF\ while reading the following
explanations. As mentioned earlier, the figure shows an example
of a family of non-intersecting up-down paths as described in the
statement of the proposition
for $k=4$, $m=3$, $n=10$,
$r_0=1$, $r_1=4$, $r_2=5$, $r_3=7$, 
$s_0=2$, $s_1=3$, $s_2=5$, and $s_3=7$.

Analogously to the proofs of Propositions~\TE\ and \TH, for all $i$ and~$j$, 
we mark all lattice points $(2i-1,2j-1)$ and all lattice
points $(2i,2j)$ in the rectangular region
$$
\{(x,y):0\le y\le 2k+2m-1,\ 0\le x\le 2n\},
$$
that are not occupied by any of the paths.
In the figure, these lattice points are marked by circles.

Also here, each marked point acquires a label. 
As earlier, the label is equal to $1$ plus the
number of paths and marked points that are below $X$ and have the
same abscissa as~$X$. Equivalently, the label
of a marked point $X=(x,y)$ equals $\cl{(y+1)/2}$. The figure 
shows the labels of these marked points. 

Finally, for $i=0,1,\dots,2n$, we read the labels of the marked
points with abscissa~$i$ and put them into a column, and then
we concatenate the columns into an array of
rectangular shape as in (\DB). In this way, for the example in Figure~\FF,
we obtain the array
$$
\matrix
6&7&6&7&6&6&6&7&7&7&7&7&7&7&6&7&7&7&7&7&6\\
3&5&4&4&3&4&4&5&4&6&5&5&5&5&5&5&5&6&5&5&4\\
2&2&1&1&1&2&1&3&3&3&2&3&1&4&4&4&3&3&1&1&1
\endmatrix
\tag\DBa$$

It is not difficult to see that (in general) the obtained array satisfies
the constraints in the statement of the proposition. 

Conversely, given an array of the form (\DB) as in the statement of
the proposition, if one reverses the previously described construction
then one sees that this array determines a unique family 
 of non-intersecting up-down paths as described in the
statement of the proposition.\quad \quad \qed
\enddemo

\proclaim{Theorem \TL}
With the assumptions and notation of Theorem~{\rm\TJ},
the number of rectangular arrays of integers of the form {\rm(\DB)}
that satisfy the constraints given in the statement of Proposition~{\rm\TK} 
is equal to
$$
(-1)^{\sum_{i=0}^{m-1}(\bar r_i+\bar s_i)}
\det\left(C_{-2n}^{(2k+2m-1)}(2\bar r_i-2\to 2\bar s_j-2)\right)_{0\le i,j\le m-1}.
\tag\DC
$$
\endproclaim

\demo{Sketch of proof of Theorem \TL}
We observe that, by Corollary~\UF, the determinant in
(\DC) can be alternatively written as
$$
\det\left(\big\vert\Cal A_{2n+1}^{(k+m)}(\bar r_i\to \bar s_j)
\big\vert\right)_{0\le i,j\le m-1}.
$$

We claim that this determinant equals the number of families
$(P_0,P_1,\dots,P_{m-1})$ of non-intersecting lattice paths in the
directed graph $\Cal G_{k+m}$ (see the proof of Theorem~\TF\ and
in particular Figure~\FDa), where $P_i$ runs from $(0,\bar r_i)$ to
$(2n+1,\bar s_i)$, and starts and ends with a horizontal step. 
Figure~\FDd\ shows an example of such a path family
for $k=4$, $m=3$, $n=10$,
$r_0=1$, $r_1=4$, $r_2=5$, $r_3=7$, 
$s_0=2$, $s_1=3$, $s_2=5$, and $s_3=7$, so that
$\bar r_0=2$, $\bar r_1=3$, $\bar r_2=6$,
$\bar s_0=1$, $\bar s_1=4$, and $\bar s_2=6$.
(The numbers should be ignored at this point.)

\midinsert
$$
\Gitter(22,8)(0,0)
\Koordinatenachsen(22,8)(0,0)
\Pfad(0,6),12161216111211111116121111161\endPfad
\Pfad(0,3),1221611612112161221611111112161161\endPfad
\Pfad(0,2),116111216122111612166122211161166111\endPfad
\DickPunkt(0,6)
\DickPunkt(0,3)
\DickPunkt(0,2)
\DickPunkt(21,6)
\DickPunkt(21,4)
\DickPunkt(21,1)
\Label\l{P_0}(0,2)
\Label\l{P_1}(0,3)
\Label\l{P_2}(0,6)
\Label\lo{\text{\it2}}(1,2)
\Label\lo{\text{\it3}}(1,3)
\Label\lo{\text{\it6}}(1,6)
\Label\ro{\text{\it2}}(1,2)
\Label\ro{\text{\it5}}(1,5)
\Label\ro{\text{\it7}}(1,7)
\Label\lo{\text{\it1}}(3,1)
\Label\lo{\text{\it4}}(3,4)
\Label\lo{\text{\it6}}(3,6)
\Label\ro{\text{\it1}}(3,1)
\Label\ro{\text{\it4}}(3,4)
\Label\ro{\text{\it7}}(3,7)
\Label\lo{\text{\it1}}(5,1)
\Label\lo{\text{\it3}}(5,3)
\Label\lo{\text{\it6}}(5,6)
\Label\ro{\text{\it2}}(5,2)
\Label\ro{\text{\it4}}(5,4)
\Label\ro{\text{\it6}}(5,6)
\Label\lo{\text{\it1}}(7,1)
\Label\lo{\text{\it4}}(7,4)
\Label\lo{\text{\it6}}(7,6)
\Label\ro{\text{\it3}}(7,3)
\Label\ro{\text{\it5}}(7,5)
\Label\ro{\text{\it7}}(7,7)
\Label\lo{\text{\it3}}(9,3)
\Label\lo{\text{\it4}}(9,4)
\Label\lo{\text{\it7}}(9,7)
\Label\ro{\text{\it3}}(9,3)
\Label\ro{\text{\it6}}(9,6)
\Label\ro{\text{\it7}}(9,7)
\Label\lo{\text{\it2}}(11,2)
\Label\lo{\text{\it5}}(11,5)
\Label\lo{\text{\it7}}(11,7)
\Label\ro{\text{\it3}}(11,3)
\Label\ro{\text{\it5}}(11,5)
\Label\ro{\text{\it7}}(11,7)
\Label\lo{\text{\it1}}(13,1)
\Label\lo{\text{\it5}}(13,5)
\Label\lo{\text{\it7}}(13,7)
\Label\ro{\text{\it4}}(13,4)
\Label\ro{\text{\it5}}(13,5)
\Label\ro{\text{\it7}}(13,7)
\Label\lo{\text{\it4}}(15,4)
\Label\lo{\text{\it5}}(15,5)
\Label\lo{\text{\it6}}(15,6)
\Label\ro{\text{\it4}}(15,4)
\Label\ro{\text{\it5}}(15,5)
\Label\ro{\text{\it7}}(15,7)
\Label\lo{\text{\it3}}(17,3)
\Label\lo{\text{\it5}}(17,5)
\Label\lo{\text{\it7}}(17,7)
\Label\ro{\text{\it3}}(17,3)
\Label\ro{\text{\it6}}(17,6)
\Label\ro{\text{\it7}}(17,7)
\Label\lo{\text{\it1}}(19,1)
\Label\lo{\text{\it5}}(19,5)
\Label\lo{\text{\it7}}(19,7)
\Label\ro{\text{\it1}}(19,1)
\Label\ro{\text{\it5}}(19,5)
\Label\ro{\text{\it7}}(19,7)
\Label\lo{\text{\it1}}(21,1)
\Label\lo{\text{\it4}}(21,4)
\Label\lo{\text{\it6}}(21,6)
\hskip10.5cm
$$
\centerline{\eightpoint A family of non-intersecting paths in the
directed graph $\Cal G_7
$}
\vskip7pt
\centerline{\eightpoint Figure \FDd}
\endinsert

From here on, everything 
just works completely analogously to the corresponding arguments in the
proof of Theorem~\TF. In particular, under the correspondence between
families of non-intersecting paths and arrays of alternating sequences
described there, the path family in Figure~\FDd\ corresponds to
the array in (\DBa).
We leave the details to the reader.\quad \quad \qed
\enddemo

The second set of results in this section concerns determinants of numbers of
paths of odd length and flagged alternating tableaux of rectangular
shape with an even number of columns. Since the proofs are completely
analogous to the previous proofs in this section, we skip them here entirely
for the sake of brevity and leave them to the reader.

\proclaim{Theorem \TM}
Let $n,k,m$ be positive integers, and let $r_0<r_1<\dots<r_{k-1}$
and $s_0<s_1<\dots<s_{k-1}$ be sequences of positive integers with
$1\le r_i,s_i\le k+m$ for all~$i$. Then
$$
\multline
\det\left(C_{2n-1}^{(2k+2m-1)}(2r_i-2\to 2s_j-1)
\right)_{0\le i,j\le k-1}\\
=(-1)^{\sum_{i=0}^{m-1}(\bar r_i+\bar s_i)}
\det\left(C_{-2n+1}^{(2k+2m-1)}(2\bar r_i-2\to 2\bar s_j-1)\right)_{0\le i,j\le m-1},
\endmultline
\tag\DD
$$
where 
$$
\align
\{\bar r_0,\bar r_1,\dots,\bar r_{m-1}\}&=
\{1,2,\dots,k+m\}\setminus\{r_0,r_1,\dots,r_{k-1}\},\\
\{\bar s_0,\bar s_1,\dots,\bar s_{m-1}\}&=
\{1,2,\dots,k+m\}\setminus\{s_0,s_1,\dots,r_{k-1}\},
\endalign
$$
and where we assume that 
$\bar r_0<\bar r_1<\dots<\bar r_{m-1}$ and
$\bar s_0<\bar s_1<\dots<\bar s_{m-1}$.
\endproclaim

By the Lindstr\"om--Gessel--Viennot theorem \cite{\LindAA, Lemma~1}, 
the determinant
on the left-hand side of (\DA) equals the number of 
families $(P_0,P_1,\dots,P_{k-1})$ of non-intersecting up-down
paths of height at most $2k+2m-1$ that do not pass below the $x$-axis, 
where $P_i$ runs from
$(0,2r_i-2)$ to $(2n-1,2s_i-1)$, $i=0,1,\dots,k-1$. 

Hence, we see that Theorem~\TM\ will follow by combining the above
observation with Proposition~\TN, Corollary~\UF\ with $k$ replaced
by~$k+m$, and Theorem~\TO.

\proclaim{Proposition \TN}
With the assumptions and notation of Theorem~{\rm\TM},
there is a bijection between 
families $(P_0,P_1,\dots,P_{k-1})$ of non-intersecting up-down
paths of height at most $2k+2m-1$ that do not pass below the $x$-axis, 
where $P_i$ runs from
$(0,2r_i-2)$ to $(2n-1,2s_i-1)$ and rectangular arrays of integers of the form
$$
\matrix 
a_{1,1}&a_{1,2}&\dots&a_{1,2n}\\
a_{2,1}&a_{2,2}&\dots&a_{2,2n}\\
\hdotsfor4\\
a_{m,1}&a_{m,2}&\dots&a_{m,2n}\\
\endmatrix
\tag\DE
$$
where $a_{i,1}=\bar r_{m-i}$ and $a_{i,2n}=\bar s_{m-i}$ for all~$i$,
in which each row is alternating, that is,
$$
a_{i,1}\le a_{i,2}\ge a_{i,3}\le
a_{i,4}\ge\dots\le a_{i,2n}
$$
for all $i$,
and in which we have 
$$
1\le a_{i,j}\le k+m
$$
and
$$
a_{i+1,2j}<a_{i,2j+1}>a_{i+1,2j+2}
$$
for all $i$ and $j$.
\endproclaim

\proclaim{Theorem \TO}
With the assumptions and notation of Theorem~{\rm\TM},
the number of rectangular arrays of integers of the form {\rm(\DE)}
that satisfy the constraints given in the statement of Proposition~{\rm\TN} 
is equal to
$$
(-1)^{\sum_{i=0}^{m-1}(\bar r_i+\bar s_i)}
\det\left(C_{-2n+1}^{(2k+2m-1)}(2\bar r_i-2\to 2\bar s_j-1)\right)_{0\le i,j\le m-1}.
\tag\DF
$$
\endproclaim

\subhead 8. Weighted enumeration\endsubhead
We now embark on the extension of
the results from Sections~2--7 to a weighted setting.
In this section we begin with the extensions of the enumeration results for 
paths and alternating sequences in Sections~2 and~3.
In doing so, the Chebyshev polynomials of the second kind get replaced
by two more general sequences of orthogonal polynomials. We define
the sequence $(P_n(x))_{n\ge0}$ of polynomials by the two-term recurrence
$$
P_{n+1}(x)=xP_n(x)-B_{n}P_{n-1}(x),\quad \text{for }n\ge1,
\tag\EA
$$
with initial conditions $P_0(x)=1$ and $P_1(x)=x$. Here, the $B_i$'s
are indeterminates. The sequence $(P_n(x))_{n\ge0}$ is a sequence
of orthogonal polynomials (in the formal sense) due to Favard's
theorem (cf\. \cite{\FavaAA} and \cite{\VienAE, Cor.~19 on p.~I-15}).
Comparison with (\BD) shows that $P_n(x)=U_n(x/2)$ if all $B_i$'s
are equal to~$1$. Later we shall use two $q$-analogues of the Chebyshev
polynomials~(\BAb) (with $x$ replaced by~$x/2$), namely
$$
P_n(x)\Big\vert_{B_i=q^{i-1}}=
\sum_{j\ge0}(-1)^jq^{j(j-1)}\bmatrix {n-j}\\{j}\endbmatrix_q x^{n-2j}
\tag\EAa$$
and
$$
P_n(x)\Big\vert_{B_i=q^{i}}=
\sum_{j\ge0}(-1)^jq^{j^2}\bmatrix {n-j}\\{j}\endbmatrix_q x^{n-2j},
\tag\EAb$$
where the $q$-binomial coefficient is defined by
$$
\bmatrix
N\\j\endbmatrix_q =\frac {(1-q^N)(1-q^{N-1})\cdots (1-q^{N-j+1})} 
{(1-q^j)(1-q^{j-1})\cdots (1-q)}.
\tag\EAc$$
The validity of (\EAa) and (\EAb) is easily verified by induction
on~$n$, using the recurrence~(\EA).

\medskip
Our second sequence, $(Q_n(x))_{n\ge0}$ is defined by
$$
Q_{n+1}(x)=\cases
V_{{(n+2)/2}}xQ_{n}(x)-Q_{n-1}(x),&\text{if $n$ is even},\\
A_{{(n+1)/2}}xQ_{n}(x)-Q_{n-1}(x),&\text{if $n$ is odd},\\
\endcases
\quad \text{for }n\ge1,
\tag\EB
$$
with initial conditions $Q_0(x)=1$ and $Q_1(x)=V_1x$. Here, the $A_i$'s
and $V_i$'s are indeterminates. 
The sequence $(Q_n(x))_{n\ge0}$ is also a sequence
of orthogonal polynomials (in the formal sense), for the same reason
as above. If all $A_i$'s and $V_i$'s equal~$1$, 
then $Q_n(x)$ reduces to $U_n(x/2)$.
Otherwise, however, the polynomials
$Q_n(x)$ are not monic, as opposed to the polynomials $P_n(x)$.
Namely, the leading coefficient of $Q_{2n}(x)$ is
$A_1A_2\cdots A_nV_1V_2\cdots V_n$, while the leading coefficient of
$Q_{2n+1}(x)$ is $A_1A_2\cdots A_nV_1V_2\cdots V_{n+1}$.
Let $Q^*_n(x)$ denote the monic form of $Q_n(x)$, that is,
$Q_n(x)$ divided by its leading coefficient. Then,
rewriting~(\EB), we have
$$
\align
Q^*_{2n}(x)
&=xQ^*_{2n-1}(x)
-\frac {1} {A_nV_n}Q^*_{2n-2}(x),\\
Q^*_{2n+1}(x)
&=xQ^*_{2n}(x)
-\frac {1} {A_nV_{n+1}}Q^*_{2n-1}(x).
\endalign
$$
If we compare these relations with (\EA), then we see how the
polynomials $Q_n(x)$ are related to the polynomials $P_n(x)$:
we have
$$
P_n(x)\Bigg\vert_{\smallmatrix B_{2i-1}=A_i^{-1}V_i^{-1}\kern.3cm
\\B_{2i}=A_i^{-1}V_{i+1}^{-1}
\endsmallmatrix}=Q^*_n(x), \quad \text{for }n\ge0.
\tag\EBa
$$

\medskip
We are now going to extend Theorem~\VA\ to a weighted setting.
We introduce the following weight for up-down paths, which we are
going to denote by $w_B$. 
To all up-steps we assign a weight of~$1$, while to a 
down-step from height~$h$ to height~$h-1$ we assign a weight of~$B_h$.
The weight $w_B(P)$ of an up-down path~$P$ is defined to be the
product of the weights of all its steps. Thus, the weight of the
second path from top in Figure~\FD\ (the path from $(-4,0)$ to
$(18,0)$) is $B_5B_7B_6B_5B_6B_6B_6B_5B_4B_3B_2B_1=
B_1B_2B_3B_4B_5^3B_6^4B_7$.

Here, and in the following, given a set $\Cal O$ of combinatorial objects,
we write $\GF(\Cal O;w)$ for the
generating function $\sum_{t\in \Cal O}w(t)$.

With the above notations and definitions, the weighted extension of
Theorem~\VA\ is the following. (As indicated earlier, a further
generalisation to Motzkin paths is possible, see 
\cite{\VienAE, Ch.~V, Eq.~(27)}, \cite{\KratBP, proof of Theorem~A2},
\cite{\KratCL, proof of Theorem~10.11.1}. While it is not relevant
here, we shall come back to it in Section~12.(3).)

\proclaim{Theorem \TP}
For all non-negative integers $r,s,k$ with $0\le r,s\le k$, we have
$$
\sum_{n\ge0}\GF\left(\Cal C_n^{(k)}(r\to s);w_B\right)\, x^{n}
=\cases \dfrac {P_r(1/x)\,T_B^{s+1}P_{k-s}(1/x)} {x\,P_{k+1}(1/x)},
&\text{if }r\le s,\\
B_{s+1}B_{s+2}\cdots B_r
\dfrac {P_s(1/x)\,T_B^{r+1}P_{k-r}(1/x)}
{x\,P_{k+1}(1/x)},
&\text{if }r\ge s,
\endcases
\tag\EC
$$
where the operator $T_B$,  when applied to a polynomial in the $B_i$'s.
replaces each $B_i$ by $B_{i+1}$, for all~$i$.
Here, as before, $\Cal C_n^{(k)}(r\to s)$ stands for the set of all up-down paths
from $(0,r)$ to $(n,s)$ of height at most~$k$ that do not pass below
the $x$-axis.
\endproclaim

\demo{Sketch of proof}
The proof of Theorem~\VA\ using the theory of heaps also proves
the present theorem. The only thing that has to be done is to extend
Lemma~\VC\ to our weighted setting. Indeed, we claim that
$$
\sum_T (-1)^{\vert T\vert}w_B(T)\,x^{2\vert T\vert}=x^{k+1}P_{k+1}(1/x),
$$
where the sum is over all trivial heaps $T$ of dimers on $[0,k]$, and
in which $w_B(T)$ is the obvious transfer of the weight $w_B$ to heaps
of dimers where the weight of a dimer $d_i$ is $B_{i+1}$.\footnote{%
We alert the reader that our definition of weight
of a dimer here is not consistent with the definition of weight of a dimer
in Section~2 in that, here, we do not
include $x^2$ --- which keeps track of the length of a path,
respectively of the size of a heap --- in the weight. This choice makes the exposition
in Section~10 more convenient. We hope that the reader bears with us.}
The proof of the claim works in exactly the same way as the proof of
Lemma~\VC. In particular, the extension of the computation (\BGc) is
$$
(-B_{k}x^2)x^{k-1}P_{k-1}(1/x)+x^kP_{k}(1/x)=x^{k+1}P_{k+1}(1/x),
$$
here following from the two-term recurrence~(\EA).

The appearance of the shift $T_B^{s+1}$ in the first case of~(\EC)
is explained by the fact that the second factor in the numerator
represents the contribution of trivial heaps of dimers on $[s+1,k]$,
that is, in comparison to heaps on $[0,k-s-1]$ there is a shift of
$s+1$ in the $B_i$'s that appears in the weights of heaps.
There is an analogous explanation for the appearance of $T_B^{r+1}$
in the second case of~(\EC). In that case, we also have to insert
the factor $B_{s+1}B_{s+2}\cdots B_r$ because the reflection argument
used in the proof of Theorem~\VA\ maps down-steps to up-steps and
vice versa, which is weight-preserving only for {\it pairs}
consisting of an up- and a matching down-step;\footnote{``Matching" 
up- and down-steps in up-down paths refers to the following procedure:
for an up-step look to the right; if there is
a down-step on the same height as the up-step then this is the
``matching" down-step. Conversely, 
for a down-step look to the left; if there is
an up-step on the same height as the down-step then this is the
``matching" down-step.} 
indeed, if $r>s$ then there are
unmatched down-steps from height~$r$ to height~$r-1$, 
from height~$r-1$ to height~$r-2$, \dots,
from height~$s+1$ to height~$s$, which altogether contribute a
weight of $B_rB_{r-1}\cdots B_{s+1}$.\quad \quad \qed
\enddemo

If we specialise $B_i$ to $q^{i-1}$ then the weight $w_B(P)$ of a Dyck
path~$P$ reduces to the classical {\it area weight\/} $q^{a(P)}$ of~$P$, with
$a(P)$ denoting the number of full squares that fit between the path
and the zigzag lower bound on the Dyck path caused by the $x$-axis. 
This is illustrated in Figure~\FDe. The full squares between path and
$x$-axis are indicated by the dotted line segments. Thus, for the
path~$P_0$ in the figure we have $a(P_0)=11$.

\midinsert
$$
\Gitter(18,6)(-1,-1)
\Koordinatenachsen(18,6)(-1,-1)
\Pfad(0,0),3344333433443444\endPfad
\SPfad(1,1),43\endSPfad
\SPfad(5,1),4343434343\endSPfad
\SPfad(6,2),43434343\endSPfad
\SPfad(9,3),43\endSPfad
\DickPunkt(0,0)
\DickPunkt(16,0)
\Label\ro{P_0}(7,3)
\hskip9cm
$$
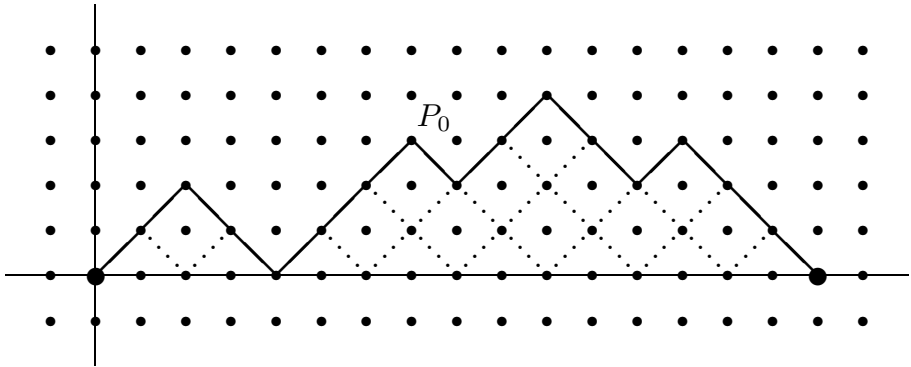
\centerline{\eightpoint Area below a Dyck  path}
\vskip7pt
\centerline{\eightpoint Figure \FDe}
\endinsert

Theorem~\TP\ with $r=s=0$ and $B_i=q^{i-1}$ for all~$i$ 
together with~(\EAa) and~(\EAb) then implies
that the area generating function for Dyck paths is given by
(cf\. \cite{\OwPrAA, Cor.~2})
$$
\sum_{n\ge0}\GF\left(\Cal C_n^{(k)};q^{a(\,.\,)}\right)\, x^{n}
=\dfrac {\dsize\sum_{j\ge0}(-1)^jq^{j^2}\bmatrix {k-j}\\{j}\endbmatrix_q x^{2j}} 
{\dsize\sum_{j\ge0}(-1)^jq^{j(j-1)}\bmatrix {k+1-j}\\{j}\endbmatrix_q x^{2j}}.
\tag\ECa$$

\medskip
Next we turn to alternating sequences. We define the weight $w_{AV}(S)$
of an alternating sequence $S=
a_1\le a_2\ge a_3\le a_4\ge \cdots {}\diamond\,a_{n-1}{}\square\,a_n,
$
where $(\diamond,\square)=(\ge,\le)$ if $n$ is even and
$(\diamond,\square)=(\le,\ge)$ if $n$ is odd, by 
$\Big(\prod _{i=1} ^{\fl{n/2}}A_{a_{2i}}\Big)
\Big(\prod _{i=1} ^{\cl{n/2}}V_{a_{2i-1}}\Big)$. 
In colloquial terms, the ``top entries" in the alternating sequence
are assigned the $A$-variables as weights, 
and the ``bottom entries" the $V$-variables.\footnote{As the reader
may have guessed, the choice of
letters for the weights has its explanation in the interpretation of the
letter~$A$ as an arrow-like symbol pointing up and of the letter~$V$ as
arrow-like symbol pointing down.}
Thus, the weight of the
alternating sequence on the top of Figure~\FA\ is
$$
\multline
V_4 A_4 V_3 A_3 V_1 A_1 V_1 A_3 V_2 A_6 V_6 A_8 V_4 A_7 V_4 A_4 V_2
 A_3 V_2 A_2 V_2 A_5 V_5 A_6 V_3 A_6 V_5 A_8 V_6\\
=A_1A_2A_3^3A_4^2A_5A_6^3A_7A_8^2V_1^2V_2^4V_3^2V_4^3V_5^2V_6^2.
\endmultline
$$

The weighted extension of Theorem~\UA\ then reads as follows.

\proclaim{Theorem \TQ}
For all positive integers $r,s,k$ with $1\le r,s\le k$, we have
$$
\multline
\sum_{n\ge0}\GF\left(\Cal A_{2n+1}^{(k)}(r\to s);w_{AV}\right)\, x^{2n}
\\=
\cases 
\dsize
(-1)^{r+s+1}V_rV_s\frac {xQ_{2r-2}(x)\,T_{AV}^{s-1}R_{AV}^{(k-s+1)}Q_{2k+1-2s}(x)} 
{Q_{2k}(x)},&\text{if }r< s,\\
\dsize
V_r-V_r^2\frac {xQ_{2r-2}(x)\,T_{AV}^{r-1}R_{AV}^{(k-r+1)}Q_{2k+1-2r}(x)} 
{Q_{2k}(x)},&\text{if }r=s,\\
\dsize
(-1)^{r+s+1}V_rV_s\frac {xQ_{2s-2}(x)\,T_{AV}^{r-1}R_{AV}^{(k-r+1)}Q_{2k+1-2r}(x)} 
{Q_{2k}(x)},&\text{if }r> s,
\endcases
\endmultline
\tag\ED
$$
and
$$
\multline
\sum_{n\ge0}\GF\left(\Cal A_{2n+2}^{(k)}(r\to s);w_{AV}\right)\, x^{2n+1}
\\=
\cases 
\dsize
(-1)^{r+s+1}V_rA_s\frac {xQ_{2r-2}(x)\,T_{AV}^{s}Q_{2k-2s}(x)}
{Q_{2k}(x)},&\text{if }r\le s,\\
\dsize
(-1)^{r+s+1}V_rA_s\frac {xQ_{2s-1}(x)\,T_{AV}^{r-1}R_{AV}^{(k-r+1)}Q_{2k+1-2r}(x)}
{Q_{2k}(x)},&\text{if }r> s,
\endcases
\endmultline
\tag\EE
$$
where the operator $T_{AV}$, when applied
to a polynomial in the $A_i$'s and $V_i$'s, 
replaces each $A_i$ by $A_{i+1}$ and each $V_i$ by $V_{i+1}$,
for all~$i$, while the operator $R_{AV}^{(j)}$ 
replaces $A_i$ by $V_{j+1-i}$ and $V_i$ by $A_{j+1-i}$ for all~$i$.\footnote{%
While it seems that such a substitution may produce non-positively indexed 
$A_{j+1-i}$ and $V_{j+1-i}$, this is in fact not so since, in all situations
where the operator $R_{AV}^{(j)}$ is applied, it is applied to~$A_i$'s and~$V_i$'s
with $i\le j$.}
Here, as before,
$\Cal A_n^{(k)}(r\to s)$ is the set of alternating sequences
$r\le a_2\ge a_3\le a_4\ge \cdots {}\diamond\,a_{n-1}{}\square\,s,$
where $\diamond=\,\ge$ and $\square=\,\le$ if $n$ is even and
$\diamond=\,\le$ and $\square=\,\ge$ if $n$ is odd,
in which all $a_i$'s are integers between~$1$ and~$k$.
\endproclaim

\demo{Sketch of proof}
The proof of Theorem~\UA\ via the bijection of Lemma~\UC\ between
alternating sequences and heaps of segments also works here.
In particular, this bijection will be weight-preserving with respect
to $w_{AV}$ if we declare the weight\footnote{Here also, our definition of weight
of a segment is not consistent with the definition of weight of a segment
in Section~3 in that, here, we do not
include $x^2$ --- which keeps track of the size of an alternating sequence,
respectively of a heap --- in the weight. This choice makes the exposition
in Section~10 more convenient. Again, we hope that the reader bears with us.}
of a segment $j${}$\pmb-${}$i$ to be 
$w_{AV}(j${}$\pmb-${}$i):=V_iA_j$, and the weight $w_{AV}(H)$ of a heap $H$
with marked vertex $(0,r)$ (corresponding to the first element $r$
of the associated alternating sequence; 
see Figure~\FB) to be $V_r$ times the 
product of all the weights $w_{AV}(s)$ of segments~$s$ of~$H$.

Here we have to extend Lemmas~\UD\ and \UE\ to our weighted setting.
In generalisation of Lemma~\UD, we claim that
$$\sum_{T\text{ trivial heap on }[1,k]} (-1)^{\vert T\vert}w_{AV}(T)\,x^{2\vert T\vert}=
(-1)^{k}Q_{2k}(x).
\tag\EF$$
Proceeding again by induction, and recalling the three possibilities
(a)--(c) in the proof of Lemma~\UD, the sum of the contributions of
these three possible sets of trivial heaps is
$$
\align
V_{k+1}A_{k+1}(-x^2)&(-1)^kQ_{2k}(x)+(-1)^kQ_{2k}(x)
+\frac {A_{k+1}} {A_k}\left((-1)^kQ_{2k}(x)-(-1)^{k-1}Q_{2k-2}(x)\right)\\
&=
(-1)^{k+1}\big(
V_{k+1}A_{k+1}x^2Q_{2k}(x)-Q_{2k}(x)-A_{k+1}xQ_{2k-1}(x)
\big)\\
&=
(-1)^{k+1}\big(
A_{k+1}xQ_{2k+1}(x)-Q_{2k}(x)
\big)\\
&=
(-1)^{k+1}Q_{2k+2}(x),
\endalign
$$
where, during the computation, we applied (\EB) for $n=2k-1,2k,2k+1$,
in that order. This establishes the induction step.

In order to extend Lemma~\UE\ to our weighted setting, we start with
the observation that
$$
R_{AV}^{(j)}Q_{2j}(x)=Q_{2j}(x),\quad \text{for }j\ge0.
\tag\EG$$ 
This is not obvious at all from the definition~(\EB) of the
polynomials $Q_n(x)$, but it is obvious from~(\EF) by vertically reflecting the heaps.
The same reasoning yields
$$
R_{AV}^{(j+1)}T_{AV}Q_{2j}(x)=Q_{2j}(x),\quad \text{for }j\ge0.
\tag\EH$$ 
Recalling (\BId), we have to compute the sum
$$
-\underset T\not\subseteq[1,r-1]\cup[s+1,k]
\to{\underset T\subseteq[1,r-1]\cup[s,k]
\to{\sum_{T\text{ trivial heap}}}}
(-1)^{\vert T\vert}w_{AV}(T)\,x^{2\vert T\vert}.
$$
By the analogue of the computation (\BIa) in the weighted setting,
this sum is equal to
$$\align
-(-1)^{r-1}&Q_{2r-2}(x)\left((-1)^{k-s+1}T_{AV}^{s-1}Q_{2k-2s+2}(x)
-(-1)^{k-s}T_{AV}^sQ_{2k-2s}(x)\right)\\
&=(-1)^{k+r+s+1}Q_{2r-2}(x)\,T_{AV}^{s-1}\big(Q_{2k-2s+2}(x)
+T_{AV}Q_{2k-2s}(x)\big)\\
&=(-1)^{k+r+s+1}Q_{2r-2}(x)\\
&\kern1.5cm
\times
T_{AV}^{s-1}R_{AV}^{(k-s+1)}
\big(R_{AV}^{(k-s+1)}Q_{2k-2s+2}(x)
+R_{AV}^{(k-s+1)}T_{AV}Q_{2k-2s}(x)\big)\\
&=(-1)^{k+r+s+1}Q_{2r-2}(x)\,T_{AV}^{s-1}R_{AV}^{(k-s+1)}
\big(Q_{2k-2s+2}(x)+Q_{2k-2s}(x)\big)\\
&=(-1)^{k+r+s+1}Q_{2r-2}(x)\,T_{AV}^{s-1}R_{AV}^{(k-s+1)}
A_{k-s+1}x\,Q_{2k-2s+1}(x)\\
&=(-1)^{k+r+s+1}V_sx\,Q_{2r-2}(x)\,T_{AV}^{s-1}R_{AV}^{(k-s+1)}
Q_{2k-2s+1}(x).
\tag\EHa
\endalign
$$
Here we used (\EG) and (\EH) to obtain the third equality and the
defining two-term recurrence (\EB) to obtain the fourth equality.

These arguments prove (\ED).

\medskip
For the proof of the first case in (\EE), we use the weighted version of the
recurrence~(\BIc), 
$$
\GF\left(\Cal A_{2n+2}^{(k)}(r\to s);w_{AV}\right)
=\sum_{j=1}^sA_s\GF\left(\Cal A_{2n+1}^{(k)}(r\to j);w_{AV}\right).
$$
The recurrence implies
$$
\sum_{n\ge0}\GF\left(\Cal A_{2n+2}^{(k)}(r\to s);w_{AV}\right)\, x^{2n+1}=
A_s\sum_{j=1}^s\sum_{n\ge0}\GF\left(\Cal A_{2n+1}^{(k)}(r\to j);w_{AV}\right)
x^{2n+1}.
$$
Let $r>s$. 
By the second and third case of (\ED), we obtain
$$
\align
\sum_{n\ge0}\GF&\left(\Cal A_{2n+2}^{(k)}(r\to s);w_{AV}\right)\, x^{2n+1}\\
&=
A_s\sum_{j=1}^sxV_rV_j(-1)^{r+j+1}
\frac {xQ_{2j-2}(x)\,T_{AV}^{r-1}R_{AV}^{(k-r+1)}Q_{2k+1-2r}(x)} 
{Q_{2k}(x)}\\
&=(-1)^{r}V_rA_s\frac {xT_{AV}^{r-1}R_{AV}^{(k-r+1)}Q_{2k+1-2r}(x)} 
{Q_{2k}(x)}
\sum_{j=1}^s(-1)^{j-1}V_jxQ_{2j-2}(x)\\
&=(-1)^{r}V_rA_s\frac {xT_{AV}^{r-1}R_{AV}^{(k-r+1)}Q_{2k+1-2r}(x)} 
{Q_{2k}(x)}
\sum_{j=1}^s(-1)^{j-1}\big(Q_{2j-1}(x)+Q_{2j-3}(x)\big)\\
&=(-1)^{r}V_rA_s\frac {xT_{AV}^{r-1}R_{AV}^{(k-r+1)}Q_{2k+1-2r}(x)} 
{Q_{2k}(x)}
(-1)^{s-1}Q_{2s-1}(x),
\endalign$$
where we used (\EB) to obtain the next-to-last equality.
Thus, we obtain exactly the second expression on the right-hand side
of~(\EE). 

\medskip
Now let $r\le s$. We proceed as in the proof of Theorem~\UA, where
we mapped the alternating sequences~$S$ in $\Cal A_{2n+2}^ {(k)}(r\to s)$
bijectively to heaps~$H$ of $n$~segments on $[1,k]$ whose maximal segments
are not contained in $[1,r-1]\cup[s+1,k]$. This bijection is weight-preserving
in the sense that $w_{AV}(S)=V_rA_sw_{AV}(H)$. Then, by (\BIe), with the weight $w$
being defined by $w_{AV}(T)x^{2\vert T\vert}$, and by~(\EF), we obtain the first
expression on the right-hand side of~(\EE) straightforwardly.

\medskip
This completes the proof of the theorem.\quad \quad \qed
\enddemo

As consequences of the above theorem, we can now 
give the weighted generalisations of Corollaries~\UG\ and \UH.
We omit the proofs since they are identical with those of these
corollaries, except that they use Theorem~\TQ\ instead of Theorem~\UA.

\proclaim{Corollary \TR}
For all positive integers $k$, we have
$$
\sum_{n\ge1}\GF\left(\Cal A_{2n-1}^{(k)};w_{AV}\right)\, x^{2n}
=-\frac {xQ_{2k-1}(x)} {Q_{2k}(x)}.
\tag\EI
$$
\endproclaim

\proclaim{Corollary \TS}
For all positive integers $k$, we have
$$
\sum_{n\ge0}\GF\left(\Cal A_{2n}^{(k)};w_{AV}\right)\, x^{2n+1}
=(-1)^k\frac {x} {Q_{2k}(x)}.
\tag\EJ
$$
\endproclaim

For later use,
we work out one more special case of Theorem~\TQ. Namely, in
Section~12.(5) we shall need the special case of~(\ED) in which
$r=s=1$,
$$
\sum_{n\ge0}\GF\left(\Cal A_{2n+1}^{(k)}(1\to 1);w_{AV}\right)\, x^{2n}
=
-V_1^2\frac {x\,R_{AV}^{(k)}Q_{2k-1}(x)} 
{Q_{2k}(x)},
\tag\EJa$$
in which $A_i=(yq)^i$ and $V_i=y^{-i}$ for $i\ge1$.
We claim that we have
$$
Q_{2k}(x)\bigg\vert_{\smallmatrix A_{i}=(yq)^i
\\V_{i}=y^{-i}\kern.2cm\endsmallmatrix}
=1+\sum_{j=1}^k(-1)^{k-j}x^{2j}q^{\binom {j+1}2}
\sum_{i=0}^{k-j}(yq)^i\bmatrix k-i\\j\endbmatrix_q
\bmatrix i+j-1\\j-1\endbmatrix_q,
\tag\EJb$$
and
$$
Q_{2k-1}(x)\bigg\vert_{\smallmatrix A_{i}=(yq)^i
\\V_{i}=y^{-i}\kern.2cm\endsmallmatrix}
=
\sum_{j=1}^k(-1)^{k-j}x^{2j-1}q^{\binom {j}2}
\sum_{i=0}^{k-j}y^{-k+i}\bmatrix k-i-1\\j-1\endbmatrix_q
\bmatrix i+j-1\\j-1\endbmatrix_q.
$$
This is easily verified by induction on~$k$ using the defining
recurrence~(\EB). However, what we need in~(\EJa) is not the above
specialisation of $Q_{2k-1}(x)$, but rather that specialisation
of $R_{AV}^{(k)}Q_{2k-1}(x)$. The key to obtain the latter
specialisation is the identity
$$
V_1xR_{AV}^{(k)}Q_{2k-1}(x)=Q_{2k}(x)+T_{AV}Q_{2k-2}(x)
$$
that we (implicitly) established in~(\EHa) (choose $r=s=1$ there).
Since the total degree in the~$A_i$'s in the coefficient of~$x^{2j}$
in $Q_{2k}(x)$ is~$j$ for all~$j$, and the same is true for the total degree in
the~$V_i$'s, the specialisation
$$
T_{AV}Q_{2k-2}(x)\bigg\vert_{\smallmatrix A_{i}=(yq)^i
\\V_{i}=y^{-i}\kern.2cm\endsmallmatrix}
=
Q_{2k-2}(x)\bigg\vert_{\smallmatrix A_{i}=(yq)^{i+1}
\\V_{i}=y^{-i-1}\kern.2cm\endsmallmatrix}
$$
can be easily obtained from~(\EJb) by replacing~$x$ by $xq^{1/2}$.
After simplification, we obtain
$$
R_{AV}^{(k)}Q_{2k-1}(x)\bigg\vert_{\smallmatrix A_{i}=(yq)^i
\\V_{i}=y^{-i}\kern.2cm\endsmallmatrix}
=
y\sum_{j=1}^k(-1)^{k-j}x^{2j-1}q^{\binom {j+1}2}
\sum_{i=0}^{k-j}(yq)^i\bmatrix k-i-1\\j-1\endbmatrix_q
\bmatrix i+j-1\\j-1\endbmatrix_q.
\tag\EJc$$

\subhead 9. Generating functions for bounded paths ``with negative
length" and bounded alternating sequences\endsubhead
Our next goal is to find the weighted generalisations of
Corollaries~\UF--\TC\ on ``numbers of paths with negative length".
We shall abuse notation here: when, for a {\it positive} integer~$n$ we write
$$
\GF\left(\Cal C_{-n}^{(k)}(r\to s);w_B\right) ,
\tag\EK
$$
then we mean the value of the corresponding term in the continuation
of the sequence
$$
\left(\GF\left(\Cal C_{n}^{(k)}(r\to s);w_B\right)\right)_{n\ge0} ,
$$
to negative integers, using the linear recurrence that it satisfies in
the backward direction. (The above sequence satisfies indeed a linear
recurrence since its generating function is rational; cf\. Theorem~\TP.)

In view of (\EBa), we introduce the weight $w_{BAV}$ on up-down paths
by what we get from $w_B$ when we do the replacements
$B_{2i-1}=A_i^{-1}V_i^{-1}$ and $B_{2i}=A_i^{-1}V_{i+1}^{-1}$ for all~$i$.
Explicitly, all up-steps are assigned a weight of~$1$, while a 
down-step from height~$2i-1$ to height~$2i-2$ is assigned a weight of 
$A_i^{-1}V_i^{-1}$, and a 
down-step from height~$2i$ to height~$2i-1$ is assigned a weight of 
$A_i^{-1}V_{i+1}^{-1}$.

We then have the following weighted extension of Corollary~\UF.

\proclaim{Corollary \TT}
Let $n,k,r,s$ be positive integers with $1\le r,s\le k$.
We have
$$\multline
(-1)^{r+s}\Bigg(\prod _{i=r} ^{s-1}A_i^ {-1}\Bigg)\Bigg(
\prod _{i=r+1} ^{s-1}V_i^{-1}\Bigg)
\GF\left(\Cal C_{-2n}^{(2k-1)}(2r-2\to 2s-2);w_{BAV}\right)
\\=
\GF\left(\Cal A_{2n+1}^{(k)}(r\to s);w_{AV}\right) .
\endmultline
\tag\EL
$$
Furthermore, we have
$$\multline
(-1)^{r+s}\Bigg(\prod _{i=r} ^{s-1}A_i^{-1}\Bigg)\Bigg(
\prod _{i=r+1} ^{s}V_i^{-1}\Bigg)
\GF\left(\Cal C_{-2n+1}^{(2k-1)}(2r-2\to
2s-1);w_{BAV}\right)\\
=
\GF\left(\Cal A_{2n}^{(k)}(r\to s);w_{AV}\right).
\endmultline
\tag\EM$$ 
Here and in the sequel, the
products have to be interpreted in an extended sense by
$$\prod _{k=R} ^{S-1}\text {\rm Expr}(k)=\cases \hphantom{-}
\prod _{k=R} ^{S-1} \text {\rm Expr}(k),&R<S,\\
\hphantom{-}1,&R=S,\\
\prod _{k=S} ^{R-1}\big(\text {\rm Expr}(k)\big)^{-1},&R>S.\endcases
\tag\EMa$$
\endproclaim

\demo{Proof}
We do the replacements $k\to 2k-1$, $r\to 2r-2$, $s\to 2s-2$, and we set
$B_{2i-1}=A_i^{-1}V_i^{-1}$ and $B_{2i}=A_i^{-1}V_{i+1}^{-1}$ for all~$i$
in Theorem~\TP. 
Subsequently, in order to obtain the generating function for
the quantities (\EK), we apply (\BDa). By (\EBa) and 
the easily verified relations
$$
T_B^2P_n(x)\Bigg\vert_{\smallmatrix B_{2i-1}=A_i^{-1}V_i^{-1}
\kern.3cm
\\B_{2i}=A_i^{-1}V_{i+1}^{-1}
\endsmallmatrix}=
T_{AV}\left(P_n(x)\Bigg\vert_{\smallmatrix B_{2i-1}=A_i^{-1}V_i^{-1}
\kern.3cm
\\B_{2i}=A_i^{-1}V_{i+1}^{-1}
\endsmallmatrix}\right)
\tag\EBb
$$
and
$$
T_BP_{2n+1}(x)\Bigg\vert_{\smallmatrix B_{2i-1}=A_i^{-1}V_i^{-1}
\kern.3cm
\\B_{2i}=A_i^{-1}V_{i+1}^{-1}
\endsmallmatrix}=
R_{AV}^{(n+1)}\left(P_{2n+1}(x)
\Bigg\vert_{\smallmatrix 
B_{2i-1}=A_i^{-1}V_i^{-1}\kern.3cm
\\B_{2i}=A_i^{-1}V_{i+1}^{-1}
\endsmallmatrix}\right),
\tag\EN
$$
we obtain after some simplification
$$\align
\sum_{n\ge1}\GF&\left(\Cal C_{-2n}^{(2k-1)}(2r-2\to 2s-2);w_{BAV}\right)\, x^{2n}\\
&=\cases -\dfrac {
A_r\cdots A_{s-1}
V_r\cdots V_s
x\,Q_{2r-2}(x)\,T_{AV}^{s-1}R_{AV}^{(k-s+1)}Q_{2k-2s+1}(x)}
{Q_{2k}(x)},
&\text{if }r\le s,\\
-\dfrac {x\,Q_{2s-2}(x)\,T_{AV}^{r-1}R_{AV}^{(k-r+1)}Q_{2k-2r+1}(x)}
{A_s\cdots A_{r-1}V_{s+1}\cdots V_{r-1}
Q_{2k}(x)},
&\text{if }r\ge s,
\endcases\\
&=
(-1)^{r+s}\Bigg(\prod _{i=r} ^{s-1}A_i\Bigg)\Bigg(
\prod _{i=r+1} ^{s-1}V_i\Bigg)
\sum_{n\ge1}\GF\left(\Cal A_{2n+1}^{(k)}(r\to s);w_{AV}\right)\, x^{2n},
\endalign
$$
where we appealed to Theorem~\TQ\ to get the last line. 
Comparison of coefficients of $x^{2n}$ then establishes the
relation~(\EL).

For proving (\EM), we proceed in the same way. 
Here, we do the replacements $k\to 2k-1$, $r\to 2r-2$, $s\to 2s-1$, and
we set $B_{2i-1}=A_i^{-2}$ and $B_{2i}=A_i^{-1}A_{i+1}^{-1}$ for all~$i$
in Theorem~\TP. 
The generating function identity that we obtain here is
$$\align
\sum_{n\ge1}\GF&\left(\Cal C_{-2n+1}^{(2k-1)}(2r-2\to
2s-1);w_{BAV}\right)\, x^{2n-1}\\
&=\cases -\dfrac {A_r\cdots A_sV_r\cdots V_s 
x\,Q_{2r-2}(x)\,T_{AV}^{s}Q_{2k-2s}(x)} 
{Q_{2k}(x)},
&\text{if }r\le s,\\
-\dfrac {x\,Q_{2s-1}(x)\,T_{AV}^{r-1}R_{AV}^{(k-r+1)}Q_{2k-2r+1}(x)}
{A_{s+1}\cdots A_{r-1}V_{s+1}\cdots V_{r-1}Q_{2k}(x)},
&\text{if }r\ge s,
\endcases
\\
&=
(-1)^{r+s}\Bigg(\prod _{i=r} ^{s-1}A_i\Bigg)\Bigg(
\prod _{i=r+1} ^{s}V_i\Bigg)
\sum_{n\ge1}\GF\left(\Cal A_{2n}^{(k)}(r\to s);w_{AV}\right)\, x^{2n-1},
\endalign
$$
from which (\EM) follows immediately upon comparison of coefficients
of~$x^{2n-1}$.\quad \quad \qed
\enddemo

By specialising the above corollary appropriately, we obtain the
weighted extensions of Corollaries~\TB\ and~\TC.

\proclaim{Corollary \TU}
For positive integers $n$ and $k$, we have
$$
A_k^{-1}R_{AV}^{(k)}
\GF\left(\Cal C_{-2n}^{(2k-1)};w_{BAV}\right)
=
\GF\left(\Cal A_{2n-1}^{(k)};w_{AV}\right) .
\tag\EO
$$
\endproclaim

\demo{Proof}
We put $r=s=1$ in (\EL). Then, on the right-hand side, we have the
generating function for sequences $1\le a_2\ge a_3\le
a_4\ge \cdots \le a_{2n}\ge 1$. By skipping the $1$'s at the beginning
and at the end, and replacing each $a_i$ by
$a_{k+1-i}$, we obtain a sequence in $\Cal A_{2n-1}^{(k)}$.
However, in terms of the weight $w_{AV}$, this replacement amounts to
an application of the operator $R_{AV}^{(k)}$.\quad \quad \qed
\enddemo

\proclaim{Corollary \TV}
Let $n$ and $k$ be positive integers. Furthermore, denote the set of
up-down paths from $(0,0)$ to $(2n+2k-1,2k-1)$ of height at
most~$2k-1$ that do not pass below the $x$-axis
by $\Cal D_{2n}^{(2k-1)}$. Then we have
$$
(-1)^{k+1}\Bigg(\prod _{i=1} ^{k}A_i^{-1}V_i^{-1}\Bigg)
\GF\left(\Cal D_{-2n-2k}^{(2k-1)};w_{BAV}\right)
=
\GF\left(\Cal A_{2n}^{(k)};w_{AV}\right).
\tag\EP$$ 
\endproclaim

\demo{Proof}
In the same way as the previous corollary, by putting
$r=1$ and $s=k$ in (\EM), we get
$$
(-1)^{k+1}\Bigg(\prod _{i=1} ^{k}A_i^{-1}V_i^{-1}\Bigg)
R_{AV}^{(k)}\GF\left(\Cal D_{-2n-2k}^{(2k-1)};w_{BAV}\right)
=
\GF\left(\Cal A_{2n}^{(k)};w_{AV}\right),
$$
or, equivalently,
$$
(-1)^{k+1}\Bigg(\prod _{i=1} ^{k}A_i^{-1}V_i^{-1}\Bigg)
\GF\left(\Cal D_{-2n-2k}^{(2k-1)};w_{BAV}\right)
=R_{AV}^{(k)}
\GF\left(\Cal A_{2n}^{(k)};w_{AV}\right).
$$
Now, by (\EG), we know that the operator $R_{AV}^{(k)}$ leaves $Q_{2k}(x)$
invariant. If we use this in combination with Corollary~\TS, then we see that
$R_{AV}^{(k)}$ also leaves $\GF\left(\Cal A_{2n}^{(k)};w_{AV}\right)$ invariant.
Hence, the above identity can be simplified to~(\EP).\quad \quad \qed
\enddemo

\subhead 10. Weighted reciprocity theorems\endsubhead
We now turn to the weighted versions of the reciprocity laws
in Sections~5--7. We begin with the weighted generalisation of
Theorem~\TD.

\proclaim{Theorem \TW}
For all non-negative integers $n,k,m$, we have
$$
\multline
\det\left(\GF\left(\Cal C_{2n+2i+2j+4m-2}^{(2k+2m-1)};w_{BAV}\right)
\right)_{0\le i,j\le k-1}\\
=
\Bigg(V_1^{k}
A_{k+m}^{-m}
\prod _{i=1} ^{k+m}A_i^{-(n+2m+2k-2i-1)}V_i^{-(n+2m+2k-2i)}\Bigg)\\
\times
\det\left(R_{AV}^{(k+m)}\GF\left(\Cal C_{-2n-2i-2j}^{(2k+2m-1)};w_{BAV}\right)
\right)_{0\le i,j\le m-1}.
\endmultline
\tag\EQ
$$
\endproclaim

\remark{Remark}
A notable special case arises if we choose $A_i=V_i=q^{-i+1}$, $i\ge1$.
The reader should recall that the weight $w_{BAV}$ resulted from the
weight $w_B$ under the substitutions 
$B_{2i-1}=A_i^{-1}V_i^{-1}$ and $B_{2i}=A_i^{-1}V_{i+1}^{-1}$,
$i\ge1$. Thus, the above choice of the $A_i$'s and $V_i$'s implies
the choice $B_i=q^{i-1}$, $i\ge1$. Now we should remember
that for this choice of the $B_i$'s the
weight~$w_B$ reduces to the area weight~$q^{a(\,.\,)}$, with the
corresponding generating function for the terms
$\GF\big(\Cal C_n^{(k)};q^{a(\,.\,)}\big)$ given in~(\ECa).
After considerable simplification, we obtain the reciprocity law
$$
\multline
\det\left(\GF\left(\Cal C_{2n+2i+2j+4m-2}^{(2k+2m-1)};q^{a(\,.\,)}\right)
\right)_{0\le i,j\le k-1}\\
=
q^{2\binom k2\left(\frac {4k+1} {6}+n+2m-2\right)
-2\binom m2\left(n-1+\frac {8m-1} {6}\right)}
\det\left(\GF\left(\Cal C_{-2n-2i-2j}^{(2k+2m-1)};q^{-a(\,.\,)}\right)
\right)_{0\le i,j\le m-1}.
\endmultline
\tag\EQa
$$
\endremark

For the proof of Theorem~\TW, we follow the arguments of the proof of
Theorem~\TD. Namely, the
Lindstr\"om--Gessel--Viennot theorem \cite{\LindAA, Lemma~1} shows
that the determinant on the left-hand side of~(\EQ) is the generating
function for families $\Cal P=(P_0,P_1,\dots,P_{k-1})$ of non-intersecting Dyck
paths of height at most $2k+2m-1$, where $P_i$ runs from
$(-2i,0)$ to $(2n+4m+2i-2,0)$, $i=0,1,\dots,k-1$, with respect to the
weight~$w_{BAV}$. Here, the weight $w_{BAV}(\Cal P)$ is defined as
$\prod _{i=0} ^{k-1}w_{BAV}(P_i)$. 
By the bijection in Proposition~\TE, the families $\Cal P$ of
non-intersecting Dyck paths are related to trapezoidal arrays of
alternating sequences. How the weights are related under this
bijection is spelled out in Lemma~\TX\ below. The weight $w_{BAV}$
for families of Dyck paths has been just explained, while the
weight $w_{AV}(\Cal A)$ of a trapezoidal array $\Cal A=(a_{i,j})$
is defined by $\prod _{i,j} ^{}A_{a_{i,2j}}V_{a_{i,2j-1}}$, 
that is, entries in
even-indexed columns contribute $A$-variables, while entries in
odd-indexed columns contribute $V$-variables. Finally,
the arguments in the proof of Theorem~\TF\ together with 
Corollary~\TU, the weighted version of Corollary~\TB, show that,
up to some scaling, the determinant on the right-hand side is
the generating function for these trapezoidal arrays of
alternating sequences with respect to the weight~$w_{AV}$,
see Theorem~\TY.

\proclaim{Lemma \TX}
Let $n$ be a non-negative integer and $k,m$ be positive integers.
Furthermore, let $\Cal P=(P_0,P_1,\dots,P_{k-1})$
be a family of non-intersecting Dyck
paths of height at most $2k+2m-1$, where $P_i$ runs from
$(-2i,0)$ to $(2n+4m+2i-2,0)$, $i=0,1,\dots,k-1$, and $\Cal A$ the 
trapezoidal array of integers of the form~{\rm(\AB)} that
corresponds to~$\Cal P$ under the bijection of Proposition~{\rm\TE}.
Then
$$
\Bigg(V_1^{-k}
\prod _{i=1} ^{k+m}A_i^{-(n+2m+2k-2i-1)}V_i^{n+2m+2k-2i}\Bigg)
w_{BAV}(\Cal P)=w_{AV}(\Cal A).
\tag\ER$$
\endproclaim

\demo{Sketch of proof}
The reader should recall that the bijection of Proposition~\TE\
worked as follows: in the region~(\AEa) all the lattice points
$(x,y)$ with $x$ and $y$ of the same parity, which are not
occupied by any of the Dyck paths of~$\Cal P$,
contribute an entry
to the trapezoidal array~$\Cal A$; more precisely, such a lattice point
contributes the entry $\cl{(y+1)/2}$. If we now look at the weight $w_{AV}$ of
arrays, then we see that such a lattice
point contributes $V_{{(y+2)/2}}$ to the weight $w_{AV}(\Cal A)$
of the array~$\Cal A$ if $y$ is even, while it contributes
$A_{{(y+1)/2}}$ if $y$ is odd.

In order to see how the weights $w_{BAV}(\Cal P)$ and $w_{AV}(\Cal A)$ are
related, we change $w_{BAV}$ to the equivalent weight~$\bar w_{BAV}$.
This weight~$\bar w_{BAV}$ is defined by assigning a weight of
$A_{y}^{-1}$ to up-steps from height~$2y-2$ to height~$2y-1$
and to down-steps from height~$2y$ to height~$2y-1$, and assigning
a weight of $V_{y+1}^{-1}$ to up-steps from height~$2y-1$ to height~$2y$
and to down-steps from height~$2y+1$ to height~$2y$.
For a Dyck path~$P$
we have indeed $w_{BAV}(P)=\bar w_{BAV}(P)$ since up- and down-steps
in Dyck paths can be paired up (cf\. 
Footnote~6), and a matching
pair of steps $(x_1,y-1)\to(x_1+1,y)$ and $(x_2,y)\to(x_2+1,y-1)$
contributes a weight of 
$A_{\fl{(y+1)/2}}^{-1}V_{\cl{(y+1)/2}}^{-1}$ to $\bar w_{BAV}(P)$
as well as to $w_{BAV}(P)$. 

We should think of the weight $A_{{(y+1)/2}}^{-1}$,
respectively $V_{{(y+2)/2}}^{-1}$, as attached to
the end point $(x,y)$ of a step (up or down). More precisely,
due to the fact that $x$~and~$y$ must have the same parity, the weight
$A_{{(y+1)/2}}^{-1}$ is attached to $(x,y)$ if $x$~and~$y$~are odd, and if
$x$~and~$y$~are even then it is $V_{{(y+2)/2}}^{-1}$ which is attached to it. 
It then becomes
apparent that the weighting $\bar w_{BAV}$ is chosen so as to ``kill"
the contribution of $(x,y)$ to the array weight $w_{AV}$. 
Hence, to verify the relation~(\ER), all that remains is to determine
the total contribution of lattice points $(x,y)$ with $x$ and $y$ of
the same parity in the region~(\AEa). A small detail is that, in this
count, we have to leave out the points $(-2i,0)$, $i=0,\dots,k-1$,
since all these points are occupied by paths (being actually the
starting points of the paths) but are not end points of
steps.\quad \quad \qed
\enddemo

\proclaim{Theorem \TY}
Let $n$ be a non-negative integer and $k,m$ be positive integers.
The generating function $\sum_{\Cal A}w_{AV}(\Cal A)$, where the sum is
over all trapezoidal arrays of integers of the form {\rm(\AB)}
that satisfy {\rm(\AC)--(\AE)}, and the additional restrictions
mentioned in the statement of Proposition~{\rm\TE} in the case where
$n=0$, is equal to
$$
A_{k+m}^{-m}
\det\left(R_{AV}^{(k+m)}\GF\left(\Cal C_{-2n-2i-2j}^{(2k+2m-1)};w_{BAV}\right)
\right)_{0\le i,j\le m-1}.
\tag\ES
$$
\endproclaim

\demo{Proof}
It should be noted that, by Corollary~\TU, the determinant in
(\AF) can be alternatively written as
$$
\det\left(\GF\left(\Cal A_{2n+2i+2j-1}^{(k+m)};w_{AV}\right)
\right)_{0\le i,j\le m-1}.
$$
The proof of Theorem~\TF\ can then be used verbatim, the only
difference being that here the paths in the application of
the Lindstr\"om--Gessel--Viennot theorem carry weights.\quad \quad \qed
\enddemo

Following the same line of argument, we obtain the following 
weighted generalisation of Theorem~\TG. We content ourselves with
stating the result together with the auxiliary result addressing
the weight property that the bijection in Proposition~\TH\ satisfies
and the relevant generating function result for arrays of alternating
sequences, see Lemma~\TTA\ and Theorem~\TTB\ below.
We point out a subtlety in the proof of Lemma~\TTA\ though:
since the paths are not Dyck paths but rather paths starting at
height~$0$ and ending at height~$2k+2m-1$, not every up-step
has a matching down-step, and therefore the weights $w_{BAV}$
and $\bar w_{BAV}$ do not agree on these paths. More precisely,
for such a path~$P$ we have
$$
w_{BAV}(P)=\Bigg(V_1^{-1}
\prod _{i=1} ^{k+m}A_iV_i\Bigg)\bar w_{BAV}(P).
\tag\ESa
$$
A second subtlety in this context is that the invariance of
$\GF\left(\Cal A_{2n}^{(k+m)};w_{AV}\right)$
under the action of the operator $R_{AV}^{(k+m)}$, which was already used
in the proof of Corollary~\TV, is also applied here.
We leave the remaining details of the proofs to the reader.

\proclaim{Theorem \TZ}
With the sets $\Cal D_{2n}^{(2k+2m-1)}$ explained in Corollary~{\rm\TV},
for all non-negative integers $n$ and positive integers $k,m$, we have
$$
\multline
\det\left(\GF\left(\Cal D_{2n+2j-2i}^{(2k+2m-1)};w_{BAV}\right)
\right)_{0\le i,j\le k-1}
=
(-1)^{km}\Bigg(
\prod _{i=1} ^{k+m}A_i^{-n-m}V_i^{-n-m}\Bigg)\\
\times
\det\left(
\GF\left(\Cal D_{-2n-2j+2i-2k-2m}^{(2k+2m-1)};w_{BAV}\right)
\right)_{0\le i,j\le m-1}.
\endmultline
\tag\ET
$$
\endproclaim

\proclaim{Lemma \TTA}
Let $n$ be a non-negative integer and $k,m$ be positive integers.
Furthermore, let $\Cal P=(P_0,P_1,\dots,P_{k-1})$
be a family of non-intersecting up-down
paths of height at most $2k+2m-1$ that do not pass below the $x$-axis, 
where $P_i$ runs from
$(2i,0)$ to $(2n+2m+2k+2i-1,2k+2m-1)$, $i=0,1,\dots, k-1$, 
and $\Cal A$ the 
rhomboidal array of integers of the form~{\rm(\CB)} that
corresponds to~$\Cal P$ under the bijection of Proposition~{\rm\TH}.
Then
$$
\Bigg(
\prod _{i=1} ^{k+m}A_i^{n}V_i^{n}\Bigg)
w_{BAV}(\Cal P)=w_{AV}(\Cal A).
\tag\EU$$
\endproclaim

\proclaim{Theorem \TTB}
Let $n$ be a non-negative integer and $k,m$ be positive integers.
The generating function $\sum_{\Cal A}w_{AV}(\Cal A)$, where the sum is
over all rhomboidal arrays of integers of the form {\rm(\CB)}
that satisfy {\rm(\CC)--(\CE)} is equal to
$$
(-1)^{km}\Bigg(\prod _{i=1} ^{k+m}A_i^{-m}V_i^{-m}\Bigg)
\det\left(\GF\left(\Cal D_{-2n-2j+2i-2k-2m}^{(2k+2m-1)}
;w_{BAV}\right)
\right)_{0\le i,j\le m-1}.
\tag\EV
$$
\endproclaim

It is interesting to note that Theorem~\TZ\ implies a (seemingly) more
general reciprocity result for Hankel determinants of sequences
which satisfy a linear recurrence with constant coefficients, which we state below.

\proclaim{Theorem \TTTTTT}
Let $k$ and $m$ be positive integers, and let
$a_0,a_1,\dots,a_{k+m}$ be given constants. Furthermore,
let the sequence $\big(c(n)\big)_{n\in\Bbb Z}$
{\rm(}with $\Bbb Z$ denoting the set of all integers{\rm)} be defined by the recurrence
$$
\sum_{i=0}^{k+m}a_i\,c(n-i)=0,
\tag\EVa
$$
with initial conditions
$c(n)=0$ for $-k-m+1\le n\le -1$. It is understood that the recurrence is
read both ways, that is, in the positive direction as well as in the
negative direction. Then, for all non-negative integers~$n$, we have
$$\multline
\det\big(c(n+i+j)\big)_{0\le i,j\le k-1}\\
=(-1)^{\binom {n+k+m}2+\binom n2}c^{k-m}(0)\left(\frac {a(k+m)} {a(0)}\right)^{n+k+m-1}
\det\big(c(-n-i-j-2k)\big)_{0\le i,j\le m-1}.\kern.5cm
\endmultline
\tag\EVb
$$
\endproclaim

\demo{Sketch of proof}
We rewrite (\ET) by replacing $i$ by $k-1-i$ on the left-hand side and by replacing
$i$ by $m-1-i$ on the right-hand side, both replacements leaving the corresponding
determinants invariant up to a sign
since they amount to reversals of the order of rows of the matrices of which
the determinants are taken. Subsequently, we replace $n$ by~$n+k-1$. The result is
$$
\multline
\det\left(\GF\left(\Cal D_{2n+2i+2j}^{(2k+2m-1)};w_{BAV}\right)
\right)_{0\le i,j\le k-1}
=
(-1)^{\binom {k+m}2}\Bigg(
\prod _{i=1} ^{k+m}A_i^{-n-m}V_i^{-n-m}\Bigg)\\
\times
\det\left(
\GF\left(\Cal D_{-2n-2i-2j-4k}^{(2k+2m-1)};w_{BAV}\right)
\right)_{0\le i,j\le m-1}.
\endmultline
\tag\EVc
$$
By (\EJ) and (\EU), we know that the generating function 
$$\sum_{n\ge0}\GF\left(\Cal D_{2n}^{(2k+2m-1)};w_{BAV}\right)x^{2n}$$
has the form
$C/Q_{2k+2m}(x)$, where $C$ is some constant and $Q_{2k+2m}(x)$ is the
polynomial defined in~(\EB). We know that it is a polynomial in~$x^2$
of degree~$k+m$. In particular, these considerations imply that the
quantities $\GF\left(\Cal D_{2n}^{(2k+2m-1)};w_{BAV}\right)$ satisfy
a linear recurrence with constant coefficients of the form~(\EVa).
Now we argue --- and this is not difficult to see ---
that the polynomials $Q_{2N}(\sqrt x)$, $N=0,1,\dots$
(these are polynomials in~$x$ of degree~$N$) are {\it generic} polynomials,
by which we mean that any polynomial of degree~$N$ is equal to
$Q_{2N}(\sqrt x)$ with an appropriate choice of the $A_i$'s and $V_i$'s.
However, then~(\EVc) describes a reciprocity relation for general sequences
$\big(c(n)\big)_{n\in\Bbb Z}$ that satisfy~(\EVa). It only remains to
connect the $A_i$'s and $V_i$'s to the constants $a_i$ in the statement of
the theorem to see that~(\EVc) and~(\EVb) are equivalent.
We leave these details to the reader.\quad \quad \qed
\enddemo

Finally, we present the weighted generalisations of the reciprocity
laws in Theorems~\TJ\ and~\TM, see Theorems~\TTC\ and~\TTF\ below, 
together with the relevant facts that are needed in their proofs.
Since the arguments for their proofs are completely analogous to
previous ones, we leave it to the reader to fill in the details.

\proclaim{Theorem \TTC}
Let $n$ be a non-negative integer and $k,m$ be positive integers,
and let $r_0<r_1<\dots<r_{k-1}$
and $s_0<s_1<\dots<s_{k-1}$ be sequences of positive integers with
$1\le r_i,s_i\le k+m$ for all~$i$. Then
$$
\align
\det&\left(\GF\left(\Cal C_{2n}^{(2k+2m-1)}(2r_i-2\to 2s_j-2);w_{BAV}\right)
\right)_{0\le i,j\le k-1}\\
&=(-1)^{\sum_{i=0}^{m-1}(r_i+s_i)}
\Bigg(
\prod _{i=1} ^{k+m}A_i^{-n}V_i^{-n}\Bigg)
\Bigg(
\prod _{j=0} ^{k-1}V_{s_j}V_{r_j}^{-1}
\prod _{i=r_j} ^{s_j-1}A_i^{2}V_i^{2}\Bigg)
\\
&\kern1cm
\times
\det\left(\GF\left(\Cal C_{-2n}^{(2k+2m-1)}(2\bar r_i-2\to 2\bar
s_j-2);w_{BAV}\right)\right)_{0\le i,j\le m-1}
\\
&=(-1)^{\sum_{i=0}^{m-1}(\bar r_i+\bar s_i)}
\Bigg(
\prod _{i=1} ^{k+m}A_i^{-n}V_i^{-n}\Bigg)
\Bigg(
\prod _{j=0} ^{m-1}
V_{\bar r_j}V_{\bar s_j}^{-1}\prod _{i=\bar r_j} ^{\bar s_j-1}A_i^{-2}V_i^{-2}
\Bigg)
\\
&\kern1cm
\times
\det\left(\GF\left(\Cal C_{-2n}^{(2k+2m-1)}(2\bar r_i-2\to 2\bar
s_j-2);w_{BAV}\right)\right)_{0\le i,j\le m-1},
\tag\EW
\endalign
$$
where 
$$
\align
\{\bar r_0,\bar r_1,\dots,\bar r_{m-1}\}&=
\{1,2,\dots,k+m\}\setminus\{r_0,r_1,\dots,r_{k-1}\},\\
\{\bar s_0,\bar s_1,\dots,\bar s_{m-1}\}&=
\{1,2,\dots,k+m\}\setminus\{s_0,s_1,\dots,r_{k-1}\},
\endalign
$$
and where we assume that 
$\bar r_0<\bar r_1<\dots<\bar r_{m-1}$ and
$\bar s_0<\bar s_1<\dots<\bar s_{m-1}$.
As earlier, the products have to be interpreted according to~{\rm(\EMa)}.
\endproclaim

\proclaim{Lemma \TTD}
Let $n$ be a non-negative integer and $k,m$ be positive integers.
Furthermore, let $\Cal P=(P_0,P_1,\dots,P_{k-1})$
be a family of non-intersecting up-down
paths of height at most $2k+2m-1$ that do not pass below the $x$-axis, 
where $P_i$ runs from
$(0,2r_i-2)$ to $(2n,2s_i-2)$, $i=0,1,\dots,k-1$, 
and $\Cal A$ the 
rectangular array of integers of the form~{\rm(\DB)} that
corresponds to~$\Cal P$ under the bijection of Proposition~{\rm\TK}.
Then
$$
\Bigg(
\prod _{i=1} ^{k+m}A_i^{n}V_i^{n+1}\Bigg)
\Bigg(
\prod _{j=0} ^{k-1}V_{s_j}^{-1}
\prod _{i=r_j} ^{s_j-1}A_i^{-1}V_i^{-1}\Bigg)
w_{BAV}(\Cal P)=w_{AV}(\Cal A).
\tag\EX$$
\endproclaim

Given an up-down path $P$ from height $2r-2$ to height $2s-2$, 
the weighted analogue of~(\ESa) that is relevant here is
$$
w_{BAV}(P)=\Bigg(
\prod _{i=r+1} ^{s}A_{i-1}V_i
\Bigg)\bar w_{BAV}(P),
\tag\EXa
$$
where, again, the product has to be understood as explained in~(\EMa).

\proclaim{Theorem \TTE}
Let $n$ be a non-negative integer and $k,m$ be positive integers,
and let $r_0<r_1<\dots<r_{k-1}$
and $s_0<s_1<\dots<s_{k-1}$ be sequences of positive integers with
$1\le r_i,s_i\le k+m$ for all~$i$.
The generating function $\sum_{\Cal A}w_{AV}(\Cal A)$, where the sum is
over all rectangular arrays of integers of the form~{\rm(\DB)}
that satisfy the constraints given in the statement of Proposition~{\rm\TK} 
is equal to
$$
\multline
(-1)^{\sum_{i=0}^{m-1}(\bar r_i+\bar s_i)}
\Bigg(
\prod _{j=0} ^{m-1}
V_{\bar r_j}\prod _{i=\bar r_j} ^{\bar s_j-1}A_i^{-1}V_i^{-1}
\Bigg)\\
\times
\det\left(\GF\left(\Cal C_{-2n}^{(2k+2m-1)}(2\bar r_i-2\to 2\bar
s_j-2);w_{BAV}\right)\right)_{0\le i,j\le m-1}.
\endmultline
\tag\EY
$$
\endproclaim

\proclaim{Theorem \TTF}
Let $n,k,m$ be positive integers, and let $r_0<r_1<\dots<r_{k-1}$
and $s_0<s_1<\dots<s_{k-1}$ be sequences of positive integers with
$1\le r_i,s_i\le k+m$ for all~$i$. Then
$$
\align
\det&\left(\GF\left(\Cal C_{2n-1}^{(2k+2m-1)}(2r_i-2\to 2s_j-1);w_{BAV}\right)
\right)_{0\le i,j\le k-1}\\
&=(-1)^{\sum_{i=0}^{m-1}(r_i+s_i)}
\Bigg(
\prod _{i=1} ^{k+m}A_i^{-n}V_i^{-n}\Bigg)
\Bigg(
\prod _{j=0} ^{k-1}A_{r_j}^2V_{r_j}A_{s_j}^{-1}
\prod _{i=r_j+1} ^{s_j}A_i^{2}V_i^{2}\Bigg)
\\
&\kern1cm
\times
\det\left(\GF\left(\Cal C_{-2n+1}^{(2k+2m-1)}(2\bar r_i-2\to 2\bar
s_j-1);w_{BAV}\right)\right)_{0\le i,j\le m-1}
\\
&=(-1)^{\sum_{i=0}^{m-1}(\bar r_i+\bar s_i)}
\Bigg(
\prod _{i=1} ^{k+m}A_i^{-n+1}V_i^{-n+1}\Bigg)
\Bigg(
\prod _{j=0} ^{m-1}A_{\bar r_j}^{-2}V_{\bar r_j}^{-1}A_{\bar s_j}
\prod _{i=\bar r_j+1} ^{\bar s_j}A_i^{-2}V_i^{-2}
\Bigg)\\
&\kern1cm
\times
\det\left(\GF\left(\Cal C_{-2n+1}^{(2k+2m-1)}(2\bar r_i-2\to 2\bar
s_j-1);w_{BAV}\right)\right)_{0\le i,j\le m-1},
\tag\EZ
\endalign
$$
where 
$$
\align
\{\bar r_0,\bar r_1,\dots,\bar r_{m-1}\}&=
\{1,2,\dots,k+m\}\setminus\{r_0,r_1,\dots,r_{k-1}\},\\
\{\bar s_0,\bar s_1,\dots,\bar s_{m-1}\}&=
\{1,2,\dots,k+m\}\setminus\{s_0,s_1,\dots,r_{k-1}\},
\endalign
$$
and where we assume that 
$\bar r_0<\bar r_1<\dots<\bar r_{m-1}$ and
$\bar s_0<\bar s_1<\dots<\bar s_{m-1}$. Once more,
the products have to be interpreted according to~{\rm(\EMa)}.
\endproclaim

\proclaim{Lemma \TTG}
Let $n,k,m$ be positive integers.
Furthermore, let $\Cal P=(P_0,P_1,\dots,P_{k-1})$
be a family of non-intersecting up-down
paths of height at most $2k+2m-1$ that do not pass below the $x$-axis, 
where $P_i$ runs from
$(0,2r_i-2)$ to $(2n-1,2s_i-1)$, $i=0,1,\dots,k-1$, 
and $\Cal A$ the 
rectangular array of integers of the form~{\rm(\DE)} that
corresponds to~$\Cal P$ under the bijection of Proposition~{\rm\TN}.
Then
$$
\Bigg(
\prod _{i=1} ^{k+m}A_i^{n}V_i^{n}\Bigg)
\Bigg(
\prod _{j=0} ^{k-1}
\prod _{i=r_j} ^{s_j}A_i^{-1}V_i^{-1}\Bigg)
w_{BAV}(\Cal P)=w_{AV}(\Cal A).
\tag\EZa$$
\endproclaim

Given an up-down path $P$ from height $2r-2$ to height $2s-1$, 
the weighted generalisation of~(\ESa) that is relevant here is
$$
w_{BAV}(P)=\Bigg(A_r
\prod _{i=r+1} ^{s}A_iV_i
\Bigg)\bar w_{BAV}(P),
\tag\EZb
$$
where, again, the product has to be understood as explained in~(\EMa).

\proclaim{Theorem \TTH}
Let $n,k,m$ be positive integers, and let $r_0<r_1<\dots<r_{k-1}$
and $s_0<s_1<\dots<s_{k-1}$ be sequences of positive integers with
$1\le r_i,s_i\le k+m$ for all~$i$.
The generating function $\sum_{\Cal A}w_{AV}(\Cal A)$, where the sum is
over all rectangular arrays of integers of the form~{\rm(\DE)}
that satisfy the constraints given in the statement of Proposition~{\rm\TN} 
is equal to
$$
\multline
(-1)^{\sum_{i=0}^{m-1}(\bar r_i+\bar s_i)}
\Bigg(
\prod _{j=0} ^{m-1}
\prod _{i=\bar r_j+1} ^{\bar s_j}A_{i-1}^{-1}V_i^{-1}
\Bigg)\\
\times
\det\left(\GF\left(\Cal C_{-2n+1}^{(2k+2m-1)}(2\bar r_i-2\to 2\bar
s_j-1);w_{BAV}\right)\right)_{0\le i,j\le m-1}.
\endmultline
\tag\EZc
$$
\endproclaim

\subhead 11. Enumeration of alternating tableaux\endsubhead
The key for proving our reciprocity laws in Sections~5--7 and~10 
has been the observation from Sections~4 and~9 that numbers
(respectively generating functions) of bounded up-down paths
``of negative length" are combinatorially modelled by numbers 
(respectively generating functions) of
bounded alternating sequences and to prove
enumeration results for certain arrays of alternating sequences,
see Theorems~\TF, \TI, \TL, \TO, \TY, \TTB, \TTE, and~\TTH. 
In the present section, we present two more general theorems, of which
the above theorems turn out to be special cases.

To begin with, we now officially define {\it alternating tableaux}.
Given $m$-tuples $\boldsymbol\lambda=(\la_1,\la_2,\dots,\la_{m})$
and $\boldsymbol\mu=(\mu_1,\mu_2,\dots,\mu_{m})$, where all elements
of $\boldsymbol\mu$ are even and $\mu_i\le \la_i$ for all~$i$, 
we call an array of integers of the form
$$\matrix 
&&&a_{1,\mu_1+1}&\innerhdotsfor3\after\quad &a_{1,\la_1}\\
&&a_{2,\mu_2+1}\quad \dots&a_{2,\mu_1+1}&\innerhdotsfor2\after\quad 
&a_{2,\la_2}\\
&\iddots&&\vdots&&\iddots\\
a_{m,\mu_m+1}&\innerhdotsfor3\after\quad &a_{m,\la_m}
\endmatrix
\tag\GAa$$
an {\it alternating tableau of shape
$\boldsymbol\lambda/\boldsymbol\mu$}
if the entries along each row are alternating in the sense that
$$
a_{i,2j-1}\le a_{i,2j}\quad \text{and}\quad a_{i,2j}\ge a_{i,2j+1},
\tag\GA
$$
for all $i$ and $j$, and if
$$
a_{i+1,2j}<a_{i,2j+1}>a_{i+1,2j+2}
\tag\GB
$$
for all $i$ and $j$, whenever the respective entries are defined.
Examples of such arrays can be found in (\AEb), (\CEa), and (\DBa).
As earlier, we define the {\it weight\/} $w_{AV}$ of such an array
$\Cal A$ by $w_{AV}(\Cal A)=
\prod _{i,j} ^{}A_{a_{i,2j}}V_{a_{i,2j-1}}$.
Then we have the following generating function results.

\proclaim{Theorem \TTI}
Let $k$ be a positive integer or $k=\infty$, and
let $m$ be a positive integer. Furthermore, let
$\boldsymbol\lambda=(\la_1,\la_2,\dots,\la_{m})$
and $\boldsymbol\mu=(\mu_1,\mu_2,\dots,\mu_{m})$ be
sequences of integers with $\mu_i\le\la_i$ for all~$i$,
where $\boldsymbol\mu$ is decreasing and all $\mu_i$'s are even. 

If $\boldsymbol\lambda$ is increasing and all $\la_i$'s are odd,
then the generating function $\sum_{\Cal A}w_{AV}(\Cal A)$, where the sum is
over all alternating tableaux of shape
$\boldsymbol\lambda/\boldsymbol\mu$ whose entries are at least~$1$
and at most~$k$, is equal to
$$
\det\left(\GF\left(\Cal A_{\la_j-\mu_i}^{(k)};w_{AV}\right)
\right)_{1\le i,j\le m}.
\tag\GC
$$
If $\boldsymbol\lambda$ is decreasing and all $\la_i$'s are even,
then the generating function $\sum_{\Cal A}w_{AV}(\Cal A)$, where the sum is
over all alternating tableaux of shape
$\boldsymbol\lambda/\boldsymbol\mu$ whose entries are between~$1$
and~$k$, is also given by~{\rm(\GC)}.
\endproclaim

\remark{Remarks}(1)
If $k=\infty$, then the generating functions 
$\sum_{\Cal A}w_{AV}(\Cal A)$ and
$\GF\big(\Cal A_{\la_j-\mu_i}^{(\infty)};\mathbreak w_{AV}\big)$ 
are not polynomials
anymore, but instead formal power series in the $A_i$'s and $V_i$'s. 
If $A_i=V_i$ for all~$i$, then it is known that 
$\GF\big(\Cal A_{s}^{(\infty)};w_{AV}\big)$ 
is a quasi-symmetric function (cf\. \cite{\StanBI, paragraph
containing Eq.~(7.92)}, where the underlying poset~$P$ is the
zigzag poset $Z_s$ from \cite{\StanAP, Ex.~3.66}). If we leave
the variables~$A_i$ and~$V_i$ unrelated, then the functions
$\GF\big(\Cal A_{s}^{(\infty)};w_{AV}\big)$ may stand at the beginning
of a theory of quasi-symmetric functions in two sets of variables
that would have to be developed.

\medskip
(2)
The special case of this theorem 
where $\la_i=2n+2m+2i-5$ and $\mu_i=2m-2i$, $i=1,2,\dots,m$,
and $k$ is replaced by $k+m$ is at the heart of the proofs of
Theorems~\TF\ and~\TY, while 
the special 
where $\la_i=2n+2m-2i$ and $\mu_i=2m-2i$, $i=1,2,\dots,m$,
and $k$ is replaced by $k+m$ is at the heart of the proofs of
Theorems~\TI\ and~\TTB.
\endremark

\demo{Sketch of proof of Theorem \TTI}
The theorem can be proved in the same way as Theorem~\TF.
Here, in the first case the alternating tableaux under consideration
are in bijection with families
$(P_1,P_2,\dots,P_{m})$ of non-intersecting lattice paths in the
directed graph $\Cal G_{k}$, where $P_i$ runs from $(\mu_i,k)$ to
$(\la_i,k)$, $i=1,2,\dots,m$. In the second case, 
the alternating tableaux under consideration
are in bijection with families
$(P_1,P_2,\dots,P_{m})$ of non-intersecting lattice paths in the
directed graph $\Cal G_{k}$, where $P_i$ runs from $(\mu_i,k)$ to
$(\la_i,0)$, $i=1,2,\dots,m$.
We leave it to the reader to fill in the
details.\quad \quad \qed
\enddemo

The finer ``flagged" version of Theorem~\TTI\ is the following.

\proclaim{Theorem \TTJ}
Let $k$ be a positive integer or $k=\infty$, 
let $m$ be a positive integer, and let
$\boldsymbol\lambda=(\la_1,\la_2,\dots,\la_{m})$
and $\boldsymbol\mu=(\mu_1,\mu_2,\dots,\mu_{m})$ be
sequences of integers with $\mu_i\le\la_i$ for all~$i$,
where $\boldsymbol\mu$ is non-increasing and all $\mu_i$'s are even. 
Furthermore, let $r_1>r_2>\dots>r_{m}$
and $s_1>s_2>\dots>s_{k}$ be sequences of positive integers with
$1\le r_i,s_i\le k$ for all~$i$. 

If $\boldsymbol\lambda$ is non-decreasing and all $\la_i$'s are odd,
then the generating function $\sum_{\Cal A}w_{AV}(\Cal A)$, where the sum is
over all alternating tableaux of shape
$\boldsymbol\lambda/\boldsymbol\mu$ whose entries are between~$1$
and~$k$, and in which the first entry in row~$i$ is $r_i$ and the
last entry is $s_i$, is equal to
$$
\det\left(\GF\left(\Cal A_{\la_j-\mu_i}^{(k)}(r_i\to s_j);w_{AV}\right)
\right)_{1\le i,j\le m}.
\tag\GD
$$
If $\boldsymbol\lambda$ is non-increasing and all $\la_i$'s are even,
then the generating function $\sum_{\Cal A}w_{AV}(\Cal A)$, where the sum is
over all alternating tableaux of shape
$\boldsymbol\lambda/\boldsymbol\mu$ whose entries are between~$1$
and~$k$, and in which the first entry in row~$i$ is $r_i$ and the
last entry is $s_i$, is also given by~{\rm(\GD)}.
\endproclaim

\remark{Remark}
The special case of this theorem 
where $\la_i=2n(+1)$, $\mu_i=0$, 
$r_i=\bar r_{m-i}$, and
$s_i=\bar s_{m-i}$,
$i=1,2,\dots,m$,
and $k$ is replaced by $k+m$ is at the heart of the proofs of
Theorems~\TL, \TO, \TTE, and~\TTH.
\endremark

\demo{Sketch of proof of Theorem \TTJ}
The theorem can be proved in the same way as Theorem~\TL.
Here, the alternating tableaux under consideration
are in bijection with families
$(P_1,P_2,\dots,P_{m})$ of non-intersecting lattice paths in the
directed graph $\Cal G_{k}$, where $P_i$ runs from $(\mu_i,r_i)$ to
$(\la_i,s_i)$, $i=1,2,\dots,m$. 
We leave it to the reader to fill in the
details.\quad \quad \qed
\enddemo

We close this section by pointing out that our alternating tableaux
are plane partitions in disguise. A {\it plane partition} is an
array of integers of the form~(\GAa) such that entries along rows
and columns are non-increasing (cf\. e.g\. \cite{\KratCM} or 
\cite{\StanBI, Sec.~7.20}). 

Instead of a formal description of the correspondence between
alternating tableaux and plane partitions in full generality,
we content ourselves with illustrating the correspondence
by considering the array in~(\AEb). 
If one subtracts $i$ from row~$i$ (counted from  bottom),
then one obtains
$$
\matrix 
 & & & &4&4&3&4&4\\
 & &4&4&3&3&3&3&3&4&3\\
3&4&1&2&0&3&3&3&2&3&2&3&0
\endmatrix
$$
We now rearrange this array, so that each row becomes a
zigzag ``strip", see the array on the left of Figure~\FI.

\midinsert
$$
\smatrix \format\sa\c\s\c\s\c\s\c\s\c\s\c\s\c\se\\
\omit&\omit&\hlinefor{13}\\
\omit& &&\hbox to10pt{\hss 4\hss}&\omit&\hbox to10pt{\hss 4\hss}&&\hbox to10pt{\hss 4\hss}&\omit&\hbox to10pt{\hss 3\hss}&&\hbox to10pt{\hss 3\hss}&\omit&\hbox to10pt{\hss 0\hss}&\\
\hlinefor{3}&\omit&\hlinefor{3}&\omit&\hlinefor{3}&\omit&\hlinefor{3}\\
&4&\omit&3&&3&\omit&3&&3&\omit&2&\\
&\omit&\hlinefor{3}&\omit&\hlinefor{3}&\omit&\hlinefor{3}\\
&4&&3&\omit&3&&3&\omit&2&\\
\hlinefor{3}&\omit&\hlinefor{3}&\omit&\hlinefor{3}\\
&4&\omit&3&&3&\omit&3&\\
&\omit&\hlinefor{3}&\omit&\hlinefor{3}\\
&4&&2&\omit&0&\\
\hlinefor{3}&\omit&\hlinefor{3}\\
&4&\omit&1&\\
\omit&\omit&\hlinefor{3}\\
&3&\\
\hlinefor{3}
\endsmatrix
\kern2cm
\matrix 
 &4&4&4&3&3&0\\
4&3&3&3&3&2\\
4&3&3&3&2\\
4&3&3&3\\
4&2&0\\
4&1\\
3
\endmatrix$$
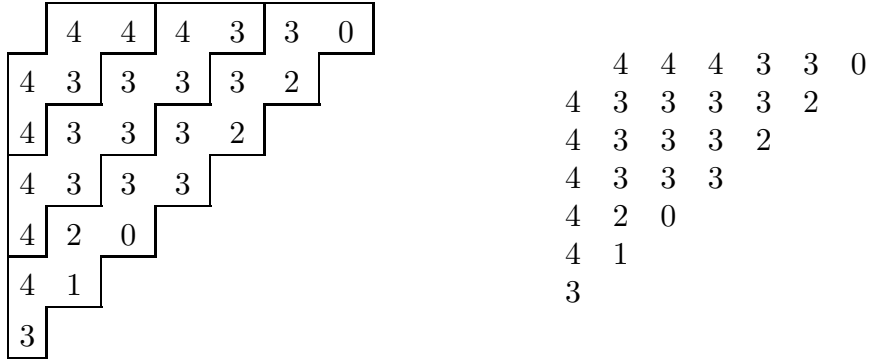
\centerline{\eightpoint From an alternating tableau of trapezoidal
 shape to a plane partition}
\vskip7pt
\centerline{\eightpoint Figure \FI}
\endinsert

On the right, we have erased the line segments separating the
zigzag strips.
This new array has the property that entries along rows and columns
are non-increasing. In other words, we have obtained a plane partition.

If one applies the same transformation to the rhomboidal array
in~(\CEa), then one obtains the plane partition shown in Figure~\FJ.

\midinsert
$$
\smatrix \format\sa\c\s\c\s\c\s\c\s\c\s\c\s\c\se\\
\omit&\omit&\omit&\omit&\omit&\omit&\omit&\omit&\hlinefor{3}&\omit&\hlinefor{3}\\
\omit& &\omit&\hbox to10pt{\hss \hss}&\omit&\hbox to10pt{\hss \hss}&\omit&\hbox to10pt{\hss \hss}&&\hbox to10pt{\hss 3\hss}&\\
\omit&\omit&\omit&\omit&\omit&\omit&\hlinefor{3}&\omit&\hlinefor{3}&\omit&\\
\omit&\omit&\omit& &\omit& &&4&\omit&2&\\
\omit&\omit&\omit&\omit&\hlinefor{3}&\omit&\hlinefor{3}&\omit&\hlinefor{3}\\
\omit&\omit&\omit& &&4&\omit&4&&1&\\
\omit&\omit&\hlinefor{3}&\omit&\hlinefor{3}&\omit&\hlinefor{3}&\omit&\\
\omit&\omit&&4&\omit&4&&3&\omit&1&\\
\hlinefor{3}&\omit&\hlinefor{3}&\omit&\hlinefor{3}&\omit&\hlinefor{3}\\
&4&\omit&3&&3&\omit&3&&0&\\
&\omit&\hlinefor{3}&\omit&\hlinefor{3}&\omit&\hlinefor{3}&\omit&\\
&4&&3&\omit&3&&3&\omit&0&\\
\hlinefor{3}&\omit&\hlinefor{3}&\omit&\hlinefor{3}\\
&4&\omit&3&&3&\omit&3&\\
&\omit&\hlinefor{3}&\omit&\hlinefor{3}\\
&4&&2&\omit&0&\\
\hlinefor{3}&\omit&\hlinefor{3}\\
&4&\omit&1&\\
\omit&\omit&\hlinefor{3}\\
&3&\\
\hlinefor{3}
\endsmatrix
\kern2cm
\matrix 
 & & & &3\\
 & & &4&2\\
 & &4&4&1\\
 &4&4&3&1\\
4&3&3&3&0\\
4&3&3&3&0\\
4&3&3&3\\
4&2&0\\
4&1\\
3
\endmatrix$$
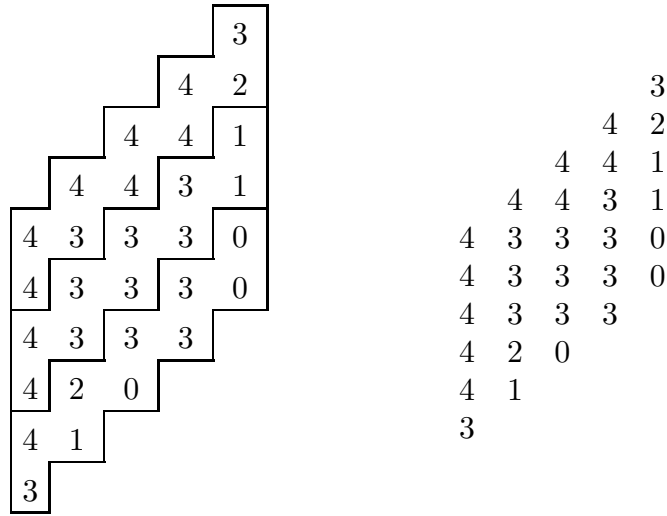
\centerline{\eightpoint From an alternating tableau of rhomboidal
 shape to a plane partition}
\vskip7pt
\centerline{\eightpoint Figure \FJ}
\endinsert

By examining in more depth the mapping from alternating tableaux to plane
partitions that we indicated above, we see that
it sets up a bijection between alternating tableaux of shape
$(2n+2m-3,2n+2m-1,\dots,2n+4m-5)/(2m-2,2m-4,\dots,0)$ 
with entries at least~$1$ and at
most~$k+m$ (see Proposition~\TE) and plane partitions of shape
$(n+2m-2,n+2m-3,\dots,1)/(n-2,n-3,\dots,1,0,\dots,0)$ 
with entries at least~$0$ and at most~$k$. (Here, the inner shape has 
$2m$~occurrences of~$0$.)
Likewise --- referring to the alternating tableaux in 
Proposition~\TH\ ---, 
this mapping sets up a bijection between alternating tableaux of shape
$(2n+2m-2,2n+2m-4,\dots,2n)/(2m-2,2m-4,\dots,0)$ 
with entries at least~$1$ and at
most~$k+m$ and plane partitions of shape
$(n,\dots,n,n-1,\dots,1)/(n-1,n-2,\dots,1,0,\dots,0)$ 
with entries at least~$0$ and at most~$k$.
Here, the outer shape has $2m$~occurrences of $n$, and the
inner shape has $2m$~occurrences of~$0$.

Using the known determinant formula \cite{\KratAH, Theorem~6.1}
for the number of plane partitions of a given shape and subject to
separate bounds on the entries in each row with $x=1$, 
$\al=\be=0$, $r=n+2m-2$, $\boldsymbol\lambda=(n+2m-2,n+2m-3,\dots,1)$, 
$\boldsymbol\mu=(n-2,n-3,\dots,1,0,\dots,0)$, $a_i=k$ and $b_i=0$
for all~$i$, we get the following alternative determinantal formula
for the determinants in Theorems~\TD\ and~\TF.

\proclaim{Proposition \TTK}
The determinants in {\rm(\AA)} and\/ {\rm(\AF)} are equal to
$$
\det\left(
\binom {n+2m-1-i-\max\{0,n-1-j\}+k} 
{k+i-j}
\right)_{1\le i,j\le n+2m-2}.
$$
Here, the binomial coefficient has to be interpreted as $0$
if its upper or lower parameter is negative.
\endproclaim

Likewise, using \cite{\KratAH, Theorem~6.1} with $x=1$,
$\al=\be=0$, $r=n+2m-1$, $\boldsymbol\lambda=(n,\dots,n,n-1,\dots,1)$
(with $2m$ occurrences of~$n$), 
$\boldsymbol\mu=(n-1,n-2,\dots,1,0,\dots,0)$ 
(with $2m$ occurrences of~$0$), $a_i=k$ and $b_i=0$
for all~$i$, we get the following alternative determinantal formula
for the determinants in Theorems~\TG\ and~\TI.

\proclaim{Proposition \TTL}
The determinants in {\rm(\CA)} and\/ {\rm(\CF)} are equal to
$$
\det\left(
\binom {\min\{n,n+2m-i\}-\max\{0,n-j\}+k} 
{k+i-j}
\right)_{1\le i,j\le n+2m-1}.
$$
Again, the binomial coefficient has to be interpreted as $0$
if its upper or lower parameter is negative.
\endproclaim

The ``strip shapes'' of the plane partitions that we obtain here (see
Figures~\FI\ and~\FJ) have been considered in \cite{\BaRoAA} and
in~\cite{\JinYAA}, however not for plane partitions, but rather for
standard Young tableaux, respectively for semistandard tableaux.
More precisely, in the terminology of \cite{\BaRoAA} our shapes are
$2m$-strip shapes, the {\it width}~$2m$ referring to the maximal length
of the columns in the ``body" of the shape. In this sense,
Figures~\FI\ and~\FJ\ contain $6$-strip shapes.

A comparison with the considerations in \cite{\LaPrAC, \GoHaAA} reveals that
our formulae in Theorems~\TTI\ and~\TTJ\ have the flavour of
the ribbon determinant formulae of Lascoux and Pragacz \cite{\LaPrAC} and their
generalisation by Hamel and  Goulden \cite{\GoHaAA}. However, our
weight $w_{AV}$ is incompatible with the monomial weight assigned to 
semistandard tableaux in order to obtain Schur functions. Hence, 
while our results overlap with those of \cite{\LaPrAC, \GoHaAA} 
in the case of plain counting (that is, if all weights are equal to~1), 
they do not follow from those of \cite{\LaPrAC, \GoHaAA}, nor do our
results imply those of \cite{\LaPrAC, \GoHaAA}.

\subhead 12. Some comments, and questions\endsubhead
In this concluding section, we discuss some points of interest and
questions that are posed by the results (and non-results) in
the present article.

\medskip
(1) It will not have escaped the attention of the reader that,
although Theorems~\VA\ and~\TP\ hold for odd {\it and\/} even upper
bounds~$k$ for the up-down paths under consideration, all our
reciprocity results require this upper bound to be {\it odd}.
The obvious question therefore is: 
$$\text{\it What happens for even upper bounds?}$$

The short answer is: there is no reciprocity law in that case.

A more elaborate answer is three-fold: first, it is a simple matter of
fact that, in a straightforward manner, the sequence 
$\big(C_{n}^{(2k)}\big)_{n\ge0}$ cannot be extended to negative
indices~$n$ since it already does not satisfy its linear recurrence
with constant coefficients {\it in the beginning}. What we mean by
this is best illustrated with the special case $k=1$: here we have
$$
\big(C_{n}^{(2)}\big)_{n\ge0}=(1,0,1,0,2,0,4,0,8,\dots).
$$
Clearly, it satisfies the linear recurrence
$C_{n+2}^{(2)}=2C_n^{(2)}$ for $n\ge1$, but not for $n=0$.
The technical explanation for this phenomenon is hidden in the generating
function for these numbers, given in (\BGa). Namely, recalling that
$C_n^{(2k)}=C_n^{(2k)}(0\to0)$, we have
$$
\sum_{n\ge0}C_n^{(2k)}x^{n}=
\frac {U_{2k}(1/2x)} {x\,U_{2k+1}(1/2x)}.
$$
The degree of $U_m(x)$ is equal to~$m$, but at the same time it is
an even polynomial for even~$m$ and an odd polynomial for odd~$m$.
Thus, if we rewrite the above right-hand side as a reduced rational
fraction, then we see that numerator degree and denominator degree are
the same. The implication is that the constant term in the series
expansion is ``special" in that it does not fit into the linear
recurrence the coefficients $C_n^{(2k)}$ satisfy otherwise.

We could remedy this by ``correcting" $C_0^{(2k)}$ to
$\frac {k} {k+1}$. Then the linear recurrence is satisfied from the
very beginning, and thus the sequence can be extended to negative
indices. The drawback here is that this extension to negative indices
does not produce integers. For example, for $k=1$ we would obtain
$C_{-2n}^{(2)}=2^{-n-1}$,
or for $k=2$ we would obtain $C_{-2n}^{(4)}=\frac {1} {2}(1+3^{-n-1})$.
In particular, we cannot expect {\it combinatorial\/} reciprocity laws in this
setting.

The second part of the more elaborate answer to the above question
mentions that, at least, there is an even analogue of the
identity~(\BAc) (the $m=1$ case of Theorem~\TD), even if it is not
a reciprocity law. Namely, for $n\ge1$ we have
$$
\det\left(C_{2n+2i+2j}^{(2k)}\right)_{0\le i,j\le k-1}
=
(k+1)^{n-1}.
\tag\HA
$$
Indeed, this cannot be interpreted as a reciprocity law (in whatever
sense) since the sequence $(C_n^{(2k)})_{n\ge0}$ satisfies a linear
recurrence that is fundamentally different from the (obvious) one
that the right-hand side of~(\HA) satisfies.

Why does (\HA) hold?
This is (again) best explained by a picture.
As before, the determinant on the left-hand side of~(\HA) equals the
number of families of non-intersecting up-down paths in a strip,
here families $(P_0,P_1,\dots,P_{k-1})$ of non-intersecting Dyck
paths of height at most $2k$, where $P_i$ runs from
$(-2i,0)$ to $(2n+2i,0)$, $i=0,1,\dots,k-1$. Figure~\FGa\ provides an
example for $k=4$ and $n=8$. It illustrates a bijection between these
families of non-intersecting Dyck paths and sequences
$(a_1,a_2,\dots,a_{n-1})$ with $1\le a_i\le k+1$ for all~$i$,
by reading the heights of the ``even" lattice points that are not occupied by
any of the paths, and finally dividing these heights by~2 and adding~1
to the result. In the
figure, these non-occupied points are indicated by circles, and the
value that is assigned to them according to the above recipe
is put as a label. Thus, for our example in Figure~\FGa\
we obtain the sequence $(4,2,1,5,2,3,1)$.
It should be emphasised that here, because of the ``tighter"
upper bound on the paths, there are non-occupied points 
only above even abscissa, and hence 
the elements~$a_i$ of these sequences are independent of each
other. In particular, there is no condition of being alternating in
this context.

\midinsert
$$
\Einheit.4cm
\Gitter(24,10)(-7,-1)
\Koordinatenachsen(24,10)(-7,-1)
\Pfad(-6,0),3333333343434433434344444444\endPfad
\Pfad(-4,0),333334334344334343444444\endPfad
\Pfad(-2,0),33343343443344334444\endPfad
\Pfad(0,0),3434334434343344\endPfad
\SPfad(1,-1),2222222222\endSPfad
\SPfad(15,-1),2222222222\endSPfad
\PfadDicke{1pt}
\Pfad(-7,8),111111111111111111111111111111\endPfad
\DickPunkt(-6,0)
\DickPunkt(-4,0)
\DickPunkt(-2,0)
\DickPunkt(0,0)
\DickPunkt(16,0)
\DickPunkt(18,0)
\DickPunkt(20,0)
\DickPunkt(22,0)
\Kreis(2,6)
\Kreis(4,2)
\Kreis(6,0)
\Kreis(8,8)
\Kreis(10,2)
\Kreis(12,4)
\Kreis(14,0)
\Label\l{4\kern-5pt}(2,6)
\Label\l{2\kern-5pt}(4,2)
\Label\l{1\kern-5pt}(6,0)
\Label\l{5\kern-5pt}(8,8)
\Label\l{2\kern-5pt}(10,2)
\Label\l{3\kern-5pt}(12,4)
\Label\l{1\kern-5pt}(14,0)
%
\hskip6.3cm
$$
\centerline{\eightpoint Another 
family of non-intersecting bounded Dyck paths}
\vskip7pt
\centerline{\eightpoint Figure \FGa}
\endinsert

This bijection
immediately proves~(\HA), and more generally
$$\multline
\det\left(\GF\left(\Cal C_{2n+2i+2j}^{(2k)};w_B\right)
\right)_{0\le i,j\le k-1}\\
=
\prod _{i=1} ^{2k}B_i^{k-\fl{i/2}}
\big(B_1B_3\cdots B_{2k-1}
+B_1B_3\cdots B_{2k-3}B_{2k}\kern3cm\\
+B_1B_3\cdots B_{2k-5}B_{2k-2}B_{2k}+\dots
+B_2B_4\cdots B_{2k}\big)^{n-1}.
\endmultline
\tag\HB
$$

\midinsert
$$
\Einheit.4cm
\Gitter(24,14)(-7,-1)
\Koordinatenachsen(24,14)(-7,-1)
\Pfad(-6,0),3333333343333344434444444444\endPfad
\Pfad(-4,0),333334333444334343444444\endPfad
\Pfad(-2,0),33343343443434334444\endPfad
\Pfad(0,0),3434334434343344\endPfad
\SPfad(1,-1),22222222222222\endSPfad
\SPfad(15,-1),22222222222222\endSPfad
\PfadDicke{1pt}
\Pfad(-7,12),111111111111111111111111111111\endPfad
\DickPunkt(-6,0)
\DickPunkt(-4,0)
\DickPunkt(-2,0)
\DickPunkt(0,0)
\DickPunkt(16,0)
\DickPunkt(18,0)
\DickPunkt(20,0)
\DickPunkt(22,0)
\Kreis(2,6)
\Kreis(3,9)
\Kreis(4,2)
\Kreis(4,10)
\Kreis(5,5)
\Kreis(5,11)
\Kreis(6,0)
\Kreis(6,8)
\Kreis(6,12)
\Kreis(7,7)
\Kreis(7,9)
\Kreis(8,6)
\Kreis(8,8)
\Kreis(8,10)
\Kreis(9,7)
\Kreis(9,9)
\Kreis(10,4)
\Kreis(10,8)
\Kreis(10,12)
\Kreis(11,7)
\Kreis(11,11)
\Kreis(12,4)
\Kreis(12,8)
\Kreis(13,7)
\Kreis(14,0)
\Label\l{4\kern-5pt}(2,6)
\Label\l{5\kern-5pt}(3,9)
\Label\l{2\kern-5pt}(4,2)
\Label\l{6\kern-5pt}(4,10)
\Label\l{3\kern-5pt}(5,5)
\Label\l{6\kern-5pt}(5,11)
\Label\l{1\kern-5pt}(6,0)
\Label\l{5\kern-5pt}(6,8)
\Label\l{7\kern-5pt}(6,12)
\Label\l{4\kern-5pt}(7,7)
\Label\l{5\kern-5pt}(7,9)
\Label\l{4\kern-5pt}(8,6)
\Label\l{5\kern-5pt}(8,8)
\Label\l{6\kern-5pt}(8,10)
\Label\l{4\kern-5pt}(9,7)
\Label\l{5\kern-5pt}(9,9)
\Label\l{3\kern-5pt}(10,4)
\Label\l{5\kern-5pt}(10,8)
\Label\l{7\kern-5pt}(10,12)
\Label\l{4\kern-5pt}(11,7)
\Label\l{6\kern-5pt}(11,11)
\Label\l{3\kern-5pt}(12,4)
\Label\l{5\kern-5pt}(12,8)
\Label\l{4\kern-5pt}(13,7)
\Label\l{1\kern-5pt}(14,0)
%
\hskip6.3cm
$$
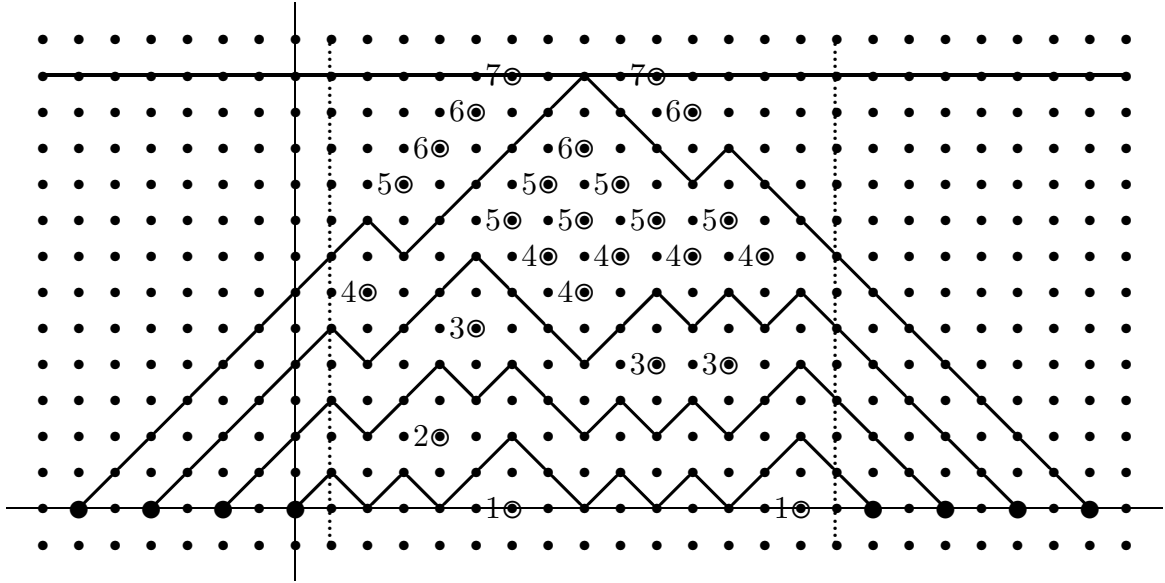
\centerline{\eightpoint Yet another 
family of non-intersecting bounded Dyck paths}
\vskip7pt
\centerline{\eightpoint Figure \FH}
\endinsert

Third, we should address the question of what happens if the bound on the
paths is lifted. Phrased differently, if we consider the determinants
$$
\det\left(C_{2n+2i+2j+4m}^{(2k+2m)}\right)_{0\le i,j\le k-1},
$$
what replaces the right-hand side of (\AA) in this ``even case"?
We may follow the line of argument in Section~5 
of interpreting the determinant in
terms of families of non-intersecting bounded Dyck paths and then
mapping them to alternating tableaux (see Proposition~\TE).
Figure~\FH\ shows an example of a family of Dyck paths for $k=4$,
$m=2$, and $n=4$. The array of integers that one obtains using the
idea of the proof of Proposition~\TE\ is
$$
\matrix 
 & & & &7& &6& &7\\
 & &6&6&5&5&5&5&5&6&5\\
4&5&2&3&1&4&4&4&3&4&3&4&1
\endmatrix
$$
One observes that, here, the top-most row has a gap at every second
position. The corresponding plane partition according to the
translation described at the end of Section~11 is shown in Figure~\FK.

\midinsert
$$
\smatrix \format\sa\c\s\c\s\c\s\c\s\c\s\c\s\c\se\\
\omit&\omit&\omit&\omit&\hlinefor{11}\\
\omit& &\omit&\hbox to10pt{\hss \hss}&&\hbox to10pt{\hss 4\hss}&&\hbox to10pt{\hss 4\hss}&\omit&\hbox to10pt{\hss 3\hss}&&\hbox to10pt{\hss 3\hss}&\omit&\hbox to10pt{\hss 0\hss}&\\
\omit&\omit&\hlinefor{5}&\omit&\hlinefor{3}&\omit&\hlinefor{3}\\
\omit& &&3&&3&\omit&3&&3&\omit&2&\\
\hlinefor{5}&\omit&\hlinefor{3}&\omit&\hlinefor{3}\\
&4&&3&\omit&3&&3&\omit&2&\\
\hlinefor{3}&\omit&\hlinefor{3}&\omit&\hlinefor{3}\\
&4&\omit&3&&3&\omit&3&\\
&\omit&\hlinefor{3}&\omit&\hlinefor{3}\\
&4&&2&\omit&0&\\
\hlinefor{3}&\omit&\hlinefor{3}\\
&4&\omit&1&\\
\omit&\omit&\hlinefor{3}\\
&3&\\
\hlinefor{3}
\endsmatrix
\kern2cm
\matrix 
 & &4&4&3&3&0\\
 &3&3&3&3&2\\
4&3&3&3&2\\
4&3&3&3\\
4&2&0\\
4&1\\
3
\endmatrix$$
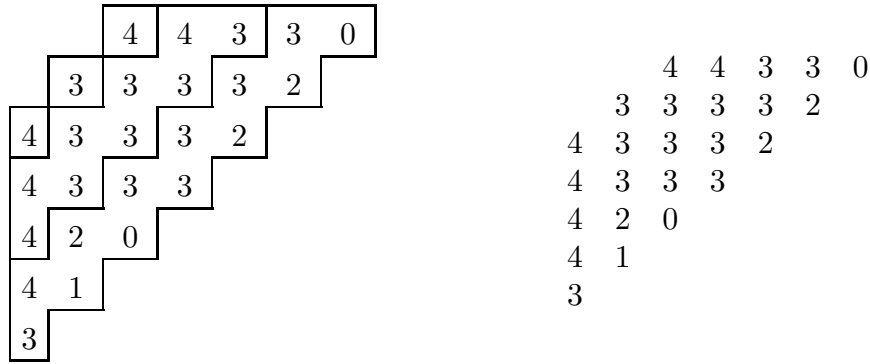
\centerline{\eightpoint From an alternating tableau of a
 shape with gaps to a plane partition}
\vskip7pt
\centerline{\eightpoint Figure \FK}
\endinsert

So, here we obtain a $5$-strip shape, that is, a strip shape 
of {\it odd\/} width~5.
Obviously, the determinant formula \cite{\KratAH, Theorem~6.1} also
applies here. However, if one wishes to find a determinant formula
that is analogous to the one in~(\AFb), then one would have to adapt
Jin's \cite{\JinYAA} extension of the ribbon determinant formulae
of Lascoux and Pragacz \cite{\LaPrAC} and of Hamel and
Goulden \cite{\GoHaAA} to thickened ribbons to the present situation.
Although this is certainly feasible, 
it would again not produce reciprocity laws.

\medskip
(2) Another question that may have occurred to the reader is whether
it is possible to ``mix alternating sequences of even and odd length"{}?
The background is that all our results (with the exception of 
Corollary~\TA) involve either alternating sequences of even length
(such as, for example, the second statement in Theorem~\UA\ or the results
in Section~6) or of odd length 
(such as, for example, the first statement in Theorem~\UA\ or the results
in Section~5), but not both together.
Here we report a conjectural identity where the two are ``mixed".
As it turns out, on the side of the up-down paths one must sum over all
possible ending heights.\footnote{See the Final Note at the end of the article.}

\proclaim{Conjecture \TTLa}
For all non-negative integers $n,k,m$, we have
$$\multline
\det\left(\sum_{s=0}^{2k+2m-1}
C_{n+i+j+2m-1}^{(2k+2m-1)}(0\to s)\right)_{0\le i,j\le k-1}\\
=(-1)^{\left(\binom k2+\binom m2\right)(n+1)}
\det\left(\big\vert\Cal A_{n+i+j}^{(k+m)}\big\vert\right)_{0\le i,j\le m-1}.
\endmultline
\tag\HBa
$$
\endproclaim

There are closed-form evaluations for the ``even cases" of the
determinants on the left-hand side of~(\HBa) if the up-down paths
are restricted to ``very tight" strips.

\proclaim{Theorem \TTLb}
Let $n$ be a non-negative integer and $k$ be a positive integer.
Then we have
$$
\det\left(\sum_{s=0}^{4k}
C_{n+i+j+1}^{(4k)}(0\to s)\right)_{0\le i,j\le 2k-1}
=(-1)^{k(n-1)}(2k+1)^n
\tag\HBb
$$
and
$$
\det\left(\sum_{s=0}^{4k-2}
C_{n+i+j}^{(4k-2)}(0\to s)\right)_{0\le i,j\le 2k-1}
=(-1)^{kn}2^n.
\tag\HBd
$$
\endproclaim

\demo{Proof}
We start with (\HBb).
Write $a(n,k)$ for $\sum_{s=0}^{k}
C_{n}^{(k)}(0\to s)$. We claim that we have
$$
\sum_{n\ge0}a(n,k)x^n=\frac {U_{\fl{k/2}}(1/2x)
\big(U_{\fl{(k+1)/2}}(1/2x)+U_{\fl{(k-1)/2}}(1/2x)\big)} 
{xU_{k+1}(1/2x)}.
\tag\HBf
$$
The claim can be easily verified by using (\BGa) and the substitution
trick explained in the proof of Corollary~\TA. Alternatively, one
may consult \cite{\CiglBD, Eq.~(1.27)}.

Now let $A_k(n)$ be the matrix
$\big(a(n+i+j,4k)\big)_{0\le i,j\le 2k-1}$. Then we claim that
$$
A_k(n+1)=A_k(n)\pmatrix 
0&0&\innerhdotsfor3\after\quad &-\al_{2k}\\
1&0&\innerhdotsfor3\after\quad &-\al_{2k-1}\\
0&1&\innerhdotsfor3\after\quad &-\al_{2k-2}\\
\vdots&\vdots&\ddots&&\\
0&0&\innerhdotsfor2\after\quad &0&-\al_2\\
0&0&\innerhdotsfor2\after\quad &1&-\al_1
\endpmatrix,\quad \text{for }n\ge1,
\tag\HBg
$$ 
where 
$$\sum_{i=0}^{2k}\al_ix^i=\frac {x^{2k+1}U_{4k+1}(1/2x)} {U_{2k}(1/2x)}
=2x^{2k+1}T_{2k+1}(1/2x).
\tag\HBh
$$
Here, $T_n(x)$ denotes the {\it $n$-th Chebyshev polynomial of the first kind},
given by
$$
\align
T_n(\cos\th)&={\cos(n\th)} ,\\
  T_n(x)&=\sum_{j\ge0}(-1)^j\frac {n} {2(n-j)}\binom {n-j}{j}(2x)^{n-2j}.
\tag\HBi
\endalign
$$
The equality in (\HBh) can again be checked using the aforementioned
substitution trick. It should be noted that the last expression in~(\HBh) is 
visibly a polynomial of degree~$2k$ with constant term~1. 
By multiplication of the right-hand
sides of~(\HBf) with~$k$ replaced by~$4k$ and of~(\HBh), we obtain
$$
\sum_{n\ge0}x^n\sum_{i=0}^{2k}\al_ia(n-i,4k)=
2x^{2k}\big(U_{2k}(1/2x)+U_{2k-1}(1/2x)\big).
$$
The right-hand side is a polynomial of degree~$2k$. Hence,
$$
\sum_{i=0}^{2k}\al_ia(n-i,4k)=0,\quad \text{for }n>2k,
$$
or, equivalently,
$$
a(n,4k)=-\sum_{i=1}^{2k}\al_ia(n-i,4k),\quad \text{for }n>2k.
$$
This establishes the claim in (\HBg).

Taking determinants on both sides of~(\HBg), we arrive at the recurrence
$$
\det A_k(n+1)=\al_{2k}\det A_k(n),\quad \text{for }n\ge1.
$$
By the explicit form (\HBi) of the Chebyshev polynomials of the first kind
and~(\HBh), we see that $\al_{2k}=(-1)^k(2k+1)$.

We have everything built up for an inductive proof of~(\HBb), except for
the start of the induction, the evaluation of $\det A_k(1)$. 
The $(i,j)$-entry of $A_k(1)$ is $a(1+i+j,4k)$, $0\le i,j\le 2k-1$, 
the number of up-down paths
starting at $(0,0)$ not passing below the $x$-axis and not passing
above height~$4k$. However, the paths take at most $1+(2k-1)+(2k-1)=4k-1$
steps so that the upper restriction is without relevance. 
The number of up-down paths taking~$n$ steps starting at $(0,0)$ and never 
passing below the $x$-axis is well known to be $\binom n{\fl{n/2}}$.
Therefore we have
$$
\det A_k(1)=\det\left(\binom {i+j+1}{\fl{(i+j+1)/2}}\right)=(-1)^k,
$$
where the last equality follows for instance from \cite{\KrYaAA,
Eq.~(6.3)}. This completes the proof of~(\HBb).

\medskip
Equation (\HBd) can be proved in a similar manner.
Let $B_k(n)$ be the matrix\linebreak
$\big(a(n+i+j,4k-2)\big)_{0\le i,j\le 2k-1}$. Then we claim that
$$
B_k(n+1)=B_k(n)\pmatrix 
0&0&\innerhdotsfor3\after\quad &-\be_{2k}\\
1&0&\innerhdotsfor3\after\quad &-\be_{2k-1}\\
0&1&\innerhdotsfor3\after\quad &-\be_{2k-2}\\
\vdots&\vdots&\ddots&&\\
0&0&\innerhdotsfor2\after\quad &0&-\be_2\\
0&0&\innerhdotsfor2\after\quad &1&-\be_1
\endpmatrix,\quad \text{for }n\ge1,
\tag\HBj
$$ 
where 
$$\sum_{i=0}^{2k}\be_ix^i=\frac {x^{2k}U_{4k-1}(1/2x)} {U_{2k-1}(1/2x)}
=2x^{2k}T_{2k}(1/2x).
\tag\HBk
$$
Here we obtain
$$
\sum_{n\ge0}x^n\sum_{i=0}^{2k}\be_ia(n-i,4k-2)=
2x^{2k-1}\big(U_{2k-1}(1/2x)+U_{2k-2}(1/2x)\big),
$$
and consequently
$$
a(n,4k-2)=-\sum_{i=1}^{2k}\be_ia(n-i,4k-2),\quad \text{for }n\ge 2k.
$$
This establishes the claim in (\HBj).

Taking determinants on both sides of~(\HBj), we arrive at the recurrence
$$
\det B_k(n+1)=\be_{2k}\det B_k(n),\quad \text{for }n\ge0.
$$
By the explicit form (\HBi) of the Chebyshev polynomials of the first kind
and~(\HBk), we see that $\be_{2k}=2(-1)^k$.

Finally, we have to evaluate the determinant of
$$B_k(0)=\big(a(i+j,4k-2)\big)_{0\le i,j\le 2k-1}
=\left(\binom {i+j}{\fl{(i+j)/2}}\right)_{0\le i,j\le 2k-1}.$$
This determinant turns out to equal~1 (cf\. e.g\. \cite{\KrYaAA, Eq.~(6.2)}).
This completes the proof of~(\HBb) and, thus, of the theorem.\quad \quad \qed
\enddemo

We refer the reader to
Corollary~\TTU\ for the extension of Theorem~\TTLb\ to
negative integers.

\medskip
The next two theorems present four further Hankel determinant evaluations of
matrices where the entries are ``cumulative" generating functions, respectively
numbers, of bounded up-down paths.

\proclaim{Theorem \TTTT}
For non-negative integers $n$ and $k$, let
$$
s(n,k,x)=\sum_{\ell=0}^{\fl{n/2}} x^{n-2\ell}C_{2\ell}^{(k)}.
$$
Then we have
$$
\det\big(s(2n+2i+2j,2k,x)\big)_{0\le i,j\le k}=(-1)^kx^{2n}(k+1)^nU_{2k}(x/2)
\tag\HJ
$$
and
$$
\det\big(s(n+i+j,2k,x)\big)_{0\le i,j\le 2k}
(-1)^{k(n+1)}x^{n}(k+1)^nU_{2k}(x/2).
\tag\HK
$$
\endproclaim

\demo{Sketch of proof}
In the determinant in (\HJ), we subtract $x^2$ times the $(i-1)$-st row from the
$i$-th row, $i=k,k-1,\dots,2$, in this order. As a result, we obtain the
determinant
$$\multline
\det\pmatrix 
\sum_{\ell=0}^{2n+2j} x^{2n+2j-2\ell}C_{2\ell}^{(2k)},&i=0\\
C_{2n+2i+2j}^{(2k)}\hfill,&i\ge1
\endpmatrix_{0\le i,j\le k}
\\
=
\det\pmatrix 
\sum_{\ell=0}^{2n+2j} x^{2\ell}C_{2n+2j-2\ell}^{(2k)},&i=0\\
C_{2n+2i+2j}^{(2k)}\hfill,&i\ge1
\endpmatrix_{0\le i,j\le k}.
\endmultline
\tag\HL
$$
If we use linearity of the determinant in the first row, then we get
$$
\sum_{\ell=0}^{2n+2k}x^{2\ell}
\det\pmatrix 
C_{2n+2j-2\ell}^{(2k)},&i=0\\
C_{2n+2i+2j}^{(2k)}\hfill,&i\ge1
\endpmatrix_{0\le i,j\le k-1}.
$$

By the Lindstr\"om--Gessel--Viennot theorem \cite{\LindAA, Lemma~1},
the last determinant enumerates families
$(P_0,P_1,\dots,P_k)$ of non-intersecting up-down paths which stay between
the lines $y=0$ and $y=2k$, where $P_0$ runs from
$(2\ell,0)$ to $(2n+2\si(0),0)$, while for $i\ge1$
the path $P_i$ runs from $(-2i,0)$ to $(2n+2\si(i),0)$, for a suitable
permutation~$\si$ of $\{0,1,\dots,k-1\}$.

In order to have an example at our disposal, see Figure~\FN, which is an
example for $k=4$, $n=5$, $\ell=6$,
$\si(0)=1$, 
$\si(1)=0$, 
$\si(2)=2$, 
$\si(3)=3$, 
and $\si(4)=4$.

\midinsert
$$
\Einheit.47cm
\Gitter(20,11)(-9,-1)
\Koordinatenachsen(20,11)(-9,-1)
\Pfad(-2,0),343434343344\endPfad
\Pfad(-4,0),333434343343444344\endPfad
\Pfad(-6,0),3333334434334344344444\endPfad
\Pfad(-8,0),33333333434433434434444444\endPfad
\PfadDicke{.1pt}
\Pfad(-9,8),1111111111111111111111111111\endPfad
\Pfad(10,-1),222222222222\endPfad
\SPfad(12,0),22822228\endSPfad
\DickPunkt(-2,0)
\DickPunkt(-4,0)
\DickPunkt(-6,0)
\DickPunkt(-8,0)
\DickPunkt(10,0)
\DickPunkt(12,0)
\DickPunkt(14,0)
\DickPunkt(16,0)
\DickPunkt(18,0)
\Kreis(10,8)
\hskip5cm
$$
\vskip7pt
\centerline{\eightpoint Figure \FN}
\endinsert

The path $P_0$ starts at $(2\ell,0)$.
If we ignore $P_0$ for the moment, then we see that the other paths
(in the example:  $P_1,P_2,P_3,P_4$) are uniquely determined at the
beginning (to be precise: until the
abscissa $x=-1$) and then ``fill everything" between $x=-1$ and
$x=2n-1$ so that there remains no space for $P_0$.
Hence, we must have $\ell\ge n$ in order to have a non-empty family of
non-intersecting up-down paths.

If $\ell\ge n$, then the point $(2\ell,0)$ is not only the starting
point of $P_0$, but at the same time also the ending point of some
path. In order to have a family of non-intersecting up-down paths,
the path~$P_0$ must be a path of length~0, which starts and ends in
$(2\ell,0)$.

In order to see how many possibilities there remain for the other
paths to the right of the abscissa $x=2n-1$,
we use Viennot's idea of a ,\kern-1pt,dual`` path. 
Starting from the given family of non-intersecting paths, we
construct a path from $(2\ell,0)$ to
$(2n,2k)$ that consists entirely of vertical steps $(0,2)$ and diagonal
steps $(-1,1)$. We start the path in $(2\ell,0)$. The rule is the following:
usually we proceed by diagonal steps $(-1,1)$, except if in this way
we would hit one of the paths of the family; if the latter would be
the case then, instead, we proceed by a vertical step $(0,2)$.
In our example in Figure~\FN, the dual path is indicated as dotted path.
It is indeed not difficult to see that the correspondence between
path portions to the right of $x=2n-1$ and the just contructed dual path
is a bijection.

Now everything is in place: there are
$\binom {2\ell-2n+k-\ell+n}{2\ell-2n}
=\binom {\ell-n+k}{k-\ell+n}$ such dual paths.
If, in addition, we replace $\ell$ by $n+k-\ell$, then the
binomial coefficient becomes $\binom {2k-\ell}{\ell}$, which is indeed
the corresponding coefficient in the Chebyshev polynomial~$U_{2k}(x/2)$.
The factor $x^{2n}$ arises since we must have $\ell\ge n$, and the
sign that goes with the coefficient of the Chebyshev polynomial comes
from the sign of the permutation~$\si$.
Finally, the factor $(k+1)^n$ accounts for the number of possibilities
the paths $P_1,P_2,\dots,P_k$ have in the region between the
abscissa $x=-1$ and $x=2n-1$.

\medskip
For (\HK), we may proceed similarly.
Here we subtract $x$ (!) times the $(i-1)$-st row from the
$i$-th row, $i=2k,2k-1,\dots,2$, in this order.
However, since we have
$$
xs(2m,2k,x)=s(2m+1,2k,x),
$$
after these row operations we obtain a matrix in which zeroes appear
in a checkerboard pattern. Consequently, the determinant decomposes
into a product of two determinants of approximately half the size.
For the exact computation, one has to distinguish between $n$
being even or odd.

If $n$ is even, then the first determinant is
$$
\det\pmatrix 
\sum_{\ell=0}^{n+2j} x^{2\ell}C_{n+2j-2\ell}^{(2k)},&i=0\\
C_{n+2i+2j}^{(2k)}\hfill,&i\ge1
\endpmatrix_{0\le i,j\le k},
$$
which is the one in~(\HL), which we already evaluated.
The second determinant is 
$$
\det\big(C_{n+2i+2j+2}^{(2k)}\big)_{0\le i,j\le k-1},
$$
which can be evaluated by means of~(\HA). 

If $n$ is odd, then everything is very similar. Here, the first
determinant is again the one in~(\HL) up to a factor of~$x$, and the
second is again covered by~(\HA).\quad \quad \qed
\enddemo

\proclaim{Theorem \TTTTT}
With the notation of Theorem~{\rm\TTTT}, we have
$$
\big(s(2n+2i+2j,6k+2,1)\big)_{0\le i,j\le 3k}=(-1)^k(3k+2)^n
\tag\HM
$$
and
$$
\big(s(n+i+j,6k+2,1)\big)_{0\le i,j\le 3k}=(-1)^{kn+n+k}(3k+2)^n.
\tag\HN
$$
\endproclaim

\demo{Sketch of proof}
These two identities can again be handled by the matrix approach from
the proof of Theorem~\TTLb. We start with~(\HM). Let
$$
A_k(n)=
\big(s(2n+2i+2j,6k+2,1)\big)_{0\le i,j\le 3k}.
$$
We claim that
$$
A_k(n+1)=A_k(n)
\pmatrix 
0&0&\innerhdotsfor3\after\quad &-\al_{3k+1}\\
1&0&\innerhdotsfor3\after\quad &-\al_{3k}\\
0&1&\innerhdotsfor3\after\quad &-\al_{3k-1}\\
\vdots&\vdots&\ddots&&\\
0&0&\innerhdotsfor2\after\quad &0&-\al_2\\
0&0&\innerhdotsfor2\after\quad &1&-\al_1
\endpmatrix,\quad \text{for }n\ge0,
$$
where
$$\sum_{i=0}^{3k+1}\al_ix^{2i}
=x^{6k+3}U_{6k+3}(1/2x).
$$

To prove the claim, we observe that, by Theorem~\VA, we have
$$
\sum_{n\ge0}s(n,6k+2,1)x^n=\frac {U_{6k+2}(1/2x)}
{x(1-x)U_{6k+3}(1/2x)},
$$
and consequently
$$
\align
\sum_{n\ge0}s(2n,6k+2,1)x^{2n}
&=\frac {1} {2}\left(
\frac {U_{6k+2}(1/2x)} {x(1-x)U_{6k+3}(1/2x)}
-
\frac {U_{6k+2}(-1/2x)} {x(1+x)U_{6k+3}(-1/2x)}
\right)\\
&=\frac {1} {2}\left(
\frac {U_{6k+2}(1/2x)} {x(1-x)U_{6k+3}(1/2x)}
+
\frac {U_{6k+2}(1/2x)} {x(1+x)U_{6k+3}(1/2x)}
\right)\\
&=
\frac {U_{6k+2}(1/2x)} {x(1-x^2)U_{6k+3}(1/2x)}.
\endalign
$$
It follows that
$$
\sum_{n\ge0}x^{2n}\sum_{i=0}^{3k+1}\al_is(2n-2i,6k+2,1)=
\frac {x^{6k+2}U_{6k+2}(1/2x)} {1-x^2}.
$$
The numerator on the right-hand side is a polynomial of degree~$6k+2$
that is divisible by $1-x^2$,  as is not difficult to see. Thus, the
right-hand side is a polynomial of degree~$6k$. In particular, we have
$$
\sum_{i=0}^{3k+1}\al_is(2n-2i,6k+2,1)=0,\quad \text{for }n>3k,
$$
or, equivalently,
$$
s(2n,6k+2,1)=-\sum_{i=1}^{3k+1}\al_is(2n-2i,6k+2,1),\quad \text{for }n>3k.
$$
This proves the claim.

By induction, we get
$$
\det A_k(n)=\left((-1)^{3k+1}\al_{3k+1}\right)^n\det A_k(0)
=\left(3k+2\right)^n\det A_k(0).
$$

For the start of the induction, we may use Theorem~\TTTT:
$$
\align
\det A_k(0)&=\det\big(s(2i+2j,6k+2,1)\big)_{0\le i,j\le 3k}\\
&=\det\big(s(2i+2j,6k,1)\big)_{0\le i,j\le 3k}\\
&=(-1)^{3k}U_{6k}(1/2)\\
&=(-1)^k,
\endalign
$$
since, for the range where we find~$i$ and~$j$, it is irrelevant
whether the height restriction is $6k+2$ or~$6k$.

\medskip
For the proof of (\HN), one uses as before row operations to generate
a matrix in which the zeroes form a checkerboard pattern. Then, again,
the determinant decomposes into a product of two determinants of
approximately half the size. One factor is the determinant that we
just treated, the other is again covered by~(\HA).\quad \quad \qed
\enddemo

\medskip
(3)
It was pointed out earlier that the generating function formulae
in Theorems~\VA\ and \TP\ can actually be extended to the ``Motzkin case".
To make this precise, we define a {\it three-step path}
to be a path consisting of up-steps $(1,1)$, level steps $(1,0)$, and
down-steps $(1,-1)$. 
We denote the set of all three-step paths from $(0,r)$ to $(n,s)$ of
height at most~$k$ that never pass below the $x$-axis 
by $\Cal M_n^{(k)}(r\to s)$.
Strictly speaking, {\it Motzkin paths} are the paths in $\bigcup_{n\ge0}\Cal
M_n^{(\infty)}(0\to 0)$, that is, those three-step paths that
start and end on the $x$-axis, and never pass below the $x$-axis.

Assign weights
to the steps by defining the weight of an up-step to be~$1$, the
weight of a level step at height~$h$ to be $b_h$, the weight of a
down-step from height~$h$ to height~$h-1$ to be $\la_h$, and extend
this weight, $w$ say, to three-step paths $P$ by defining $w(P)$ to be
the product of all weights of the steps of~$P$. 

For the statement of the formula, we need to introduce the
(orthogonal) polynomials
$(p_n(x))_{n\ge0}$ satisfying the three-term recurrence
$$
xp_n(x)=p_{n+1}(x)+b_np_n(x)+\la_np_{n-1}(x),\quad \text{for }n\ge1,
$$
with initial conditions $p_0(x)=1$ and $p_1(x)=x-b_0$. Finally, we
introduce the operator~$T$, which, when applied
to a polynomial in the $b_i$'s and $\la_i$'s, 
replaces each $b_i$ by $b_{i+1}$ and each $\la_i$ by $\la_{i+1}$
for all~$i$.

With these notations,
the aforementioned formula from \cite{\VienAE, Ch.~V, Eq.~(27)}
(see also \cite{\KratCL, Theorem~10.11.1})
is the following.

\proclaim{Theorem \TTM}
Let $k,r,s$ be non-negative integers with $r,s\le k$. Then
$$
\sum _{n\ge0} ^{}\GF\left(\Cal M_n^{(k)}(r\to s);w\right)\,x^n
=\cases \dfrac {p_r(1/x)T^{s+1}p_{k-s}(1/x)} {xp_{k+1}(1/x)},
&\text{if }r\le s,\\
\la_r\cdots\la_{s+1}\dfrac {p_s(1/x)T^{r+1}p_{k-r}(1/x)}
{xp_{k+1}(1/x)},
&\text{if }r\ge s.
\endcases
$$
\endproclaim

In particular, the generating function for bounded three-step paths from
height~$r$ to height~$s$ is
again a rational function, and therefore we
may think of extending the sequence
$\big(\GF\big(\Cal M_n^{(k)}(r\to
s);w\big)\big)_{n\ge0}$ to negative~$n$. 

We report here a conjecture hinting at the possibility of
reciprocity laws in the Motzkin path context.

Let us write $M_n^{(k)}$
for $\big\vert\Cal M_n^{(k)}(0\to 0)\big\vert$.
If we choose $r=s=0$ and $b_i=\la_i=1$ for all~$i$, then we obtain
$$
\sum _{n\ge0} ^{}M_n^{(k)}\,x^n
= \dfrac {U_{k}\left(\frac {1-x} {2x}\right)} 
{xU_{k+1}\left(\frac {1-x} {2x}\right)},
\tag\HC
$$
as is not too difficult to see. It turns out that 
the sequence $\big(M_{n}^{(k)}\big)_{n\ge0}$ can be extended to
negative indices~$n$ if and only if
$k\hbox{${}\not\equiv{}$}1$~(mod~3). 
This is similar to the sequence $\big(C_{n}^{(k)}\big)_{n\ge0}$,
which could be extended to 
negative indices~$n$ if and only if
$k\hbox{${}\not\equiv{}$}0$~(mod~2). 
The reason is the same as in Item~(1): for 
$k\equiv1$~(mod~3), numerator and denominator of the rational fraction
on the right-hand side of (\HC) have the same degree.
Also here, we could remedy this by ``correcting" $M_0^{(3k+1)}$
to $(2k+1)/(2k+2)$, again at the expense of obtaining non-integer
values for negative indices~$n$.

At any rate, computer experiments suggest the following
conjecture.\footnote{See the Final Note at the end of the article.}

\proclaim{Conjecture \TTN}
For all positive integers $n,k,m$, we have
$$
\multline
\det\left(M_{n+i+j+2m-2}^{(k+m-1)}\right)_{0\le i,j\le k-1}\\
=(-1)^{n\fl{(k+m)/3}}
\det\left(M_{-n-i-j}^{(k+m-1)}\right)_{0\le i,j\le m-1},
\quad \text{for }k+m\not\equiv2~\text{\rm (mod 3)}.
\endmultline
\tag\HD
$$
\endproclaim

Equation~(\HD) must be seen as the ``Motzkin analogue'' of~(\AA).
There is also a ``Motzkin analogue" of the ``even" evaluation in~(\HA). 
Since it can be proved by following the idea of the proof
of Theorem~\TTLb, we omit its proof.

\proclaim{Theorem \TTNa}
For all positive integers $n$ and $k$, we have
$$
\det\left(M_{n+i+j}^{(k)}\right)_{0\le i,j\le k-1}
=(-1)^{n\fl{(k+1)/3}}
\left(\frac {2k+4} {3}\right)^{n-1},
\quad \text{for }k\equiv1~\text{\rm (mod 3)}.
\tag\HE
$$
\endproclaim

It may be that there are also weighted versions of~(\HD) and~(\HE),
as well as there may be ``Motzkin analogues'' of our reciprocity
laws in Sections~6, 7, and~10.

At present, we do not know how to prove Conjecture~\TTN.
It is nevertheless tempting to approach a proof of~(\HD), say,
in the same way as we proved Theorems~\TD\ and~\TW. Namely, by the 
Lindstr\"om--Gessel--Viennot theorem \cite{\LindAA, Lemma~1}, the determinant
on the left-hand side of (\HD) equals $\sum_{\Cal P}\sgn\Cal
P$, where the sum runs over all
families $(P_0,P_1,\dots,P_{k-1})$ of non-intersecting Motzkin
paths\footnote{However, it is important to
note that a crossing of an up-step $(x,y)\to(x+1,y+1)$ and a down-step
$(x,y+1)\to(x+1,y)$ is {\it not\/} a violation of the condition of
being non-intersecting; according to \cite{\LindAA, \GeViAA, \GeViAB},
``non-intersecting" means that no two paths of a family of paths 
meet {\it in a vertex of the underlying graph}. In the present
context, the vertices of the underlying graph are the lattice points,
that is, the points with integer coordinates, in the plane.}
of height at most $k+m-1$, where $P_i$ runs from
$(-i,0)$ to $(n+2m+\si(i)-2,0)$, $i=0,1,\dots,k-1$, for some
permutation~$\si\in\frak S_k$, with $\frak S_k$ denoting the set of
permutations of $\{0,1,\dots,k-1\}$; 
the sign $\sgn\Cal P$ is by definition equal to
$\sgn\si$. 
An example of such a family of
non-intersecting Motzkin paths for $k=4$, $m=3$, $n=7$, 
and the permutation~$\si$ is given by 
$\si(0)=3$, $\si(1)=2$, $\si(2)=0$, $\si(3)=1$, is 
shown in Figure~\FL. (At this point, circled points and attached
labels, as well as dotted lines should be ignored.)

\midinsert
$$
\Einheit.6cm
\Gitter(16,8)(-4,-1)
\Koordinatenachsen(16,8)(-4,-1)
\Pfad(-3,0),333331314441444\endPfad
\Pfad(-2,0),3334414331441\endPfad
\Pfad(-1,0),33313331444444\endPfad
\Pfad(0,0),13331433344444\endPfad
\SPfad(11,-1),22222222\endSPfad
\PfadDicke{1pt}
\Pfad(-5,6),111111111111111111111\endPfad
\DickPunkt(-3,0)
\DickPunkt(-2,0)
\DickPunkt(-1,0)
\DickPunkt(0,0)
\DickPunkt(11,0)
\DickPunkt(12,0)
\DickPunkt(13,0)
\DickPunkt(14,0)
\Kreis(1,1)
\Kreis(2,0)
\Kreis(2,4)
\Kreis(3,0)
\Kreis(3,4)
\Kreis(3,6)
\Kreis(4,0)
\Kreis(4,2)
\Kreis(4,5)
\Kreis(5,1)
\Kreis(5,2)
\Kreis(5,4)
\Kreis(6,0)
\Kreis(6,3)
\Kreis(6,4)
\Kreis(7,0)
\Kreis(7,1)
\Kreis(7,5)
\Kreis(8,0)
\Kreis(8,1)
\Kreis(8,6)
\Kreis(9,0)
\Kreis(9,2)
\Kreis(10,1)
\Label\l{2\kern-5pt}(1,1)
\Label\l{1\kern-5pt}(2,0)
\Label\l{5\kern-5pt}(2,4)
\Label\l{1\kern-5pt}(3,0)
\Label\l{5\kern-5pt}(3,4)
\Label\l{7\kern-5pt}(3,6)
\Label\l{1\kern-5pt}(4,0)
\Label\l{3\kern-5pt}(4,2)
\Label\l{6\kern-5pt}(4,5)
\Label\l{2\kern-5pt}(5,1)
\Label\l{3\kern-5pt}(5,2)
\Label\l{5\kern-5pt}(5,4)
\Label\l{1\kern-5pt}(6,0)
\Label\l{4\kern-5pt}(6,3)
\Label\l{5\kern-5pt}(6,4)
\Label\l{1\kern-5pt}(7,0)
\Label\l{2\kern-5pt}(7,1)
\Label\l{6\kern-5pt}(7,5)
\Label\l{1\kern-5pt}(8,0)
\Label\l{2\kern-5pt}(8,1)
\Label\l{7\kern-5pt}(8,6)
\Label\l{1\kern-5pt}(9,0)
\Label\l{3\kern-5pt}(9,2)
\Label\l{2\kern-5pt}(10,1)
%
\hskip7cm
$$
\centerline{\eightpoint A family of non-intersecting Motzkin paths}
\vskip7pt
\centerline{\eightpoint Figure \FL}
\endinsert

Again, we may now use the non-occupied lattice points to map the
family of Motzkin paths to certain arrays of numbers. In Figure~\FL,
we have indicated the (relevant) 
non-occupied lattice points by circles, and the
labels given equal the heights of the points plus~1. These  labels are
then assembled in an array, similar to what we did in the proof of
Proposition~\TE. The array which corresponds to the path family in
Figure~\FL\ is
$$
\matrix 
 & &7&6&5&5&6&7\\
 &5&5&3&3&4&2&2&3\\
2&1&1&1&2&1&1&1&1&2
\endmatrix$$

There are two important differences to the situation in the proof of
Proposition~\TE: first, the arrays that we obtain have the property
that entries along columns are decreasing (from top to
bottom), but there is {\it no} condition for the entries along rows.
Second, the map from path families to arrays
is {\it not\/} a bijection, but rather a many-to-one mapping. Namely,
if in a path family we replace a crossing\footnote{See again Footnote~12.} of two paths 
by two parallel horizontal
edges, or vice versa, then the corresponding array remains the same.
Thus, these arrays have to be given weights taking into account
all contributions of path families that map to the same array. 
While, given an array, it is not difficult to compute the
corresponding weight in an ad hoc fashion, at present we do not know how to do
this systematically (except for the case $m=1$).

\medskip
(4) One may think of variations and generalisations of the enumeration
results for alternating tableaux in Section~11. We mention two such
possibilities only, but will not develop them here. 

In \cite{\GeViAB, Section~7}, Gessel and Viennot describe a
generalisation of the Jacobi--Trudi formula for skew Schur functions,
the latter being the generating function for semistandard tableaux of
a given skew shape, where the relation $\le$ that holds between
entries along rows in the semistandard tableaux gets replaced by a
semitransitive relation~$R$, and the strict relation $<$ along columns
gets replaced by the ``reflection" $\bar R$ of the negation of~$R$,
$\bar R=\{(a,b):b\hbox{${}\not \kern-3ptR{}\kern3pt$}a\}$.
We could generalise the results in Section~11 to this setting.

A variation that one may think of is to consider ``weak/strict"
alternating sequences, that is, sequences $a_1\le a_2>a_3\le
a_4>a_5\le\cdots$. Analogues of Theorems~\UA\ and~\TQ\ would be easy
to obtain by modifying the corresponding proofs. In the mapping from
alternating sequences to heaps in the  proof of Theorem~\UA, the only
difference that arises is that one forbids segments $i${}$\pmb-${}$i$,
$i\ge1$. 
Everything else carries through. There are also versions of
the results in Section~11 in this setting.

\medskip
(5) When we introduced heaps of segments in Section~3 (see the
paragraphs after Theorem~\UA), we mentioned that these already
appeared earlier in \cite{\BoViAA} in a completely different
context. To be precise, in Proposition~3.4 of \cite{\BoViAA},
Bousquet-M\'elou and Viennot prove (phrased in our terminology) 
that there is a bijection
between heaps of segments with a unique maximal segment of the form
$j${}$\pmb-${}$1$ for some $j\ge1$ and (so-called) {\it parallelogram
polyominoes}. (We refer the reader to \cite{\BoViAA} for the definition
of those.) Via this bijection, the number of segments of a heap corresponds to
the width $w(P)$ of the corresponding 
polyomino~$P$ (the number of columns of the polyomino),
the height $h(P)$ of~$P$ corresponds to the sum of the lengths
of the segments plus~1, and the area $a(P)$ of~$P$ corresponds to the
sum of the $y$-coordinates of the segments. 
Moreover, if the heap is restricted to the
interval $[1,k]$ then all columns of the corresponding polyomino
have a length of at most~$k$. 
Our weights allow to keep track of all these statistics.
This has the following consequence.

\proclaim{Corollary \TTO}
For a positive integer~$k$ and a non-negative integer~$n$,
we have
$$\multline
\underset \text{column lengths }\le k\to{\sum_{P\text{ polyomino}}}
q^{a(P)}y^{h(P)}x^{2w(P)}\\
=-\frac {\dsize y \sum_{j=1}^k(-1)^{j}x^{2j}q^{\binom {j+1}2}
\sum_{i=0}^{k-j}(yq)^i\bmatrix k-i-1\\j-1\endbmatrix_q
\bmatrix i+j-1\\j-1\endbmatrix_q} 
{\dsize\sum_{j=0}^k(-1)^{j}x^{2j}q^{\binom {j+1}2}
\sum_{i=0}^{k-j}(yq)^i\bmatrix k-i\\j\endbmatrix_q
\bmatrix i+j-1\\j-1\endbmatrix_q}.
\endmultline
\tag\HEd$$
\endproclaim

\demo{Proof}
Heaps of segments on $[1,k]$ with a unique maximal segment
$j${}$\pmb-${}$1$ for some $j\ge1$ are the heaps of segments that,
via the bijection in the proof of Lemma~\UB, correspond to
alternating sequences in $\Cal A_{2n+1}^{(k)}(1\to 1)$.
The corresponding generating function is given by~(\EC) with $r=s=1$,
see~(\EJa). In our context the term~$V_1^2$ can be ignored.
The specialisation of the variables $A_i$ and $V_i$ which emulates the
various statistics is $A_i=(yq)^i$ and $V_i=y^{-i}$ for $i\ge1$.
This specialisation of denominator and numerator in~(\EJa) has been
computed in~(\EJb) and~(\EJc), respectively.\quad \quad \qed
\enddemo

Bousquet-M\'elou and Viennot \cite{\BoViAA, Prop.~4.1} derive a
formula for the limiting case $k=\infty$. Their formula is easy to
obtain from~(\HEd) by letting $k\to\infty$ and then applying the
$q$-binomial theorem (cf\. \cite{\GaRaAF, Sec.~1.3})
$$
\sum_{i=0}^\infty \bmatrix i+m-1\\m-1\endbmatrix_q z^i
=\frac {1} {(z;q)_m},
$$
where  $(z;q)_m:=\prod _{i=1} ^{m}(1-zq^{i-1})$ is the standard
{\it $q$-shifted factorial}. Indeed, we have
$$\align
{\sum_{P\text{ polyomino}}}
&q^{a(P)}y^{h(P)}x^{2w(P)}\\
&=-\frac {\dsize y \sum_{j=1}^{\infty}(-1)^{j}x^{2j}q^{\binom {j+1}2}
\sum_{i=0}^{\infty}(yq)^i\frac {1} {(q;q)_{j-1}}
\bmatrix i+j-1\\j-1\endbmatrix_q} 
{\dsize\sum_{j=0}^{\infty}(-1)^{j}x^{2j}q^{\binom {j+1}2}
\sum_{i=0}^{\infty}(yq)^i\frac {1} {(q;q)_{j}}
\bmatrix i+j-1\\j-1\endbmatrix_q}\\
&=-\frac {\dsize y \sum_{j=1}^{\infty}
\frac {(-1)^{j}x^{2j}q^{\binom {j+1}2}} {(q;q)_{j-1}\,(yq;q)_j}}
{\dsize\sum_{j=0}^{\infty}
\frac {(-1)^{j}x^{2j}q^{\binom {j+1}2}} {(q;q)_{j}\,(yq;q)_j}}.
\endalign
$$

We refer the reader to \cite{\BoViAA} and \cite{\FereAA} for more
material and references on generating functions for parallelogram
polyominoes (and other classes of polyominoes). 

\medskip
(6)
Next we report
a curious fact for Hankel determinants of
numbers of up-down paths in ``very tight" strips,
implying analogous results for Hankel determinants of numbers of
bounded alternating sequences.

\proclaim{Theorem \TTP}
Let $n,k,r,s$ be non-negative integers with $0\le r,s\le 2k-1$.
Then
$$\multline
\det\left(C_{2n+2i+2j+r+s}^{(2k-1)}(r\to s)\right)_{0\le i,j\le k-1}\\
=\cases
0,
&\kern-3cm
\text{if }\gcd(r+1,2k+1)\ne1\text{ or }\gcd(s+1,2k+1)\ne1,\\
(-1)^{\sum_{i=1}^k\left(\fl{\frac {i(r+1)}{2k+1}}+\fl{\frac
{i(s+1)}{2k+1}}\right)},\\
&\kern-3cm
\text{if }\gcd(r+1,2k+1)=\gcd(s+1,2k+1)=1.
\endcases
\endmultline
\tag\HEa$$
\endproclaim

\demo{Proof}
This follows by combining Propositions~\TTQ\ and~\TTR.\quad \quad \qed
\enddemo

\proclaim{Proposition \TTQ}
With the assumptions of Theorem~{\rm\TTP}, we have
$$\multline
\det\left(C_{2n+2i+2j+r+s}^{(2k-1)}(r\to s)\right)_{0\le i,j\le k-1}\\
=
\det\left(C_{2i}^{(2k-1)}(r\to 2j+\chi(\text{$r$ odd}))\right)_{0\le i,j\le k-1}\\
\times
\det\left(C_{2j}^{(2k-1)}(2i+\chi(\text{$s$ odd})\to s)\right)_{0\le i,j\le k-1}.
\endmultline
\tag\HEb$$
Here, $\chi(\Cal A)=1$ if $\Cal A$ is true and $\chi(\Cal A)=0$ otherwise.
\endproclaim

\demo{Proof}
We show the assertion by using again non-intersecting lattice paths.
By the Lindstr\"om--Gessel--Viennot theorem \cite{\LindAA, Lemma~1},
the determinant on the left-hand side of (\HEb) equals 
$\sum_{\Cal P}\sgn\Cal P$, where the sum is over all
families $\Cal P=(P_0,P_1,\dots,P_{k-1})$ of non-intersecting up-down
paths of height at most $2k-1$ that do not pass below the $x$-axis, 
where $P_i$ runs from
$(-2i,r)$ to $(2n+2\si(i)+r+s,s)$, $i=0,1,\dots,k-1$, for some
permutation~$\si\in\frak S_k$;
again, the sign $\sgn\Cal P$ is by definition equal to
$\sgn\si$. 
Figure~\FM\ shows an example for $n=2$, $k=4$, $r=2$, $s=3$,
and the permutation~$\si$ is given by $\si(0)=0$, $\si(1)=1$,
$\si(2)=3$, $\si(3)=2$.

\midinsert
$$
\Gitter(17,9)(-7,-1)
\Koordinatenachsen(17,9)(-7,-1)
\Pfad(-6,2),3334333434343434444\endPfad
\Pfad(-4,2),4434343434343344333\endPfad
\Pfad(-2,2),3334343434344\endPfad
\Pfad(0,2),343434343\endPfad
\SPfad(9,-1),22222222\endSPfad
\PfadDicke{1pt}
\Pfad(-8,7),111111111111111111111111\endPfad
\DickPunkt(-6,2)
\DickPunkt(-4,2)
\DickPunkt(-2,2)
\DickPunkt(0,2)
\DickPunkt(9,3)
\DickPunkt(11,3)
\DickPunkt(13,3)
\DickPunkt(15,3)
\hskip4.5cm
$$
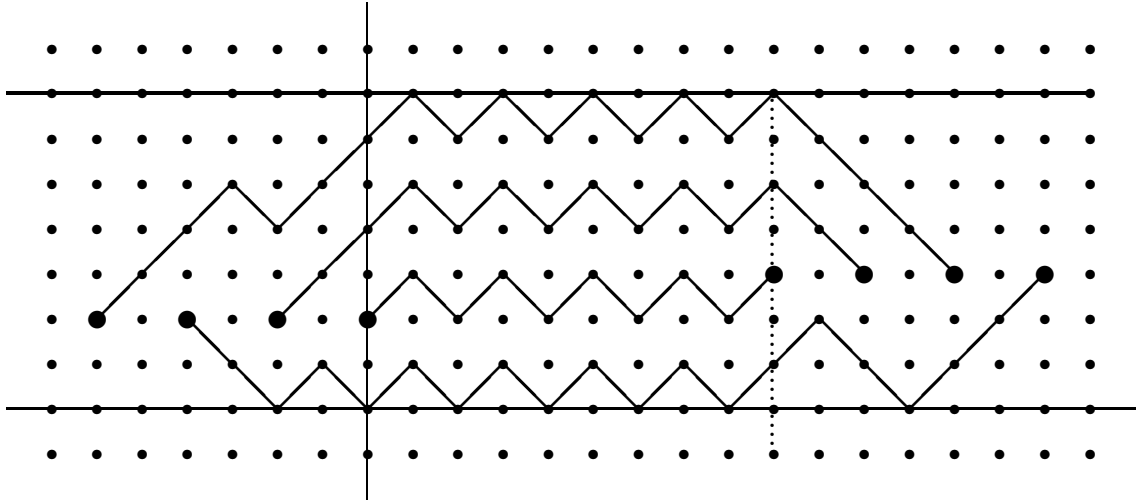
\centerline{\eightpoint A family of non-intersecting up-down paths
in a ``very tight" strip}
\vskip7pt
\centerline{\eightpoint Figure \FM}
\endinsert

Because of the upper bound on the height of the paths, they are
uniquely determined in the region 
between the abscissa $x=0$ and $x=2n+r+s$.
In the figure, this is the region between the $y$-axis and the dotted
vertical line. If we cut these path portions, then
to the left of the $y$-axis there remains a family
$\Cal Q=(Q_0,Q_1,\dots,Q_{k-1})$ of non-intersecting up-down
paths of height at most $2k-1$ that do not pass below the $x$-axis, 
where $Q_i$ runs from
$(-2i,r)$ to $(0,2\tau(i)+\chi(r\text{ odd}))$, $i=0,1,\dots,k-1$, for some
permutation~$\tau\in\frak S_k$, 
and to the right of the vertical line $x=2n+r+s$
there remains a family
$\Cal R=(R_0,R_1,\dots,R_{k-1})$ of non-intersecting up-down
paths of height at most $2k-1$ that do not pass below the $x$-axis, 
where $R_i$ runs from
$(2n+r+s,2i+\chi(s\text{ odd}))$ to $(2n+r+s+2\rho(i),s)$, 
$i=0,1,\dots,k-1$, for some permutation~$\rho\in\frak S_k$.
By the Lindstr\"om--Gessel--Viennot theorem,
the sum $\sum_{\Cal Q}\sgn \Cal Q$ is given by the first determinant
on the right-hand side of~(\HEb), and the sum $\sum_{\Cal R}\sgn \Cal
R$ is given by the second determinant on the right-hand side
of~(\HEb).\quad \quad \qed
\enddemo

\proclaim{Proposition \TTR}
Let $k$ and $s$ be non-negative integers with $0\le s\le 2k-1$.
Then
$$\multline
\det\left(C_{2j}^{(2k-1)}(2i+\chi(s\text{ odd})\to 
s)\right)_{0\le i,j\le k-1}\\
=\cases
0,&\text{if }\gcd(s+1,2k+1)\ne1,\\
(-1)^{\sum_{i=1}^k\fl{\frac {i(s+1)}{2k+1}}},&\text{if }\gcd(s+1,2k+1)=1.
\endcases
\endmultline
\tag\HEc$$
\endproclaim

\demo{Proof}
It would be desirable to find a combinatorial proof of this assertion.
Since we failed to find such proof, we base our arguments on a
formula for bounded up-down paths that goes back to Laplace, followed
by determinant manipulations.

We write $D_k(s)$ for the determinant on the left-hand side of~(\HEc).
We first treat the case where $s$~is odd, say $s=2S+1$.
We use the well-known formula (see e.g\. \cite{\KratCL, Eq.~(10.13)})
$$
\frac {2} {K+2}\sum _{k=1} ^{K+1}
\left(2\cos \frac {\pi k} {K+2}\right)^\ell\cdot
\sin\frac {\pi k(r+1)} {K+2}\cdot \sin\frac
{\pi k(s+1)} {K+2}
$$
for the number of up-down paths of height at most~$K$ from $(0,r)$ to
$(\ell,s)$ that do not pass below the $x$-axis.
The entry $C_{2j}^{(2k-1)}(2i+1\to 2S+1)$ in the determinant in~(\HEc)
(where~$s$ was replaced by~$2S+1$)
is the number of up-down paths of height at most~$2k-1$ from
$(0,2i+1)$ to $(2j,2S+1)$ that do not pass below the $x$-axis.
Therefore we have
$$\align
D_k&(2S+1)\\
&=\det\Bigg(\frac {2} {2k+1}\sum _{\ell=1} ^{2k}
\left(2\cos \frac {\pi \ell} {2k+1}\right)^{2j}\cdot
\sin\frac {\pi \ell(2i+2)} {2k+1}\cdot \sin\frac
{\pi \ell(2S+2)} {2k+1}\Bigg)_{0\le i,j\le k-1}\\
&=\frac {2^k} {(2k+1)^k}\sum _{1\le \ell_0,\dots,\ell_{k-1}\le 2k} 
\Bigg(\prod _{i=0} ^{k-1}
\sin\frac {\pi \ell_i(2i+2)} {2k+1}\cdot \sin\frac
{\pi \ell_i(2S+2)} {2k+1}\Bigg)\\
&\kern7cm
\cdot
\det\Bigg(
\left(2\cos \frac {\pi \ell_i} {2k+1}\right)^{2j}
\Bigg)_{0\le i,j\le k-1}.
\endalign$$
We take the average of the last expression over all permutations
of the $\ell_i$'s and obtain
$$\align 
D_k(2S+1)&=\frac {2^k} {(2k+1)^k}\cdot\frac {1} {k!}
\sum_{\si\in \frak S_k}\sum _{1\le \ell_0,\dots,\ell_{k-1}\le 2k} 
\Bigg(\prod _{i=0} ^{k-1}
\sin\frac {\pi \ell_{\si(i)}(2i+2)} {2k+1}\cdot \sin\frac
{\pi \ell_i(2S+2)} {2k+1}\Bigg)
\\
&\kern5cm
\cdot\sgn\si\cdot
\det\Bigg(
\left(2\cos \frac {\pi \ell_i} {2k+1}\right)^{2j}
\Bigg)_{0\le i,j\le k-1}\\
&=\frac {2^k} {(2k+1)^k}
\sum _{1\le \ell_0<\dots<\ell_{k-1}\le 2k} 
\Bigg(\prod _{i=0} ^{k-1}
 \sin\frac
{2\pi \ell_i(S+1)} {2k+1}\Bigg)
\\
&\kern1cm
\cdot\det\Bigg(\sin\frac {2\pi \ell_j(i+1)} {2k+1}\Bigg)_{0\le i,j\le k-1}
\cdot
\det\Bigg(
\left(2\cos \frac {\pi \ell_i} {2k+1}\right)^{2j}
\Bigg)_{0\le i,j\le k-1}.\\
\tag\HF
\endalign$$
The last determinant in (\HF), a Vandermonde determinant, will vanish
whenever $\ell_{i_1}=\ell_{i_2}$ for $i_1\ne i_2$, but also whenever
$\ell_{i_1}=2k+1-\ell_{i_2}$ since $\cos(\pi-x)=-\cos x$.
Hence, non-zero contributions to the multiple sum in~(\HF) arise
only if for the set $\{\ell_0,\ell_1,\dots,\ell_{k-1}\}$ of indices
we choose exactly one element out of $\{1,2k\}$, exactly one element
out of $\{2,2k-1\}$, \dots, and exactly one element out of
$\{k,k+1\}$. Moreover, each such choice produces the same summand
as is not difficult to see using the relations $\cos(\pi-x)=-\cos x$,
$\sin(2\pi+x)=\sin x$, and $\sin(-x)=-\sin x$. These arguments imply
that the sum equals $2^k$ times the summand for the choice
$\ell_i=i+1$, $i=0,1,\dots,k-1$. Thus, we have
$$\multline
D_k(2S+1)
=\frac {2^{2k}} {(2k+1)^k}
\Bigg(\prod _{i=1} ^{k}
 \sin\frac
{2\pi i(S+1)} {2k+1}\Bigg)
\\
\cdot\det\Bigg(\sin\frac {2\pi (j+1)(i+1)} {2k+1}\Bigg)_{0\le i,j\le k-1}
\cdot
\det\Bigg(
\left(2\cos \frac {\pi (i+1)} {2k+1}\right)^{2j}
\Bigg)_{0\le i,j\le k-1}.
\endmultline$$
Both these  determinants can be evaluated in closed form. Namely,
we have (cf\. \cite{\KratBT, Eq.~(5.4)})
$$\multline
\det\left(\sin\frac {\pi (i+1)\la_j}m\right)_{0\le i,j\le k-1}\\
=2^{k^2-k}\prod
_{0\le i<j\le k-1} ^{}\sin\frac {\pi(\la_i-\la_j)} {2m}\prod _{0\le i\le
j\le k-1} ^{}\sin\frac {\pi(\la_i+\la_j)} {2m},
\endmultline
$$
and in regard of the Vandermonde determinant we
have
$$
\det\Bigg(
\left(2\cos \frac {\pi \la_i} {m}\right)^{2j}
\Bigg)_{0\le i,j\le k-1}
=2^{k^2-k}
\prod _{0\le i<j\le k-1} ^{}
\sin \frac {\pi (\la_i-\la_j)} {m}
\sin \frac {\pi (\la_i+\la_j)} {m}.
$$
Substituting this back in the last expression that we obtained in our
computation, we get
$$\align 
D_k(2S+1)
&=\frac {2^{2k}} {(2k+1)^k}
\Bigg(\prod _{i=1} ^{k}
 \sin\frac
{2\pi i(S+1)} {2k+1}\Bigg)
\\
&\kern1cm
\cdot
2^{k^2-k}\Bigg(\prod
_{1\le i<j\le k} ^{}\sin\frac {\pi(i-j)} {2k+1}\prod _{1\le i\le
j\le k} ^{}\sin\frac {\pi(i+j)} {2k+1}\Bigg)\\
&\kern1cm
\cdot
2^{k^2-k}\Bigg(
\prod _{1\le i<j\le k} ^{}
\sin \frac {\pi (i-j)} {2k+1}
\sin \frac {\pi (i+j)} {2k+1}\Bigg)\\
&=\frac {2^{2k^2}} {(2k+1)^k}
\Bigg(\prod _{i=1} ^{k}
 \sin\frac
{2\pi i(S+1)} {2k+1}\Bigg)
\\
&\kern1cm
\cdot
\Bigg(\prod_{i=1} ^{k}
\bigg(\sin\frac {\pi i} {2k+1}\bigg)^{2k-2}\Bigg)
\Bigg(\prod
_{i=1} ^{k}\sin\frac {2\pi i} {2k+1}\Bigg).
\endalign$$
By Lemma~\TTS\ below, we have $\prod_{i=1} ^{k}
\sin\frac {\pi i} {2k+1}=2^{-k}\sqrt{2k+1}$. However, we also have
$$
\prod _{i=1} ^{k}\sin\frac {2\pi i} {2k+1}=
\prod _{i=1} ^{k}\sin\frac {\pi i} {2k+1}=2^{-k}\sqrt{2k+1}
$$
since $\sin\frac {2\pi i} {2k+1}=\sin\frac {\pi (2k+1-2i)} {2k+1}$
for all $i=\cl{\frac {k+1} {2}},\dots,k-1,k$. Putting these
observations together, we see that
$$
D_k(2S+1)
=\frac {2^{k}} {\sqrt{2k+1}}
\prod _{i=1} ^{k}
 \sin\frac
{2\pi i(S+1)} {2k+1}.
\tag\HG$$
If $\gcd(2S+2,2k+1)\ne1$, then for $i=(2k+1)/\gcd(2S+2,2k+1)$ we
have $ \sin\frac {2\pi i(S+1)} {2k+1}=0$, and consequently $D_k(2S+1)=0$.

On the other hand, if $\gcd(2S+2,2k+1)=1$, then we claim that,
up to a sign, the product on the right-hand side of~(\HG) 
is also equal to $\prod _{i=1} ^{k}
 \sin\frac {\pi i} {2k+1}=2^{-k}\sqrt{2k+1}$. Indeed,
for $1\le i\ne j\le k$, we have
 $2i(S+1)\hbox{${}\not\equiv{}$}\pm2j(S+1)$~(mod~$2k+1$), and
 therefore the set 
$$\{2i(S+1)\text{ mod $(2k+1)$}:i=1,2,\dots,k\}$$
contains exactly one element of the pair $\{j,2k+1-j\}$, for
$j=1,2,\dots,k$, proving our claim since $\sin(-x)=-\sin x$ and
$\sin(\pi-x)=\sin x$. The sign of $\sin\frac {2\pi i(S+1)} {2k+1}$
depends on the ``location" of $2i(S+1)$; more precisely, if
$2i(S+1)\in [m\pi,(m+1)\pi)$, then this sign is positive if $m$~is
even and this sign is negative otherwise. If these considerations are
used in~(\HG), then we arrive at
$$
D_k(2S+1)
=(-1)^{\sum_{i=1}^k\fl{2i(S+1)/(2k+1)}},
$$
provided $\gcd(2S+2,2k+1)=1$. This is equivalent to the assertion of
the proposition for $s=2S+1$.

\medskip
Now we turn to the case where $s$~is even, say $s=2S$.
The reader should recall that the determinant on the left-hand side of~(\HEc)
provides the signed enumeration of certain non-intersecting up-down paths.
Hence, by a reflection in the horizontal line $y=(2k-1)/2$, 
this case is transformed
into the odd case, where the ending height of the paths is now
$2k-1-2S$. If we substitute this in the formula on the right-hand
side of~(\HEc) and take into account that the reflection reversed
the order of the starting points, creating a sign of $(-1)^{\binom
k2}$, then after some manipulation we see that formula~(\HEc) also
holds in the even case.

\medskip
This finishes the proof of the proposition.\quad \quad \qed
\enddemo

\proclaim{Lemma \TTS}
For non-negative integers~$n$, we have
$$
\prod _{j=1} ^{n-1}2\sin\frac {j\pi} {n}=n.
$$
\endproclaim
\demo{Proof}
Let $\om=e^{\pi \bold i/n}$, where $\bold  i=\sqrt{-1}$. Then we have
$$\align 
\prod _{j=1} ^{n-1}2\sin\frac {j\pi} {n}
&=
(-\bold i)^{n-1}\prod _{j=1} ^{n-1}\left(\om^j-\om^{-j}\right)\\
&=
\bold i^{n-1}\om^{-\binom {n}2}
\prod _{j=1} ^{n-1}\left(1-\om^{2j}\right)=n,
\endalign$$
since $\om^2$ is a primitive $n$-th root of unity.\quad \quad \qed
\enddemo

\medskip
(7) It should be noted that all results of this article where this
makes sense also hold for negative values of the parameter~$n$.
In order to illustrate this: let us consider Theorem~\TD.
Our claim is that Equation~(\AA) also holds for negative values
of~$n$. The reason is simple: since the sequence
$(C_{2n}^{2k+2m-1})_{n\ge0}$ satisfies a linear recurrence with
constant coefficients (dictated by the denominator of the rational
function on the right-hand side of~(\BGa) with $r=s=0$ and $k$
replaced by $2k+2m-1$), this implies a linear recurrence with constant
coefficients for the left-hand side of~(\AA) as well as for the
right-hand side. Since~(\AA) holds for all non-negative integers~$n$,
these recurrences for the left-hand side and for the right-hand side
must be the same. Consequently, Equation~(\AA) continues to hold also
for negative values of~$n$. 

In this case, however, this is not very exciting since it turns out
that, in this way, one obtains again the same identity. Namely, let us
replace~$n$ by~$-n$ in~(\AA), to get
$$
\det\left(C_{-2n+2i+2j+4m-2}^{(2k+2m-1)}\right)_{0\le i,j\le k-1}
=
\det\left(C_{2n-2i-2j}^{(2k+2m-1)}\right)_{0\le i,j\le m-1}.
$$
Next we reverse the order of rows and columns in both determinants.
This yields
$$
\det\left(C_{-2n+4k-2i-2j+4m-4}^{(2k+2m-1)}\right)_{0\le i,j\le k-1}
=
\det\left(C_{2n-4m+2i+2j+2}^{(2k+2m-1)}\right)_{0\le i,j\le m-1}.
$$
Finally we replace $n$ by $n+2k+2m-2$. Thus, we arrive at
$$
\det\left(C_{-2n-2i-2j}^{(2k+2m-1)}\right)_{0\le i,j\le k-1}
=
\det\left(C_{2n+2i+2j+4k-2}^{(2k+2m-1)}\right)_{0\le i,j\le m-1}.
$$
Clearly, this is (\AA) in which $k$ and $m$ are interchanged.

\medskip
It is more interesting to apply this argument in, say, Theorem~\TTP.
Let us replace $n$ by $-n$ in~(\HEa). The determinant on the left-hand
side then becomes
$$
\det\left(C_{-2n+2i+2j+r+s}^{(2k-1)}(r\to s)\right)_{0\le i,j\le k-1}.
$$
In view of Corollary~\UF, this determinant can be rewritten as 
a determinant of numbers of bounded alternating sequences. More
precisely, we have
$$\multline
\det\left(C_{-2n+2i+2j+2r+2s-4}^{(2k-1)}(2r-2\to 2s-2)\right)_{0\le i,j\le
k-1}\\
=(-1)^{k(r+s)}
\det\left(\big\vert A_{2n-2i-2j-2r-2s+5}^{(k)}(r\to s)\big\vert
\right)_{0\le i,j\le k-1}
\endmultline
$$
and
$$\multline
\det\left(C_{-2n+2i+2j+2r+2s-3}^{(2k-1)}(2r-2\to 2s-1)\right)_{0\le i,j\le
k-1}\\
=(-1)^{k(r+s)}
\det\left(\big\vert A_{2n-2i-2j-2r-2s+4}^{(k)}(r\to s)\big\vert
\right)_{0\le i,j\le k-1}.
\endmultline
$$
In both determinants on the right-hand sides we reverse the order of
rows and columns and then recall that, by Proposition~\TTQ, all these
determinants are independent of~$n$. This leads to
$$\multline
\det\left(C_{-2n+2i+2j+2r+2s-4}^{(2k-1)}(2r-2\to 2s-2)\right)_{0\le i,j\le
k-1}\\
=(-1)^{k(r+s)}
\det\left(\big\vert A_{2n+2i+2j+1}^{(k)}(r\to s)\big\vert
\right)_{0\le i,j\le k-1}
\endmultline
$$
and
$$\multline
\det\left(C_{-2n+2i+2j+2r+2s-3}^{(2k-1)}(2r-2\to 2s-1)\right)_{0\le i,j\le
k-1}\\
=(-1)^{k(r+s)}
\det\left(\big\vert A_{2n+2i+2j}^{(k)}(r\to s)\big\vert
\right)_{0\le i,j\le k-1}.
\endmultline
$$
In combination with Theorem~\TTP, this proves the following result.

\proclaim{Corollary \TTT}
Let $n,k,r,s$ be non-negative integers with $0\le r,s\le 2k-1$.
Then
$$\multline
\det\left(\big\vert\Cal A_{2n+2i+2j+1}^{(k)}(r\to s)\right)_{0\le i,j\le k-1}\\
=\cases
0,
&\kern-3cm
\text{if }\gcd(2r-1,2k+1)\ne1\text{ or }\gcd(2s-1,2k+1)\ne1,\\
(-1)^{k(r+s)+\sum_{i=1}^k\left(\fl{\frac {i(2r-1)}{2k+1}}+\fl{\frac
{i(2s-1)}{2k+1}}\right)},\\
&\kern-3cm
\text{if }\gcd(2r-1,2k+1)=\gcd(2s-1,2k+1)=1.
\endcases
\endmultline
\tag\HH$$
and
$$\multline
\det\left(\big\vert\Cal A_{2n+2i+2j}^{(k)}(r\to s)\right)_{0\le i,j\le k-1}\\
=\cases
0,
&\kern-3cm
\text{if }\gcd(2r-1,2k+1)\ne1\text{ or }\gcd(2s,2k+1)\ne1,\\
(-1)^{k(r+s)+\sum_{i=1}^k\left(\fl{\frac {i(2r-1)}{2k+1}}+\fl{\frac
{2is}{2k+1}}\right)},\\
&\kern-3cm
\text{if }\gcd(2r-1,2k+1)=\gcd(2s,2k+1)=1.
\endcases
\endmultline
\tag\HI$$
\endproclaim

\medskip
We present one more application of this extension idea.
It concerns the identitites~(\HBb) and~(\HBd).
Let us again write $a(n,k)$ for $\sum_{s=0}^{k} C_{n}^{(k)}(0\to s)$.
Inspection of numerator and denominator degrees in~(\HBf)
with $k$ replaced by~$4k-2$ shows that the sequence
$(a(n,4k-2))_{n\ge0}$ can be extended to negative integers~$n$ by
using the linear recurrence with constant coefficients that the
sequence satisfies. This is not quite so for the sequence
$(a(n,4k))_{n\ge0}$ since numerator and denominator degrees in~(\HBf)
with $k$ replaced by~$4k$ are the same, causing the constant term
to violate the linear recurrence with constant coefficients that the
sequence satisfies otherwise. To ``repair" this, we
define the modified sequence $(\tilde a(n,4k))_{n\ge0}$ by
$\tilde a(n,4k)=a(n,4k)$ for $n\ge1$, and $\tilde a(0,4k)=\frac
{2k}{2k+1}$. This sequence can be extended to negative indices.
Then, using the same reasoning as above, we
obtain the following corollary from Theorem~\TTLb.

\proclaim{Corollary \TTU}
Let $n$ be a non-negative integer and $k$ be a positive integer.
Then we have
$$
\det\big(
\tilde a(-n-i-j-1,4k)\big)_{0\le i,j\le 2k-1}
=(-1)^{k(n-1)}(2k+1)^{-n-4k}
\tag\HBc
$$
and
$$
\det\big(
a(-n-i-j,4k-2)\big)_{0\le i,j\le 2k-1}
=(-1)^{kn}2^{-n-4k+2}.
\tag\HBe
$$
\endproclaim

\demo{Proof}
In (\HBb) and (\HBd), we replace $n$ by $-n$. Furthermore, we
reverse the order of rows and solumns of the matrices, that is,
we replace $i$ by $2k-1-i$ and $j$ by $2k-1-j$. To complete the proof,
$n$ has to be shifted by~$4k$ to obtain~(\HBc), and by $4k-2$
to obtain~(\HBe).\quad \quad \qed
\enddemo


\medskip
{\smc Final Note.} Since the original version of this article appeared on
the {\tt ar$\chi$iv} in 2020, several further developments related to parts
of this article have occurred that we want to mention here. More precisely:

\roster 
\item Zaimi~\cite{\ZaimAA}, by a clever use of the transfer matrix idea,
has proved a general reciprocity theorem for non-intersecting paths in periodic
planar networks, which has now appeared in paper form in the joint article~\cite{\HoZaAA}
with Hopkins. It covers our
Theorem~\TJ, and the weighted version covers Theorem~\TTC.
With some adaptation (respectively preprocessing),
one might also cover Theorems~\TM\ and~\TTF, and provide
a proof of Conjecture~\TTN.
However, the connection with (arrays of) alternating sequences is not a part
in this picture.
\item In~\cite{\JaKKAA}, Jang et al\. develop an interesting framework for
computing moments of orthogonal polynomials ``with negative indices".
In particular, this leads them to a combinatorial model for
the number of bounded Motzkin paths ``with negative indices" that
involves generalisations of alternating sequences. These
developments allow them in the end to 
give proofs of both Conjectures~\TTLa\ and~\TTN, but are in the first place
extremely interesting in their own right.
\item In our article, the ``duality" between up-down paths and alternating
sequences, expressed by the reciprocity relations in Corollaries~\UF--\TC\
and further results in subsequent sections, has been presented as an {\it observation},
by which we mean: without an intrinsic, conceptual explanation.
Viennot~\cite{\VienZZ} has found a beautiful such explanation based on the heap
pictures for up-down paths and alternating sequences that we presented in
Sections~2 and~3, respectively, via a
combinatorial reciprocity between dimers and segments.
\endroster


\Refs

\ref\no \AnBBAA\by N. Anzalone, J. Baldwin, I. Bronshtein and T. Kyle Petersen
\paper A reciprocity theorem for monomer-dimer coverings\inbook
Discrete models for complex systems, DMCS '03 (Lyon)\publ
Discrete Math. Theor. Comput. Sci. Proc., vol.~AB\publaddr Nancy\yr 2003\pages 179--193\endref

\ref\no \BaRoAA\by Y. Baryshnikov and D. Romik\paper 
Enumeration formulas for Young tableaux in a diagonal strip\jour 
Israel J. Math\.\vol 178\yr 2010\pages 157--186\endref

\ref\no \BeSaAA\by M. Beck and R. Sanyal\book Combinatorial reciprocity
theorems. An invitation to enumerative geometric
combinatorics\publ
 Graduate Studies in Mathematics, vol.~195, American
Mathematical Society\publaddr Providence, RI\yr 2018\endref

\ref\no \BoViAA\by M. Bousquet-M\'elou and X. Viennot\paper 
Empilements de segments et 
$q$-\'enum\'eration de polyominos convexes dirig\'es\jour
J. Combin\. Theory Ser.~A\vol 60\yr 1992\pages 196--224\endref

\ref\no \CaFoAA\by P. Cartier and D. Foata\book
Probl\`emes combinatoires de commutation et r\'earrangements\publ
Lecture Notes in Mathematics, no.~85,
Springer--Verlag\publaddr Berlin, New York\yr 1969\finalinfo
republished in the ``books" section of the S\'eminaire
Lotharingien de Combinatoire\endref

\ref\no \CiglBD\by J. Cigler\paper Some remarks and conjectures related to
lattice paths in strips along the $x$‐axis\paperinfo manuscript, 
{\tt ar$\chi$iv:1501.04750}\endref

\ref\no \FavaAA\by J. Favard\paper Sur les polynomes de Tchebicheff
\jour Comp\-tes Rendus Acad\. Sci\. Paris\vol 200\yr 1935\pages 
2052--2053\endref

\ref\no \FereAA\by S. Fereti\'c\paper 
An alternative method for $q$-counting directed
column-convex polyominoes\jour Discrete Math\.\vol 210\yr2000\pages
55--70\endref

\ref\no \FoStAA\by S. Formichella and A. Straub\paper 
Gaussian binomial coefficients with negative arguments\jour 
Ann\. Combin\.\vol 23\yr 2019\pages 725--748\endref

\ref\no \GaEgAA\by A. M. Garsia and \"O. E\u gecio\u glu\book
Lectures in Algebraic Combinatorics\publ Lecture Notes in Math., vol.~2277
\publaddr Springer, Cham\yr 2020\endref

\ref\no \GaRaAF\by G.    Gasper and M. Rahman \yr 2004 \book Basic
hypergeometric series\bookinfo second edition\publ Encyclopedia of
Math\. And Its Applications~96, Cambridge University Press\publaddr
Cambridge\endref 

\ref\no \GeViAA\by I. M. Gessel and X. G. Viennot \yr 1985 \paper
Binomial determinants, paths, and hook length formulae\jour Adv\. 
Math\. \vol 58\pages 300--321\endref 

\ref\no \GeViAB\by I. M. Gessel and X. G. Viennot \yr 1989 \paper
Determinants, paths, and plane partitions \paperinfo preprint,
1989\finalinfo available at {\tt
http://www.cs.brandeis.edu/\~{}ira}\endref 

\ref\no \GoHaAA\by A. M. Hamel and I. P. Goulden \yr 1995 \paper 
Planar decompositions of tableaux and Schur
function determinants\jour Europ\. J.~Combin\.
\vol 16 \pages 461--477\endref

\ref\no \HoZaAA\by S. Hopkins and G. Zaimi \yr 2023\paper 
Combinatorial reciprocity for non-intersecting paths\jour 
Enum. Combin. Appl.\vol 3\pages Article~S2R13, 13~pp\endref

\ref\no \JaKKAA\by J. Jang, D. Kim, J. S. Kim, M. Song and U-K. Song \yr \paper 
Negative moments of orthogonal polynomials\jour
preprint, {\tt ar$\chi$iv:2201.11344}\vol \pages \endref

\ref\no \JinYAA\by E. Y. Jin \yr 2018 \paper Outside nested
decompositions of skew
diagrams and Schur function determinants\jour Europ\. J.~Combin\.
\vol 67 \pages 239--267\endref

\ref\no \KratAH\by C.    Krattenthaler \yr 1989 \paper Enumeration of
lattice paths and generating functions for skew plane partitions\jour
Manuscripta Math\.\vol63\pages 129--156\endref  

\ref\no \KratBP\by C.    Krattenthaler \yr 2001 \paper Permutations
with restricted patterns and Dyck paths\jour Adv\. Appl\. Math\.\vol
27\pages 510--530\endref 

\ref\no \KratAZ\by C.    Krattenthaler \yr 2006 \paper The theory of heaps
and the Cartier--Foata monoid\inbook appendix to reference [\CaFoAA]
in the ``books" section of the S\'eminaire
Lotharingien de Combinatoire\endref  

\ref\no \KratBT\by C.    Krattenthaler \yr 2007 \paper Asymptotics for
random walks in alcoves of affine Weyl groups\jour S\'e\-mi\-naire
Lothar\-ing\-ien Combin\.\vol 52\pages Article~B52i, 72~pp\endref

\ref\no \KratCL\by C.    Krattenthaler \yr 2015 \paper Lattice path
enumeration\inbook Handbook of Enumerative Combinatorics\ed
M.~B\'ona\publ CRC Press\publaddr Boca Raton, London, New York\pages
589--678\endref 

\ref\no \KratCM\by C.    Krattenthaler \yr 2016 \paper Plane
partitions in the
work of Richard Stanley and his school\inbook The Mathematical Legacy
of Richard P. Stanley\eds P.~Hersh, T.~Lam, P.~Pylyavskyy, V.~Reiner\publ
Amer\. Math\. 
Soc\.\publaddr R.I.\yr 2016\pages 246--277\endref

\ref\no \KrYaAA\by C.    Krattenthaler and D. Yaqubi \yr 2018 \paper
Some determinants of path generating functions, II\jour
Adv\. Appl\. Math\.\vol 101\pages 232--265\endref 

\ref\no \LaPrAC\by A.    Lascoux and P. Pragacz \yr 1988 \paper Ribbon
Schur functions\jour Europ\. J. Combin\.\vol 9\pages 561--574\endref 

\ref\no \LindAA\by B.    Lindstr\"om \yr 1973 \paper On the vector
representations of induced matroids\jour Bull\. London
Math\. Soc\.\vol 5\pages 85--90\endref 

\ref\no \LoebAA\by D. Loeb\paper Sets with a negative number of 
elements\jour Adv\. Math\.\vol 91\pages 64--74\yr 1992\endref

\ref\no \OwPrAA\by A. L. Owczarek and T. Prellberg\paper Enumeration
of area-weighted Dyck paths with restricted height\jour Australas. J. 
Combin\.\vol 54\yr2012\pages 13--18\endref

\ref\no \PropXA \by J. Propp\paper
A reciprocity theorem for domino tilings\jour
Electron\. J. Combin\. \vol 8 \yr2001\pages Research Paper~18, 9~pp\endref

\ref\no \SpeyXA\by D. Speyer\paper A reciprocity sequence for perfect
matchings of linearly growing graphs\paperinfo unpublished manuscript;
available at {\tt
http://www-personal.umich.edu/\~{}speyer/TransferMatrices.pdf}\endref

\ref\no \StanAY\by R. P. Stanley\paper 
Combinatorial reciprocity theorems\jour Adv\.
 Math\.\vol 14\yr 1974\pages 194--253\endref

\ref\no \StanAZ\by R. P. Stanley\paper
Decompositions of rational convex polytopes\jour Ann\. Discrete
Math\.\vol 6\yr 1980\pages 333--342\endref

\ref\no \StanXA\by R. P. Stanley\paper
On dimer coverings of rectangles of fixed width\jour
Discrete Appl\. Math\.\vol 12 \yr 1985\pages 81--87\endref

\ref\no \StanBI\by R. P. Stanley \yr 1999 \book Enumerative
Combinatorics\bookinfo Vol.~2\publ Cambridge University Press\publaddr
Cambridge\endref  

\ref\no \StanAP\by R. P. Stanley \yr 2012 \book Enumerative
Combinatorics\bookinfo Vol.~1, second edition
\publ  Cambridge University Press\publaddr Cambridge\endref 

\ref\no \VienAE\by X.  G.  Viennot \yr 1983 \book Une th\'eorie combinatoire 
des polyn\^omes orthogonaux g\'en\'eraux\publ UQAM\publaddr Montr\'eal, 
Qu\'e\-bec\finalinfo
available at {\tt http://www.xavierviennot.org/xavier/polynomes$\underline{\ }$orthogonaux.html}\endref

\ref\no \VienAF\by X. G. Viennot\paper
Heaps of pieces. I. Basic definitions and combinatorial lemmas\inbook
Lecture Notes in Math., vol.~1234\publ Springer\publaddr Berlin
\yr 1986\pages 321--350\endref

\ref\no \VienZZ\by X. G. Viennot\paper
Lattice paths and heaps\paperinfo lecture given at the Conference on
``Lattice Paths, Combinatorics and Interactions'',
CIRM, Luminy, 2021; slides and video available at
{\tt https://viennot.org}\endref

\ref\no \ZaimAA\by G. Zaimi\paper Answer to the posting
``Reciprocity for fans of bounded Dyck paths" by Sam Hopkins
on\linebreak {\tt MathOverflow}
\paperinfo available at
{\tt https://mathoverflow.net/questions/373030/reciprocity-for-fans-of}\linebreak
{\tt -bounded-dyck-paths}
\endref

\endRefs
\enddocument